\documentclass[11pt, oneside]{article}

\usepackage{amsmath, amssymb, amsthm}
\usepackage{geometry}
\usepackage{enumitem}
\usepackage{graphicx}
\usepackage{capt-of}
\usepackage{makeidx}
\usepackage[bookmarks]{hyperref}
\usepackage{morefloats}

\newtheoremstyle{thmstyle}
  {1.2\baselineskip}         
  {.5\baselineskip}         
  {\itshape}    
  {0pt}         
  {\bfseries}   
  {}            
  {.5em}        
  {}            

\newtheoremstyle{claimstyle}
  {.8\baselineskip}         
  {0\baselineskip}         
  {\itshape}    
  {0pt}         
  {\bfseries}            
  {}            
  {.5em}        
  {}            

\theoremstyle{thmstyle}
\newtheorem{theorem}{Theorem}[section]
\newtheorem{proposition}[theorem]{Proposition}
\newtheorem{lemma}[theorem]{Lemma}
\newtheorem{corollary}[theorem]{Corollary}

\newtheorem{observation}[theorem]{Observation}

\newtheorem*{proposition*}{Proposition}
\newtheorem*{theorem*}{Theorem}
\theoremstyle{claimstyle}

\newcommand{\mcp}{\mathcal P}
\newcommand{\mcq}{\mathcal Q}
\newcommand{\mcr}{\mathcal R}
\newcommand{\mcv}{\mathcal V}
\newcommand{\mce}{\mathcal E}
\newcommand{\mcc}{\mathcal C}
\newcommand{\mcf}{\mathcal F}
\newcommand{\antipode}[0]{\genfrac{}{}{}{3}{\;\, a \;\,}{}}

\begin{document}

\title{The Structure of Critical Product Sets}

\author{Matt DeVos\thanks{Supported in part by an NSERC Discovery Grant (Canada) and a Sloan Fellowship.}\\
	mdevos@sfu.ca}

\maketitle

\begin{abstract}
Let $G$ be a multiplicative group, let $A,B \subseteq G$ be finite and nonempty, and define the product set $AB = \{ ab \mid \mbox{$a \in A$ and $b \in B$} \}$.  Two fundamental problems in combinatorial number theory are to find lower bounds on $|AB|$, and then to determine structural properties of $A$ and $B$ under the assumption that $|AB|$ is small.  We focus on the extreme case when $|AB| < |A| + |B|$, and call any such pair $(A,B)$ \emph{critical}.

In the case when $|G|$ is prime, the Cauchy-Davenport Theorem asserts that $|AB| \ge \min \{ |G| , |A| + |B| - 1 \}$, and Vosper refined this result by classifying all critical pairs in these groups.  For abelian groups, Kneser proved a natural generalization of Cauchy-Davenport by showing that there exists $H \le G$ so that $|AB| \ge |A| + |B| - |H|$ and $ABH = AB$.  Kemperman then proved a result which characterizes the structure of all critical pairs in abelian groups.

Our main result gives a classification of all critical pairs in an arbitrary group $G$.  As a consequence of this we derive the following generalization of Kneser's Theorem to arbitrary groups:  There exists $H \le G$ so that $|AB| \ge |A| + |B| - |H|$ and so that for every $y \in AB$ there exists $x \in G$ so that $y( x^{-1} H x) \subseteq AB$.

$\mbox{}$
\end{abstract}

\section*{Outline}

Our proof involves an early departure from the problem of classifying critical product sets in a group $G$ to an alternate problem of classifying certain specially structured incidence geometries on which $G$ has a natural action.  So, we will first derive a structure theorem for these incidence geometries, and then deduce our structure theorem for product sets from this.  To help illuminate the structure of this paper, we begin with a diagram showing dependencies among our sections and a table offering a brief description of the content of each section.  This is followed by a figure showing the relationships among some of our lemmas, theorems, and corollaries.

\begin{center}
\begin{figure}[htsp]
\centerline{\includegraphics[width=12.5cm]{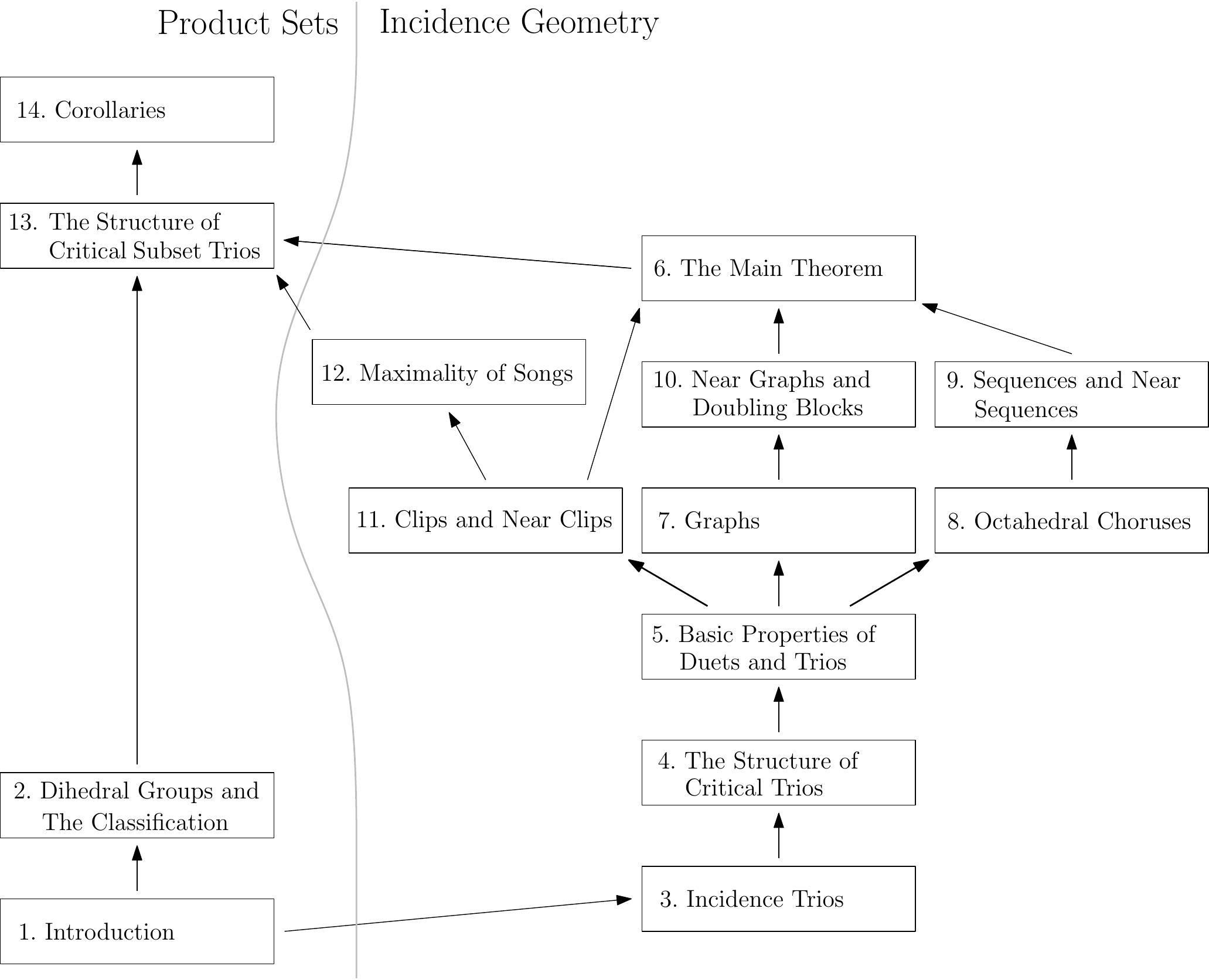}}
\label{section-dep}
\end {figure}
\end{center}

\begin{tabular}{c|p{4.8in}}
Section	&	Content	\\
\hline
1		&	A discussion of lower bounds on sizes of product sets, structure theorems for critical product sets, and a rough sketch of our proof. \\
2		&	Statement of our classification theorem for critical product sets. \\
3		&	An introduction to the incidence geometries we classify, and an example to motivate this approach. \\
4		&	Statement of our main classification theorem (incidence geometry).	\\
5		&	Numerous basic properties of our incidence geometries are established. \\
6		&	Proof of the main theorem (hard direction) assuming numerous lemmas. \\
7-11		&	Proofs of the lemmas required for the main theorem. \\
12		&	The easy direction of the main theorem is proved.\\
13		&	A derivation of our classification theorem for critical product sets from our classification theorem on incidence geometries \\
14		&	Some corollaries of the main result are established.
\end{tabular}

\begin{center}
\begin{figure}[ht]
\centerline{\includegraphics[width=13.2cm]{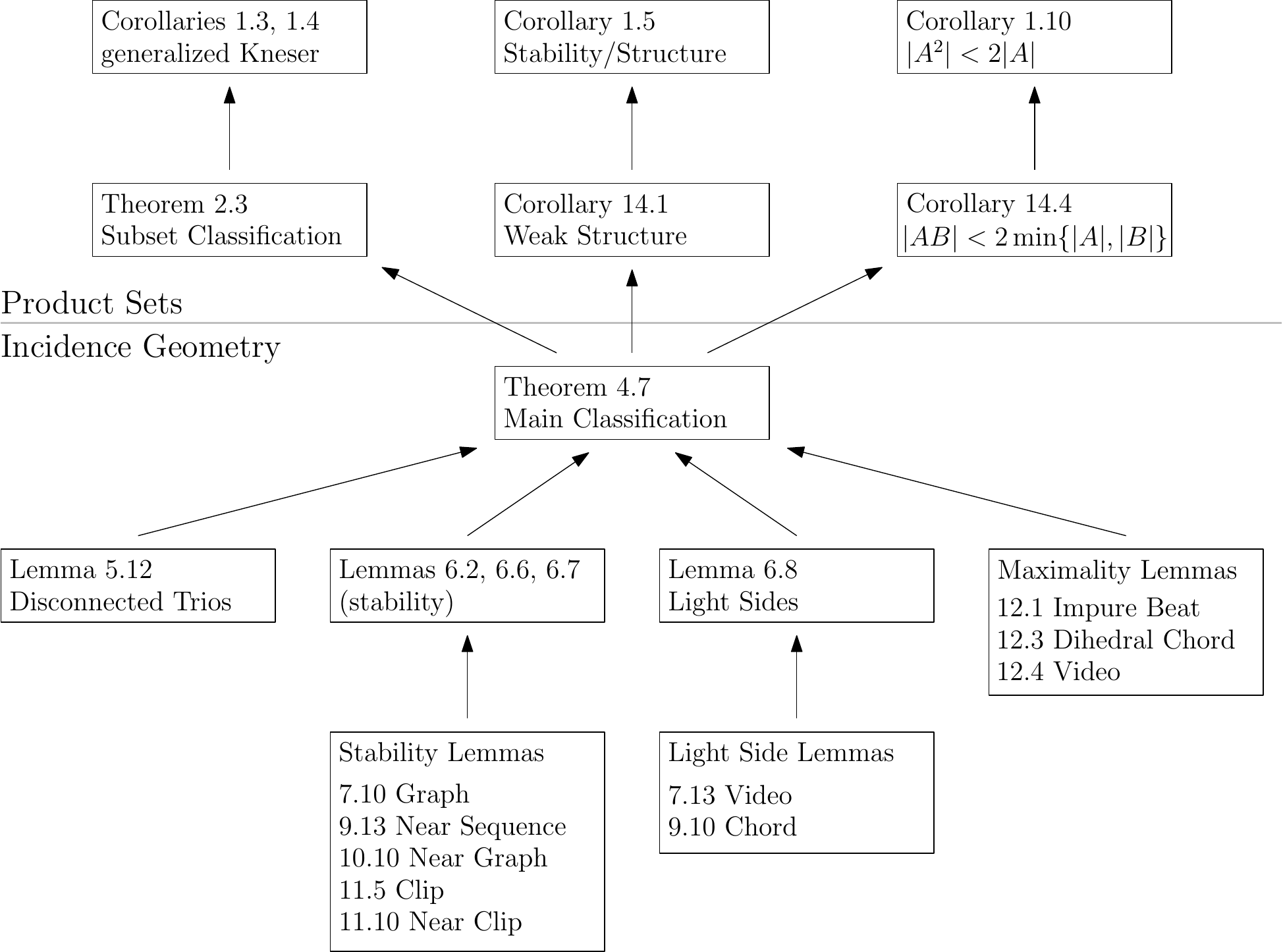}}
\label{imp_diag}
\end {figure}
\end{center}

\section{Introduction}
\label{intro_sec}

Except where noted, all groups will be written multiplicatively, and let us now fix such a group $G$.  
For a set $A \subseteq G$ we let $\overline{A} = G \setminus A$ and we define $A^{-1} =  \{ g^{-1} \mid g \in A \}$.  
If $B \subseteq G$ then we define the \emph{product set} of $A$ and $B$ to be $AB = \{ ab \mid \mbox{$a \in A$ and $b \in B$} \}$.  For an element $g \in G$ we let $gA = \{g\}A$ and let $Ag = A\{g\}$, and we call a set of either form a \emph{shift} of $A$.

One of the central goals in combinatorial number theory is to understand the structure of finite sets $A,B \subseteq G$ 
for which the product set $AB$ is small.  This program began in earnest with the pioneering work of Freiman \cite{Freiman-Thm}, and has been further developed through the contributions of Ruzsa \cite{ruzsa} and numerous others.  Quite recently, a breakthrough by Breuillard, Green, and Tao \cite{BGT} has established a powerful rough structure theorem for arbitrary groups.  

Our interest here is in the exceptional case when a pair $(A,B)$ of finite subsets of $G$ satisfy $|AB| < |A| + |B|$, and in this case we define $(A,B)$ to be \emph{critical}.  Although there are numerous constructions of critical pairs, this is a suitably restrictive notion to hope for a precise structure theorem.  Vosper \cite{vosper1} \cite{vosper2} proved such a characterization for groups of prime order, Kemperman \cite{kemperman2} extended this to abelian groups, and the goal of this paper is to prove a characterization which holds for an arbitrary group.  We will also apply this structure theorem to derive some new lower bounds on the sizes of product sets.

The remainder of this section is organized as follows.  In the first subsection we discuss old and new lower bounds on sizes of product sets, and in the second we give a brief overview of some existing structure theorems concerning critical product sets in nonabelian groups.  In the third and fourth subsections, we will transform the problem of classifying critical product sets into a new problem concerning the classification of certain triples of subsets we call maximal critical subset trios.  The fifth subsection contains a description of Vosper's Theorem and Kemperman's Theorem, and in the sixth and last subsection we give a rough sketch of our proof.  (The statement of our structure theorem for critical product sets is postponed to the following section.)

\subsection{Lower Bounds}

The starting point for the study of small product sets is the following famous result, which was first proved by Cauchy
and then rediscovered by Davenport.  For every positive integer $n$ we let ${\mathbf C}_n$ denote a cyclic group of order $n$

\begin{theorem}[Cauchy \cite{cauchy} - Davenport \cite{davenport}]
If $p$ is prime and $A,B \subseteq {\mathbf C_p}$ are nonempty, then
\[ |AB| \ge \min\{ p, |A| + |B| - 1 \}. \]
\end{theorem}

Kneser generalized this result to abelian groups as follows.

\begin{theorem}[Kneser \cite{kneser}]
If $G$ is abelian and $A,B \subseteq G$ are finite and nonempty, then there exists $H \le G$ so that 
\begin{enumerate}
\item $|AB| \ge |A| + |B| - |H|$
\item $ABH = AB$.
\end{enumerate}
\end{theorem}

The above theorem does not hold for general nonabelian groups even if we replace property 2 above with the weaker statement that $H$ must satisfy one of the conditions: $ABH = AB$ or $AHB = AB$ or $HAB = AB$.  This was demonstrated by Olson in \cite{olson} and will be discussed later in this section.  The following corollary of our main theorem generalizes Kneser's Theorem by weakening the second property in a different manner.

\begin{corollary}
\label{cor_basic}
If $A,B \subseteq G$ are finite and nonempty, there exists $H \le G$ so that 
\begin{enumerate}
\item $|AB| \ge |A| + |B| - |H|$
\item For every $y \in AB$ there exists $x \in G$ so that $ y(x^{-1}Hx) \subseteq AB$.
\end{enumerate}
\end{corollary}

There is an equivalent form of Kneser's Theorem where the first condition is replaced by the following stronger inequality:
\[ |AB| \ge |AH| + |BH|  - |H|. \]
A similar statement also holds for general groups.  For a set $A \subseteq G$ and a subgroup $H \le G$, we say that $A$ is \emph{left} (\emph{right}) $H$-\emph{stable} if $HA = A$ ($AH = A$), and we say that $A$ is $H$-\emph{stable} if it is both left and right $H$-stable.  We say that $A$ is $H$-\emph{conj-stable} if for every $y \in A$ there exists $x \in G$ so that $y(x^{-1}Hx) \subseteq A$ (note that $x^{-1}Hxy = y ( (xy)^{-1} H xy )$ so there is no distinction here between $y$ multiplied on the left or right by a conjugate of $H$).  

\begin{corollary}
\label{big_cor}
If $A,B \subseteq G$ are finite and nonempty, there exists $H \le G$ and $A^*, B^* \subseteq G$ with $A \subseteq A^*$ and $B \subseteq B^*$ so that
\begin{enumerate}
\item $A^* B^* = AB$
\item $|AB| \ge |A^*| + |B^*| - |H|$
\item $A^*$, $B^*$, $AB$, and $\overline{(AB)^{-1}}$ are $H$-conj-stable.
\end{enumerate}
\end{corollary}

Next we shall introduce an important parameter.  Define the \emph{deficiency} of a pair of finite subsets $(A,B)$ of $G$ to be 
$\delta(A,B) = |A| + |B| - |AB|$.  So, a pair has positive deficiency if and only if it is critical.  Also, note that all of the aforementioned theorems may be viewed
as upper bounds on the deficiency of a pair $(A,B)$.

Although it is not possible to replace the condition of conjugate stability in the previous two corollaries with the more desirable notions of left or right stability, by looking ``one step'' into our structure theorem, we may obtain the following result which either gives some stronger stability properties, or gives us 
some structural information about the sets $A$ and $B$.  Let us note that this result has a stronger form, Corollary \ref{weak_structure} (a weak 
version of the full structure theorem) which appears in the final section of the paper.

\begin{corollary}[\emph{Structure/Stability}]
\label{struc_or_stable}
If $A,B \subseteq G$ are finite, then one of the following holds.
\begin{enumerate}
\item One of $A$, $B$, or $\overline{AB}$ is contained in a coset of a proper subgroup.
\item There exists a finite subgroup $H \triangleleft G$ satisfying $\delta(A,B) \le \frac{1}{2}|H|$ and one of:
\begin{enumerate}
\item $G/H$ is cyclic,  $|G/H| \ge 4$, and $\max\{ |AH \setminus A| , |BH \setminus B| \} +\delta(A,B)  \le |H|$.
\item $G/H$ is dihedral, $|G/H| \ge 8$, and $\max\{ |AH \setminus A| , |BH \setminus B| \} +\delta(A,B)  \le 2|H|$.
\end{enumerate}
\item There exist $H_1, H_2, H_3 \le G$ so that  $H_1 A H_2 B H_3 = AB$ and $\delta(A,B) \le \min_{1 \le i \le 3} |H_i|$.
\end{enumerate}
\end{corollary}

There is already a wide body of existing results which provide upper bounds on the deficiency of a pair of sets in a nonabelian group under various assumptions.  Although we will not give an extensive overview of these, we will highlight a couple significant results of this type which relate closely 
with this article.  Below is  a theorem due to Kemperman which has a very similar structure to our Corollary \ref{cor_basic}.

\begin{theorem}[Kemperman \cite{kemperman1}]
\label{kemp_easy}
If $A,B \subseteq G$ are finite and nonempty and  $a \in A$ and $b \in B$,  there exists $H \le G$ so that
\begin{enumerate}
\item $\delta(A,B) \le |H|$
\item $aHb \subseteq AB$.
\end{enumerate}
\end{theorem}

The preceding theorem is also consequence of Corollary \ref{cor_basic}.  To see this, let $a \in A$ and $b \in B$ and then apply this corollary to choose $H \le G$ and $x \in G$ so that $|AB| \ge |A| + |B| - |H|$ and $(ab)(x^{-1}Hx) \subseteq AB$.  Now setting $K = bx^{-1} H xb^{-1}$ we have that $|AB| \ge |A| + |B| - |K|$ and $aKb = a ( bx^{-1} H xb^{-1})b = ab(x^{-1}Hx) \subseteq AB$.  Indeed, we may view Corollary \ref{cor_basic} as a uniform version of this theorem of Kemperman.   

We close this subsection with an important deficiency bound which will play a key early role in our proof.  In Subsection \ref{tight_subs} we will call
upon this result to reduce the overall classification problem.  

\begin{theorem}[Kemperman, Scherk\footnote{We have followed Lev \cite{lev1} in attributing this result.} \cite{kemperman1}]  
\label{rep_cor}
If $A,B \subseteq G$ are finite and nonempty, then
\[ \delta(A,B) \le \min_{z \in AB} | \{ (a,b) \in A \times B \mid ab=z \}|. \]
\end{theorem}

\subsection{Structure}

In this subsection we will briefly discuss some existing partial results on the structure of critical product sets in arbitrary groups, and we will mention one more corollary of our structure theorem.  Our focus here will be on nonabelian groups, since Kemperman's structure theorem (which gives a complete description in the abelian case) will be treated in detail later in this section.  Again, we have made no effort to give a comprehensive survey, and will present only a handful of results which connect closely with this work.

Fix a finite nonempty set $B \subset G$ and note that for every finite nonempty set $A \subseteq G$ we must have $\delta(A,B) = |A| + |B| - |AB| \le |B|$.  In light of this we may define
\[ \delta(B) = \max \{ \delta(A,B) \mid \mbox{ $A \subseteq G$ is finite and $A, \overline{AB} \neq \emptyset$} \}. \]
Below we state a theorem of Hamidoune concerning $\delta(B)$ which extends some related results of Mann \cite{mann} and Z\'{e}mor \cite{zemor}.  A 
variant of this result will play an important role in our proof.

\begin{theorem}[Hamidoune \cite{hamidoune-cayley-con}]
\label{hamidoune-sumset}
If $B \subset G$ is finite and nonempty, there exists a finite set $A \subseteq G$ so that $A, \overline{AB} \neq \emptyset$, and 
$\delta(A,B) = \delta(B)$, and either $A$ or $\overline{AB}$ is a finite subgroup.
\end{theorem}

This theorem was one of many which Hamidoune proved using the isoperimetric method (see also \cite{hamidoune3}, \cite{hamidoune-connectivity}, \cite{hamidoune-structure}, \cite{hamidoune-ks_ext}), and in fact he proved a generalization of this result which holds in the setting of vertex transitive digraphs, and is quite closely related to the precise result we will utilize.  Indeed, our incidence geometries are close relatives of Hamidoune's vertex transitive graphs, and our basic technique is a variant of the isoperimetric method.  As such, this work may be viewed as a continuation of Hamidoune's central program.  

In one of his last papers, Hamidoune \cite{hamidoune-minkowski} investigated sets $B$ for which $\delta(B) = 1$ and there exists $A \subseteq G$ with $|A|, |\overline{AB}| \ge 2$ for which $\delta(A,B) = 1$, and this work was then sharpened by Serra and Z\'{e}mor \cite{serra-zemor}.  These authors obtain a structure theorem for such sets which is a special case of our main result.  Namely, they prove that either there exists a set $A$ as in the statement of Hamidoune's Theorem for which either $A$ or $\overline{AB}$ is a nontrivial subgroup of $G$, or the set $B$ is a geometric progresion, or there exists a set $A'$ and a finite subgroup $H \le G$ so that $\delta(A',B) = 1$ and both $A'$ and $\overline{A'B}$ are unions of two right $H$-cosets.  

Next we shall introduce another recent theorem giving some structural information about critical product sets.  The essential ingredients for the proof of this theorem were provided by Hamidoune \cite{hamidoune-two_inverse} in response to a question of Tao \cite{tao-bloq_question}, but the result was proved in this form by Tao \cite{tao-small_doubling}, so we have credited both with the theorem.  This gives a generalization to nonabelian groups of a consequence of 
Kneser's Theorem.

\begin{theorem}[Hamidoune \cite{hamidoune-two_inverse}, Tao \cite{tao-small_doubling}]
Let $A,B \subseteq G$ be finite and nonempty and assume $|A| \ge |B|$ and $|AB| \le (2 - \epsilon)|B|$ for some $\epsilon > 0$.  Then one of the 
following holds.
\begin{enumerate}
\item $B$ is contained in a coset of a finite subgroup $H$ with $|H| \le \frac{2}{\epsilon} |B|$.
\item $B$ can be covered by at most $\frac{2}{\epsilon}-1$ right $H$ cosets for some $H \le G$ with $|H| \le |B|$.
\end{enumerate}
\end{theorem}

The above result may also be deduced (in somewhat stronger form) from Corollary \ref{def_vs_disp} of our structure theorem.

A case of particular interest in small product set problems, is finite subsets $A \subseteq G$ for which the set $A^2 = AA$ has small size.  In the very extreme case when $|A^2| < 1.5 |A|$, Freiman proved that the following property holds (see \cite{freiman-new1.5})

\noindent{(1)} \; There exists $H \le G$ and $x \in G$ so that $A \subseteq xH = Hx$ and $A^2 = x^2H$.

Freiman \cite{freiman-1.6} also proved a structure theorem for sets $A$ with $|A^2| < 1.6|A|$, and a modern 
translation of this appears in Husbands' Masters Thesis \cite{husbands}.  Here it is proved that such sets must either satisfy (1) above 
or the following property.

\noindent{(2)} \; There exists $K \le G$ so that the set $B = KAK$ satisfies $|B| = 2|K|$ and $|B^2| = 3|K|$. 

In fact, whenever (2) holds, there exist subgroups $L \le K \le H \le G$ and $x \in G$ so that $B \subseteq xH = Hx$ and $[K:L] \le 2$, 
and furthermore $L$ is a normal subgroup of $H$ with the property that $H/L$ is either cyclic or dihedral.  So, in effect, this critical phenomena can be found either in a cyclic or a dihedral group (more precisely, setting $C = x^{-1}B$ and $D = Bx^{-1}$, we have that the images of $C$ and $D$ under the 
canonical homomorphism from $H$ to $H/L$ form a critical pair).  A proof of the existence of these subgroups $L,H$ appears in the appendix.  It relies only upon material appearing in Section \ref{inc_sec}.

Our main theorem specializes to sets $A$ with $|A^2| < 2|A|$ to give the following result.  
This corollary also has a stronger form, Corollary \ref{def_vs_disp}, which gives added structural information in the second and third 
outcomes (and also applies more generally to pairs $(A,B)$ for which $|AB| < \min\{ |A|, |B| \}$).

\begin{corollary}
\label{a_squared}
Let $A \subseteq G$ be finite with $|A| \ge 2$ and $|A^2| < 2|A|$.  If $x \in A$ and $H$ is the unique minimal subgroup for which 
$A \subseteq xH$, then $xH = Hx$ and there exists $K < H$ satisfying one of the following.
\begin{enumerate}
\item There exists $y \in H$ so that $x^2H \setminus x^2yK \subseteq A^2$.
\item $K \triangleleft H$, and $H/K$ is cyclic with order at least $4$, and 
$\delta(A,A) + |AK \setminus A| \le |K|$.
\item $K \triangleleft H$, and $H/K$ is dihedral with order at least $8$, and 
$\delta(A,A) + |AK \setminus A| \le 2|K|$.
\end{enumerate}
\end{corollary}

\subsection{Tight Critical Pairs}
\label{tight_subs}

In this subsection we will simplify the task of classifying all critical pairs by first identifying some trivial constructions, and then reducing the general problem to that of classifying certain special critical pairs we call tight.

Let $(A,B)$ be a critical pair in $G$.  If exactly one of $A$, $B$ is empty then $(A,B)$ will be critical (as $AB = \emptyset$) and we call such a pair \emph{trivial}.  Another case of little interest is when $G$ is finite and $|A| + |B| > |G|$.  Here, it follows immediately that $(A,B)$ will be critical.  However, in this case we must also have $AB = G$ (to check this, note that for every $g \in G$ we have $B \cap A^{-1}g \neq \emptyset$ so there will exist $x \in A$ for which $x^{-1}g \in B$).  Therefore, the critical pairs $(A,B)$ with $AB = G$ are precisely those $A,B$ for which $|A| + |B| > |G|$ and accordingly, we shall also call these \emph{trivial}.  Having defined these, we now turn our attention to nontrivial critical pairs.  

We define $(A,B)$ to be \emph{tight} if the only sets $A^* \supseteq A$ and $B^* \supseteq B$ for which $A^* B^* = AB$ are given by $A^* = A$ and $B^* = B$.  Fortunately, the problem of classifying all critical pairs reduces to the problem of classifying tight critical pairs, as shown by the following result.

\begin{theorem}
The pair $(A,B)$ is critical if and only if there exists a tight critical pair $(A^*,B^*)$ for which 
$A \subseteq A^*$ and $B \subseteq B^*$ and the following equation holds:
\[ |A^* \setminus A| + | B^* \setminus B| < \delta(A^*, B^*). \]
\end{theorem}

\noindent{\it Proof:} First suppose that $(A,B)$ is a critical pair and choose $A^* \supseteq A$ and $B^* \supseteq B$ maximal so that $A^*B^* = AB$ (it follows immediately from the fact that $AB$ is finite that such maximal sets exist).  This gives us 
\begin{align*} 
\delta( A^*, B^*) 
	&= |A^* \setminus A| + |B^* \setminus B| + \delta(A,B)	\\
	&> |A^* \setminus A| + |B^* \setminus B|.
\end{align*}

Next let $(A^*,B^*)$ be a tight critical pair and let $A \subseteq A^*$ and $B \subseteq B^*$ satisfy $ |A^* \setminus A| + |B^* \setminus B| < \delta(A^*, B^*)$.  Theorem \ref{rep_cor} implies that every element in $A^*B^*$ has at least $\delta(A^*,B^*)$ representations as a product of an element in $A^*$ times an element in $B^*$.  It follows from this and $|A^* \setminus A| + |B^* \setminus B| < \delta(A^*,B^*)$ that $AB = A^* B^*$.  Therefore, the pair $(A,B)$ satisfies $\delta(A,B) = |A| + |B| - |AB| > |A^*| + |B^*| - \delta(A^*,B^*) - |A^*B^*| = 0$ so it is critical. 
\quad\quad$\Box$

\subsection{Subset Trios}

In the study of critical pairs $(A,B)$, there is a third set of significant interest which deserves to considered together with $A$ and $B$, namely $C = \overline{(AB)^{-1}}$.  
For simplicity, first suppose that $G$ is finite, and that $(A,B)$ is a nontrivial critical pair, and observe the following two properties:
\begin{itemize}
\item $|A| + |B| + |C| > |G|$,
\item $1 \not\in ABC$.
\end{itemize}
Given the above, it follows that the pair $(B,C)$ is also critical (to see this, note that  $BC$ must be disjoint from $A^{-1}$, so it will have size at most $|G| - |A| < |B| + |C|$).  Similarly, we have $1 \not\in CAB$, 
so  $(C,A)$ is also a critical pair.  In light of this, we shall now extend our definitions to triples.  To allow for the possibility that $G$ is infinite, we will permit sets in our triples to be infinite if they are cofinite.

\bigskip

\noindent{\bf Definition:} If $A,B,C \subseteq G$ are all finite or cofinite and $1 \not\in ABC$ then we define $(A,B,C)$ to be a \emph{subset trio} (\emph{in} $G$).

\bigskip

Next we extend some of our earlier terminology from pairs to subset trios.  We call the subset trio $(A,B,C)$ \emph{trivial} if one of $A$, $B$, or $C$ is empty.  Note that for a nontrivial subset trio, at most one of the three sets may be infinite.  Set $n$ to be the size of the complement of a largest set in $(A,B,C)$ and let $\ell, m$ be the sizes of the other two sets in $(A,B,C)$.  Then we define the \emph{deficiency} of $(A,B,C)$ to be $\delta_G(A,B,C) = \ell + m - n$ (we permit $\delta_G(A,B,C) = \infty$ if two of $A,B,C$ are infinite and $\delta_G(A,B,C) = -\infty$ if $A,B,C$ are finite but $G$ is infinite), and we say that $(A,B,C)$ is \emph{critical} if $\delta(A,B,C) > 0$.  As usual, we drop the subscript when the group is clear from context.  
A subset trio $(A,B,C)$ is \emph{maximal} if the only subset trio $(A^*,B^*,C^*)$ with $A \subseteq A^*$ and $B \subseteq B^*$ and $C \subseteq C^*$ is given by $A = A^*$ and $B = B^*$ and $C = C^*$.  Equivalently, $(A,B,C)$ is maximal if $C = \overline{(AB)^{-1}}$ and $A = \overline{(BC)^{-1}}$ and $B = \overline{(CA)^{-1}}$.  The natural motivation for these definitions is the following.

\bigskip

If $A,B$ are finite subsets of $G$ and $C = \overline{(AB)^{-1}}$ then:
\begin{itemize}
\item $(A,B)$ is trivial if and only if $(A,B,C)$ is trivial,
\item $(A,B)$ is critical if and only if $(A,B,C)$ is critical,
\item $\delta(A,B) = \delta(A,B,C)$,
\item $(A,B)$ is tight if and only if $(A,B,C)$ is maximal.
\end{itemize}

It is immediate that whenever $(A,B,C)$ is a subset trio, the triples $(B,C,A)$, $(C,A,B)$, $(C^{-1}, B^{-1}, A^{-1})$,  $(B^{-1}, A^{-1},C^{-1})$, and  $( A^{-1}, C^{-1}, B^{-1})$ will also be subset trios.  Similarly, for every $g \in G$ the triples $(Ag,g^{-1}B,C)$, $(A,Bg,g^{-1}C)$, and $(g^{-1}A,B,Cg)$ will be subset trios.  All of these transformations preserve the properties of interest to us (nontrivial, critical, maximal, deficiency).  Accordingly, we will say that any subset trio obtained from $(A,B,C)$ by a sequence of these transformations is \emph{similar} to $(A,B,C)$.  

At this point, we have transformed the original problem of classifying critical pairs to the new problem of classifying nontrivial maximal critical subset trios.

\subsection{Classical Critical Subset Trios}

In this subsection we will describe some classical constructions of critical subset trios, and state the theorems of Vosper and Kemperman.  Although Vosper's Theorem has a fairly simple form, Kemperman's Theorem is somewhat complicated to state.  Indeed, the difficulty of his theorem and his paper lead at one point to some confusion about the precise nature of his results.  Fortunately, this situation has been remedied by recent work of Grynkiewicz \cite{grynkiewicz-qpk} \cite{grynkiewicz-step}, Hamidoune \cite{hamidoune-kemperman} \cite{hamidoune-structure},  and Lev \cite{lev-kemp}.  Notably, Lev \cite{lev-kemp} gives a top-down version of Kemperman's Theorem which is considerably easier to work with.  We shall adopt this framework throughout.

We begin with two easy constructions for critical pairs.  The first construction is to take one of the sets to be a singleton.  For instance, if $A = \{a\}$ then $|AB| = |B|$ and so $(A,B)$ will be critical.  For our second construction, we require some additional notation.  We define any set of the form $A = \{ g^i \mid 0 \le i \le \ell\}$ for $g \in G$ to be a \emph{basic geometric progression with ratio} $g$.  We say that this progression is \emph{nontrivial} if it has size at least two and \emph{proper} if $g^{\ell + 1} \not\in A$.  More generally, we define a \emph{geometric progression}\footnote{Since we will only consider geometric progressions in abelian groups, we have not concerned ourselves with right vs. left multiplication here.} \emph{with ratio} $g$ to be any set of the form $hA$ where $h \in G$ and $A$ is a basic geometric progression with ratio $g$, and we say that this geometric progression is \emph{nontrivial} (\emph{proper}) if $A$ has this property.  
Now suppose that $A$ and $B$ are proper basic geometric progressions with ratio $g$, say $A = \{ g^i \mid 0 \le i \le \ell \}$ and $B = \{g^i \mid 0 \le i \le m \}$.  In this case $|A| = \ell + 1$ and $|B| = m + 1$ and $AB = \{ g^i \mid 0 \le i \le \ell + m \}$ satisfies $|AB| \le \ell + m + 1$, so $(A,B)$ is critical.  We may extend this pair to a critical subset trio by defining 
$C = \overline{(AB)^{-1}}$.  Up to our equivalence relation, this gives all of the possibilities for nontrivial critical trios in groups of prime order.

\begin{theorem}[Vosper \cite{vosper1}, \cite{vosper2}]
If $(A,B,C)$ is a nontrivial critical subset trio in ${\mathbf C}_p$ for $p$ prime, then one of the following holds.
\begin{enumerate}
\item $\min\{ |A|, |B|, |C| \} = 1$
\item $A,B,C$ are geometric progressions with a common ratio.
\end{enumerate}
\end{theorem}

For groups with nontrivial subgroups, there are natural generalizations of these critical subset trios which we shall now define.  For a set $A \subseteq G$, the \emph{left stabilizer} of $A$ is ${\mathit stab}_L(A) = \{ g \in G \mid gA = A \}$ and the
\emph{right stabilizer} of $A$ is ${\mathit stab}_R(A) = \{ g \in G \mid Ag = A \}$.  Note that ${\mathit stab}_L(A)$ and ${\mathit stab}_R(A)$ are both subgroups of $G$ and that $A$ is left (right) $H$-stable if and only if $H \le {\mathit stab}_L(A)$ ($H \le {\mathit stab}_R(A)$).

\bigskip

\noindent{\bf Definition:} Let $H < G$ be finite.  A subset trio is a \emph{pure beat} \emph{relative to} $H$  if it is similar to a subset trio $(A,B,C)$ satisfying:
\begin{enumerate}
\item $A = H$
\item ${\mathit stab}_L(B) = H$
\item $C = \overline{(AB)^{-1}}$.  
\end{enumerate}

Observe that every pure beat relative to $H$ is a maximal nontrivial critical trio (with deficiency $|H|$).  

If $H \triangleleft G$, we let $\varphi_{G/H} : G \rightarrow G / H$ denote the canonical homomorphism.  

\bigskip

\noindent{\bf Definition:} 
Let $H \triangleleft \, G$ be finite with $G/H$ is cyclic and $|G/H| \ge 5$.  A subset trio is a \emph{pure cyclic chord} \emph{relative to}
$H$ if there exists a similar subset trio $(A,B,C)$ and $S \in G/H$ which generates $G/H$ satisfying:
\begin{enumerate}
\item $A,B,C$ are $H$-stable.
\item $\varphi_{G/H} (A)$ and $\varphi_{G/H}(B)$ are nontrivial basic geometric progressions with ratio $S$.
\item $C = \overline{(AB)^{-1}}$ is not contained in a single $H$-coset.
\end{enumerate}
Note that every pure cyclic chord relative to $H$ is a maximal critical subset trio with deficiency $|H|$.

\bigskip

Next we shall introduce two structures which are similar to pure beats and pure chords, but allow for recursive constructions of critical pairs and critical subset trios.  

\bigskip

\begin{figure}[ht]
\centerline{\includegraphics[width=5cm]{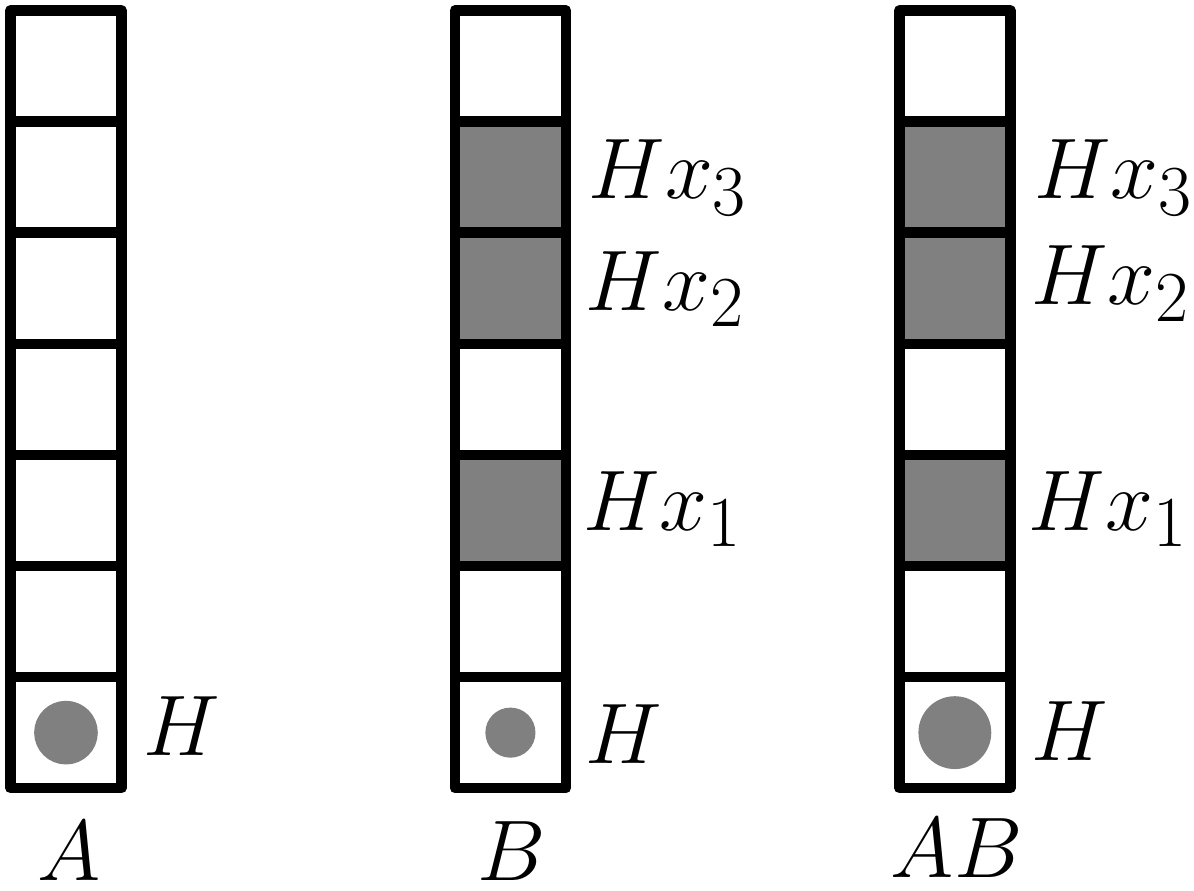}}
\caption{An impure beat}
\label{ssbeat}
\end {figure}

\noindent{\bf Definition:} 
Let $H < G$.  A subset trio $\Phi$ is an \emph{impure beat} \emph{relative to} $H$  if there exists an similar subset trio $(A,B,C)$ which satisfies the following properties.
\begin{enumerate}
\item $A \subseteq H$
\item $H \le {\mathit stab}_L(B \setminus H)$
\item $C \setminus H = \overline{(AB)^{-1}} \setminus H$
\item $B \cap H$ and $C \cap H$ are nonempty.
\end{enumerate}
In this case we define the (nontrivial) subset trio $\Phi' = (A, B \cap H, C \cap H)$ in the group $H$, to be a \emph{continuation} of $\Phi$.  Observe that 
$\delta_H(\Phi') = \delta_G(\Phi)$ and that $\Phi'$ is maximal whenever $\Phi$ is maximal.

In our definition of pure beats, we insisted that the subgroup $H$ be finite, however this restriction was not included in the definition of 
impure beat.  This restriction is actually implied for pure beats, as otherwise the definitions would force all of $A$, $B$, and $C$ to be infinite.  
For impure beats, if $H$ is infinite, then either $B$ or $C$ is contained in a single $H$-coset.  This is a natural possibility which is straightforward to construct.  To do so, let $H$ be an infinite subgroup of $G$, let $(A,B)$ be a nontrivial critical pair from $H$, and set $C = \overline{(AB)^{-1}}$.

Let us pause here to make another comment concerning the above construction.  One might have 
hoped to generalize Kneser's Theorem to nonabelian groups by showing that for every pair 
$A,B$ of finite nonempty subsets of a group there exist subgroups $H_1,H_2,H_3 \le G$ for which 
$H_1 A H_2 B H_3 = AB$ and $|AB| \ge |A| + |B| - \max \{ |H_1|, |H_2|, |H_3| \}$.  However, this 
statement is false.  In \cite{olson}, Olson constructs sets $A,B$ for which $|AB| < |A| + |B| - 1$ and the 
only possible subgroups $H_1,H_2,H_3$ with $H_1 A H_2 B H_3 = AB$ are given by
$H_1 = H_2 = H_3 = \{1 \}$.  In our language, his construction uses impure beats, and it is not difficult to indicate 
the key principle involved in this construction.  First note that in the language of trios,  
we wish to construct a maximal subset trio $(A,B,C)$ for which $\delta(A,B,C) \ge 2$ and each of $A$, $B$, 
and $C$ has trivial right and left stabilizers.  Suppose now that $(A,B,C)$ is an impure beat relative to $H$ 
where $A \subseteq H$ and $H \le {\mathit stab}_L(B \setminus H)$ and that 
$(A', B', C')$ is a continuation as given in the definition of impure beat.  In this case, it is quite possible that there exists a nontrivial subgroup $H_3 < H$ so that 
$B'H_3 = B'$ and $H_3 C' = C'$ but $BH_3 \neq B$ and $H_3 C \neq C$.  So, in other words, the outer structure of the 
impure beat $(A,B,C)$ may spoil the stabilizer behaviour found in the continuation.  Olson's construction 
uses a nesting of impure beats to spoil all of the stabilizers.

\begin{figure}[ht]
\centerline{\includegraphics[width=5cm]{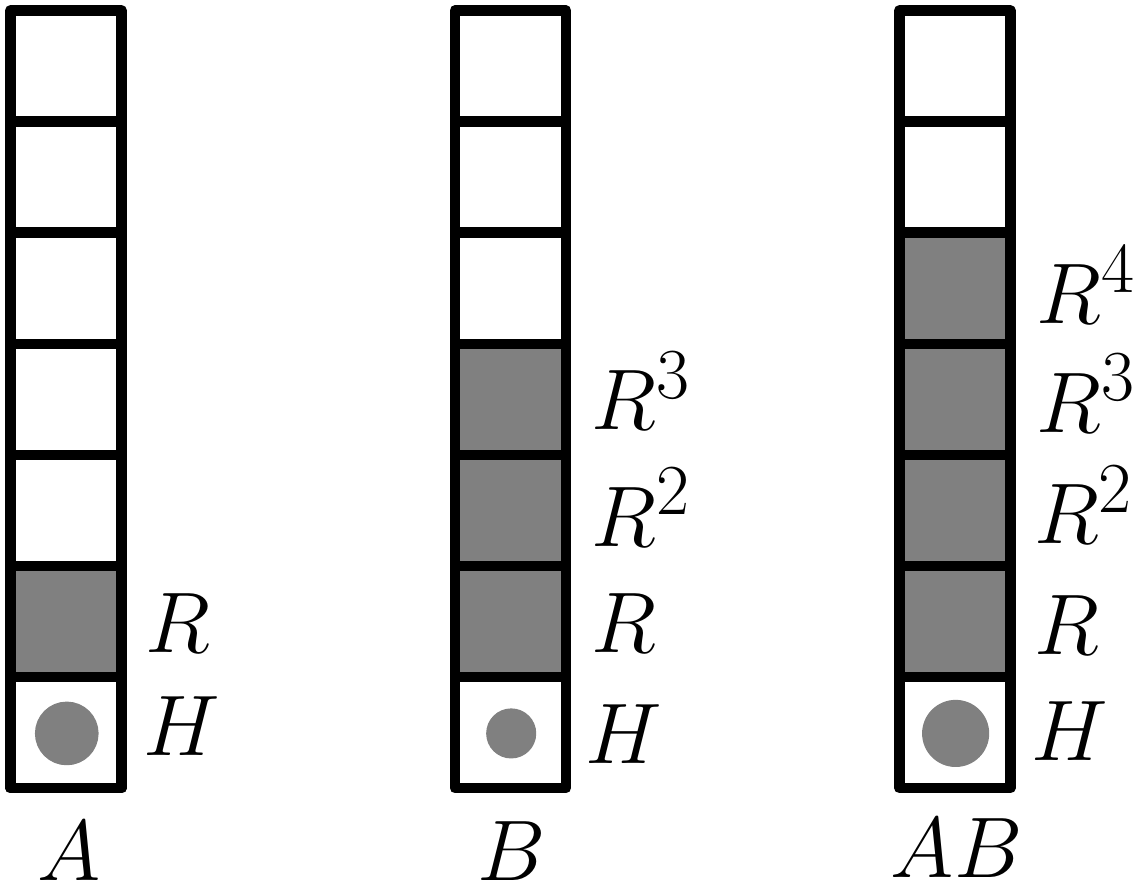}}
\caption{An impure cyclic chord}
\label{sschord}
\end {figure}

\bigskip

\noindent{\bf Definition:}
Let $H \triangleleft G$ be finite with $G/H$ is cyclic and $|G/H| \ge 4$.  A subset trio $\Phi$ is an \emph{impure cyclic chord} 
\emph{relative to} $H$ if  there exists a similar trio $(A,B,C)$ and $S \in G/H$ which generates $G/H$ satisfying
\begin{enumerate}
\item $A \cup H$ and $B \cup H$ are $H$-stable.
\item $\varphi_{G/H}(A)$ and $\varphi_{G/H}(B)$ are  nontrivial basic geometric progressions with ratio $S$.
\item $C \setminus H = \overline{(AB)^{-1}} \setminus H \neq \emptyset$ and $C \cap H \neq \emptyset$.
\end{enumerate}
We call the (nontrivial) subset trio $\Phi' = (A \cap H, B \cap H, C \cap H)$ a \emph{continuation} of $\Phi$.   Observe that 
$\delta_H(\Phi') = \delta_G(\Phi)$, and that $\Phi'$ is maximal whenever $\Phi$ is maximal.  With this terminology in place, we are ready to state Kemperman's Theorem.

\begin{theorem}[Kemperman \cite{kemperman2}]
If $\Phi$ is a nontrivial maximal critical subset trio in the abelian group $G$, there exists a sequence of subgroups 
$G = G_1 > G_2 \ldots > G_m$ and a sequence of subset trios $\Phi = \Phi_1, \Phi_2, \ldots, \Phi_m$ so that the following hold
\begin{enumerate}
\item $\Phi_i$ is a subset trio in the group $G_i$ for $1 \le i \le m$.
\item $\Phi_i$ is an impure beat or impure cyclic chord with continuation $\Phi_{i+1}$ for $1 \le i < m$.
\item $\Phi_m$ is either a pure beat or a pure cyclic chord.
\end{enumerate}
\end{theorem}

Our methods do specialize to give an alternative proof of Kemperman's Theorem.  However, there is little need to pass to the setting of incidence geometry when working with abelian groups, and certain aspects of our argument can be simplified significantly in this setting.  In a paper together with Boothby and Montejano \cite{bdm}, we have carried through these simplifications to obtain a new proof of Kemperman's Theorem.

\subsection{Proof Sketch}

In order to prove our result, we will move from the world of product sets and subsets of groups to the world of incidence geometries 
equipped with group actions.  Following the next section (where our characterization theorem for critical subset trios is stated), we will turn
our attention to these incidence geometries, and this is the setting in which all of the heavy lifting of the proof will be done.  
Although this departure is well precedented, it has the side effect that the partial results we prove have a rather different appearance
from those found in earlier works.  Despite this, all of the essential ingredients in our proof are very standard combinatorial intersection-union 
and induction arguments which have been well established by numerous authors.  To help account for some of these important contributions
we will give a brief overview of these ideas next.

The practice of exchanging a pair of sets $A,B$ from a common ground set by $A \cap B, A \cup B$ we will casually refer to as uncrossing, 
and it has a surprisingly rich history in the world of combinatorics.  From additive number theory, there are numerous transforms which are 
based on uncrossing together with shifting sets.  These transforms are key ingredients in the proofs of the Cauchy-Davenport Theorem and 
Kneser's Theorem (see \cite{nathanson}), the results from the first section due to Kemperman and Scherk, and numerous other theorems.  In 
the world of graph theory, 
Mader \cite{mader-edge} \cite{mader-vertex} and Watkins \cite{watkins} 
exploited uncrossing arguments to achieve a number of important theorems concerning connectivity properties of vertex transitive graphs, and
these tools have been further developed by many other authors.  In fact, connectivity of vertex transitive graphs is intimately linked with the study of small product sets.  Hamidoune was among the first to recognize the importance of this connection, and his work helped to establish a rich bridge between these worlds.  Our program of working in incidence geometries is a natural extension of Hamidoune's methodology.

Although our work will be done in a different setting, it is possible to describe the analogue of our procedure in the framework 
of product sets, so we shall give a rough sketch of this here.  The basic property upon which all of our results are based is the 
following uncrossing property:  If $(A,B,C)$ and $(A',B',C)$ are subset trios, then $(A \cap A', B \cup B', C)$ and $(A \cup A', B \cap B', C)$ are 
also subset trios and 
\begin{equation}
\label{subset_supermodular}
\delta(A \cap A', B \cup B', C) + \delta(A \cup A', B \cap B', C) = \delta(A,B,C) + \delta(A',B',C).
\end{equation}
We will take advantage of this in numerous ways, but let us highlight one here.  Suppose that $(A,B,C)$ is critical, let $g \in G$, 
and set $A' = Ag$ and $B' = g^{-1}B$.  Then $(A',B',C)$ is another subset trio, and it follows from equation \ref{subset_supermodular} that 
either one of $(A' \cap A, B' \cup B, C)$ or $(A' \cup A, B' \cap B, C)$ has greater deficiency than $(A,B,C)$ or they both have the same 
deficiency as $(A,B,C)$.  The first usage of this precise transformation we know of appears in a paper of Kemperman \cite{kemperman1} but 
it is closely related to other transformations used earlier by Sherck, and Davenport and others.  
The traditional method of application is to set up an inductive framework in which one of these two new subset trios has improved with regard to the inductive parameter, and then to obtain the result by applying induction to this trio.  Of course, if one of these two new subset trios is trivial, then it has rather little structure, so in general, when using these transformations one restricts attention to only those which produce nontrivial  subset trios.  

For our purposes, we will need to maintain a stronger condition than nontriviality.  A cursory glance at the structures found in Kemperman's Theorem reveals the following essential property: For both pure and impure beats, one of the three sets is contained in a coset of a proper subgroup, but the other two have relatively little individual structure; in contrast 
to this, for both pure and impure chords, all three sets in the subset trio have strong structure individually.  Roughly speaking, our plan will be to begin with a critical trio $(A,B,C)$ for which none of the three sets is contained in a coset of a proper subgroup (the other case is dealt with separately), and then try to apply an uncrossing argument to obtain two new trios with the same property.  Assuming we have improved a suitable inductive parameter, we will be able to apply the theorem inductively to one of these new subset trios, say $(A^*,B^*,C)$.  Since $(A^*,B^*,C)$ will not be a pure or impure beat (and cannot be contained in one), we may apply our result to gain useful structural information about the set $C$.  With this in hand, we will then return to try to understand the structure of $(A,B,C)$.  

There are a number of complications in following this plan.  First off, we will need to use critical trios of the form $(A^*,B^*,C)$ not constructed in the aforementioned manner of Kemperman.  For this purpose, we will develop a technique called purification  based on Hamidoune's Theorem \ref{hamidoune-sumset} which will allow for more subtle modifications.  Another complication is that we are only proving a structure theorem for maximal subset trios.  So, after our uncrossing operation we will need to pass from the new trio we have found to a maximal one.  As a consequence of this, the structural description of $C$ we obtain will be somewhat rough.  To remedy this situation, we will need to prove a number of stability results which show that whenever a maximal critical subset trio $(A,B,C)$ has one set which is suitably close to one found in our structures, then $(A,B,C)$ is one of our familiar types.  These stability lemmas will constitute the lion share of the technical work.

\section{Dihedral Groups and the Classification}
\label{dihedral_sec}

The main goal of this section is to state our structure theorem for critical subset trios.  However, before doing so, we will need to introduce 
some additional types of critical behaviour.

\subsection{Pure Dihedral Chords}

In this subsection we will introduce a simple type of critical trio which occurs in dihedral groups.  This critical behaviour has been observed 
previously.  For instance, it appears (implicitly) in a paper of Eliahou and Kervaire \cite{ek} which investigates small product sets in dihedral groups.

For every positive integer $n$ we let ${\mathbf D}_n$ denote a dihedral group of order $2n$, and we let ${\mathbf D}_{\infty}$ denote a countably infinite dihedral group.  Let us now fix a dihedral group ${\mathbf D}_n$ (here we permit $n = \infty$) with $n \ge 3$.  Then ${\mathbf D}_n$ has a unique cyclic subgroup of index two, which we call the \emph{rotation subgroup}.  We call the elements of this subgroup \emph{rotations} and all other elements \emph{flips}.  Let $r$ be a rotation, let $f$ be a flip and let $k$ be a nonnegative integer.  We define any set of the following form to be a \emph{basic dihedral progression with ratio} $r$.

\begin{equation}
\label{defa}
A =  \{ 1, r, r^2, \ldots, r^k \} \{1,f\}.
\end{equation}

We define a \emph{dihedral progression with ratio} $r$ to be any set of the form $sA$ where $s \in {\mathbf D_n}$, but we will keep our focus on basic progressions.  We say that the basic dihedral progression $A$ is \emph{nontrivial} if $k \ge 1$ and \emph{proper} if $r^{k+1} \not\in A$.  Assuming $A$ is proper, it follows immediately that $|A| = 2(k+1)$.  Since the conjugate of a rotation by a flip yields the inverse rotation we may also express $A$ as follows.

\[ A = \{1, fr^{-k} \} \{1, r, r^2, \ldots, r^k \}. \]

If $A$ is proper, there does not exist a nontrivial rotation in ${\mathit stab}_L(A) \cup {\mathit stab}_R(A)$ (such an element would have to stabilize $\{1,r,\ldots,r^k\}$ which is not possible).  It then follows that ${\mathit stab}_L(A) = \{1, f r^{-k} \}$ and ${\mathit stab}_R(A) = \{1, f\}$.  Next suppose that $\ell$ is a positive integer and let $B$ be a basic dihedral progression with ratio $r$ given as follows.

\begin{equation}
\label{defb}
B = \{1,f\} \{1, r, r^2, \ldots, r^{\ell} \} = \{1, r, \ldots, r^{\ell} \} \{ 1, f r ^{\ell} \}.
\end{equation}

This gives us 
\begin{align*}
AB 	&= \{1, r, \ldots, r^k\}  \{1,f\} \{1, r, \ldots, r^{\ell} \}	\\
	&=  \{1, f r^{-k} \} \{1, r, \ldots, r^k\}  \{1, r, \ldots, r^{\ell} \}	\\
	&=  \{1, f r^{-k} \} \{ 1, r, \ldots, r^{k + \ell} \} = \{1, r, \ldots, r^{k + \ell} \} \{ 1, fr^{\ell} \}.
\end{align*}

So, $AB$ will be another basic dihedral progression with ratio $r$.  If $AB$ is proper, then $|AB| = 2(k+\ell+1) = |A| + |B| - 2$ so $(A,B)$ is a nontrivial critical pair.  Note that in this case ${\mathit stab}_L(A) = {\mathit stab}_L(AB) = \{1, fr^{-k} \}$ and ${\mathit stab}_R(AB) = {\mathit stab}_R(B) = \{1, fr^{\ell} \}$. We now define a type of subset trio which captures this behaviour.

\bigskip

\noindent{\bf Definition:}
Let $H \triangleleft G$ be finite with $G/H$ dihedral and $|G/H| \ge 10$.  A subset trio is a \emph{pure dihedral chord} \emph{relative to} $H$
if there is a similar subset trio $(A,B,C)$ and $S \in G/H$ generating the rotation subgroup so that
\begin{enumerate}
\item $A,B,C$ are $H$-stable.
\item $\varphi_{G/H}(A)$ and $\varphi_{G/H}(B)$ are nontrivial basic dihedral progressions with ratio $S$ and 
${\mathit stab}_R(A) = {\mathit stab}_L(B)$.
\item $C = \overline{(AB)^{-1}}$ satisfies $|C| > 2|H|$.
\end{enumerate}
Observe that every pure dihedral chord is a maximal nontrivial critical trio, and setting $H_1 = {\mathit stab}_L(A) = {\mathit stab}_R(C)$ and $H_2 = {\mathit stab}_R(A) = {\mathit stab}_L(B)$ and $H_3 = {\mathit stab}_R(B) = {\mathit stab}_L(C)$ we have that $\delta(A,B,C) = |H_1| = |H_2| = |H_3|$.

\subsection{Dihedral Prechords}

As was the case with cyclic chords, there is also an impure version of dihedral chords.  Analogously, the group $G$ will have a normal subgroup $H$ so that the three sets $A$ and $B$ and $C$ are close to $H$-stable sets with quotient a dihedral progression, but with some impurity at the ends.  However, before defining these we will first need to introduce a type of triple of subsets (in a dihedral group) which is not a subset trio, but has only a small number of different combinations which multiply to give the identity.  This is the purpose of this subsection.  For a basic dihedral progression $A$ as given by equation \ref{defa} we define the \emph{ends} of $A$ to be $(1, f, r^k, fr^{-k})$.  

\bigskip

\noindent{\bf Definition:}
Let $G$ be a dihedral group of order at least $8$ with rotation subgroup generated by $r$, let $A,B$ be basic nontrivial dihedral progressions with ratio $r$ 
and assume that ${\mathit stab}_R(A) = {\mathit stab}_L(B)$ and that $X = AB \neq G$.  Then $X$ is another basic dihedral progression with ratio $r$,  we let $(x_1,x_2,x_3,x_4)$ denote its ends, and put $C =\overline{X^{-1}} \cup \{x_1^{-1},x_2^{-1},x_3^{-1},x_4^{-1}\}$.  We define the triple $(A,B,C)$ to be a \emph{dihedral prechord}.

\bigskip

One key property of a dihedral prechord $(A,B,C)$ which is an immediate consequence of our definition is the following 
\begin{equation}
\label{prechord_excess}
 |A| + |B| - |\overline{C}| = 6. 
\end{equation}
Consider a dihedral prechord as in the above definition where $A$ and $B$ are given by equations \ref{defa} and \ref{defb}.  Label the ends of $A$ and $B$ and $X$ as follows.
\begin{align*}
&\mbox{$A$ has ends $(a_1,a_2,a_3,a_4) = (1,f,r^k, fr^{-k})$},		\\
&\mbox{$B$ has ends $(b_1,b_2,b_3,b_4) = (1,fr^{\ell}, r^{\ell}, f)$},		\\
&\mbox{$X$ has ends $(x_1,x_2,x_3,x_4) = (1, fr^{\ell}, r^{k+ \ell},  fr^{-k})$}.		
\end{align*}
Now let us consider the different ways that the ends of the progression $X$ can be produced as a product of an element in $A$ and an element in $B$.  For this purpose, define the representation function ${\mathit Rep}_{A,B}(g) = \{ (a,b) \in A \times B \mid ab = g \}$.  It follows immediately from the assumption ${\mathit stab}_R(A) =  {\mathit stab}_L(B) = \{1,f\}$ that every element $x \in X$ will have $| {\mathit Rep}_{A,B}(x)| \ge 2$.  For the ends of $X$, this bound will be met with equality as indicated below.
\[
\begin{array}{lclclcl}
{\mathit Rep}_{A,B}( x_1 )		&=& {\mathit Rep}_{A,B}( 1 )		&=&\{ (1,1), (f,f) \}				&=& \{ (a_1, b_1), (a_2, b_4) \},	\\
{\mathit Rep}_{A,B}( x_2 )		&=&{\mathit Rep}_{A,B}( fr^{\ell}) 	&=&\{ (1,fr^{\ell}) , (f,r^{\ell}) \}		&=& \{ (a_1, b_2), (a_2,b_3) \}, 	\\
{\mathit Rep}_{A,B}( x_3 )		&=&{\mathit Rep}_{A,B}( r^{k+\ell})	&=&\{ (r^k,r^{\ell}), (fr^{-k}, fr^{\ell})\}	&=& \{ (a_3,b_3), (a_4,b_2) \}, \\
{\mathit Rep}_{A,B}( x_4 )		&=& {\mathit Rep}_{A,B}( fr^{-k})	&=&\{ (r^k,f), (fr^{-k},1) \}			&=& \{ (a_3, b_4), (a_4, b_1) \}.
\end{array}
\]
Next we define the following group elements (note the intentional flip in indexing yielding $c_2 = x_4^{-1}$ and $c_4 = x_2^{-1}$)\footnote{The cause of this slightly strange flip is that the object of interest is not $X$, but $C = \overline{X^{-1}} \cup \{x_1^{-1},x_2^{-1},x_3^{-1},x_4^{-1}\}$.  If $C$ is finite, then it is a basic dihedral progression with ends $(c_1,c_2,c_3,c_4)$.}.
\[ (c_1,c_2,c_3,c_4) = (x_1^{-1}, x_4^{-1}, x_3^{-1}, x_2^{-1}) = (1,fr^{-k},  r^{-k- \ell}, fr^{\ell},). \]

\begin{figure}[ht]
\centerline{\includegraphics[height=4.2cm]{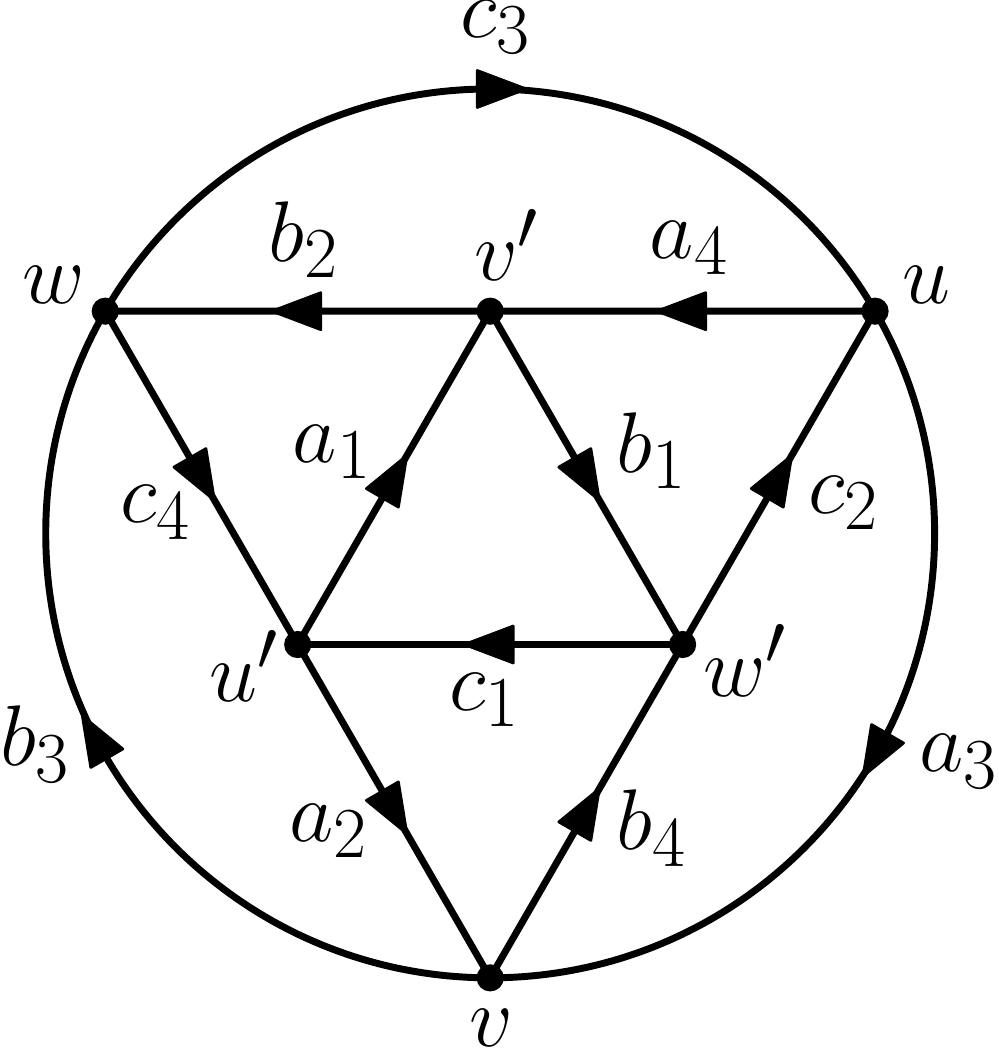}}
\caption{A labelled octahedron}
\label{lab_oct1}
\end {figure}

The key point is that the three sets $A,B,C$ have the property that there are exactly 8 elements $(a,b,c)$ in $A \times B \times C$ for which $abc = 1$ and these 8 products are precisely those labelling the edges around each face of the octahedron pictured in Figure \ref{lab_oct1}.  We call this figure a \emph{labelled octahedron associated with} $(A,B,C)$, and it will prove helpful since it conveniently encodes all trivial triple products from our sets.  Next we shall explore the structure of this labelling a little further.

Assign the vertex labels from Figure \ref{lab_oct1} as follows
\begin{align*}
	u' = v' &= w' = 1,	\\
	u &= c_2,	\\
	v &= a_2,	\\
	w &= b_2.
\end{align*}
Now consider an edge $e$ in our octahedron which goes from a vertex with label $x$ to one with label $y$.  A check of our assignments reveals that the label on the edge $e$ is given by $x^{-1} y$.  In fact, this is a special case of a more general phenomena.  For an arbitrary directed graph $\Gamma$, a map $\psi : E(\Gamma) \rightarrow G$ is called a \emph{tension} if for every closed walk in the underlying graph with edge sequence $e_1, \ldots, e_n$ we have $ \prod_{i=1}^n (\psi(e_i))^{\eta_i}  = 1 $ where $\eta_i = 1$ if $e_i$ is traversed in the forward direction and $\eta_i = -1$ if $e_i$ is traversed backward.  A basic property of every tension $\psi$ is that there exists a map $\tau : V( \Gamma ) \rightarrow G$ called a \emph{potential} function so that every edge $e = (x,y)$ satisfies $\psi(e) = \tau(x)^{-1} \tau(y)$.  In our case, the product of the edge labels on each of the closed walks bounding a face of our octahedron is trivial, and it then follows by an easy inductive argument that these labels form a tension.  The vertex labelling given above is a potential function associated with this tension.

\subsection{Impure Dihedral Chords}

In this subsection we will use the notion of a dihedral prechord to introduce a type of critical behaviour which occurs in groups which have dihedral quotients.  

\bigskip

\begin{figure}[ht]
\centerline{\includegraphics[height=4.2cm]{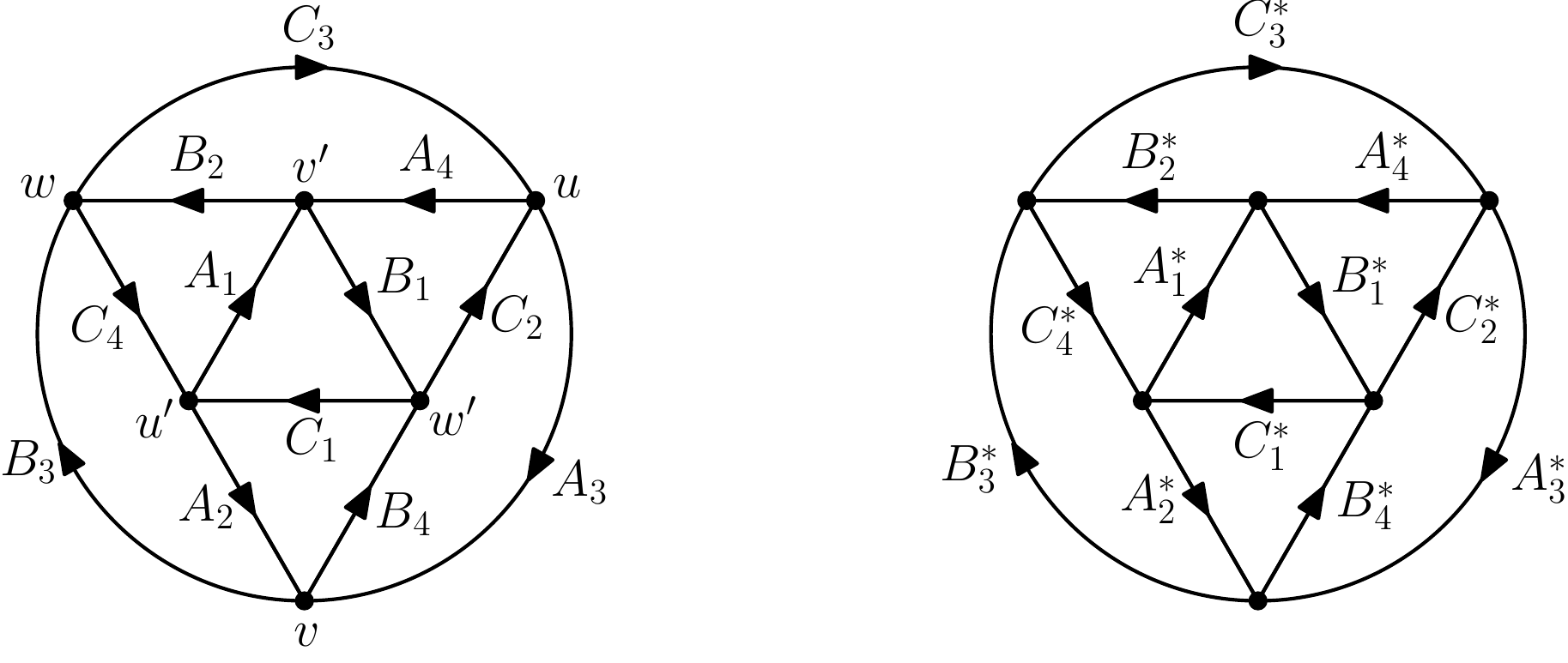}}
\caption{labelled octahedra}
\label{oct_rem}
\end {figure}

\noindent{\bf Definition:}
Let $H \triangleleft G$ be finite, with $G/H$ is dihedral and $|G/H| \ge 8$.  A subset trio is an \emph{impure dihedral chord} (\emph{relative to} $H$) if there exists a similar subset trio $(A,B,C)$ and $A_1,\ldots,A_4$, $B_1,\ldots,B_4$, $C_1, \ldots, C_4 \in G/H$ so that $A^+ = A \cup (\cup_{i =1}^4 A_i)$ and $B^+ = B \cup (\cup_{i =1}^4 B_i)$ and $C^+ = C \cup (\cup_{i=1}^4 C_i)$ satisfy
\begin{enumerate}
\item $A^+$ and $B^+$ and $C^+$ are $H$-stable.
\item $(\varphi_{G/H} (A^+), \varphi_{G/H}(B^+), \varphi_{G/H}(C^+))$ is a dihedral prechord with labelled octahedron indicated on the left in Figure \ref{oct_rem}.
\item $0 < |A^+ \setminus A| < 4 |H|$ and $0 < |B^+ \setminus B| < 4|H|$ and $0 < |C^+ \setminus C| < 4|H|$.  
\end{enumerate}

Let us pause to make a technical comment concerning the bound $|G/H| \ge 8$ in the above definition.  Although it would be possible to extend the notion of an impure dihedral chord to the case when $G/H$ is a dihedral group of order 6, it follows from the analysis in Section \ref{sequence_sec} that a maximal critical subset trio of this form would also be either a pure or impure beat.

Our next task will be to investigate the possible behaviour of the sets $A \cap (\cup_{i=1}^4 A_i)$ and $B \cap (\cup_{i=1}^4 B_i)$ and 
$C \cap (\cup_{i=1}^4 C_i)$ (which is intentionally obscured in the above definition).  To do so, we shall follow the same procedure as before to choose 
$u,u',v,v',w,w'  \in G$ so that in Figure \ref{oct_rem} the label on the edge $(x,y)$ is given by $x^{-1} H y$.  Now, for $1 \le i \le 4$ we define a new set $A_i^*$ for as follows:  If the initial vertex of the $A_i$ edge is labelled with $x$ and the terminal vertex labelled with $y$ then we set $A_i^* = x ( A_i \cap A ) y^{-1}$.  We define the sets $B_i^*$ and $C_i^*$ for $1 \le i \le 4$ similarly.  

It follows from our definitions that $A_1^* \ldots A_4^*, B_1^*, \ldots, B_4^*, C_1^* \ldots C_4^* \subseteq H$.  Furthermore, in the octahedron on the right in Figure \ref{oct_rem} we have the property that the product of the three sets (in order) labelling the edges of any oriented triangle does not contain the identity.  To see this last claim, suppose (for a contradiction) that $1 \in A_3^* B_3^* C_3^*$ and observe that
\[ 1 \in u (A \cap A_3) v^{-1} v (B \cap B_3) w^{-1} w (C \cap C_3) u^{-1} = u (A \cap A_3)(B \cap B_3)(C \cap C_3) u^{-1} \]
which contradicts the assumption that $(A,B,C)$ is a subset trio.  It follows from equation \ref{prechord_excess} that the 
deficiency of $(A,B,C)$ is given by the following formula
\[ \delta(A,B,C) = \sum_{i=1}^4 \big( |A_i^*| + |B_i^*| + |C_i^*| \big) - 6 |H|. \]
This prompts a new definition.

\bigskip

\noindent{\bf Definition:}
For a finite group $H$, an \emph{octahedral configuration} (\emph{in} $H$) consists of an oriented octahedron in which each edge is labelled with a subset of $H$ as in Figure \ref{oct_rem}, with the property that the product of the labels (in order) on any three edges forming a directed triangle does not contain $1$.  

\bigskip

An octahedral configuration \emph{maximal} if  it is not possible to replace any of the sets with a proper superset while maintaining an octahedral configuration.  We define the \emph{deficiency} of an octahedral configuration in $H$ as given on the right in Figure \ref{oct_rem} to be
$\sum_{i=1}^4 \big( |A_i^*| + |B_i^*| + |C_i^*| \big) - 6 |H|$ and we say that such a configuration is \emph{critical} if it has positive deficiency.

We say that an octahedral configuration derived as above from an impure dihedral chord $(A,B,C)$ is \emph{associated with} 
$(A,B,C)$.  Our definitions of maximal and deficiency imply that the deficiency of an associated octahedral configuration is equal to the deficiency of $(A,B,C)$, and this associated octahedral configuration is maximal whenever $(A,B,C)$ is maximal.

To complete our classification we shall determine all maximal critical octahedral configurations in a finite group $H$.  In our next subsection we turn our attention to this problem.

\subsection{Critical Octahedral Configurations}

In a maximal critical octahedral configuration in the group $H$, many of the sets involved will be either $\emptyset$ or the entire group $H$.  Next we will introduce some notation to help work with this situation.  In the following figures, we have depicted octahedra with three types of edges: a dotted line indicates that the label on the associated edge is empty, a solid dark line indicates it is the full group $H$, and a grey line indicates that its label is a proper nonempty subset of $H$.  

\bigskip

\begin{figure}[ht]
\centerline{\includegraphics[width=12cm]{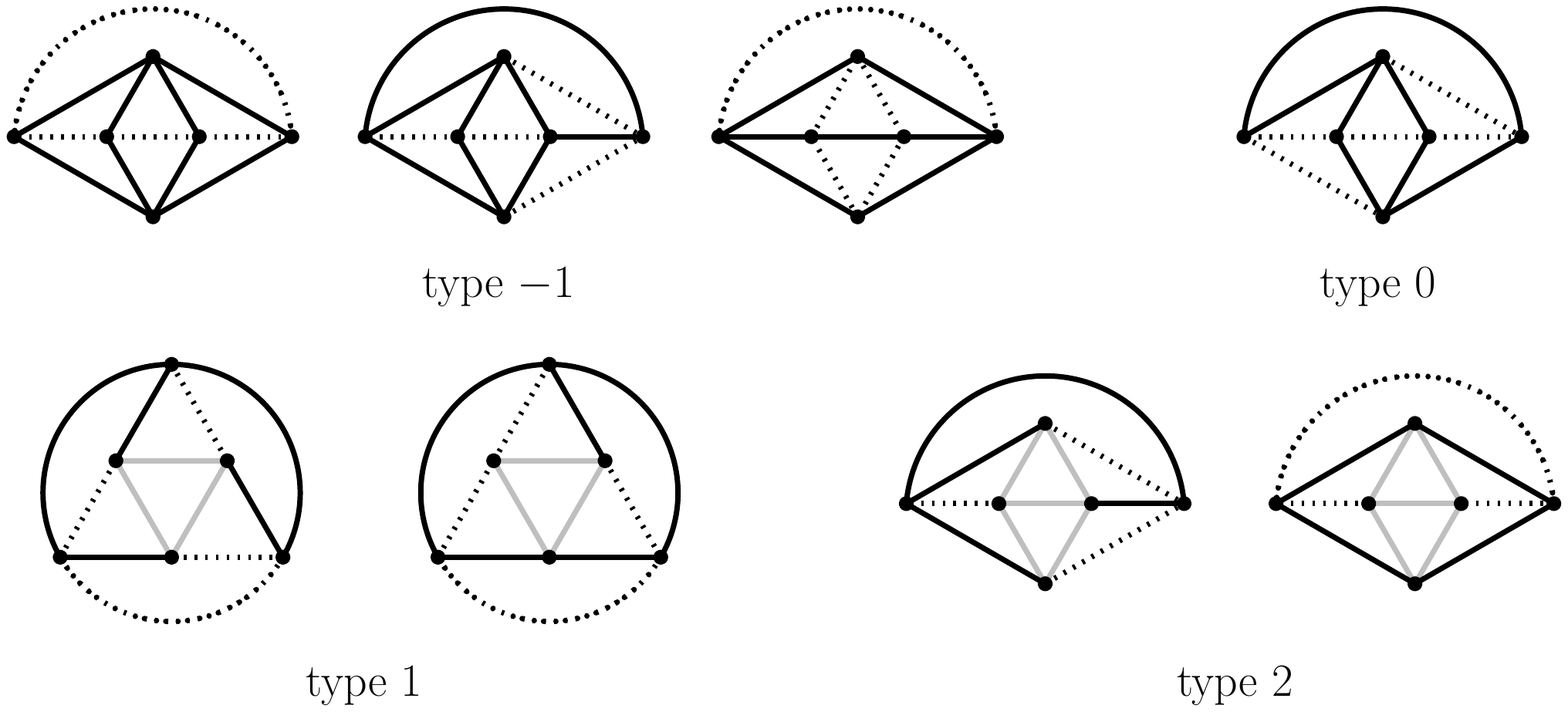}}
\caption{Octahedral configurations}
\label{all_oct}
\end {figure}

\noindent{\bf Definition:} We say that an octahedral configuration has \emph{type} $-1$, $0$, $1$, or $2$ if its associated graph appears in Figure \ref{all_oct} with the corresponding label.\footnote{To help explain this numbering, let us comment that configurations of type $-1$ will be inessential, and configurations of type $0$, $1$, and $2$ have the corresponding number of triangles of grey edges.}.  

In Section \ref{prod_again} we will prove the following theorem which classifies maximal critical octahedral configurations.

\begin{theorem}
\label{oct_config_class}
Every maximal critical octahedral configuration has type $-1$, $0$, $1$, or $2$.
\end{theorem} 

If $(A,B,C)$ is an impure dihedral chord with an associated octahedral configuration of type $n$, then we say that $(A,B,C)$ has \emph{type} $n$.  We shall prove that an impure dihedral chord of type $-1$ cannot be a maximal subset trio, so this type of impure dihedral chord will not appear in our structure theorem.  An impure dihedral chord of type $0$ is a fully described maximal critical trio which will appear as a possible outcome in our structure theorem.  Next consider an impure dihedral chord $(A,B,C)$ of type $1$, and assume that it has an associated octahedral configuration in $H$ for which the labels on the unique directed triangle of grey edges are $A'$, $B'$, $C'$ (in this order).  Then $(A',B',C')$ is a nontrivial subset trio in $H$ which we call a \emph{continuation} of $(A,B,C)$.  As with our other structures, $\delta(A,B,C) = \delta(A',B',C')$ and $(A',B',C')$ is maximal whenever $(A,B,C)$ is maximal.  Maximal critical impure dihedral chords of type $2$ have some additional structure, and we shall now introduce some definitions to handle this.

\bigskip

\noindent{\bf Definition:} Consider an octahedral configuration in $H$ of type $2$ in which the labels on the two directed triangles of grey edges are $A',B',C'$ and $A',B'',C''$, and let $x \in A'$ and $K',K'' < H$.
\begin{enumerate}
\item[(2A)] We say this configuration has \emph{type} $2A$ if $A' \subseteq xK' \subseteq xK''$ and $(A',B',C')$ is an impure beat relative to $K'$ and $(xK'',B'',C'')$ is a pure 
beat relative to $K''$.  In this case we define a \emph{continuation} of this configuration to be a continuation of $(A',B',C')$.  
\item[(2B)] We say this configuration has \emph{type} $2B$ if $A' = xK' \cap xK''$ and both 
$(xK', B', C')$ and $(xK'', B'', C'')$ are pure beats (relative to $K'$, $K''$).  
\end{enumerate}

\bigskip

Extending our notation in the obvious manner, we will say that a dihedral chord has \emph{type} $2A$ or $2B$ if it has an associated octahedral configuration of this type.  A dihedral configuration of type $2B$ is a fully described maximal critical trio, so these will appear as a possible outcome in our structure theorem.  For a dihedral chord $(A,B,C)$ of type $2A$, we define a \emph{continuation} of $(A,B,C)$ to be a continuation of its associated octahedral configuration.  As in the other cases, if $(A',B',C')$ is such a continuation, then $\delta(A',B',C') = \delta(A,B,C)$ and $(A',B',C')$ will be maximal whenever $(A,B,C)$ is maximal.  In Section \ref{prod_again} we shall prove the following result.

\begin{theorem}
\label{oct_set_classify}
Every type $2$ maximal critical octahedral configuration has type $2A$ or $2B$.
\end{theorem}

In summary, every maximal critical subset trio which is an impure dihedral chord will have type $0$, $1$, $2A$, or $2B$.  Those of types $0$ or $2B$ are fully described maximal critical trios, while those of types $1$ and $2A$ have continuations which are maximal critical subset trios in a proper (finite) subgroup.

\subsection{The Classification Theorem}

There is another simple and flexible construction of a critical pair (or subset trio) which is well known, and appears for instance as the exceptional case in Serra and Zemor's paper \cite{serra-zemor} discussed earlier.  
Namely, choose a finite subgroup $H \le G$ and $x \in G \setminus H$ and set $A = H \cup xH$ and $B = H \cup Hx$.  In this case $|A| = |B| = 2|H|$ but 
$AB = H \cup xH \cup Hx \cup xHx$ will have size less than $4|H|$ (since $xH \cap Hx \neq \emptyset$), so $(A,B)$ will be critical.  
Of course, in the special case when $G$ is abelian, this is nothing new since in this case $A = B$ is a geometric progression in the quotient group $G/H$.  However, for general groups this does yield new critical pairs.  

There are additional constructions of critical pairs $(A,B)$, where $A$ is a union of a small number of left $H$-cosets and $B$ is a union of a small number of right $H$-cosets.  More precisely, there are two infinite families for which $AH = A$ and $HB = B$ and $\big\{ \tfrac{|A|}{|H|}, \tfrac{|B|}{|H|} \big\} = \{2,3\}$, and there are a small number of exceptional examples which are based on the Platonic solids and the regular maps on the projective plane.  We shall postpone discussion of these since they are easier to understand within the framework developed in the following section.  We now state our main classification theorem for critical subset trios.  

For every positive integer $n$ we let ${\mathbf S}_n$ and ${\mathbf A}_n$ denote the symmetric and alternating groups of permutations of $\{1,2,\ldots,n\}$.

\begin{theorem}
\label{main_subset}
If $\Phi$ is a nontrival maximal critical subset trio in $G$, then there exists a sequence of subgroups 
$G = G_1 > G_2 \ldots G_m$ and a sequence of subset trios $\Phi = \Phi_1, \Phi_2, \ldots, \Phi_m$ 
so that the following hold:
\begin{enumerate}
\item $\Phi_i$ is a subset trio in the group $G_i$ for $1 \le i \le m$.
\item $\Phi_i$ is either an impure beat, an impure cyclic chord, or an impure dihedral chord of type $1$ or $2A$  with continuation
$\Phi_{i+1}$ for $1 \le i < m$.
\item  $\Phi_m$ is either a pure beat, a pure cyclic chord, a pure dihedral chord, an impure dihedral chord of type $0$ or $2B$, or there exists a subset trio $(A,B,C)$ similar to $\Phi_m$ and a (possibly trivial) subgroup
$H \triangleleft G_m$ together with subgroups $H \triangleleft H_1,H_2,H_3 \le G_m$ so that $H_1 = {\mathit stab}_L(A) = {\mathit stab}_R(C)$,  
$H_2 = {\mathit stab}_R(A) = {\mathit stab}_L(B)$, $H_3 = {\mathit stab}_R(B) = {\mathit stab}_L(C)$ and one of the following holds:
\begin{enumerate}
\item $|A| = 2|H_2| = |B|$ and $\delta(A,B,C) =  |H_1| =  |H_3| < |H_2|$.
\item $9|H_1| = |A| = 3|H_2|$ and $2|H_2| = |B| = 3|H_3|$ and $\delta(A,B,C) = |H_1|$.  
Further, if $G_m$ is finite, then 
$[H_1 : H]$ divides $16$.
\item $4|H_1| = |A| = 2|H_2|$ and $3|H_2| = |B| = 2|H_3|$ and $\delta(A,B,C) = |H_1|$.  
Further, if $G_m$ is finite, then $[H_1 : H]$ divides $16$.
\item $G_m/H$ and $A,B,C$ and $H_1,H_2,H_3$ satisfy one of the following. \\
\begin{tabular}{|c|c|c|c|c|c|c|}
\hline
\begin{picture}(0,19)\end{picture} \raisebox{.95ex}{$G_m/H$}	
	&\raisebox{.95ex}[0pt]{$\tfrac{|G_m|}{|H_1|}$}	&\raisebox{.95ex}[0pt]{$\tfrac{|G_m|}{|H_2|}$}	
	&\raisebox{.95ex}[0pt]{$\tfrac{|G_m|}{|H_3|}$}	&\raisebox{.95ex}[0pt]{$\tfrac{|G_m|}{|A|}$}	
	&\raisebox{.95ex}[0pt]{$\tfrac{|G_m|}{|B|}$}		&\raisebox{.95ex}[0pt]{$\tfrac{|G_m|}{|C|}$} 	\\
\hline
${\mathbf C_2} \times {\mathbf S_4}$, ${\mathbf C_2} \times {\mathbf A_4}$, or ${\mathbf S_4}$ 
		& $8$	& $12$	& $6$	& $4$	& $3$	& $2$	\\	
${\mathbf C_2} \times {\mathbf A_5}$ or ${\mathbf A_5}$
		& $12$	& $20$	& $30$	& $4$	& $10$	& $3/2$	\\	
${\mathbf C_2} \times {\mathbf A_5}$ or ${\mathbf A_5}$
		& $20$	& $30$	& $12$	& $10$	& $6$	& $4/3$	\\	
${\mathbf C_2} \times {\mathbf A_5}$ or ${\mathbf A_5}$
		& $20$	& $12$	& $30$	& $4$	& $6$	& $5/3$	\\	
${\mathbf S_5}$
		& $12$ & $15$ & $10$ & $3$ & $5$ & $2$ \\
${\mathbf A_5}$
		& $6$ & $15$ & $10$ & $3$ & $5$ & $2$ \\
${\mathbf S_6}$, or ${\mathbf A_6}$, or ${\mathbf S}_5$
		& $20$ & $15$ & $6$ & $5$ & $3$ & $2$ \\
\hline
\end{tabular}
\end{enumerate}
\end{enumerate}
\end{theorem}

This paper gives an almost entirely self-contained proof of this theorem.   Indeed, the only significant result we rely on which is not proved here is a theorem of Tutte concerning arc-transitive cubic graphs which we use to get the bounds on $[H_1 : H]$ in (b) and (c) of the third part of the outcome.
This theorem as stated is not a proper characterization since we have not supplied a converse.  However, in Section \ref{prod_again} we will remedy this by establishing a variant of the above result for which we do have a converse.

\section{Incidence Trios}
\label{inc_sec}

As indicated earlier, our approach to classifying critical subset trios in the group $G$ will involve an incidence geometry on which $G$ has a natural action.  This step brings us from groups to group actions, and has the advantage that it permits us to first classify the extremal structures, and then to determine which groups act on them.  The purpose of this section is to introduce the incidence geometries we will work with, and to motivate this approach with a natural example.  This example will also give us opportunity to introduce some of the exceptional critical phenomena.

\subsection{Incidence Geometry}

In this subsection we will introduce the basic terminology we will use for working with incidence geometries.

\bigskip

\noindent{\bf Definition:} An \emph{incidence geometry} $\Lambda = (X_1, \ldots, X_n; \sim)$ consists of a finite sequence of pairwise disjoint sets $(X_1,\ldots,X_n)$ together with a symmetric binary relation $\sim$ on $\cup_{i=1}^n X_i$, called the \emph{incidence} relation, which satisfies the following rule: For every $1 \le i \le n$,  if $x,y \in X_i$ then $x \not\sim y$.  We define $\cup_{i=1}^n X_i$ to be the \emph{ground set} of $\Lambda$ and define the \emph{rank} of $\Lambda$ to be $n$.  If $x,y \in X_i$ for some $1 \le i \le n$ we say that $x$ and $y$ have the \emph{same type}, and otherwise they have \emph{different type}.  

\bigskip

For every point $x \in X_i$, the \emph{neighbours} of $x$ are $N_{\Lambda}(x) = \{ y \in \cup_{i=1}^n X_i \mid y \sim x \}$ and for a subset of points $A$ we let $N_{\Lambda}(A) = \cup_{x \in A} N_{\Lambda}(x)$.  When the incidence geometry is clear from context, we will drop these subscripts. 

For brevity, we will frequently write $\Lambda = (X_1, \ldots, X_n)$ to specify an incidence geometry, and treat the incidence relation $\sim$ as implied.  
Also for convenience, when $X_1, \ldots, X_n$ are disjoint sets and $\sim$ is a relation defined on a superset of $\cup_{i=1}^n X_i$ satisfying $x \not\sim y$ 
whenever $x,y \in X_i$, we will let $(X_1, \ldots, X_n ; \sim )$ denote the incidence geometry given by $(X_1, \ldots, X_n; \sim')$ where $\sim'$ is the restriction of $\sim$ to $\cup_{i=1}^n X_i$.   

Next we introduce a few ways to construct a new incidence geometry from an existing one.   If $\sigma \in {\mathbf S}_n$ then $\Lambda' = (X_{\sigma(1)}, \ldots, X_{\sigma(n)} ; \sim )$ is also an incidence geometry, and we say that $\Lambda$ and $\Lambda'$ are \emph{equivalent} and write $\Lambda \equiv \Lambda'$.  
If $i_1, i_2, \ldots, i_{\ell} \in \{1,2,\ldots,n\}$ are distinct, and for every $1 \le j \le \ell$ we have a set $X'_{i_j} \subseteq X_{i_j}$, then 
$\Lambda'' = (X'_{j_1}, \ldots, X'_{j_{\ell}} ; \sim)$ is an incidence geometry which we say is \emph{contained} in $\Lambda$, and we write $\Lambda'' \subseteq \Lambda$.  In this case we will also call $\Lambda''$ a \emph{subgeometry} of $\Lambda$.  Finally, if $1 \le i,j \le n$ and $i \neq j$, then $(X_i,X_j ; \sim)$ is a subgeometry of $\Lambda$ which we call a \emph{side of} $\Lambda$.  

For incidence geometries defined on the same sequence of sets, there is a natural partial order which will be of interest to us.  If $\Lambda^* = (X_1, \ldots, X_n, \sim^*)$ is another incidence geometry we write $\Lambda \le \Lambda^*$ if for every $x,y \in \cup_{i=1}^n X_i$ we have $x \sim y$ only if $x \sim^* y$.  If $\Lambda \le \Lambda^*$ and $\Lambda \neq \Lambda^*$ we write $\Lambda < \Lambda^*$.

Isomorphisms will play an important role for us, so let us introduce this concept here.
If $\Pi = (Y_1, \ldots, Y_n ; \smile)$ is another incidence geometry, an \emph{isomorphism} from $\Lambda$ to $\Pi$ is a map $\phi: \cup_{i=1}^n X_i \rightarrow \cup_{i=1}^n Y_i$ so that $\phi$ maps $X_i$ bijectively to $Y_i$ for every $1 \le i \le n$, and so that for all $x,x' \in \cup_{i=1}^n X_i$ we have $x \sim x'$ if and only if $\phi(x) \smile \phi(x')$.  If such a map exists we say that $\Lambda$ and $\Pi$ are \emph{isomorphic} and write $\Lambda \cong \Pi$.  An isomorphism from $\Lambda$ to itself is called an \emph{automorphism} and we let $Aut(\Lambda)$ denote the group of automorphisms (under composition).  

Next we introduce a concept which plays an important role in our proof.  We say that an incidence geometry $\Delta$ of rank $2$ is \emph{disconnected} if there exists a partition of the ground set of $\Delta$ into $\{U,W\}$ so that no element in $U$ is incident with an element in $W$.  We say that $\Delta$ is \emph{connected} if it is not disconnected, and we define a \emph{connected component} of $\Delta$ to be a maximal connected subgeometry.  
For an incidence geometry $\Lambda$ of rank greater than $2$, we define $\Lambda$ to be \emph{connected}\footnote{Although this will prove to be the meaningful type of connectivity for us, it should be noted that this is not the standard use of the term in incidence geometry.} if every side of $\Lambda$ is connected.   
 
 Finally, let us close this subsection by introducing a particularly important family of rank 2 incidence geometries.  We 
define a \emph{graph} to be a rank 2 incidence geometry $(V,E)$ with the property that $|N(e)| = 2$ for every $e \in E$.  This is the standard definition of a graph where multiple edges are permitted, but loops are forbidden.  We say that a graph is \emph{simple} if it has no parallel edges (i.e. there do not exist distinct $e,e' \in E$ with $N(e) = N(e')$), and for convenience, we will say that an incidence geometry is \emph{graphic} if it is a graph.

\subsection{Choruses}

Our main focus will be on incidence geometries which are equipped with a natural group action.
Throughout this paper our groups will always act on the left.  So, if $G$ acts on the set $X$ and $x \in X$ and $g \in G$ we let $gx$ denote the image of $x$ under the permutation associated with $g$.  For a subset $S \subseteq X$ we let $gS = \{ gs \mid s \in S \}$ and define $G_{S} = \{ g \in G \mid gS = S \}$ and $G_{ (S) } = \{ g \in G \mid \mbox{$gs= s$ for every $s \in S$} \}$.  We will shortcut this notation for singletons and write $G_x = G_{ \{ x \} }$.  The action is \emph{transitive} if for every $x,y \in X$ there exists $g \in G$ for which $gx = y$ and it is \emph{regular} if there is a unique such $g$.  A \emph{block of imprimitivity} (sometimes abbreviated \emph{block}) is a set $S \subseteq X$ so that for every $g \in G$ either $gS = S$ or $gS \cap S = \emptyset$.  If $S$ is a block, and the action of $G$ is transitive, then $\{ gS \mid g \in G \}$ is a partition of $X$ which we call a \emph{system of imprimitivity generated by} $S$.  

\bigskip

\noindent{\bf Definition:} For a group $G$, a $G$-\emph{chorus} consists of an incidence geometry 
$\Lambda =(X_1,\ldots,X_n)$ equipped with an action of $G$ on $\Lambda$ (i.e. a homomorphism $G \rightarrow Aut(\Lambda)$) 
so that $G$ acts transitively on $X_i$ for every $1 \le i \le n$.

\bigskip

Next we define a natural family of choruses.  

\bigskip

\noindent{\bf Definition:} Let $G$ be a group and let $M$ be an $n \times n$ matrix of subsets of $G$ with the property that $M_{ij} = M_{ji}^{-1}$ and $M_{ii} = \emptyset$ for every $1 \le i,j \le n$.  We define the \emph{Cayley Chorus}, denoted ${\mathit CayleyC}(G;M)$, to be the incidence geometry 
$(G \times \{1\}, G \times \{2\}, \ldots, G \times \{n\})$ where $(x,i)$ and $(y,j)$ are incident if and only if $y \in x M_{ij}$.  

\bigskip

Cayley choruses of rank 2 and 3 will be of frequent interest for us, so we add some shortcuts to the above notation.  For sets $A,B,C \subseteq G$ we let
\[ M = \left[ \begin{array}{cc}
			\emptyset		&A	\\
			A^{-1}	&\emptyset
			\end{array}  \right]
		\quad\quad\quad
 M' = \left[ \begin{array}{ccc}
  			\emptyset		&A		&C^{-1}	\\
			A^{-1}	&\emptyset	&B		\\
			C		&B^{-1}	&\emptyset
			\end{array} \right]
\]
and define 
\begin{align*}
{\mathit CayleyC}(G;A) &= {\mathit CayleyC}(G;M), \\
{\mathit CayleyC}(G; A,B,C) &= {\mathit CayleyC}(G;M'). 
\end{align*}
Note that in the above definition of a Cayley chorus, we have defined incidence in terms of right multiplication.  It follows from this that the group $G$ has a natural left action on the structure (sending each $(x,i)$ to $(gx,i)$) which shows that Cayley choruses are indeed choruses.  The following observation follows immediately from our definitions, and determines when a Cayley chorus is connected.  

\begin{observation} If $A \subseteq G$ then ${\mathit CayleyC}(G; A)$ is connected if and only if $A$ is not contained in a coset of a proper subgroup.
\end{observation}

The central goal of the remainder of this subsection is to make the connection between general choruses and Cayley choruses explicit, but to do so will require a couple of concepts.  Let $\Lambda = (X_1,\ldots,X_n)$ be an incidence geometry and let $\mcp$ be a partition of $X_j$.  We define $\Lambda' = (X_1,\ldots,X_{j-1}, \mcp, X_{j+1}, \ldots, X_n)$ to be the incidence geometry with incidence relation given by the rule that $x \in X_i$ and $y \in X_k$ are incident in $\Lambda'$ if they are incident in $\Lambda$, and $x \in X_i$ and $P \in \mcp$ are incident in $\Lambda'$ if there exists $y \in P$ so that $x$ and $y$ are incident in $\Lambda$.  We call any geometry formed from $\Lambda$ by a sequence of these operations a \emph{quotient} of $\Lambda$.  Note that if $\Lambda$ is a $G$-chorus and $\mcp$ is a system of imprimitivity on $X_j$, then $\Lambda'$ is also a $G$-chorus.

Next we introduce an important relation.  We define two points $x,y$ in the ground set of an incidence geometry $\Lambda$ to be \emph{clones} if they are of the same type and $N(x) = N(y)$.  We say that $\Lambda$ is \emph{clone free} if it does not contain distinct points which are clones.  It follows immediately that the clone relation is an equivalence relation.  For every $1 \le i \le n$ let $\mcp_i$ be the partition of $X_i$ given by the clone relation.  We say that $\mcp_i$ is the \emph{clone partition} of $X_i$ and we define $\Lambda^{\bullet} = (\mcp_1, \ldots, \mcp_n)$.  Again we note that if 
$\Lambda$ is a $G$-chorus, then $\Lambda^{\bullet}$ is also a $G$-chorus.  More generally, if $\mcq_j$ is a system of imprimitivity on $X_j$ for some $1 \le j \le n$ so that any two points in the same block of $\mcq_j$ are clones, we call $\mcq_j$ a \emph{system of clones}.  If $\mcq_i$ is a system of clones for every $1 \le i \le n$ then we say that $(\mcq_1, \ldots, \mcq_n)$ is obtained from $\Lambda$ by \emph{identifying clones}.

Next we state an observation which follows immediately from our definitions.  If $H$ is a subgroup of the group $G$ then we let $G/H$ denote 
the set of all left $H$-cosets in $G$.

\begin{observation}
Let $A \subseteq G$ and assume that $H_1 = {\mathit stab}_L(A)$ and $H_2 = {\mathit stab}_R(A)$.  Then in the chorus 
${\mathit CayleyC}(G; A)$, the clone partition of $G \times \{i\}$ is $G/H_i \times \{i\}$.  
\end{observation}

If $\Lambda = (X_1, \ldots, X_n)$ and $\Pi = (Y_1, \ldots, Y_n)$ are two $G$-choruses, we say that a function $\phi : \cup_{i=1}^n X_i \rightarrow \cup_{i=1}^n Y_i$ is a \emph{strong isomorphism}
if it is an isomorphism, and furthermore, for every $g \in G$ and $x \in \cup_{i=1}^n X_i$ we have that $\phi(gx) = g \phi(x) $.  If such a map exists, we say that $\Lambda$ and $\Pi$ are \emph{strongly isomorphic}.

The next theorem shows that every chorus is strongly isomorphic to one obtained from a Cayley chorus by identifying clones.  
We give a particularly detailed statement of this result since this will be required in the final section.  This relationship given by this result mirrors the connection 
between Cayley graphs and general vertex transitive graphs, so we have attributed it to Sabidussi who proved the corresponding statement in that setting.

\begin{theorem}[Sabidussi \cite{sabidussi}]
\label{sabidussi}
Let $\Lambda = (X_1,\ldots,X_n)$ be a $G$-chorus and let $x_i \in X_i$ for every $1 \le i \le n$.  For every $1 \le i \le n$ 
associate each $x \in X_i$ with $S_x = \{ g \in G \mid g x_i = x \}$, and for every $1 \le i, j \le n$ define $A_{ij} = \cup_{y \in X_j \cap N(x_i)} S_y$.
Then we have:
\begin{enumerate}
\item $S_x$ is a left $G_{x_i}$ coset for every $x \in X_i$.
\item if $x,y \in X_i$ and $g \in G$, then $gx = y$ if and only if $gS_x = S_y$.
\item $G_{x_i} A_{ij} G_{x_j} = A_{ij}$ for every $1 \le i,j \le n$.
\item if $y \in X_i$ and $z \in X_j$, then $y \sim z$ if and only if $(S_y)^{-1} S_z \subseteq A_{ij}$.
\item if $M$ is the $n \times n$ matrix with $ij$ entry $A_{ij}$, then ${\mathit CayleyC}(G; M)$ 
has the property that $G / G_{x_i} \times \{i\}$ is a system of clones for $1 \le i \le n$, and associating each $x \in X_i$ with $S_x \times \{i\}$
is a strong isomorphism from $\Lambda$ to $(G/G_{x_1} \times \{1\}, \ldots, G/G_{x_n} \times \{n\})$.
\end{enumerate}
\end{theorem}

\noindent{\it Proof:} Associating every  $x \in X_i$ with $S_x$ yields the standard representation of the action of $G$ on $X_i$ by left $G_{x_i}$ cosets with ``base point'' $x_i$.  
The first two parts of the result are basic features of this representation, and follow immediately from the definitions.
For the third part, observe that we must have $A_{ij} G_{x_j} = A_{ij}$ since $A_{ij}$ is a union of left $G_{x_j}$ cosets.  To see that $G_{x_i} A_{ij} = A_{ij}$, note that $N(x_i) \cap X_j$ is setwise stabilized by $G_{x_i}$.  
For the fourth part, choose $g \in G$ so that $gx_i = y$ and observe that $y \sim z$ if and only if $x_i = g^{-1}y \sim g^{-1}z$ which, by part 2, holds 
if and only if $g^{-1} S_z \subseteq A_{ij}$.  However, in light of 3, this last condition is equivalent to $(S_y)^{-1} S_z = (g G_{x_i})^{-1} S_z \subseteq A_{ij}$.  
For the last part, observe that from our construction we have $A_{ii} = \emptyset$ and from 4 we have $A_{ij} = A_{ji}^{-1}$ for every $1 \le i, j \le n$.  So, we may define the Cayley chorus ${\mathit CayleyC}(G; M)$.  It follows from 3 that  $G/G_{x_i} \times \{i\}$ is a system of clones for every $1 \le i \le n$, and it follows from parts 2 and 4 that associating each $x \in X_i$ to $S_x \times \{i\}$ gives a strong isomorphism.
\quad\quad$\Box$

\bigskip

We close this subsection by introducing some important incidence geometries.  Though we have worked only with multiplicative groups so far, for indexing purposes and for defining Cayley incidence geometries, it will be helpful for us to occasionally use the additive groups ${\mathbb Z}/n{\mathbb Z}$ and ${\mathbb Z}$.  For $J = {\mathbb Z}$ or $J = {\mathbb Z}/n {\mathbb Z}$ we call an incidence geometry by the name below if it is isomorphic to an incidence geometry in the corresponding description.
\begin{center}
\begin{tabular}{|l|l|}
\hline
Name	&	Description	\\
\hline
\emph{polygon}	&	${\mathit CayleyC}( J ; \{0,1\} )$ with $|J| \ge 2$	\\
\emph{sequence}	&	${\mathit CayleyC}( J ; \{0,1,\ldots, \ell-1 \} )$ with $|J| \ge 4$ and $2 \le \ell \le |J|-2$ \\
\hline
\end{tabular}
\end{center}
We define a \emph{triangle} to be a polygon $(X,Y)$ with $|X| = |Y| = 3$.  Note that every polygon $(X,Y)$ with $|X| = |Y| \ge 4$ 
is also a sequence.

\subsection{Trios}
\label{finite_trios}

In this subsection we will restrict our attention to finite groups and finite incidence geometries since these are somewhat simpler to work with, and 
this is sufficient to motivate our approach.  The purpose of this subsection will be to introduce the type of chorus which will be our 
principle object of study (for the finite case).   In Subsection \ref{weights} we will extend this to the general case.  

\begin{figure}[ht]
\centerline{\includegraphics[height=4cm]{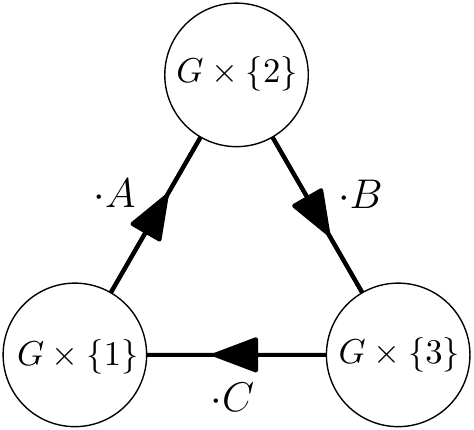}}
\caption{A Cayley chorus based on a subset trio}
\label{trio_basic}
\end {figure}

We say that an incidence geometry $\Lambda$ is \emph{collision free} if there do not exist distinct points $x,y,z$ in the ground set of $\Lambda$ with $x \sim y$ and $y \sim z$ and $z \sim x$.
If $(X_1,X_2)$ is an incidence geometry with $X_1$ and $X_2$ finite, we define the \emph{density} of $(X_1,X_2)$ to be 
\[ q(X_1,X_2) = \frac{ |\{ (x,y) \in X_1 \times X_2 \mid x \sim y \} |}{ |X_1| \cdot |X_2|}. \]  
Note that if there exists an integer $d$ so that every point in $X_1$ has exactly $d$ neighbours in $X_2$, then $q(\Lambda) = \frac{d}{|X_2|}$.  
Next we have a key observation concerning the Cayley chorus depicted in Figure \ref{trio_basic}.

\begin{observation}
Let $(A,B,C)$ be a critical subset trio in the finite group $G$.  Then  
$\Theta = {\mathit CayleyC}(G; A, B, C)$ satisfies the following properties:
\begin{itemize}
\item $\Theta$ is a $G$-chorus,
\item $\Theta$ is collision free,
\item $\sum_{1 \le i < j \le 3} \, q( G \times \{i\}, G \times \{j\} ) =  \frac{|A|}{|G|} + \frac{|B|}{|G|} + \frac{|C|}{|G|} > 1$.
\end{itemize}
\end{observation}

This prompts the following definition, which we make here only for finite groups.  

\bigskip

\noindent{\bf Definition:} If $G$ is a finite group, a $G$-\emph{(incidence) trio} $(X,Y,Z)$ is a collision free $G$-chorus of 
rank 3.  We say that $(X,Y,Z)$ is \emph{critical} if  $q(X,Y) + q(Y,Z) + q(Z,X) > 1$.

\bigskip

A key property to note here is that our notion of critical depends only on the underlying structure of the geometry and not on the particular group acting on it.  This feature permits us to first classify critical incidence trios and then to determine which groups act on them. 

We have already defined a number of important properties of subset trios, and now we shall introduce analogous properties for (incidence) trios.  We say that a trio $\Theta = (X,Y,Z)$ is \emph{trivial} if it has a side which does not contain a pair of incident points, and we say that $\Theta$ is \emph{maximal} if there does not exist a $G$-trio $\Theta'$ with $\Theta < \Theta'$.  Note that if $\Theta$ is a maximal $G$-trio and $x,y$ are points of distinct type in the ground set with $x \not\sim y$ then there must exist a point $z$ with $x \sim z$ and $y \sim z$ (otherwise we could increase the incidence relation by adding all incidences between pairs of points of the form $(gx,gy)$ for $g \in G$ to obtain a new trio contradicting the maximality of $\Theta$).  If $(A,B,C)$ is a subset trio in $G$, we call $\Theta = {\mathit CayleyC}(G; A, B, C)$  a \emph{Cayley (incidence) trio}.  Note that in this case we have the following equivalences:
\[ \mbox{$\Theta$ is } \left\{ \begin{array}{c} \mbox{critical} \\ \mbox{maximal} \\ \mbox{trivial} \end{array} \right\} \mbox{if and only if $(A,B,C)$ is} 
		\left\{ \begin{array}{c} \mbox{critical} \\ \mbox{maximal} \\ \mbox{trivial} \end{array} \right\}. \]

Next let us consider the notion of similarity for subset trios.  
Let $(A,B,C)$ be a subset trio in $G$ and define $\Theta = {\mathit CayleyC}(G; A, B, C)$.  Now let $x,y,z \in G$ and apply Theorem \ref{sabidussi} to 
$\Theta$ for the points $(x,1)$, $(y,2)$, and $(z,3)$ to produce a new trio $\Theta'$.  A quick check reveals that  
$\Theta' = {\mathit CayleyC}(G; A', B', C')$ where $A' = x A y^{-1}$, $B' = y B z^{-1}$, and $C' = z C x^{-1}$.  So, we see that 
our notion of similarity for subset trios corresponds to switching ``base points'' and equivalence of incidence trios.

\subsection{An Example}

In this subsection, we offer an example of a natural (incidence) trio as motivation for the introduction of incidence geometry.  This will also permit us to introduce some of the exceptional structures which appear in our classification.  

\begin{figure}[ht]
\centerline{\includegraphics[width=12cm]{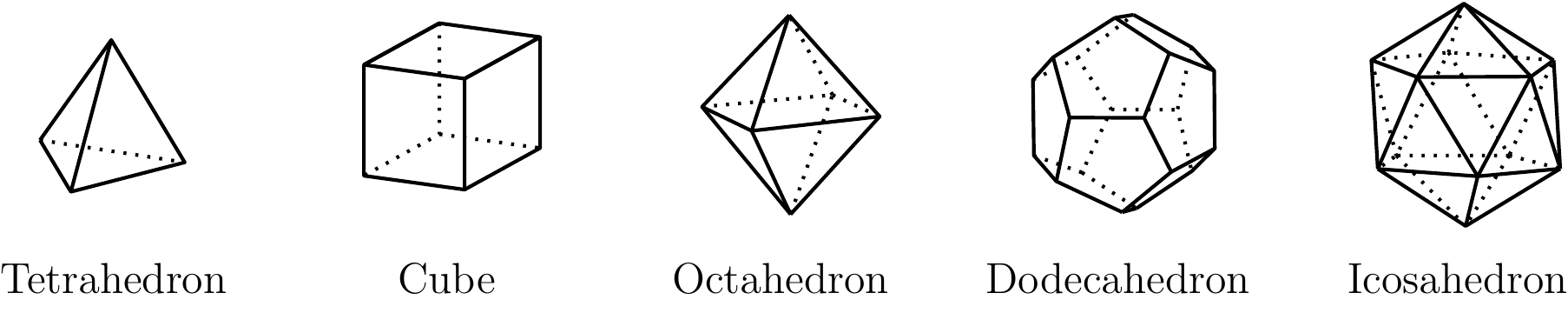}}
\caption{The Platonic Solids}
\label{platonic}
\end {figure}

Consider a map $\Xi$ on a surface with vertex set $V$, edge set $E$, and face set $F$ and assume that $\Xi$ is vertex, edge, and face transitive and that $V,E,F$ are all finite.

Now, define an incidence geometry $\Lambda$ on $(V,E,F)$ by the rule that an edge $e \in E$ is incident with an element in $V \cup F$ in $\Lambda$ if they are incident in $\Xi$, and $v \in V$ and $f \in F$ are incident in $\Lambda$ if they are not incident in $\Xi$.  It follows that $\Lambda$ is a $G$-chorus where $G$ is the automorphism group of $\Xi$.  Further, it is collision free by construction, so $\Lambda$ is a $G$-trio.  To determine if this trio is critical, we need to determine the densities of the sides.  First note that the density of $(V,E)$ is $\frac{2}{|V|}$ and the density of $(E,F)$ is $\frac{2}{|F|}$.  For the side $(V,F)$, observe that the number of vertex-face incidences in $\Xi$ is exactly $2 |E|$ since for every $v \in V$ and $f \in F$ which are incident, there are exactly two edges which are incident with both $v$ and $f$, and each edge contributes to four such incident pairs.  It follows from this that the sum of the densities of 
$(F,E)$, $(V,F)$, and $(V,E)$ is given by
\[ \tfrac{2}{|F|} + \left( 1 - \tfrac{2 |E|}{|V| \cdot |F|} \right) + \tfrac{2}{|V|} = 1 + \tfrac{2}{|V| \cdot |F|} \left( |V|  - |E| + |F| \right). \]
So, we see that this construction gives a critical trio if and only if the Euler characteristic is positive, which holds only when the surface is either a sphere or projective plane.  Applying this construction to any of the dual pairs Cube/Octahedron, Dodecahedron/Icosahedron, Petersen/$K_6$ results in an exceptional critical trio (for smaller maps the resulting trio is degenerate and falls into another class).  

\begin{figure}[ht]
\centerline{\includegraphics[width=3cm]{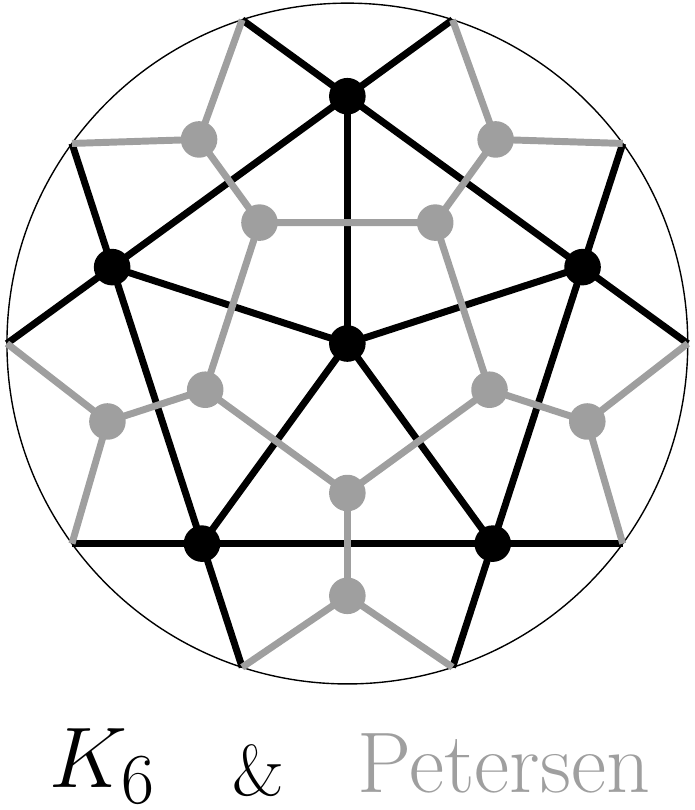}}
\caption{Two dual maps on the projective plane}
\label{k6pete}
\end {figure}

To see how these structures will appear in the context of subset trios, let us now assume that $(A,B,C)$ is a critical subset trio in $G$ and assume that the (incidence) trio 
$\Theta = {\mathit CayleyC}(G; A, B, C)$ has the property that $\Theta^{\bullet}$ is isomorphic to the exceptional structure just described for the cube.  Now our group $G$ 
must have a homomorphism to the automorphism group of the cube, so we have 
\[ G \rightarrow {\mathit Aut}(\Theta^{\bullet}) \cong {\mathbf C}_2 \times {\mathbf S}_4. \]
The kernel of this map is a normal subgroup $H \triangleleft G$ and the image of $G$ is a subgroup of ${\mathit Aut}(\Theta^{\bullet})$ which acts transitively on the vertices, edges, and faces of the cube.  The only proper subgroups which have this property are isomorphic to ${\mathbf C}_2 \times {\mathbf A}_4$ and ${\mathbf S}_4$, so $G$ must have a normal subgroup $H$ with $G/H$ isomorphic to one of ${\mathbf C}_2 \times {\mathbf S}_4$, ${\mathbf C}_2 \times {\mathbf A}_4$, or ${\mathbf S}_4$.  Furthermore, $A,B,C$ must all be $H$-stable, and with some further analysis, we can determine the precise structure of these sets.

\subsection{Weights and General Trios}
\label{weights}

So far we have only defined (incidence) trios in the special case when the group is finite.  In this subsection we will extend this  notion to allow for infinite groups.  To do this, we  will introduce weights, and insist that the geometries of interest to us all have either finite or cofinite weight.  In fact, these weights are nothing more than the sizes of sets involved were we to model our incidence geometry as a Cayley chorus using Theorem \ref{sabidussi}.

Let $G$ act transitively on the set $X$.  Then $|G_x| = |G_y|$ for every $x,y \in X$, and we call this number (or $\infty$) the \emph{point weight} of $X$ and denote it by $w^{\circ}_G(X)$.  Note that this group action is regular if and only if $w^{\circ}_G(X) = 1$.  For $A \subseteq X$ we define the \emph{weight} of $A$ to be $w_G(A) = w_G^{\circ}(X) |A|$ and the \emph{coweight} of $A$ to be $\overline{w}_G(A) = w_G^{\circ}(X) |X \setminus A|$.  When the group is clear from context, we shall drop these subscripts.  

Let $(X,Y)$ be a $G$-chorus.   Then $|N(x)|$ is independent of the choice of $x \in X$ and we call this number (or $\infty$) the \emph{degree} of $(X,Y)$ and denote it by $d(X,Y)$.  Note that for a graphic duet, this aligns with the standard definition of degree.  Similarly, $| Y \setminus N(x)|$ is independent of the choice of $x$ and we denote this quantity by $\overline{d}(X,Y)$.  We define $w_G(X,Y) = w_G(N(x)) = d(X,Y) w_G^{\circ}(Y)$ and $\overline{w}_G(X,Y) = \overline{w}_G(N(x)) = \overline{d}(X,Y) w_G^{\circ}(Y)$, and call these quantities the \emph{weight} and \emph{coweight} of $(X,Y)$.  As usual, when the group is clear from context, we drop the subscript. 

\bigskip

\noindent{\bf Definition:} A $G$-\emph{duet} is a $G$-chorus $\Delta = (X,Y)$ for which $w^{\circ}(X)$ and $w^{\circ}(Y)$ are finite and either $w(X,Y)$ or $\overline{w}(X,Y)$ is finite.

\bigskip

Next we establish an easy but essential property of duets.  

\begin{proposition}
\label{duet_identity}
For every $G$-duet $(X,Y)$ we have 
\begin{align*}
	w_G(X,Y) &= w_G(Y,X),	\\
	\overline{w}_G(X,Y) &= \overline{w}_G(Y,X).
\end{align*}
\end{proposition}

\noindent{\it Proof:} Apply Theorem \ref{sabidussi} to choose a Cayley chorus ${\mathit CayleyC}(G; A)$ with incidence relation $\sim$ and subgroups $H_1,H_2 \le G$ so that $G/H_1 \times \{1\}$ and $G/H_2 \times \{2\}$ are systems of clones, and $(X,Y)$ is strongly isomorphic to 
$(G/H_1 \times \{1\}, G/H_2 \times \{2\}; \sim)$.  Now we have $w(X,Y) = |A| = |A^{-1}| = w(Y,X)$ and 
$\overline{w}(X,Y) = \overline{A} = \overline{A^{-1}} = \overline{w}(Y,X)$. 
\quad\quad$\Box$

\bigskip

Next we generalize our earlier notion of finite (incidence) trios, to allow for infinite groups.  This will be our principle object of study.

\bigskip

\noindent{\bf Definition:} A $G$-\emph{trio} is a collision-free $G$-chorus for which each side is a duet.

\bigskip

Extending our definitions from the finite case, we will say that a $G$-trio is \emph{trivial} if it has a side which is empty and \emph{maximal} if there does not exist a $G$-trio $\Theta'$ with $\Theta < \Theta'$.  Next we establish the notion of deficiency for (incidence) trios.

\bigskip

\noindent{\bf Definition:} Let $\Theta = (X,Y,Z)$ be a $G$-chorus and assume (without loss) $w(X,Y) \le w(Y,Z) \le w(X,Z)$.  We define the \emph{deficiency} of 
$\Theta$ to be $\delta_G(\Theta) = w(X,Y) + w(Y,Z) - \overline{w}(X,Z)$, and we say that $\Theta$ is \emph{critical} if $\delta_G(\Theta)> 0$.  As usual, we drop the subscript when the group is clear from context.

\bigskip

Note that for $A,B,C \subseteq G$ we have that $(A,B,C)$ is a subset trio if and only if $\Theta = {\mathit CayleyC}(G; A,B,C)$ is an (incidence) trio, and $\delta(A,B,C) = \delta(\Theta)$.  As was the case with finite trios, the next observation shows that our notion of critical for a general trio depends only on the structure of the underlying incidence geometry and not on the particular group acting on it.  

\begin{observation}
If $(X,Y,Z)$ is a $G$-trio with $\overline{d}(X,Z) < \infty$ and $d(X,Y) > 0$, then $(X,Y,Z)$ is critical if and only if
\[ 1 > \frac{ \overline{d}(Z,X) }{ d(Y,X) } - \frac{ d(Z,Y) }{ d(X,Y) }. \]
\end{observation}

\noindent{\it Proof:}
This follows from the equation
\[ 1 + \frac{d(Z,Y)}{d(X,Y)} - \frac{ \overline{d}(Z,X)}{d(Y,X)}  
	= \frac{1}{ w(X,Y) }\left( w(X,Y) + w(Z,Y) - \overline{w}(Z,X) \right). \quad\quad\Box \]

\section{The Structure of Critical Trios}
\label{structure}

The main goal of this section is to state our structure theorem for critical (incidence) trios.  Although we have already stated our theorem for subset trios, 
the setting of incidence geometry provides us with a different, and sometimes more beneficial, perspective on these structures.

\subsection{Videos}
\label{video_subsec}

Consider a maximal subset trio $(A,B,C)$ for which $H = {\mathit stab}_R(A) = {\mathit stab}_L(B)$ satisfies $|AH| = 2|H|$ and define 
$\Theta  = {\mathit CayleyC}(G; A,B,C)$.  It follows from our assumptions that whenever $x,y \in G$ lie in the same left $H$-coset, the points $(x,2)$ and $(y,2)$  will be clones.  It then follows that in the trio $\Theta^{\bullet}$, the side corresponding to $(G \times \{1\},G \times \{2\})$ will be a graph.  So, to get a handle on subset trios of this form, we will investigate (incidence) trios which have a graphic side.  Conveniently, the next lemma shows that a natural graph parameter determines when such a trio is critical.  For a set of vertices $A$ in a graph $\Gamma$, we let $\partial(A)$ denote the set of all edges with exactly one end in $A$, and we call any set of this form an \emph{edge cut}.  We define $c(A) = |\partial(A)|$.  

\begin{observation}
\label{2d_obs}
Let $\Theta = (V,E,Z)$ be a nontrivial maximal trio for which $(V,E)$ is a graph of degree $d$.  If $z \in Z$ and $A = N(z) \cap V$, then $\Theta$ is critical if and only if $c(A) < 2d$.
\end{observation}

\noindent{\it Proof:}  Set $S = N(z) \cap E$.  First suppose that $\overline{w}(Z,V)$ is finite, and observe that by maximality, $S$ is precisely the set of edges with both endpoints in $V \setminus A$.  Now summing the degrees of the vertices in $V \setminus A$ counts every edge in $S$ twice and every edge in $\partial(A)$ once.  This gives us $d |V \setminus A| = 2 |S| + c(A)$.  So $\Theta$ is critical if and only if 
\[ 1 	>	\tfrac{\overline{d}(Z,V)}{d(E,V)} - \tfrac{d(Z,E)}{d(V,E)} 	
	=	 \tfrac{|V \setminus A|}{2}	-  \tfrac{|S|}{d}	
	=	\tfrac{c(A)}{2d} \]
which is precisely the desired condition.

Next consider the case when $\overline{w}(Z,E)$ is finite and set $R = E \setminus (S \cup \partial(A))$.  It follows from maximality that $R$ is precisely the set of edges with both ends in $A$.  Summing the degrees of the vertices in $A$ counts every edge in $R$ twice and every edge in $\partial(A)$ once, so $d |A| = 2 |R| + c(A)$.  Now, $\Theta$ will be critical if and only if 
\[ 1 	>	\tfrac{\overline{d}(Z,E)}{d(V,E)} - \tfrac{d(Z,V)}{d(E,V)} 	
	= 	\tfrac{|R| + c(A)}{d} - \tfrac{|A|}{2}
	=	\tfrac{c(A)}{2d} \]
which is again the desired condition.  Since one of $\overline{w}(Z,V)$ or $\overline{w}(Z,E)$ must be finite, this completes the proof.
\quad\quad$\Box$

\bigskip

So, to understand critical trios which have a graphic side, we will need to classify edge-cuts of size at 
most $2d$ in vertex and edge-transitive graphs of degree $d$.  This will be a major focus of Section \ref{graph_sec}.  
For now we will introduce some notation to help us describe the trios which result from this construction.

Let $\Gamma$ be a graph, let $m$ be a positive integer, and define 
$C_m(\Gamma)$ ($P_m(\Gamma)$) to be the set of all cycles (paths) of $\Gamma$ 
with exactly $m$ vertices.  Note that if $\Gamma$ is a $G$-duet, then $G$ has a natural action on both 
$C_m(\Gamma)$ and $P_m(\Gamma)$.  If $\Gamma'$ is a subgraph of $\Gamma$, 
$v \in V(\Gamma)$, and $e \in E(\Gamma)$ then we define
$v \sim \Gamma'$ if $v \in V(\Gamma')$ and $e \sim \Gamma'$ if $e \in E(\Gamma')$.  
Now, for any graph $\Gamma$ and positive integer $m$, we have operators $V$, $E$, 
$P_m$, and $C_m$ which take $\Gamma$ and return a set.  Using these operators 
together with the various relations $\sim$ defined here permits the construction of numerous
incidence geometries based on $\Gamma$.  If $A,B,C$ are three such operators and we have 
define a relation $\sim$ between $A$ and $B$ and between $B$ and $C$ then we let 
$\Gamma : A \sim B \sim C$ denote the incidence geometry $( A(\Gamma), B(\Gamma), C(\Gamma) )$ 
with incidence relation given by the rule that $y \in B(\Gamma)$ is incident 
with another point $w$ of different type if $w \sim y$, and two points $x \in A(\Gamma)$ and 
$z \in C(\Gamma)$ are incident if and only if there does not exist $y \in B(\Gamma)$
with $x \sim y \sim z$.  It follows from this that $\Gamma : A \sim B \sim C$
is always collision free.  We shall permit $A$ and $C$ to be the same operator and in this case for the 
incidence geometry we just use two different copies of the set $A(\Gamma) = C(\Gamma)$.  
With this, we can introduce three families of critical incidence trios.

\bigskip

\noindent{\bf Definition:} 
A trio $\Theta \equiv (X,Y,Z)$ is a {\it standard video} if there exists a graphic duet $\Gamma$ of degree $d$ so that $(X,Y,Z)^{\bullet}$ is 
isomorphic to one of the following.
\begin{enumerate}
\item $\Gamma : E \sim V \sim E$ where $d \ge 3$, and $|V(\Gamma)| \ge 5$
\item $\Gamma : P_3 \sim V \sim E$ where $d = 3$, and $|V(\Gamma)| \ge 6$
\item $\Gamma : P_3 \sim E \sim V$ where $d=3$, and $|V(\Gamma)| \ge 6$.
\end{enumerate}

As indicated earlier, the six graphs Cube, Octahedron, Dodecahedron, Icosahedron, Petersen, and $K_6$
will play a special role in our theorem.  These graphs give rise to exceptional critical trios because they possess exceptional small edge cuts.  
The first four of these graphs embed in the sphere and the last two embed naturally in the projective plane 
(as can be seen in Figures \ref{platonic} and \ref{k6pete}).   
For any one of these graphs $G$, we define $F(G)$, to be the 
{\it faces} of $G$ in this embedding and we use $\sim$ to denote
the incidence relation between vertices and faces and between
edges and faces.  We treat the operator $F$ as $V,E,P_m,C_m$ above.
With this, we can introduce the exceptional videos.  

\bigskip

\noindent{\bf Definition:} 
A trio $\Theta \equiv (X,Y,Z)$ is an {\it exceptional video} if $(X,Y,Z)^{\bullet}$ is 
isomorphic to one of the following.
\begin{enumerate}
\item Cube/Octahedron : $V \sim E \sim F$	
\item Dodecahedron : $F \sim V \sim E$
\item Dodecahedron/Icosahedron : $V \sim E \sim F$
\item Icosahedron : $F \sim V \sim E$	
\item Petersen : $C_5 \sim E \sim V$	
\item Petersen/$K_6$ : $F \sim E \sim V$
\item $K_6$ : $C_3 \sim E \sim V$.
\end{enumerate}

\bigskip

\noindent{\bf Definition:} A trio is a \emph{video} if it is either a standard video or an exceptional video.  

\bigskip

\subsection{Beats and Chords}

In this subsection we will introduce the remaining types of critical (incidence) trios.  We have already encountered these structures in the context of subset trios, and our names for them will be similar to those we used in that setting.  

We say that an incidence geometry $(X,Y)$ is a \emph{matching} if $|X| \ge 2$ and there is a bijection $f : X \rightarrow Y$ so that for every $x \in X$ and $y \in Y$ we have $x \sim y$ if and only if $f(x) = y$.  

\bigskip

\noindent{\bf Definition:} We define a trio $\Theta \equiv (X,Y,Z)$ to be a \emph{pure beat relative to} $(X,Y)$ if $(X,Y)^{\bullet}$ is a matching, and $\Theta$ is maximal and nontrivial.

\bigskip

Consider a pure beat as in the above definition for which $\Delta = (X,Y)$ and $\Theta \equiv (X,Y,Z)$, and let $x \in X$ and $y \in Y$ satisfy 
$x \sim y$.  It then follows from the structure of $\Delta$ and the maximality of $\Theta$ that every point $z \in Z$ will satisfy either $z \sim x$ or 
$z \sim y$.  It follows from this that $\Theta$ is critical and $\delta(\Theta) = w(X,Y)$.  

\bigskip

\noindent{\bf Definition:} A trio $\Theta \equiv (X,Y,Z)$ is a \emph{pure chord} if there exist $\ell,m \ge 1$ and a group $J = {\mathbb Z}$ or $J = {\mathbb Z}/n {\mathbb Z}$ with $|J| \ge \ell+m+3$ so that  $(X,Y,Z)^{\bullet}$ is isomorphic to 
\[ {\mathit Cayley}( J ; \{0,\ldots,\ell\}, \{0,\ldots,m\}, 
	J \setminus \{-\ell-m,\ldots,0\} ). \]

It is straightforward to verify that pure chords are maximal and critical.  Furthermore, by the bounds on the parameters, pure chords will also 
be connected.

We say that an incidence geometry $\Lambda$ is \emph{empty} if no two points in $\Lambda$ are incident, and we say that $\Lambda$ is \emph{full} if 
$x \sim y$ whenever $x$ and $y$ are points of different type.  Finally, we call $\Lambda$ \emph{partial} if it is neither empty nor full.  

If $\Delta = (X,Y)$ is a $G$-duet which is nonempty and disconnected, then there exist systems of imprimitivity $\mcp = \{ P_i \mid i \in I \}$ of $X$ and $\mcq = \{ Q_i \mid i \in I \}$ of $Y$ so that $\{ (P_i,Q_i ) \mid i \in I \}$ is the set of connected components of $\Delta$, and we call $( \mcp, \mcq )$ the 
\emph{component quotient} of $(X,Y)$.  Note that $( \mcp, \mcq )$ is a matching which is a $G$-chorus, but $(\mcp,\mcq)$ will be a $G$-duet if 
and only if $w^{\circ}(\mcp) = w^{\circ}(\mcq) < \infty$.

\bigskip

\begin{figure}[ht]
\centerline{\includegraphics[height=4cm]{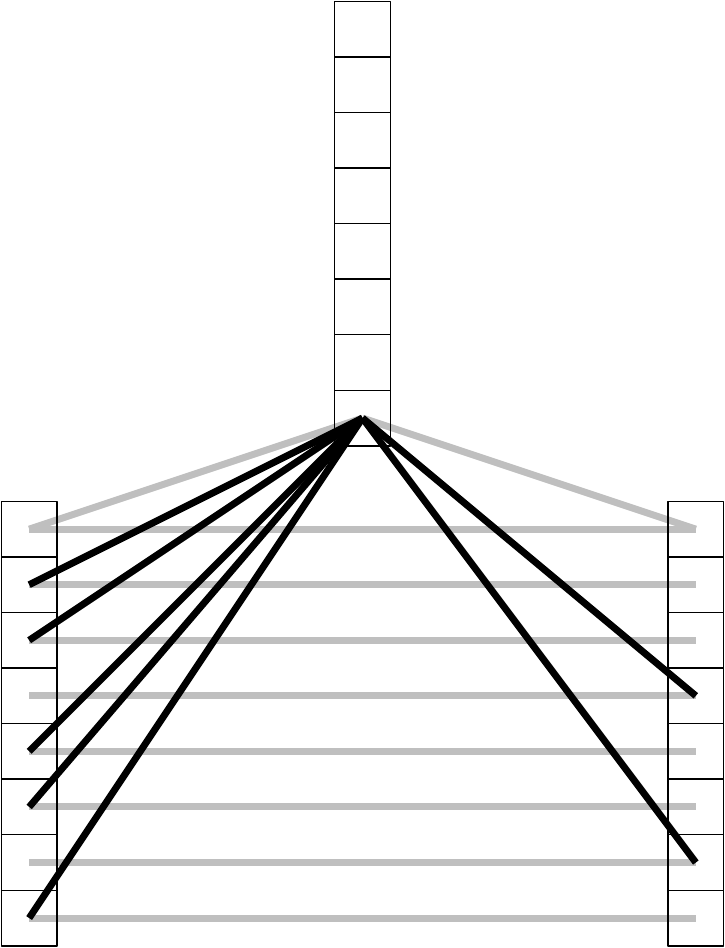}}
\caption{An impure beat}
\label{beat_fig}
\end {figure}

\noindent{\bf Definition:} We define a trio $\Theta \equiv (X,Y,Z)$ to be an \emph{impure beat relative to} $(X,Y)$ if 
$(X,Y)$ is disconnected with component quotient $(\mcp, \mcq$), and for every $z \in Z$ we have
\begin{enumerate}
\item there is exactly one $P_z \in \mcp$ with $(z,P_z)$ partial, exactly one $Q_z \in \mcq$ with $(z,Q_z)$ partial, and these blocks satisfy $P_z \sim Q_z$.
\item for every $P \in \mcp \setminus \{P_z\}$ and $Q \in \mcq \setminus \{Q_z\}$ with $P \sim Q$, one of $(z,P)$, $(z,Q)$ is full and the other is empty.
\end{enumerate}

Let $P \in {\mathcal P}$ and $Q \in {\mathcal Q}$ satisfy $P \sim Q$ and define 
$R = \{ z \in Z \mid \mbox{$( \{z\}, P )$ is partial} \}$.  It follows immediately that $R$ is a block of imprimitivity, 
and we say that the $\mcp, \mcq$, and the system of imprimitivity generated by $R$ are \emph{associated} with this impure beat.  
We call $(P,Q,R)$ a \emph{continuation} of $\Theta$, and note that any two such continuations are isomorphic.
Setting $H = G_P = G_Q = G_R$ we have that $(P,Q,R)$ is an $H$-trio and $\delta_H(P,Q,R) = \delta_G(\Theta)$.  

\bigskip

\begin{figure}[ht]
\centerline{\includegraphics[height=4cm]{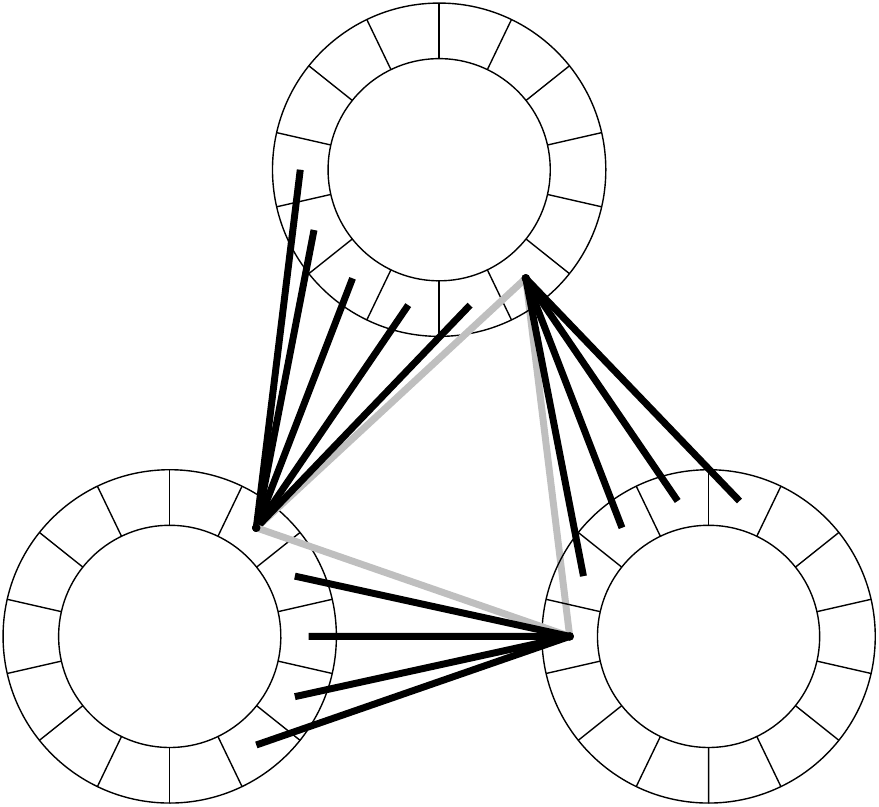}}
\caption{A cyclic chord}
\label{cyc_fig}
\end {figure}

\noindent{\bf Definition:} We say that $\Theta \equiv (X,Y,Z)$ is an (\emph{impure}) \emph{cyclic chord} if there exist $\ell,m \ge 1$, a group $J = {\mathbb Z}$ or $J = {\mathbb Z}/n{\mathbb Z}$ 
with $|J| \ge \ell + m +2$ and systems of imprimitivity ${\mathcal P} = \{ P_j : j \in J \}$, ${\mathcal Q} = \{ Q_j : j \in J \}$, ${\mathcal R} = \{ R_j : j \in J \}$ of
$X,Y,Z$ so that the following hold
\begin{align*}
(P_i,Q_j) \mbox{ is } &\left\{\begin{array}{ll}
		\mbox{partial}		&	\mbox{ if $j=i$}	\\
		\mbox{full}	&	\mbox{ if $i+1 \le j \le i+\ell$}	\\
		\mbox{empty}		&	\mbox{ otherwise}
		\end{array} \right. \\
(Q_i,R_j) \mbox{ is } &\left\{\begin{array}{ll}
		\mbox{partial}		&	\mbox{ if $j=i$}	\\
		\mbox{full}	&	\mbox{ if $i+1 \le j \le i+m$}	\\
		\mbox{empty}		&	\mbox{ otherwise}
		\end{array} \right. 	\\
(P_i,R_j) \mbox{ is } &\left\{\begin{array}{ll}
		\mbox{partial}		&	\mbox{ if $j=i$}	\\
		\mbox{empty}		&	\mbox{ if $i+1 \le j \le i+\ell + m$}	\\
		\mbox{full}	&	\mbox{ otherwise}
		\end{array} \right. 
\end{align*}
In this case we will say that $\Theta$ is a cyclic chord \emph{relative to} $\mcp, \mcq, \mcr$.  
It follows immediately that the subgroup $H = G_{ (\mcp) } = G_{ (\mcq) } = G_{ (\mcr )}$ has $G/H$ cyclic, and furthermore 
$H$ is the stabilizer of each block in $\mcp \cup \mcq \cup \mcr$.  
Let $P \in {\mathcal P}$, $Q \in {\mathcal Q}$, and $R \in {\mathcal R}$ satisfy $(P,Q)$, $(Q,R)$, and $(R,P)$ partial.  Then
$(P,Q,R)$ is an $H$-trio, and we call $(P,Q,R)$ a \emph{continuation} of $\Theta$.  It follows easily 
that any two continuations are isomorphic, and as usual, we have $\delta_H(P,Q,R) = \delta_G(\Theta) $.  

We need one added bit of notation before introducing our next structure.  If a group $G$ acts on a set $X$, we say the action is \emph{dihedral} if $G$ is dihedral, and the rotation subgroup of $G$ acts regularly on $X$.  So, the standard action of the dihedral group ${\mathbf D}_n$ on the vertices of a regular $n$-gon is dihedral.  

\bigskip

\begin{figure}[ht]
\centerline{\includegraphics[height=4cm]{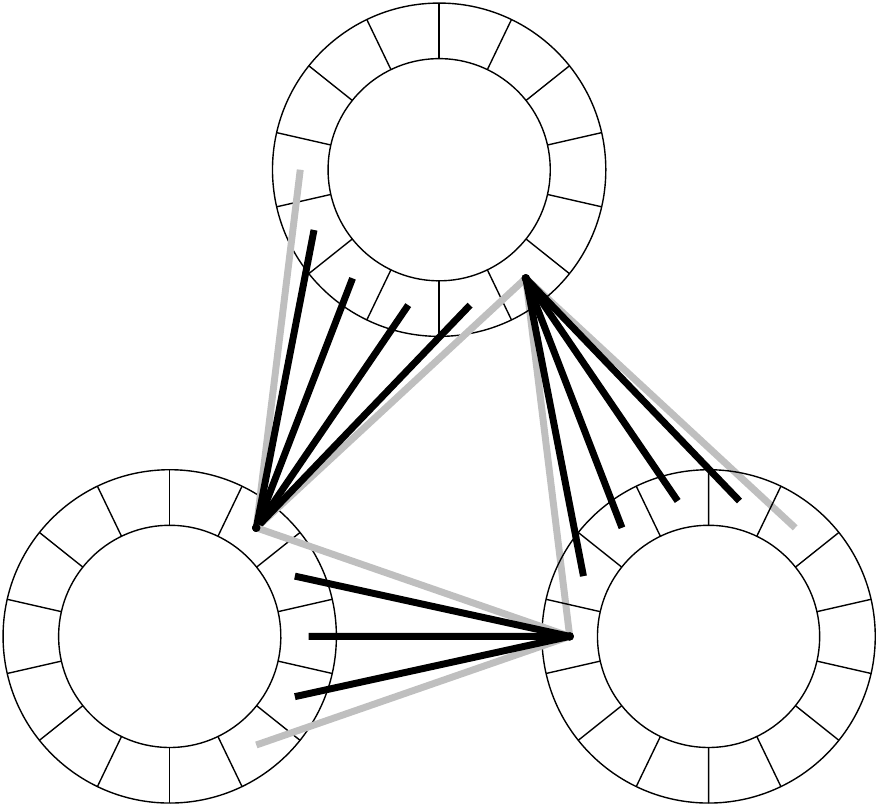}}
\caption{A  dihedral chord}
\label{dih_fig}
\end {figure}

\noindent{\bf Definition:} We define a connected trio $\Theta \equiv (X,Y,Z)$ to be an (\emph{impure}) \emph{dihedral chord} 
if there exist $\ell,m \ge 1$, a group $J = {\mathbb Z}$ or $J = {\mathbb Z}/n{\mathbb Z}$
with $|J| \ge \ell+ m + 2$, and systems of imprimitivity ${\mathcal P} = \{ P_j : j \in J \}$, ${\mathcal Q} = \{ Q_j : j \in J \}$, 
${\mathcal R} = \{ R_j : j \in J \}$ of $X,Y,Z$ so that
\begin{align*}
(P_i,Q_j) \mbox{ is } &\left\{\begin{array}{ll}
		\mbox{partial}		&	\mbox{ if $j=i$ or $j=i+\ell$}	\\
		\mbox{full}	&	\mbox{ if $i+1 \le j \le i+\ell-1$}	\\
		\mbox{empty}		&	\mbox{ otherwise}
		\end{array} \right. \\
(Q_i,R_j) \mbox{ is } &\left\{\begin{array}{ll}
		\mbox{partial}		&	\mbox{ if $j=i$ or $j=i+m$}	\\
		\mbox{full}	&	\mbox{ if $i+1 \le j \le i+m-1$}	\\
		\mbox{empty}		&	\mbox{ otherwise}
		\end{array} \right. 	\\
(P_i,R_j) \mbox{ is } &\left\{\begin{array}{ll}
		\mbox{partial}		&	\mbox{ if $j=i$ or $j=i + \ell + m$}	\\
		\mbox{empty}		&	\mbox{ if $i+1 \le j \le i+ \ell + m-1$}	\\
		\mbox{full}	&	\mbox{ otherwise}
		\end{array} \right. 
\end{align*}
and in addition, for the normal subgroup $H = G_{ ( {\mathcal P} ) } = G_{ ( {\mathcal Q} ) } = G_{ ( {\mathcal R} ) }$, the action of $G/H$ 
on each of ${\mathcal P}$, ${\mathcal Q}$, ${\mathcal R}$ is dihedral.  We will also say that $\Theta$ is a dihedral chord \emph{relative to}
$\mcp, \mcq, \mcr$.  

In the next two subsections we will develop the notion of a continuation of a dihedral chord.  However, we are now  
able to state a significant corollary of our structure theorem, so we end the present subsection
with this.

\begin{corollary}
\label{one_step_cor}
Every nontrivial maximal critical trio is one of the following: a pure beat, a pure chord, a 
video, an impure beat, a cyclic chord, or a dihedral chord.
\end{corollary}

\subsection{Square and Octahedral Choruses}

In this subsection we will develop a framework to investigate the structure of (impure) dihedral chords.  We begin by defining 
a type of duet which appears as a side of a dihedral chord of finite weight.

\bigskip
\begin{figure}[ht]
\centerline{\includegraphics[height=2cm]{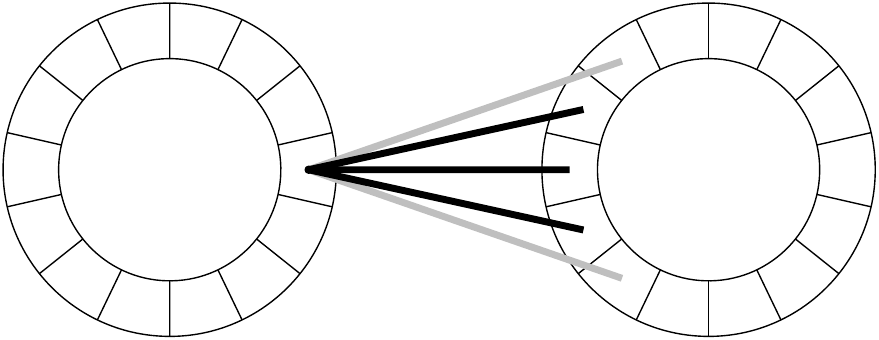}}
\caption{A fringed sequence}
\label{fringed_fig}
\end {figure}

\noindent{\bf Definition:} For a positive integer $\ell$, we say that a duet $(X,Y)$ is a \emph{fringed sequence} of \emph{length} $\ell+1$ if there exists a group $J = {\mathbb Z}$ or $J = {\mathbb Z}/n {\mathbb Z}$ with $n \ge \min\{4, \ell+2\}$, 
and systems of imprimitivity ${\mathcal P} = \{ P_j \mid j \in J \}$ of $X$ and 
${\mathcal Q} = \{ Q_j \mid j \in J \}$ of $Y$ so that 
\[ (P_i,Q_j) \mbox{ is } \left\{\begin{array}{ll}
		\mbox{partial}		&	\mbox{ if $j=i$ or $j=i+\ell$}	\\
		\mbox{full}	&	\mbox{ if $i+1 \le j \le i+\ell-1$}	\\
		\mbox{empty}		&	\mbox{ otherwise}
		\end{array}
		\right. \]
and in addition $H = G_{ ( {\mathcal P} ) } = G_{ ( {\mathcal Q} ) }$ has the property that the action of $G/H$ on each of 
${\mathcal P}$ and ${\mathcal Q}$ is dihedral.  In this case we will also say that $(X,Y)$ is a fringed sequence \emph{relative to} $(\mcp, \mcq)$.

\bigskip

Let us pause to note a couple of consequences of the dihedral action of $G/H$ on ${\mathcal P}$ and ${\mathcal Q}$ (these will be discussed in detail in Section \ref{sequence_sec}).  First, it follows from this that whenever $P,P' \in {\mathcal P}$ and $Q, Q' \in {\mathcal Q}$ satisfy $(P,Q)$ and $(P',Q')$ partial, there exists $g \in G$ for which $gP = P'$ and $gQ = Q'$.  It also follows from this dihedral action that the stabilizers of $P$ and $Q$ satisfy $[G_P : H] = 2$ and $[G_Q : H] = 2$.  So, $P$ ($Q$) will consist of either one or two orbits under the action of $H$.  Next we will introduce a type of chorus equipped with some additional structure to work with this situation.

\bigskip

\noindent{\bf Definition:} For a finite group $H$, an $H$-\emph{square chorus} is an 
$H$-chorus $\Lambda = (X_1,\ldots,X_n)$ equipped with a partition $\{ \Lambda_1, \Lambda_2 \}$ of $\{ X_1, \ldots, X_n \}$ so that $|\Lambda_i| = 1,2$ for $i = 1,2$ and so that $(X,Y)$ is empty whenever $X,Y \in \Lambda_i$ are distinct.

\bigskip

Continuing with $\mcp, \mcq$ from above, let $P \in \mcp$ and $Q \in \mcq$ satisfy $(P,Q)$ partial and define $\Lambda_1 = \{ X_1, \ldots, X_m \}$ to be the set of $H$-orbits in $P$ (so $m=1$ or $m=2$) and $\Lambda_2 = \{ X_{m+1}, \ldots X_n \}$ to be the set of $H$-orbits in $Q$.  Now $\Lambda = (X_1, \ldots, X_n)$ is an $H$-chorus and together with $\{ \Lambda_1, \Lambda_2\}$ it forms an $H$-square chorus.  We will say that $\Lambda$ is 
\emph{associated with} $(X,Y)$.  

We will need to relate the weight of $(X,Y)$ to the weights of the sides of $\Lambda$.  To do so, put  
$m_i = 3 - |\Lambda_i|$ for $i=1,2$ and define the \emph{adjusted weight} of $\Lambda$ to be 
\[ \hat{w}_H ( \Lambda ) = \sum_{S \in \Lambda_1, T \in \Lambda_2} m_1 m_2 w_H(S,T) \]

\begin{proposition}
\label{fringe_sq_weight}
Let $(X,Y)$ be a fringed sequence of length $\ell+1$ and let $\Lambda$ be an associated $H$-square  chorus.  Then
\[ w(X,Y) = 2(\ell-1)|H| + \hat{w}_H( \Lambda ). \]
\end{proposition}

\noindent{\it Proof:} 
Let ${\mathcal P}$ and ${\mathcal Q}$ and $P,Q$ and $\Lambda_1, \Lambda_2$ be as in the above discussion
and let $Q' \in \mcq \setminus \{Q\}$ satisfy $(P,Q')$ partial.  
Choose $x \in P$ and if $|\Lambda_1| = 2$ choose $x' \in P$ so that $x$ and $x'$ are in distinct $H$-orbits.  Now 
$w^{\circ}({\mathcal P}) = 2|H| = w^{\circ}({\mathcal Q})$ so we have
\[ w(X,Y)  = w(N(x)) = 2(\ell-1)|H| + w( N(x) \cap (Q \cup Q')). \]
The following properties (all of which follow immediately) permit us to express $w( N(x) \cap (Q \cup Q'))$ in terms of weights of sides of $\Lambda$.
\begin{itemize}
\item	If $|\Lambda_1| = 1$ then $w_G(N(x) \cap (Q \cup Q')) = 2 w_G(N(x) \cap Q)$.
\item If $|\Lambda_1| = 2$ then $w_G(N(x) \cap (Q \cup Q')) = w_G(N(x) \cap Q) + w_G(N(x') \cap Q)$.
\item If $\Lambda_2 = \{Q\}$ then for every $D \subseteq Q$ we have $w_G(D) = 2 w_H(D)$.
\item If $\Lambda_2 = \{Q_1,Q_2\}$ then for $i=1,2$ and every $D \subseteq Q_i$ we have $w_G(D) = w_H(D)$.
\end{itemize}
For instance, if $\Lambda_1 = \{P_1,P_2\}$ and $\Lambda_2 = \{Q\}$ then
\begin{align*}
w_G(N(x) \cap (Q \cup Q'))	
	&=	w_G( N(x) \cap Q ) + w_G( N(x') \cap Q )	\\
	&=	2w_H(N(x) \cap Q) + 2 w_H( N(x') \cap Q)		\\
	&= 	2 w_H(P_1,Q) + 2 w_H(P_2,Q)		\\
	&=	\hat{w}_H(\Lambda).
\end{align*}
A similar analysis in the remaining cases reveals that 
$w_G(N(x) \cap (Q \cup Q')) = \hat{w}(\Lambda)$
holds in general, and this completes the proof.
\quad\quad$\Box$

\bigskip

All of our analysis extends naturally from fringed sequences to dihedral chords, but for this, we need to step up from square choruses to the following.

\bigskip

\noindent{\bf Definition:} For a finite group $H$, an $H$-\emph{octahedral chorus} is a collision free
$H$-chorus $\Lambda = (X_1,\ldots,X_n)$ equipped with a partition $\{ \Lambda_1, \Lambda_2, \Lambda_3 \}$ of $\{ X_1, \ldots, X_n \}$ 
so that $|\Lambda_i| = 1,2$ for $i = 1,2,3$ and so $(X,Y)$ is empty whenever $X,Y \in \Lambda_i$ are distinct.
We set  $m_i = 3 - |\Lambda_i|$ for $i=1,2,3$ and define the \emph{adjusted weight} of $\Lambda$ to be 
\[ \hat{w}_H ( \Lambda ) = \sum_{1 \le i < j \le 3} \; \sum_{S \in \Lambda_i, T \in \Lambda_j} m_i m_j w_H(S,T). \]

\bigskip

If $\Lambda = (X_1, \ldots, X_n)$ together with the partition $\{ \Lambda_1, \Lambda_2, \Lambda_3 \}$ is an octahedral chorus, then
we may obtain a square chorus from $\Lambda$ by choosing $1 \le i \le 3$ and removing all members of $\Lambda_i$ from the sequence of sets, 
and using the partition $\{ \Lambda_1, \Lambda_2, \Lambda_3 \} \setminus \Lambda_i$.  We say that such a square chorus is \emph{contained} 
in $\Lambda$.

Now let us return to the dihedral chord $\Theta$ defined in the previous subsection, and choose $P \in {\mathcal P}$ and 
$Q \in {\mathcal Q}$ and $R \in {\mathcal R}$ so that $(P,Q)$ and $(Q,R)$ and $(P,R)$ are partial.  Then define 
$\Lambda_1, \Lambda_2, \Lambda_3$ to be the respective partitions of $P,Q,R$ into $H$-orbits.  These orbits 
together with the partition $\{\Lambda_1,\Lambda_2,\Lambda_3\}$ form an octahedral chorus, $\Lambda$ which we say is
\emph{associated with} $\Theta$.  Let $\Lambda_{12}$, $\Lambda_{23}$, and $\Lambda_{13}$ be the three square choruses
contained in $\Lambda$.  Using our earlier analysis for fringed sequences (and noting $w^{\circ}({\mathcal P}) = w^{\circ}({\mathcal Q}) = w^{\circ}({\mathcal R}) = 2|H|$)
we have
\begin{align*}
\delta(\Theta) 
	&= w_G(X,Y) + w_G(Y,Z) - \overline{w}_G(X,Z) 	\\
	&= (\ell-1) 2 |H| + \hat{w}_H(\Lambda_{12}) + (m-1) 2 |H| + \hat{w}_H(\Lambda_{23}) 
			- (\ell+m+1) 2 |H| + \hat{w}_H(\Lambda_{13})		\\
	&= \hat{w}_H(\Lambda) - 6 |H|.
\end{align*}
Accordingly, we now define the \emph{deficiency} of an octahedral chorus $\Lambda$ to be $\delta(\Lambda) = \hat{w}_H(\Lambda) - 6 |H|$ and call $\Lambda$ \emph{critical} if it has positive deficiency.  The following observation records the obvious purpose for this definition.

\begin{observation}
If $\Theta$ is a dihedral chord and $\Lambda$ is an associated octahedral chorus 
\[ \delta(\Theta) = \delta(\Lambda). \]
\end{observation}

\subsection{Critical Octahedral Choruses}

It follows from the developments of the previous subsection that every maximal critical trio which is a dihedral chord has an associated octahedral chorus which will be maximal and critical.  The goal of this subsection is to introduce terminology to describe all such octahedral choruses.  

As we did with octahedral configurations in Section \ref{dihedral_sec}, we shall associate each octahedral chorus $\Lambda = (X_1, \ldots, X_n)$ 
with a graph $\Gamma$ with three types of edges.  The vertex set of $\Gamma$ is $\{X_1, \ldots, X_n\}$, and if $X_i,X_j$ are in distinct members of the distinguished partition, there is an edge between them in $\Gamma$.  This edge will be dotted if $(X_i,X_j)$ is empty, grey if $(X_i,X_j)$ is partial, and solid black if $(X_i,X_j)$ is full.  

\begin{figure}[ht]
\centerline{\includegraphics[width=12cm]{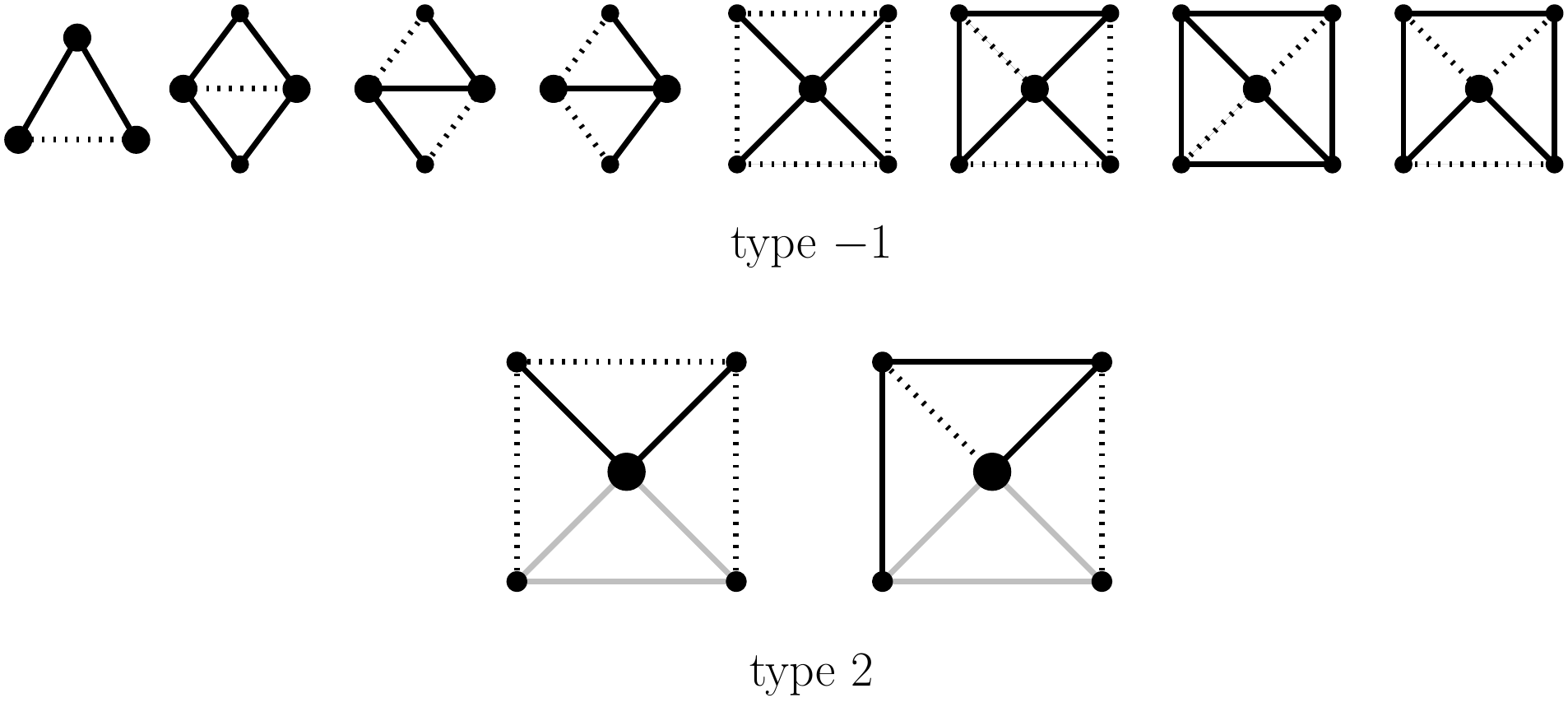}}
\caption{Octahedral choruses of low rank}
\label{deg_oct}
\end {figure}

\bigskip

\noindent{\bf Definition:} We define an octahedral chorus to be of \emph{type} $-1$, $0$, $1$, or $2$ if its associated graph appears in either Figure \ref{all_oct} or in Figure \ref{deg_oct} with the corresponding label.  

\bigskip

In Section \ref{oct_sec} we will prove the following theorem, which classifies maximal critical octahedral choruses.  

\begin{theorem}[Classification of critical octahedral choruses]
\label{oct_type}
Every maximal critical octahedral chorus has type $-1$, $0$, $1$, or $2$.
\end{theorem}

We will say that a dihedral chord has \emph{type} $n$ if it has an associated octahedral chorus of type $n$.  We will prove that dihedral chords 
of type $-1$ cannot be maximal trios, so they will not appear in our structure theorem.  We will prove that a dihedral chord of type $0$ is also a pure chord (with a special type of group action), and hence these structures will appear elsewhere in our structure theorem.  Next consider a $G$-trio $\Theta$ which is a dihedral chord of type $1$, let $\Lambda$ be an associated $H$-octahedral chorus, and assume that $(X,Y)$, $(Y,Z)$, and $(X,Z)$ are distinct sides of $\Lambda$ which are partial.  Then $(X,Y,Z)$ is an $H$-trio for which $\delta_H(X,Y,Z) = \delta_H(\Lambda) = \delta_G(\Theta)$ and we call $(X,Y,Z)$ a \emph{continuation} of $\Theta$.  Maximal critical trios which are type 2 have additional structure, and we now introduce some definitions to capture this.  For 
a symmetric relation $\sim$ on $X$ and $Y \subseteq X$ we let $\sim|_Y$ denote the restriction of $\sim$ to $Y$.  

\bigskip

\noindent{\bf Definition:} Let $\Lambda$ be an $H$-octahedral chorus of type 2 with incidence relation $\sim$ and first assume that $(X,Y)$, $(Y,Z_1)$, $(X,Z_1)$, $(X,Z_2)$, $(Y,Z_2)$ are distinct sides of $\Lambda$ which are partial (so the associated graph is one from Figure \ref{all_oct}).  For $i=1,2$ let $\Upsilon_i = (X,Y,Z_i ; \sim_i)$ be an $H$-trio, and assume that $\sim|_{X \cup Z_i} =\; \sim_i|_{X \cup Z_i}$ and  $\sim|_{Y \cup Z_i} =\; \sim_i|_{Y \cup Z_i}$ for $i = 1,2$.
\begin{enumerate}
\item[(2A)] We say that $\Lambda$ has \emph{type} $2A$ if $\sim|_{X \cup Y} = \; \sim_1|_{X \cup Y} \subseteq \; \sim_2|_{X \cup Y}$ and $\Upsilon_1$ is an impure beat relative to $(X,Y)$ and $\Upsilon_2$ is a pure beat relative to $(X,Y)$.  In this case we define a \emph{continuation} of $\Lambda$ to be a continuation of $\Upsilon_1$.  
\item[(2B)] We say that $\Lambda$ has \emph{type} $2B$ if $\sim|_{X \cup Y} = \; \big( \sim_1|_{X \cup Y} \big) \cap \big( \sim_2|_{X \cup Y} \big)$ and both $\Upsilon_1$ and $\Upsilon_2$ 
are pure beats relative to $(X,Y)$.
\end{enumerate}
Next assume that $\Lambda$ is an $H$-octahedral chorus with incidence relation $\sim$ and assume that $(X,Y)$, $(Y,Z)$, $(X,Z)$ are sides of $\Lambda$ which are partial and that $\{Z\}$ is a member of the distinguished partition (so the associated graph is one from Figure \ref{deg_oct}).  
\begin{enumerate}
\item[(2C)] We say that $\Lambda$ has \emph{type} $2C$ if $(X,Y,Z)$ is a pure beat relative to $(X,Y)$.
\end{enumerate}

\bigskip

As before, we will say that a dihedral chord has \emph{type} $2A$, $2B$, or $2C$ if it has an associated octahedral chorus of this type.  Note that 
a dihedral chord of type $2B$ or $2C$ is a fully described critical trio.  For a dihedral chord $\Theta$ of type $2A$, we define a \emph{continuation} of 
$\Theta$ to be a continuation of its associated octahedral chorus.  If $\Theta'$ is such a continuation, then $\delta(\Theta') = \delta(\Theta)$ and 
$\Theta'$ will be maximal whenever $\Theta$ is maximal.  Our next result gives us a classification of the type $2$ choruses of interest.

\begin{theorem}[Classification of type 2 octahedral choruses]
\label{classify_oct2}
Every maximal critical octahedral chorus of type $2$ has type $2A$, $2B$, or $2C$.
\end{theorem}

In summary, every maximal critical trio which is a dihedral chord has type $0$, $1$, $2A$, $2B$, or $2C$.  Those of type $0$ appear elsewhere in our classification theorem.  Dihedral chords of type $2B$, or $2C$ are fully described critical trios.  Finally, dihedral chords of types $1$ and $2A$ have continuations which are maximal critical trios in a proper (finite) subgroup.

\subsection{The Classification Theorem}

In this subsection, we will state our structure theorem for critical trios.  We begin by naming the trios of interest.

\bigskip

\noindent{\bf Definition:}  A trio $\Theta$ is a \emph{song}\footnote{To help explain our choice of terminology, observe that a song is recursively composed of beats and chords, possibly ending with a video} if there exists a sequence of trios 
$\Theta = \Theta_1, \ldots, \Theta_m$ so that the following hold.
\begin{enumerate}
\item 	For $1 \le i < m$ the trio $\Theta_i$ is an impure beat, a cyclic chord, or a dihedral 
		chord of type $1$ or $2A$, and has continuation $\Theta_{i+1}$.
\item 	$\Theta_m$ is a video, a pure beat, a pure chord, or a dihedral chord of type $2B$ or $2C$.
\end{enumerate}

We are now ready to state our main theorem.

\begin{theorem}
\label{main_incidence}
A nontrivial trio is maximal and critical if and only if it is a song.
\end{theorem}

We will give an elementary proof of this theorem which is entirely contained in this article.  In the next section we will develop
some basic bounds and tools for working with critical trios, and in Section \ref{proof_sec} we will give a proof of the ``only if'' direction
of this theorem under the assumption of a number of lemmas.  These lemmas will then be proved in Sections 7-11.  The ``if'' direction
of the proof is given in Section 12.

\section{Basic Properties of Duets and Trios}

In this section we will develop some basic tools for working with critical trios, and establish some simple bounds on 
deficiency.  Most if not all of these results are straightforward generalizations of existing theorems concerning product sets.

\subsection{Uncrossing}

The goal of this subsection is to develop a basic uncrossing procedure which can be used to move 
from one trio to another one.  We will use this to prove a variant of a theorem of Hamidoune which will 
play a key role in the proof of our main theorem.

\bigskip

\noindent{\bf Definition:} Let $(X,Y)$ be a duet,  let $A \subseteq X$ and $B \subseteq Y$, 
and assume that each of $A,B$ is either finite or cofinite.  We say that $(A,B)$ is a \emph{cross} 
\emph{(of} $(X,Y)$) if $x \not\sim y$ for every $x \in A$ and $y \in B$.

\bigskip

Next we introduce two definitions which identify a key correspondence between $G$-trios
and crosses in $G$-duets.

\bigskip

\noindent{\bf Definition:} Let $(X,Y,Z)$ be a trio and let $z \in Z$.  Then 
$(N(z) \cap X, N(z) \cap Y)$ is a cross of $(X,Y)$ which we call the \emph{cross 
associated with} $z$.  

\bigskip

\noindent{\bf Definition:} Let $(X,Y)$ be a $G$-duet and let $(A,B)$ be a cross of $(X,Y)$.  
Define $Z = \{ (gA, gB) \mid g \in G \}$ and extend the incidence relation $\sim$ from 
$X \cup Y$ to $X \cup Y \cup Z$ by the rule that $(A',B') \in Z$ is incident with 
$A' \cup B'$.  It follows that $(X,Y,Z)$ is a $G$-trio and we call it the trio \emph{associated with the cross} $(A,B)$.

\bigskip

Consider a $G$-trio $(X,Y,Z)$, choose $z \in Z$ and let $(A,B)$ be the cross associated 
with $z$.  Then starting from the $G$-duet $(X,Y)$ let $(X,Y,Z')$ be the trio associated with
the cross $(A,B)$.  It follows that $(X,Y,Z')$ is strongly isomorphic to 
the trio obtained from $(X,Y,Z)$ by identifying clones in $Z$ (i.e. the trio $(X,Y,\mcr)$ where 
$\mcr$ is the clone partition of $Z$).  

So (up to clones), we see that trios correspond to duets with distinguished crosses. 
Next we will introduce some natural terminology for crosses which correspond
with already established definitions for trios.  Let $(A,B)$ be a cross of the $G$-duet $(X,Y)$ 
and let $(X,Y,Z)$ be the trio associated with $(A,B)$.  
We say that $(A,B)$ is \emph{trivial} if $A = \emptyset$ or $B = \emptyset$.  We define 
$\delta(A,B) = \delta(X,Y,Z)$ and we say that $(A,B)$ is \emph{critical} if 
$\delta(A,B) > 0$.  Naturally, we will say that $(A,B)$ is maximal if the only cross 
$(A^*,B^*)$ with $A^* \supseteq A$ and $B^* \supseteq B$ is given by $A^* = A$ and 
$B^* = B$.  

The natural motivation for these definitions are the following equivalences which hold whenever 
$\Theta = (X,Y,Z)$ is the trio associated with the cross $(A,B)$ of the duet $(X,Y)$, and also whenever 
$(A,B)$ is the cross of $(X,Y)$ associated with a point $z \in Z$ of the trio $\Theta = (X,Y,Z)$.  
\begin{itemize}
\item $(X,Y,Z)$ is nontrivial if and only if $(X,Y)$ is nonempty and $(A,B)$ is nontrivial;
\item $(X,Y,Z)$ is critical if and only if $(A,B)$ is critical;
\item $w(A) = w(Z,X)$, $\overline{w}(A) = \overline{w}(Z,X)$, $w(B) = w(Z,Y)$, and
		$\overline{w}(B) = \overline{w}(Z,Y)$;
\item If $B$ is cofinite, $\delta(A,B) = w(A) + w(X,Y) - \overline{w}(B)$;
\item If $\overline{d}(X,Y)$ is finite, $\delta(A,B) = w(A) + w(B) - \overline{w}(X,Y)$;
\item $(A,B)$ is maximal if and only if every trio $\Theta' = (X,Y,Z; \sim')$ for which
$\Theta \le \Theta'$ satisfies $\sim'|_{X \cup Z} = \sim|_{X \cup Z}$ and $\sim'|_{Y \cup Z} = \sim|_{Y \cup Z}$.
\end{itemize}

With this, we are ready for the following observation which is our basic uncrossing argument.  The first 
part of this observation follows immediately from the definition of cross, while the second follows from the above 
definitions (here we need to split into cases depending on which of $\overline{w}(X,Y)$ and $\overline{w}(A)$ and $\overline{w}(B)$ is 
finite).  

\begin{observation}[Uncrossing]
\label{uncrossing}
If $(A,B)$ and $(A',B')$ are crosses of a duet, then
\begin{enumerate}
\item $(A \cap A', B \cup B')$ and $(A \cup A', B \cap B')$ are crosses.
\item $\delta(A,B) + \delta(A',B') = \delta(A \cap A', B \cup B') + \delta(A \cup A', B \cap B')$.
\end{enumerate}
\end{observation}

Since connectivity will play an important role in our arguments, we will be interested 
in when the trio $(X,Y,Z)$ which is associated with the cross $(A,B)$ of $(X,Y)$ is 
connected.  For this purpose, let us introduce a little further notation.  If $G$ acts transitively 
on the set $X$ and $A \subseteq X$, the \emph{closure} of $A$, denoted $[A]$, is the 
minimal block of imprimitivity which contains $A$ (since blocks of imprimitivity are closed under intersection, 
there is always a unique such minimal block).

\begin{observation}
\label{cross_con}
Let $(A,B)$ be a cross of the duet $(X,Y)$ and let $(X,Y,Z)$ be the associated trio.  Then 
$(X,Z)$ is connected if and only if $[A] = X$.
\end{observation}

\noindent{\it Proof:}
First suppose that $(Z,X)$ is disconnected and let $\{ (Z_i,X_i) \}_{i \in I}$ be its connected components.  
Then each $X_i$ is a block of imprimitivity, and one must contain $A$, so $[A] \neq X$.
Next suppose that $[A] \neq X$ and let ${\mathcal P}$ be the associated 
system of imprimitivity.  Then every $z \in Z$ has $N(z) \subseteq P$ for some $P \in {\mathcal P}$, so 
$(Z,X)$ is disconnected.  
\quad\quad$\Box$

\bigskip

Our next goal will be to define a notion of deficiency for a duet so that the deficiency of a duet is 
equal to the maximum deficiency of a nontrivial cross.  However, we will require the following easy 
lemma to ensure that our upcoming definition makes sense.  

\begin{lemma} 
\label{nontriv_bound}
For every nontrivial $G$-trio $(X,Y,Z)$ we have
\begin{enumerate}
\item $\delta(X,Y,Z) \le  w(X,Y) \le \overline{w}(Y,Z)$
\item if $G$ is finite, then $\delta(X,Y,Z) \le \frac{1}{2}|G|$.
\end{enumerate}
\end{lemma}

\noindent{\it Proof:}
For the second inequality in part 1, choose $x \in X$ and $z \in Z$ with $x \sim z$ and observe
\[ w(X,Y) = w( N(x) \cap Y ) \le w( Y \setminus N(z) ) = \overline{w}(Y,Z). \]
To prove the other inequality in part 1, we may assume that $\overline{w}(Y,Z) < \infty$.  By an 
argument similar to the above we find that  $w(X,Z) \le \overline{w}(Y,Z)$ and then we have
\[ \delta(X,Y,Z) = w(X,Y) + w(X,Z) - \overline{w}(Y,Z) \le w(X,Y). \]
For part 2, assume $G$ is finite and suppose (without loss) that 
$w(X,Y) \le w(Y,Z)$.  Then part 1 implies $w(X,Y) \le \overline{w}(Y,Z)$ so $w(X,Y) \le \frac{1}{2}|G|$.  
Using part 1 again we have $\delta(X,Y,Z) \le w(X,Y) \le \frac{1}{2}|G|$ as desired.
\quad\quad$\Box$

\bigskip

Let $\Delta$ be a nonempty duet and let $(A,B)$ be a nontrivial cross of $\Delta$.  Then by applying the
previous lemma to the trio associated with the cross $(A,B)$ we deduce that $\delta(A,B) \le \min \{ w( \Delta), \overline{w}(\Delta) \}$.  So, in particular, there do not exist nontrivial crosses of $\Delta$ with arbitrarily large deficiency.  
This permits the following definition.  

\bigskip

\noindent{\bf Definition:} For a partial duet $(X,Y)$, the \emph{deficiency} of $(X,Y)$ is  
\[ \delta(X,Y) = \max \{  \delta(A,B) \mid \mbox{ $(A,B)$ is a nontrivial cross of $(X,Y)$ } \}. \]

\smallskip

If $(X,Y)$ is a duet and $A \subseteq X$ is finite or cofinite, then $(A, Y \setminus N(A))$ is a cross.  
Similarly, if $B \subseteq Y$ is finite or cofinite, then $(X \setminus N(B), B)$ is a cross.  It will be 
convenient for us to further extend our notion of deficiency to subsets of $X$ and $Y$ in accordance 
with this observation.

\bigskip

\noindent{\bf Definition:} If $(X,Y)$ is a duet and $A \subseteq X$ and $B \subseteq Y$ are finite or 
cofinite, we define the \emph{deficiency} of $A$ and $B$ as follows:
\begin{align*}
	\delta(A) &= \delta(A, Y \setminus N(A)),	\\
	\delta(B) &= \delta(X \setminus N(B), B).
\end{align*}

As usual, we say that a set is \emph{critical} if it has positive deficiency.  We say that a nonempty subset 
$A \subseteq X$ (or $B \subseteq Y$) is \emph{nontrivially critical} if $\delta(A) > 0$
 ($\delta(B) > 0$) and $N(A) \neq Y$ ($N(B) \neq X$).  Since it will be frequently used, let us 
 highlight the following observation.  
 
\begin{observation}
\label{nbr_obs}
If $(X,Y)$ is a duet with $w(X,Y) < \infty$ and $A \subseteq X$ is finite then 
\begin{enumerate}
\item $\delta(A) = w(A) + w(X,Y) - w(N(A))$
\item if $(X,Y)$ is nonempty, $w(N(A)) \ge w(A)$.
\end{enumerate}
\end{observation}

\noindent{\it Proof:} Part 1 is a restatement of our definitions.  For part 2 we let 
$(X,Y,Z)$ be the trio associated with the cross $(A, Y \setminus N(A))$.  If $A = \emptyset$ or $N(A) = Y$ the 
result holds immediately.  Otherwise, Lemma \ref{nontriv_bound} gives
\[ w(N(A)) = \overline{w}(Y,Z) \ge w(Z,X) = w(A). \quad\quad\Box \]

\smallskip

Our next result reveals an important property of duets.  We have credited Hamidoune with this theorem 
since he proved the analogue for digraphs.  To see this connection between duets and digraphs, note that if $\Gamma = (V,E)$ is a digraph and 
$G \le Aut(\Gamma)$ acts transitively on $V$ with finite vertex stabilizers, then we may form a $G$-duet on 
$(V,V')$ where $V'$ is a copy of $V$, by the rule that $x \in V$ and $y \in V'$ satisfy $x \sim y$ if $(x,y) \in E$.  

\begin{theorem}[Hamidoune]
\label{hamidoune}
If $\Delta = (X,Y)$ is a partial $G$-duet, 
there exists a finite nontrivially critical block of imprimitivity $T \subseteq X$ or 
$T \subseteq Y$ so that $\delta(T) = \delta(\Delta)$.
\end{theorem}

\noindent{\it Proof:}
Choose a nontrivial cross $(A,B)$ subject to the constraints
\begin{enumerate}
\item[$\alpha$)] $\delta(A,B)$ is maximum.
\item[$\beta$)] $\min\{ w(A), w(B) \}$ is minimum (subject to $\alpha$).
\end{enumerate}
Assume (without loss) that $w(A) \le w(B)$, and note that this implies that $A$ is finite.  
Now suppose (for a contradiction) that $A$ is not a block of imprimitivity.  
Then we may choose $g \in G$ so that 
$\emptyset \neq A \cap gA \neq A$.  Set $(A',B') = (A \cap gA, B \cup gB)$ and 
$(A'',B'') = (A \cup gA, B \cap gB)$ and note that these are both crosses.  If 
$(A'',B'')$ is trivial then $B \cap gB = \emptyset$, so $B$ is finite, but then
\begin{align*}
\delta(A',B') 
	&=	w(A') + w(B') - \overline{w}(\Delta)	\\
	&>	2 w(B) - \overline{w}(\Delta)	\\
	&\ge w(A) + w(B) - \overline{w}(\Delta)	\\
	&= \delta(A,B)
\end{align*}
so $(A',B')$ contradicts our choice of $(A,B)$ for $\alpha$.  It follows that both $(A',B')$ and $(A'',B'')$ are nontrivial.  
However, then by Observation \ref{uncrossing} (applied to $(A,B)$ and $(gA,gB)$), 
either one of $(A',B')$ or $(A'',B'')$ has greater deficiency than $(A,B)$ which again contradicts 
the choice of $(A,B)$ for $\alpha$, or $\delta(A',B') = \delta(A,B)$ and then $(A',B')$ contradicts the 
choice of $(A,B)$ for $\beta$.  
\quad\quad$\Box$

\subsection{Purification}

In this subsection we will use Theorem \ref{hamidoune} from the previous subsection to develop a process 
called purification which we will later use to simplify critical crosses.  The results in this subsection have been 
proved only for duets of finite weight since this is all that we require, but they can easily be extended to duets of 
infinite (but cofinite) weight.  

\begin{proposition}
\label{block_mult}
Let $\Delta = (X,Y)$ be a duet with $w(\Delta) < \infty$, let ${\mathcal P}$ be a system of imprimitivity on $X$ with finite blocks, and let $T \in {\mathcal P}$.  Setting $\Delta' = ( {\mathcal P}, Y )$ we have
\begin{enumerate}
\item $w(N_{\Delta}(T)) = w(T) d( Y, {\mathcal P} )$.
\item If $T \neq X$ and $w(N(T)) = w(T)$, then $\Delta$ is disconnected.
\item If $T$ is critical, then $w(N(T)) = w(T) \Big\lceil \frac{w(\Delta)}{w(T)} \Big\rceil$.
\end{enumerate}
\end{proposition}

\noindent{\it Proof:}
Part 1 follows from the equation
\[ w(N_{\Delta}(T)) = w(N_{\Delta'}(T)) = w( {\mathcal P}, Y) 
	= w(T) d( Y, {\mathcal P}).	\]
For part 2,  note that if $w(N(T)) = w(T)$ then by part 1 we have
 $d( Y, {\mathcal P} ) = 1$.  Since $T \neq X$ we find that $\Delta'$ is disconnected, 
but then $\Delta$ is also disconnected.  For part 3, choose $x \in T$ and note that $w(N(T)) \ge w(N(x)) = w(\Delta)$.  Since $T$ is critical we also have $w(T) + w(\Delta) > w(N(T))$.  However, Part 1 implies that $w(N(T))$ is a multiple of $w(T)$, and this yields the desired result.
\quad\quad$\Box$

\bigskip

As already noted, the first part of the above proposition implies that $w(N(T))$ will always be a multiple of $w(T)$.  Our next 
lemma uses this to establish a nice property of subsets of critical blocks of imprimitivity.

\begin{lemma}
\label{block_best}
Let $\Delta = (X,Y)$ be a duet with $w(\Delta) < \infty$ and let $T \subseteq X$ be a finite critical block of imprimitivity.  
If $S \subseteq T$, then $\delta(S) \le \delta(T)$.
\end{lemma}

\noindent{\it Proof:}
Choose $S \subseteq T$ subject to
\begin{enumerate}
\item[$\alpha$)] $\delta(S)$ maximum.
\item[$\beta$)] $|S|$ minimum (subject to $\alpha$).
\end{enumerate}
As in the proof of Theorem \ref{hamidoune}, it follows from our choice that $S$ must be a block of imprimitivity 
(otherwise we may choose $g \in G$ with $\emptyset \neq S \cap gS \neq S$ and either one of $gS \cap S$ or 
$gS \cup S$ has greater deficiency than $S$, or both have the same deficiency, and either case contradicts our choice 
of $S$).   

Therefore $k = \frac{w(T)}{w(S)}$ is an integer.  Now, choose a point $x \in S$ and observe 
the following (the last inequality follows from the assumption $T$ is critical).
\[ w(N(T) \setminus N(S)) \le w(N(T) \setminus N(x)) = w(N(T)) - w(\Delta) < w(T) = k w(S) \]
It follows from the previous proposition that $w(N(T)) - w(N(S))$ is a multiple of $w(S)$ so we conclude that this 
quantity is at most $(k-1) w(S)$.  But then we have 
\begin{align*}
\delta(T) - \delta(S) 
	&= w(T) - w(N(T)) - w(S) + w(N(S))	\\
	&= (k-1) w(S) - \big( w(N(T)) - w(N(S)) \big) \ge 0.
\end{align*}
It follows that every $S \subseteq T$ satisfies $\delta(S) \le \delta(T)$ as desired.\quad\quad$\Box$

\bigskip

Throughout much of this paper, we will be considering a group acting transitively on a set $X$ with a distinguished 
system of imprimitivity ${\mathcal P}$.  For a set $A \subseteq X$ and a block $P \in {\mathcal P}$ we say that 
\[ \mbox{$P$ is } \left\{ \begin{array}{c}
				\mbox{\emph{inner}}	\\
				\mbox{\emph{outer}}	\\
				\mbox{\emph{boundary}}
				\end{array}	\right\}
\;\;
\mbox{(\emph{with respect to} $A$) if } \; \left\{ \begin{array}{c}
				P \subseteq A	\\
				P \cap A = \emptyset	\\
				otherwise.
				\end{array} \right\}.
\]	
We say that a set $A$ is \emph{pure} with respect to ${\mathcal P}$ if there are no boundary blocks, and \emph{impure} otherwise.
Similarly, if $(X,Y)$ is a duet and $\mcp, \mcq$ are systems of imprimitivity of $X,Y$ then we say that a cross $(A,B)$ of $(X,Y)$ is 
\emph{pure with respect to} $(\mcp,\mcq)$ if $A$ is pure with respect to $\mcp$ and $B$ is pure with respect to $\mcq$.  

\bigskip

\noindent{\bf Definition:} If $(A,B)$ is a cross of $(X,Y)$ and $P \subseteq X$ is 
a finite critical boundary block, we define the cross $(A \cup P, B \setminus N(P))$ to be 
the \emph{weak purification of} $(A,B)$ \emph{at} $P$.  We define weak 
purification similarly for blocks of $Y$.

\bigskip

A disadvantage of weak purification is that purifying a maximal cross may result 
in a cross which is no longer maximal.  Next we introduce a stronger process 
which does not have this defect.

\bigskip

\noindent{\bf Definition:} If $(A,B)$ is a cross of $(X,Y)$ and $P \subseteq X$ is 
a finite critical boundary block, we set $B' = B \setminus N(P)$ and $A' = X \setminus N(B')$ 
and define $(A',B')$ to be the \emph{purification of} $(A,B)$ \emph{at} $P$.  We define
purification similarly for blocks of $Y$.

\bigskip

We close this subsection with a theorem which encapsulates the essential properties of purification.

\begin{theorem}[Purification]
\label{purify_thm}
Let $(X,Y)$ be a duet with $w(X,Y) < \infty$, let $(A,B)$ be a nontrivial critical cross, and let $P \subseteq X$ be a finite critical 
boundary block.  If $(A',B')$ and $(A'',B'')$ are the weak purification and the purification of $(A,B)$ at $P$ then 
\begin{enumerate}
\item $\delta(A'',B'') \ge \delta(A',B') \ge \delta(A,B)$
\item If $(A,B)$ is maximal, then $(A'',B'')$ is maximal.
\item $w(B \setminus B') = w(B \setminus B'') < w(P)$.
\end{enumerate}
\end{theorem}

\noindent{\it Proof:} Part 1 follows from $\delta(A'',B'') \ge \delta(A',B')$ and 
the following equation (which calls on Lemma \ref{block_best}).
\begin{eqnarray*}
\delta(A',B') - \delta(A,B) 
	& = &		w( P \setminus A ) - w( B \cap N(P) )	\\
	& \ge &	w( P \setminus A ) - \big( w(N(P)) - w(N(P \cap A)) \big)		\\
	& = &		\delta(P) - \delta(P \cap A)	\\
	& \ge &	0.	\quad\quad\Box
\end{eqnarray*}
Part 2 is an immediate consequence of our definitions, and part 3 follows from $B' = B''$ and  the inequality
\[ 0 \le \delta(A',B') - \delta(A,B) = w(A' \setminus A) - w(B \setminus B') < w(P) - w(B \setminus B'). \quad\quad\Box \]

\subsection{Connectivity}
\label{discon_sec}

In this subsection we will prove a couple of results concerning connectivity.  The first is a variant of a theorem of Olson \cite{olson2}
which asserts that whenever $A,B$ are finite subsets of the group $G$ and $B$ is a generating set containing $1$, either
$AB= G$ or $|AB| \ge |A| + \frac{1}{2}|B|$.  Equivalently, whenever $B \subset G$ is a generating set containing $1$, the duet
$\Delta = {\mathit CayleyC}(G; B)$ satisfies $\delta(\Delta) \le \frac{1}{2}|B|$.

\begin{theorem}[Olson]
\label{connected_bound}
If $\Delta = (X,Y)$ is a connected partial $G$-duet, 
\[ \delta(\Delta) \le \min \left\{ \tfrac{1}{2} w(\Delta), \tfrac{1}{3}|G| \right\}.\]
\end{theorem}

\noindent{\it Proof:} First we establish an easy inequality for integers.  Let $d,t$ be nonnegative integers
and choose integers $a,k$ so that $d = at + k$ where $1 \le k \le t$.  Then we have
\[ d + t - t \lceil \tfrac{d}{t} \rceil = (at+k) + t - t(a+1) = k \le (at+k)/(a+1) = d / \lceil \tfrac{d}{t} \rceil \]
If $w(X,Y) = \infty$ then we have nothing to prove, so we may assume $w(X,Y) < \infty$.  
Now, apply Theorem \ref{hamidoune} to choose a block of imprimitivity $T$ so that 
$T$ is nontrivally critical and $\delta(T) = \delta(\Delta)$.  We assume (without loss) 
that $T \subseteq X$.  Then by Proposition \ref{block_mult} (and the above inequality) we have
\begin{align*}
\delta(\Delta) 
	&= w(\Delta) + w(T) - w(N(T))	\\ 
	&= w(\Delta) + w(T) - w(T) \left\lceil \tfrac{w(\Delta)}{w(T)} \right\rceil		\\
	&\le w(\Delta) \big/ \left\lceil \tfrac{w(\Delta)}{w(T)} \right\rceil.
\end{align*}
If $w(\Delta) \le w(T)$ then Proposition \ref{block_mult} implies $w(N(T)) = w(T)$ and that $\Delta$ is disconnected
which is contradictory.  Thus $w(\Delta) > w(T)$ and this implies $\delta(\Delta) \le \tfrac{1}{2} w(\Delta)$.  
Next observe that if $G$ is finite then $w(N(T))$ and $|G|$ are both multiples of $w(T)$ and $w(T) < w(N(T)) < |G|$.  Thus $|G| \ge 3 w(T)$ and this gives $\delta(\Delta) \le w(T) \le \frac{1}{3}|G|$.  
\quad\quad$\Box$

\begin{corollary}
\label{connected_set_bound}
If $(X,Y)$ is a nonempty $G$-duet and $A \subseteq X$ satisfies $[A] = X$ and $N(A) \neq Y$, then 
\[\delta(A) \le \min\left\{ \tfrac{1}{2}w(A), \tfrac{1}{3}|G| \right\}.\]
\end{corollary}

\noindent{\it Proof:} 
Let $(X,Y,Z)$ be the trio associated with the cross $(A,Y \setminus N(A))$.  Then, by Observation \ref{cross_con} we have that $(X,Z)$ is connected.  So, by the previous theorem gives
\[ \delta(A) = \delta(X,Y,Z) \le \delta(X,Z) \le \min\{ \tfrac{1}{2}w(X,Z), \tfrac{1}{3}|G| \} 
=   \min \{ \tfrac{1}{2}w(A), \tfrac{1}{3}|G| \} \]
as desired. \quad\quad$\Box$

\bigskip

Next we prove a key lemma concerning critical crosses in disconnected duets.

\begin{lemma} 
\label{cross_in_discon}
Let $(X,Y)$ be a disconnected duet with component quotient $(\mcp$, $\mcq)$.  
If $(A,B)$ is a maximal critical cross of $(X,Y)$  then either $(A,B)$ is pure with respect to $(\mcp,\mcq)$, or 
there is exactly one $P \in {\mathcal P}$ $(Q \in {\mathcal Q})$ which is a boundary block and further
\[ \delta(A,B) \le \min\{ \tfrac{1}{2}w(X,Y), \tfrac{1}{3}w^{\circ}(\mcp) \}. \]
\end{lemma}

\noindent{\it Proof:} Assume that $\Delta = (X,Y)$ is a $G$-duet, and note that $w(X,Y) < \infty$.  Let $\mcp = \{ P_i \mid i \in I \}$
and $\mcq = \{ Q_i \mid i \in I \}$ and assume $P_i \sim Q_i$ for every $i \in I$.  We may assume that $w(A) < \infty$, 
and then for every $i \in I$ set $A_i = A \cap P_i$.  This gives us
\begin{align}
\label{discon_eqn}
	0	&<	\delta_G(A)	 \nonumber	\\
		&=	w_G(\Delta) -w_G(N(A)) +  w_G(A)	\nonumber	\\
		&=	w_G(\Delta) - \sum_{i \in I} \big( w_G(N(A_i)) - w_G(A_i) \big).
\end{align}
If $N(A_i) = Q_i$ then $w_G(A_i) \le w_G(P_i) = w_G(Q_i) = w_G(N(A_i))$.  
Next suppose that $i \in I$ satisfies $A_i \neq \emptyset$ and $N(A_i) \neq Q_i$.  Then 
$(A_i, Q_i \setminus N(A_i))$ is a nontrivial cross of the connected duet $\Delta_i = (P_i,Q_i)$.  
If we set $H = G_{P_i} = G_{Q_i}$ then $w_G^{\circ}(X) = w_H^{\circ}(P_i)$ and 
$w_G^{\circ}(Y) = w_H^{\circ}(Q_i)$ (since every point $x \in P_i$ satisfies $G_x = H_x$ 
and similarly every $y \in Q_i$ satisfies $G_y = H_y$).  It follows that 
$w_G(\Delta) = w_H(\Delta_i)$.  Since $\Delta_i$ is connected, 
Theorem \ref{connected_bound} gives us 
\begin{align*}
\tfrac{1}{2} w_G( \Delta)	
	&=	\tfrac{1}{2} w_H(\Delta_i) \\
	&\ge	\delta_H(\Delta_i) \\
	&\ge w_H(\Delta_i)  - w_H( N(A_i))  + w_H(A_i)\\
	&= w_G(\Delta) - w_G(N(A_i)) + w_G(A_i) 
\end{align*}
so $w_G(N(A_i)) - w_G(A_i) \ge \frac{1}{2}w_G(\Delta)$.  It follows from this, maximality, and equation \ref{discon_eqn} that there is 
at most one $i \in I$ so that $P_i$ is a boundary block with respect to $A$ (and similarly for $\mcq$).  If there does exist such a boundary block 
$P_i$, then $\delta(A,B) \le \frac{1}{2}w(\Delta)$, and by applying Theorem \ref{connected_bound} to the $H$-duet $\Delta_i = (P_i,Q_i)$
we find $\delta_{\Delta}(A,B) = \delta_{\Delta_i}(A_i) \le \frac{1}{3}|H| = \frac{1}{3} w^{\circ}(\mcp)$
which completes the proof.
\quad\quad$\Box$

\subsection{Disconnected Trios and Closure}

In this section we prove a lemma which reveals the structure of maximal critical disconnected trios, and a theorem concerning
the closure of a critical set in a duet.

\begin{lemma}[Disconnected Trios]
\label{disconnected}
Let $\Theta = (X,Y,Z)$ be a maximal nontrivial critical trio in $G$ with $(X,Y)$  disconnected.  Then $\Theta$ is either a pure or 
impure beat relative to $(X,Y)$.  In the latter case, if $\Theta' = (X',Y',Z')$ is an $H$-trio which is a continuation of $\Theta$, 
one of the following holds:
\begin{enumerate}
\item $H$ is finite and $|H| < |G|$.
\item $H$ is infinite, so there are sides $\Delta$ of $\Theta$ and $\Delta'$ of $\Theta'$ with $w_G(\Delta) = \infty = w_H(\Delta')$, and 
	these sides satisfy $\bar{w}_G(\Delta) = \bar{w}_H(\Delta')$.  Furthermore, $\Theta'$ has a side of finite weight which is connected, 
	but all sides of $\Theta$ of finite weight are disconnected.
\end{enumerate}
\end{lemma}

\noindent{\it Proof:}
Let $\Theta = (X,Y,Z)$ be a nontrivial maximal critical $G$-trio with $(X,Y)$ disconnected, and let $(\mcp, \mcq)$ be the component quotient of $(X,Y)$.
Choose $z \in Z$, let $(A,B)$ be the cross of $(X,Y)$ associated with $z$ and assume (without loss) that $A$ is finite.  If $(A,B)$ is pure with respect to $(\mcp,\mcq)$, then $A$ contains a block of $\mcp$ so we have $w^{\circ}( \mcp ) = w^{\circ}( \mcq ) < \infty$.  It then follows from the maximality of $\Theta$ 
that any two points in the same block of $\mcp$ $(\mcq)$ are clones, so $\Theta$ is a pure beat relative to $(X,Y)$, and we are finished.
Otherwise, it follows from the previous lemma that there is exactly one $P \in \mcp$ and one $Q \in \mcq$ for which 
$P$ and $Q$ are boundary blocks.  It then follows from the maximality of $\Theta$ that $\Theta$ is an impure beat relative to 
$(X,Y)$.  

Using these blocks $P \in \mcp$ and $Q \in \mcq$ we may now choose $R \subseteq Z$ so that $\Theta' = (P,Q,R)$ is a continuation of $\Theta$ 
and assume $\Theta'$ is an $H$-trio.  If $H$ is finite, then $|H| < |G|$ and there is nothing left to prove.  Otherwise, the blocks in 
$\mcp$ and $\mcq$ are infinite, and it follows that $A$ is contained in a single block of $\mcp$ (since $A$ is finite).  In this case, 
$w(Y,Z) = \infty = w(Q,R)$ and $\bar{w}(Y,Z) = \bar{w}(Q,R)$.  Furthermore, the sides $(X,Y)$ and $(X,Z)$ of $\Theta$ have finite weight 
and are disconnected, but the side $(P,Q)$ of $\Theta'$ has finite weight and is connected.
\quad\quad$\Box$

\smallskip

Our next result reveals an important property of the closure of a critical set.  Note that in this 
theorem, the block of imprimitivity $T$ may be trivially critical.

\begin{theorem}[Critical Closure]
\label{critical_closure}
Let $\Delta = (X,Y)$ be a connected $G$-duet with $w(\Delta) < \infty$, let $A \subseteq X$ be finite and critical 
and assume that $T = [A] \neq X$.  Then $T$ is finite and critical. 
Further, if $N(A) \neq N(T)$ then $\delta_{\Delta}(A) \le \frac{1}{3}w(T)$.
\end{theorem}

\noindent{\it Proof:} Let $(X,Y,Z; \sim)$ be the trio associated with $(X,Y)$ and the cross $(A, Y \setminus N_{\Delta}(A))$.  
The duet $\Upsilon = (X,Z; \sim)$ is disconnected, so it has a nontrivial component quotient $(\mcp, \mcq)$, and our assumptions 
imply that $T \in \mcp$.  Now, let $(B,C)$ be a cross of $(X,Z)$ associated with some point in $Y$.  It follows from the construction 
of $(X,Y,Z; \sim)$ that $N_{\Upsilon}(B) = Z \setminus C$, so we may choose $B \subseteq B^* \subseteq X$ so that $(B^*,C)$ is
a maximal cross of $\Upsilon$.  Now by Lemma \ref{cross_in_discon}, the set $B^*$ is either pure with respect 
to $\mcp$, or it has exactly one boundary block.  It follows from the connectivity of $(X,Y)$ that $[B] = X$ so $B^*$ is not contained
in a block of~$\mcp$.  Thus, $B^*$ contains a block of $\mcp$, and since $B^*$ is finite, we have that $T$ is finite.

Now define $\Upsilon^* = (X,Z, \sim^*)$ by the rule that $x \in X$ and $z \in Z$ satisfy $x \sim^* z$ if the blocks $P \in \mcp$, $Q \in \mcq$ with 
$x \in P$ and $z \in Q$ satisfy $P \sim Q$.  Suppose first that $C$ is pure relative to $\mcq$ and note that this implies
$N_{\Upsilon}(B)$ is pure with respect to $\mcq$.  In this case, for every $y \in Y$ and $z \in Z$ we have that there exists $x \in X$ with 
$y \sim x \sim^* z$ if and only if there exists $x' \in X$ with $y \sim x' \sim z$.  It follows from this that $N_{\Delta}(A) = N_{\Delta}(T)$, 
which implies that $T$ is critical (and completes the proof of this case).  So, we may now assume that $C$ is not pure relative to $\mcq$, which implies
$N_{\Delta}(A) \neq N_{\Delta}(T)$.   Lemma \ref{cross_in_discon} implies
\[ \delta_{\Delta}(A) = \delta(X,Y,Z; \sim) = \delta_{\Upsilon}(B,C) \le \delta_{\Upsilon}(B^*,C) \le \tfrac{1}{3} w^{\circ}(\mcp) = w(T). \]
Furthermore, since $N_{\Upsilon}(B)$ has a single boundary block, we have 
$w(N_{\Upsilon^*}(B)) - w(N_{\Upsilon}(B)) \le w(\Upsilon^*) - w(\Upsilon)$.  This gives
$\delta_{\Delta}(T) = \delta_{\Upsilon^*}(B) \ge \delta_{\Upsilon}(B) = \delta_{\Delta}(A) > 0$ as desired.
\quad\quad$\Box$

\bigskip

Next we record a simple corollary of this for future use.

\begin{corollary}
\label{block_subset_bound2}
Let $\Delta = (X,Y)$ be a connected $G$-duet with $w(X,Y) < \infty$, let $T \subset X$ be a finite block of imprimitivity and let
$A \subseteq T$ satisfy $N(A) \neq N(T)$.  Then one of the following holds.
\begin{enumerate}
\item $w(N(A)) \ge w(\Delta) + w(A) - \frac{1}{3} w(T)$
\item There is a block $A \subseteq S \subseteq T$ so that $N(A) = N(S)$ 
	and $w(S) = \frac{1}{2} w(T)$.
\end{enumerate}
\noindent{In particular, $w(N(A)) \ge w(\Delta) + w(A) - \frac{1}{2}w(T)$
is always satisfied.  }
\end{corollary}

\noindent{\it Proof:} We may assume that $A$ is critical, as otherwise the first outcome holds trivially.
Let $S = [A]$.  If $N(A) \neq N(S)$, then the previous lemma (and the observation that $w(S) \le w(T)$) 
implies that 1 holds.  Thus we may assume $N(A) = N(S)$.  If $w(S) \le \frac{1}{3} w(T)$ then 
$w(N(A)) \ge w(\Delta) \ge w(S) + w(\Delta) - \frac{1}{3}w(T)$ and 1 holds.  Otherwise 
$w(S) = \frac{1}{2} w(T)$ so $w(N(A)) \ge w(\Delta) \ge w(\Delta) + w(A) - \frac{1}{2}w(T)$ and 2 holds.
\quad\quad$\Box$

\section{The Main Theorem}
\label{proof_sec}

The goal of this section is to give a proof of the main theorem assuming the validity of several lemmas 
which will be stated in this section but proved later in the paper.

\subsection{Sidon Duets}

We say that a rank 2 incidence geometry $(X,Y)$ is \emph{Sidon}\footnote{Rank 2 incidence geometries with this property are sometimes 
called linear spaces, but the word ``linear'' is rather heavily used in mathematical writing, so we have adopted the term Sidon instead.} if 
there do not exist distinct points $x_1,x_2 \in X$ and $y_1,y_2 \in Y$ with $x_i \sim y_j$ for every $i,j \in \{1,2\}$.  Note that that 
${\mathit CayleyC}(G; A)$ is Sidon if and only if $A$ is a Sidon set.  Our goal in this subsection will be to prove the following lemma.

\begin{lemma}[Sidon stability]
\label{sidon_stability}
Every connected critical trio with a Sidon side is either a video or 
a pure chord.
\end{lemma}

Our proof of this lemma will rely upon a stability lemma stated below.  The proof of the first part of this lemma is established by Lemma \ref{graph_stability}
and the second part follows from Lemma \ref{clip_stability}.  

If $\Gamma$ is a graph and $A$ and $B$ are operators given by $V$, $E$, $C_k$, or $P_k$ as in Subsection \ref{video_subsec} for which we have defined an incidence relation between $A(\Gamma)$ and $B(\Gamma)$, then we let $\Gamma(A,B)$ denote the incidence geometry $(A(\Gamma), B(\Gamma) )$ given by this relation.  

\begin{lemma}
\label{graph_and_video}
Let $\Theta$ be a connected critical trio and assume one of the following 
\begin{enumerate}
\item $\Theta$ has a side which is a simple graph.
\item $\Theta$ has a side isomorphic to $K_5(E,C_3)$\footnote{This is more commonly known as Desargues configuration.} or $K_6(E,C_3)$.
\end{enumerate}
Then $\Theta$ is either a video or a pure chord.
\end{lemma}

It is worth noting that in both of the lemmas just stated, we have \emph{not} assumed that the trio is maximal.  
In sharp contrast to most of the other theorems in this paper, we do not need to assume maximality 
here, since the property Sidon has extremely strong forcing.  Simple graphs are Sidon by definition, 
and they will play an important role for us as the primary instances of Sidon duets appearing 
in critical trios (the only others are $K_5(E,C_3)$ and $K_6(E,C_3)$).  The next lemma gives an indication of how they will emerge.   

\begin{lemma}
\label{two_sidon}
Let $\Theta = (X,Y,Z)$ be a connected critical trio.  If both $(Z,X)$ and
$(Z,Y)$ are Sidon, then one of $(X,Z)$, $(Z,X)$, $(Z,Y)$, or
$(Y,Z)$ is a simple graph.
\end{lemma}

\noindent{\it Proof:}
Since $(Z,X)$ and $(Z,Y)$ are Sidon, both $d(Z,X)$ 
and $d(Z,Y)$ are finite. Choose $z \in Z$ and let $(A,B)$ be
the cross of $(X,Y)$ corresponding to $z$.  We assume (without loss)
that $w(A) \le w(B)$.  Since $\Theta$ is connected, we may choose
$z' \in Z \setminus \{z\}$ so that the cross $(A',B')$ of $(X,Y)$
which corresponds to $z'$ satisfies $A' \cap A \neq \emptyset$.  
It follows from this (and the assumption that $(X,Z)$ is Sidon) 
that $A \cap A' = \{x\}$ for some $x \in X$.

Suppose (for a contradiction) that $B \cap B' = \emptyset$.  In this case
$x$ is not incident with any point in $B \cup B'$ so we have that 
we have $\overline{w}(X,Y) = w( Y \setminus N(x) ) \ge  2w(B) \ge w(A) + w(B)$, 
but this contradicts the assumption that $\Theta$ is critical.  Therefore $B \cap B'$ is nonempty 
and since $(Y,Z)$ is Sidon it follows that $B \cap B' = \{y\}$ for some $y \in Y$.
Furthermore, we must have $B \cup B' = Y \setminus N(x)$ as otherwise we 
would again have $w(Y \setminus N(x)) \ge 2 w(B)$ giving us the same contradiction.

If the only points in $Z$ incident with $x$ are $z$ and $z'$, then $(Z,X)$ is a graph and
we are done.  Otherwise we may choose $z'' \in Z \setminus \{z,z'\}$ 
which is incident with $x$.  Let $(A'',B'')$ be the cross associated 
with $z''$.  By our assumptions, we must have $B'' \subseteq B \cup B'$, 
$|B'' \cap B'| \le 1$, and $|B'' \cap B| \le 1$.  It follows that
$|B''| = 2$, so $(Y,Z)$ is a simple graph.  
\quad\quad$\Box$

\bigskip

\begin{figure}[ht]
\centerline{\includegraphics[height=3cm]{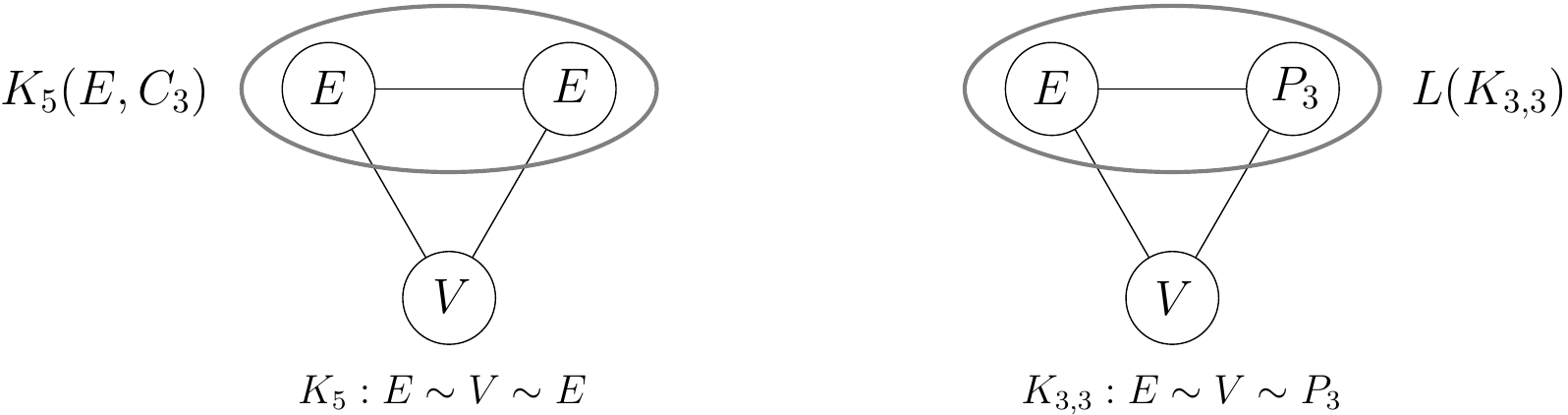}}
\caption{Two small videos}
\label{two_videos}
\end {figure}

Next we require an observation concerning some basic properties of videos.  Here the two small videos shown in Figure \ref{two_videos}
will play a special role.  As indicated by this figure, the side of $K_5 : E \sim V \sim E$ on the pair of copies of $E(K_5)$ is isomorphic to the 
duet $K_5(E,C_3)$, and the side of $K_{3,3} : E \sim V \sim P_3$ on the pair $(E(K_{3,3}), P_3(K_{3,3}))$ is isomorphic to the line graph 
of $K_{3,3}$.

\begin{observation}
\label{sidon_video}
If $\Delta$ is a Sidon side of a video, or a pure chord, then either $\Delta$ or its reverse is a
simple connected graph, or is isomorphic to one of $K_5(E,C_3)$ or $K_6(E,C_3)$.  
\end{observation}

\noindent{\it Proof:} This is immediate for pure chords, and can be verified for the 
exceptional videos by a straightforward check.  Consider a standard video given by
$(E_1, V, E_2) = \Gamma : E \sim V \sim E$.  The sides $(V,E_1)$ and $(V,E_2)$ are both graphs.  
If $|V| = 5$ then $\Gamma \cong K_5$ and then $(E_1, E_2) \cong K_5(E,C_3)$.  Otherwise $|V| \ge 6$, 
so there are two vertex disjoint subgraphs in $P_3(\Gamma)$ and the four edges contained in these 
two paths show that $(E_1,E_2)$ is not Sidon.   

Next consider a standard video given by $(P_3, V, E) = \Gamma : P_3 \sim V \sim E$.  The side $(V,E)$ is graphic
and the side $(P_3,V)$ is not Sidon.  If $\Gamma \cong K_{3,3}$, then $(E,P_3)$ is a graph.  Otherwise, $|V| \ge 8$ and 
we may choose a vertex $v$ and edges $e,e'$ so that neither $v$ nor any vertex adjacent to $v$ is incident with either $e$ or 
$e'$.  Considering $e,e'$ and the three two edge paths with middle vertex $v$ we see that $(P_3,E)$ is not Sidon.

Finally, consider a standard video $(P_3, E, V) = \Gamma : P_3 \sim E \sim V$.  In this case both $(V,E)$ and $(E,P_3)$ 
are graphs.  Choose a vertex $v$ and vertices $w,w'$ not adjacent with $v$.  By considering the vertices 
$w,w'$ and the three two edge paths with middle vertex $v$ we find that $(P_3,V)$ is not Sidon.  
\quad\quad$\Box$

\bigskip

\noindent{\it Proof of Lemma \ref{sidon_stability}:}
Let $\Theta = (X,Y,Z)$ be a connected critical trio and assume that $(X,Y)$ is Sidon.  We begin 
by proving a key claim (note that this claim applies even to trivially critical blocks of imprimitivity).  

\bigskip

\noindent{\it Claim:} There does not exist a finite proper nontrivial block of imprimitivity of $X$ or $Y$ 
which is critical in $(X,Y)$. 

\smallskip

Suppose (for a contradiction) that $T \subset X$ is a finite proper critical block of imprimitivity with $|T| \ge 2$ 
(the case when $T \subset Y$ is similar).  Choose distinct points $x_1,x_2 \in T$.  Then $w(N(T)) \ge 2w(T)$ (by Proposition \ref{block_mult})
and $w(N(T)) \setminus N(x_i)) = w(N(T)) - w(X,Y) < w(T)$ (since $T$ is critical) so for $i=1,2$ we have that 
$x_i$ is incident with more than half of the points in $N(T)$.  Since $(X,Y)$ is 
Sidon it follows that $|N(x_1) \cap N(x_2)| = 1$ and $N(T) = N(x_1) \cup N(x_2)$.  
If there exists 
$x_3 \in T \setminus \{x_1,x_2\}$, then $N(x_3)$ is subset of $N(T)$ of size $d = d(X,Y)$ 
which meets $N(x_1)$ and $N(x_2)$ in at most one point.  In this
case, we must have $d=2$ so $(Y,X)$ is a simple graph.  However, then Lemma \ref{graph_and_video} implies
that $\Theta$ is either a video or a pure chord, and we are finished.  Thus
we may assume $T = \{x_1,x_2\}$.  Now $w(N(T))$ is a multiple of 
$2w^{\circ}(X)$ and $w(N(x_1)) + w(N(x_2))$ is a multiple of $2 w^{\circ}(X)$, and 
it follows from this that $w^{\circ}(Y)$ is a multiple of $2 w^{\circ} (X)$.  
However, now we have 
$w(N(T)) \ge w(N(x_1)) + w^{\circ}(Y) \ge w(X,Y) + w(T)$ and this contradicts
the assumption that $T$ is critical.  

\bigskip

Now, let $(A,B)$ be a critical cross of $(X,Y)$ with $|A|, |B| \ge 2$ (note that one 
exists by our assumptions) chosen so that
\begin{enumerate}
\item[$\alpha$)] $\delta(A,B)$ is maximum.
\item[$\beta$)] $\min\{ w(B) + w^{\circ}(X), w(A)+ w^{\circ}(Y)  \}$ is minimum (subject to $\alpha$).
\end{enumerate}
Assume (without loss) that $w(A) + w^{\circ}(Y) \le w(B) + w^{\circ}(X)$, and note that 
this implies $A$ is finite.  Now suppose (for a contradiction) that there exists $g \in G$ with $2 \le |gA \cap A| < |A|$.  
If $|gB \cap B| \le 1$ then $B$ is also finite and we find
\begin{align*}
w(gA \cap A) + w(gB \cup B) 
	&\ge 2 w^{\circ}(X) + 2 w(B) - w^{\circ}(Y) 	\\
	&\ge w^{\circ}(X) + w(A) + w(B).
\end{align*}
However, then $\delta(gA \cap A, gB \cup B) > \delta(A,B)$ which contradicts our choice 
of $(A,B)$.  Therefore, we must have $|gB \cap B| \ge 2$.  But then, either one of $(gA \cap A, gB \cup B)$ or 
$(gA \cup A, gB \cap B)$ has deficiency greater than $(A,B)$, which contradicts $\alpha$, or both 
are equally critical and then $(gA \cap A, gB \cup B)$ contradicts the choice of $(A,B)$ for $\beta$.  
It follows from this that every $g \in G$ with $|gA \cap A| < |A|$ satisfies $|gA \cap A| \le 1$.  

Our claim together with Theorem \ref{critical_closure} imply that  $[A] = X$ and $[B] = Y$,  and it follows from
this that the trio $\Theta' = (X,Y,Z')$ associated with the cross $(A,B)$ is connected.  Furthermore, by the 
above property of $A$ we must have that $(X,Z')$ is Sidon.  Therefore $\Theta'$ has two Sidon 
sides, so by Lemma \ref{two_sidon} we must have that one of $(X,Y)$, $(Y,X)$, $(X,Z)$, or $(Z,X)$ is a simple connected graph.  
It now follows from Lemma \ref{graph_and_video} that $\Theta'$ is either a video or a pure chord.  
But then, Observation \ref{sidon_video} implies that one of $(X,Y)$ or $(Y,X)$ is either a simple connected 
graph or is isomorphic to one of $K_5(E,C_3)$ or $K_6(E,C_3)$.  Now by applying Lemma \ref{graph_and_video} to $\Theta$ and
either $(X,Y)$ or $(Y,X)$ we deduce that $\Theta$ is either a video or pure chord as desired.
\quad\quad$\Box$

\subsection{Nearness and Constituent Lemmas}

In this subsection we will introduce the important concept of near, and then state some lemmas concerning this property 
which will be used to prove our main theorem.  The proofs of these lemmas is postponed to later in the article.  

\bigskip

\noindent{\bf Definition:} Let $\Delta = (X,Y)$ and $\Lambda$ be duets and let 
${\mathcal P}$, ${\mathcal Q}$ be systems of imprimitivity on $X$, $Y$ with finite blocks.   
We say that $\Delta$ is \emph{near} $\Lambda$ \emph{relative to} $({\mathcal P}$, ${\mathcal Q})$ 
if $({\mathcal P}, {\mathcal Q}) \cong \Lambda$ and the following equation is satisfied
\[ w( {\mathcal P}, {\mathcal Q} ) - w(\Delta) < \min \{ w^{\circ}( {\mathcal P} ), w^{\circ} ( {\mathcal Q} ) \}.  \]
Next we  record a useful observation about the notion of near for future 
use.  We define a system of imprimitivity to be \emph{critical} if its blocks are critical.  

\begin{observation}
\label{near_obs}
If $\Delta = (X,Y)$ is near $\Lambda$ relative to $( {\mathcal P}, {\mathcal Q} )$ and $w(\Delta) < \infty$ then
\begin{enumerate}
\item ${ \mathcal P}$ and ${\mathcal Q}$ are critical.
\item For every $P \in {\mathcal P}$ $(Q \in {\mathcal Q})$ the set $N(P)$ $(N(Q))$ is pure with respect to ${\mathcal Q}$
	$({\mathcal P})$. 
\end{enumerate}
\end{observation}

\noindent{\it Proof:} Let $\Delta' = (\mcp, \mcq)$.  For part 2, let $P \in {\mathcal P}$ and observe that 
\[ w(\Delta') - w^{\circ}(P) < w(\Delta) \le w(N_{\Delta}(P)) \le w(N_{\Delta'}(P)) = w(\Delta').\]  
It follows from Proposition \ref{block_mult} that $w(N_{\Delta}(P))$ and $w(\Delta') = w(N_{\Delta'}(P))$ must be a 
multiples of $w^{\circ}(P)$.  
But then we must have $w(N_{\Delta}(P)) = w(N_{\Delta'}(P))$, so $N(P)$ is pure with respect to ${\mathcal Q}$.  A similar 
argument shows that $N(Q)$ is pure with respect to ${\mathcal P}$ for every $Q \in {\mathcal Q}$. 

For part 1 we have
\[ \delta_{\Delta}(P) = w(P) + w(\Delta) - w( N_{\Delta}(P)) = w(P) + w(\Delta) - w(\Delta') > 0 \]
so ${\mathcal P}$ is critical, and by a similar argument ${\mathcal Q}$ is critical.
\quad\quad$\Box$

\bigskip

Before introducing our next stability lemma, we need to define a particular set of duets which 
will play an essential role in our proof.  A graph is \emph{cubic} if every vertex is incident with $3$ edges.

\bigskip

\noindent{\bf Definition:} 
A duet $\Delta$ is a \emph{clip} if it is isomorphic to one of the following 
\begin{enumerate}
\item		$\Gamma( V, P_3 )$ for a cubic graph $\Gamma$ with $|V(\Gamma)| \ge 6$
\item 	$K_5 ( E, C_3 )$
\item		$K_6 (E, C_3)$
\item		Petersen$( E, C_5 )$
\item		Dodecahedron$( V, F )$.
\end{enumerate}

We call a graph $(V,E)$ \emph{true} if it is simple, connected, and all vertices have degree at least $3$.  We will say that 
a duet $\Delta$ is a \emph{near true graph}, \emph{near sequence}, or \emph{near clip} if it is near a duet of this type.  
The next lemma summarizes three stability lemmas, the first of which follows from Lemma \ref{near_graph_stability}, 
the second from Lemma \ref{near_clip_stability}, and the third from Lemma \ref{near_sequence_stability}.

\begin{lemma}[Near True Graph, Near Clip and Near Sequence Stability]
\label{total_stability}
Let $\Delta$ be a side of the maximal connected critical trio $\Theta$.
\begin{enumerate}
\item If $\Delta$ is a near true graph, then $\Theta$ is a video.
\item If $\Delta$ is a near clip, then $\Theta$ is a video.
\item If $\Delta$ is a near sequence, then $\Theta$ is either a pure chord, a cyclic chord, 
or a dihedral chord of type $1$, $2A$, $2B$, or $2C$.
\end{enumerate}
\end{lemma}

For a somewhat technical purpose in the proof, we require one additional stability result.  In a duet 
$(X,Y)$ we say that a block of imprimitivity $T \subseteq X$ or $T \subseteq Y$ is \emph{doubling} if $T$ is critical and $w(N(T)) = 2 w(T)$.  Doubling blocks naturally give rise to graphs, and the following lemma follows from Lemma \ref{near_graph_stability}.  

\begin{lemma}[Doubling Block Stability]
\label{doublingblock}
Let $\Delta$ be a side of the maximal connected critical trio $\Theta$.  If $\Delta$ has a doubling block, then either $\Theta$ is a 
video or $\Delta$ is a near polygon.
\end{lemma}

In addition to these stability lemmas, we will require some lemmas which give some control over the sides of a critical trio $\Theta$ with the property that $\Theta \le \Theta^*$ for a song $\Theta^*$.
The following lemma fits this need.  It will be proved by Lemmas  \ref{light_sides_video} and Lemma \ref{light_sides_chord}.  We say that a side $(X,Y)$ of a trio $\Theta$ is \emph{light} if there is a side of $\Theta$ other than $(X,Y)$ and $(Y,X)$ which has weight at least $w(X,Y)$.  

\begin{lemma}[Light Sides]
\label{light_sides}
Let $\Theta$ be a connected critical trio and assume that $\Theta \le \Theta^*$ for a song $\Theta^*$.  
If $(X,Y)$ is a light side of $\Theta$, then either $(X,Y)$ or $(Y,X)$ is a near true graph, a near clip, or a near sequence. 
\end{lemma}

\subsection{The Proof}

In this subsection we prove our main theorem.  We will say that a cross $(A,B)$ of a duet $(X,Y)$ is \emph{connected} if $[A] = X$ and 
$[B] = Y$.  Note that whenever $(X,Y)$ is connected, and $(A,B)$ is connected, the trio associated with this cross will also be connected.

\bigskip

\noindent{\it Proof of the ``only if'' direction of Theorem \ref{main_incidence}:} 
Let $\Theta = (X,Y,Z)$ be a nontrivial maximal critical $G$-trio, let $\Delta_1 = (X,Y)$, $\Delta_2 = (Z,X)$, and $\Delta_3 = (Z,Y)$ and assume that $w(\Delta_1) \le w(\Delta_2) \le w(\Delta_3)$.  We shall assume that $\Theta$ is a counterexample to Theorem \ref{main_incidence} which is chosen according to the following criteria.
\begin{enumerate}
\item if there is a finite counterexample, then $|G|$ is minimum.
\item $\overline{w}(\Delta_3)$ is minimum (subject to 1).
\item $|\{ i \in \{1,2,3\} \mid \mbox{$\Delta_i$ is connected} \}|$ is maximum (subject to 1 and 2).
\end{enumerate}

By passing from $\Theta$ to $\Theta^{\bullet}$ we may assume that $\Theta$ is clone free.  Our first task will 
be to establish a sequence of claims.

\bigskip

\noindent{(i)} $\Theta$ is connected.

\smallskip

Suppose (for a contradiction) that $\Theta$ is disconnected, and apply Lemma \ref{disconnected}.  If $\Theta$ is a pure beat, then we have an immediate contradiction.  Otherwise, $\Theta$ is an impure beat, and we let $\Theta'$ be an $H$-trio which is a continuation of $\Theta$.  If $H$ is finite, then 
$\Theta'$ contradicts the choice of $\Theta$ for the first criteria.  Otherwise, $\Theta'$ and $\Theta$ will be the same with regard to the first two criteria, but 
$\Theta'$ contradicts the choice of $\Theta$ for the third.  

\bigskip

\noindent{(ii)} $\Delta_1$ is not a near true graph, a near clip, or a near sequence.  

If $\Delta_1$ is a near true graph or a near clip, then it follows from Lemma \ref{total_stability} that $\Theta$ is a 
video which is a contradiction.  If $\Delta_1$ is a near sequence, then by the same lemma, $\Theta$
is either a pure chord, cyclic chord, or dihedral chord of type $1$, $2A$, $2B$, or $2C$.  For the pure chord 
and dihedral chords of type $2B$, and $2C$ we have an immediate contradiction.  In the remaining cases, 
a continuation of $\Theta$ must be a smaller counterexample to our theorem, which again contradicts our choice of 
$\Theta$.

\bigskip

\noindent{(iii)} $\Delta_1$ is not a near polygon.

\smallskip

Suppose (for a contradiction) that $(X,Y)$ is a  near polygon relative to $( {\mathcal P}, {\mathcal Q} )$.  It follows
from (ii) that $|{\mathcal P}| \le 3$ (all other polygons are also sequences).  If $|{\mathcal P}| = 2$ then we have $w(\Delta_3) \ge w(\Delta_2) \ge w(\Delta_1) > \frac{1}{2}|G|$ which contradicts part 2 of Lemma \ref{nontriv_bound}.  Thus, we may assume that 
${\mathcal P} = \{P_1,P_2,P_3\}$ and ${\mathcal Q} = \{Q_1, Q_2, Q_3\}$ where 
$P_i \not\sim Q_i$ for $i=1,2,3$ (so then $P_i \sim Q_j$ whenever $1 \le i,j \le 3$ satisfy $i \neq j$).  Next, choose 
$z \in Z$ and let $(A,B)$ be the cross of $\Delta_1$ associated with $z$.  Since $[A] = X$ and $[B] = Y$ the set $A$ cannot be contained in a block of ${\mathcal P}$ and $B$ cannot be contained in a block of ${\mathcal Q}$.  It follows that there must exist $1 \le i \le 3$ so that $A' = A \cap P_i \neq \emptyset$ and $B' = B \cap Q_i \neq \emptyset$.  Set 
$A'' = A \setminus A'$ and $B'' = B \setminus B'$ and observe that $N(A')$ and $B''$ are disjoint subsets of $Y \setminus Q_i$.  
So, setting $k = w^{\circ}({\mathcal P}) = w^{\circ}({\mathcal Q})$ we have $w(N(A')) \le 2k - w(B'')$.  The next inequality uses this together with Lemma \ref{block_best}.
\begin{align*}
0	&\le	\delta(P_i) - \delta(A')	\\
	&=	w(P_i) - w(N(P_i)) - w(A') + w(N(A'))	\\
	&\le	-k -w(A') + (2k - w(B''))	\\
	&=	k - w(A') - w(B'').
\end{align*}
So $w(A') + w(B'') \le k$ and similarly $w(B') + w(A'') \le k$.  However, then $w(\Delta_2) + w(\Delta_3) = w(A) + w(B) \le 2k < 2w(\Delta_1)$ and this is a contradiction.

\bigskip

\noindent{(iv)} If a block of imprimitivity $T \subset X$ or $T \subset Y$ is critical in $\Delta_1$, then $w(T) < \frac{1}{2}w(\Delta_1)$.

\smallskip

It follows from Proposition \ref{block_mult} and the assumption that $\Delta_1$ is connected that every proper critical block of imprimitivity $T$ of $X$ or $Y$ must satisfy $w(T) < w(\Delta_1)$.  If there exists such a block with $w(T) \ge \frac{1}{2}w(\Delta_1)$ then $T$ is a doubling block in $\Delta_1$.  Then by Lemma \ref{doublingblock}, we find that either $\Theta$ is a video, which is a contradiction, or $\Theta$ is a near polygon, which contradicts (iii).

\bigskip

\noindent{(v)}  If $(A,B)$ is a critical cross of $\Delta_1$ and $\min\{ w(A), w(B) \} \ge \frac{1}{2}w(\Delta_1)$ then $(A,B)$ is connected.

\smallskip

If we suppose (for a contradiction) that $A$ is finite and $[A] \neq X$, then Theorem \ref{critical_closure} implies that 
$[A]$ is critical in $(X,Y)$ and then we have $w([A]) \ge w(A) \ge \frac{1}{2}w(\Delta_1)$ which contradicts (iv).

\bigskip

With these claims in hand, we are ready for the heart of the argument.  
Let ${\mathcal P}$ and ${\mathcal Q}$ be maximal proper systems of imprimitivity on $X$ and $Y$ for which the blocks are finite and critical (they are permitted to be trivially critical) in $\Delta_1$.  (such systems of imprimitivity exist since the discrete partition is one, and there is a maximal one by (iv)).

Let $z \in Z$, let $(A,B)$ be the cross of $\Delta_1$ associated with $z$ and suppose (for a contradiction) that there is a boundary block $Q \in {\mathcal Q}$ with respect to $B$.  Let $(A',B')$ be the new cross obtained from purifying $(A,B)$ at 
$Q$.  It follows from Theorem \ref{purify_thm} that $w(A') > w(A) - w^{\circ}({\mathcal Q}) \ge \frac{1}{2}w(\Delta_1)$ and 
trivially, $w(B') \ge w(B) \ge w(\Delta_1)$.  Therefore, by (v) we have that the trio $\Theta'$ which is associated with 
$(A',B')$ is also connected.  Let $\Theta''$ be a maximal trio with $\Theta' \le \Theta''$.  Now, it follows from our choice of $\Theta$ that $\Theta''$ is not a counterexample, so it is a connected song.  It now follows from applying Lemma \ref{light_sides} to $\Theta'$ that $\Delta_1$ is either a near true graph, a near clip, or a near sequence, but this contradicts (ii).  It follows that there is no boundary block of ${\mathcal Q}$ with respect to $B$.  However, then by the maximality of $\Theta$ we must have that $\mcq$ is a system of clones.  Since, by assumption, $\Theta$ is clone free, we conclude that ${\mathcal Q}$ is the discrete partition.  

Ignoring our earlier assignment, let us now choose $(A,B)$ to be a connected cross of $\Delta_1$ so that
\begin{enumerate}
\item[$\alpha$)] $\delta(A,B)$ is maximum.
\item[$\beta$)] the number of boundary blocks in ${\mathcal P}$ w.r.t. $A$ is minimum (subject to $\alpha$).
\item[$\gamma$)] $\min\{ w(B) + w^{\circ}({\mathcal P}), w(A)+ w^{\circ}(Y)  \}$ is minimum (subject to $\alpha$ and $\beta$).
\end{enumerate}
Note that by using a cross associated with some point $z \in Z$ we must have
$\delta(A,B) \ge \delta(\Theta)$.  Furthermore, it follows immediately from $\alpha$ that $(A,B)$ is maximal.  

Suppose (for a contradiction) that $w(A) < w(\Delta_1)$.  Then we have
\[ \overline{w}(B) =  w(\Delta_1) + w(A) - \delta(A,B) < w(\Delta_1) + w(\Delta_2) - \delta(\Theta) = \overline{w}(\Delta_3). \]
Now, let $\Theta'$ be the trio associated with $(A,B)$ and let $\Theta''$ be a maximal trio with $\Theta' \le \Theta''$.  It 
follows from the above equation and our choice of $\Theta$ that $\Theta''$ is a song.  But then applying Lemma \ref{light_sides} to $\Theta'$ shows that $\Delta_1$ is either a near true graph, a near clip, or a near sequence, and this contradicts (ii).  We obtain a similar contradiction under the assumption that $w(B) < w(\Delta_1)$, so we conclude that
$w(A), w(B) \ge w(\Delta_1)$.

Next, suppose (for a contradiction) that there is a boundary block $P \in {\mathcal P}$ with respect to $A$.  Let 
$(A',B')$ denote the cross obtained from $(A,B)$ by purifying at $P$.  It follows from Theorem \ref{purify_thm}, the above argument and (iv) that $w(B') \ge w(B) - w^{\circ}({\mathcal P}) \ge w(\Delta_1) - \tfrac{1}{2}w(\Delta_1)$.  It now follows from (v) that $[B'] = Y$.  Since $A' \supseteq A$ we also have $[A'] = X$ and now $(A',B')$ contradicts the choice of $(A,B)$.  
It follows that $A$ is pure with respect to ${\mathcal P}$.  We now prove two properties which will be used to handle
the two possibilities for which term in optimization criterion $\gamma$ is minimum.  Since these arguments are quite similar, we only sketch the second.

\bigskip

\noindent{$(\star)$} \; If  $w(A) + w^{\circ}(Y) \le w(B) + w^{\circ}({\mathcal P})$ then 
	$w(gA \cap A) \le w^{\circ}({\mathcal P})$ for every $g \in G$ with $gA \neq A$.

\smallskip

Note that it follows from the above inequality that $A$ is finite.  Suppose (for a contradiction) that there exists 
$g \in G$ so that $2 w^{\circ}({\mathcal P}) \le w(gA \cap A) < w(A)$.  Let $(A',B') = (gA \cap A, gB \cup B)$ and 
$(A'', B'') = (gA \cup A, gB \cap B)$.  If $|B''| \le 1$ then $B$ is also finite and we have
\begin{align*}
w(A') + w(B') 
	&\ge 2 w^{\circ}({\mathcal P}) + 2 w(B) - w^{\circ}(Y) 	\\
	&\ge w^{\circ}({\mathcal P}) + w(A) + w(B).
\end{align*}
This implies that $\delta(A',B') > \delta(A,B)$.  It now follows from the maximality of 
${\mathcal P}$ and Theorem \ref{critical_closure} that $[A'] = X$.  However, 
then $(A',B')$ is connected and contradicts our choice.  It follows that we must have $|B''| \ge 2$.  Now, since 
$Y$ has no nontrivial critical blocks of imprimitivity, Theorem \ref{critical_closure} implies that $[B''] = Y$ so 
$(A'',B'')$ will also be a connected cross.  However, then one of $(A',B')$ or $(A'',B'')$ contradicts our choice of $(A,B)$.  

\bigskip

\noindent{$(\star \star)$} \; If $w(B) + w^{\circ}({\mathcal P}) \le w(A) + w^{\circ}(Y)$ then 
	$w(gB \cap B) \le w^{\circ}(Y)$ for every $g \in G$ with $gB \neq B$.

\smallskip

Note that it follows form the above inequality that $B$ is finite.  Arguing as in the previous case, we suppose (for a contradiction) that there exists $g \in G$ with $2 w^{\circ}(Y) \le w(gB \cap B) < w(B)$.  Let $(A',B') = (gA \cup A, gB \cap B)$ and $(A'',B'') = (gA \cap A, gB \cup B)$.  If $w(A'') \le w^{\circ}({\mathcal P})$ then $A$ is also finite and 
\begin{align*}
w(A') + w(B') 
	&\ge 2 w^{\circ}(Y) + 2 w(A) - w^{\circ}({\mathcal P}) 	\\
	&\ge w^{\circ}(Y) + w(A) + w(B).
\end{align*}
so as before we have $\delta(A',B') > \delta(A,B)$.  By Theorem \ref{critical_closure} we find that $(A',B')$ is connected and this contradicts our choice.  So, it must be that $w(A'') \ge 2 w^{\circ}({\mathcal P})$.  But then both 
$(A',B')$ and $(A'',B'')$ are connected and one contradicts our choice of $(A,B)$.  It follows that every $g \in G$ for which 
$gB \neq B$ must satisfy $w( gB \cap B) \le w^{\circ}(Y)$.  

\bigskip

Now, in either case we let $\Theta' = (X,Y,Z' ; \sim')$ be the trio associated with $(A,B)$ and let $\Theta'' = (X,Y,Z' ; \sim'')$ be a maximal trio with $\Theta' \le \Theta''$.  It follows from the maximality of $(A,B)$ that 
$\sim'|_{X \cup Z'} = \sim''|_{X \cup Z'}$ and $\sim'|_{Y \cup Z'} = \sim''|_{Y \cup Z'}$.  If $x,x' \in P$ for some $P \in {\mathcal P}$ then $N_{\Theta'}(x) \cap Z' = N_{\Theta'}(x') \cap Z'$ by construction, so this also holds in the trio $\Theta''$, but it follows from this that $x$ and $x'$ will be clones in $\Theta''$.  Therefore, $\mcp$ is a system of clones in $\Theta''$.  Let $\Theta^*$ be the trio obtained from $\Theta''$ by taking the quotient 
$( \mcp, Y, Z')$ and note that by property $(\star)$ or $(\star \star)$, the trio $\Theta^*$ has a side which is Sidon.  It now follows from Lemma \ref{sidon_stability} that $\Theta^*$ is either a video or a pure chord, and thus $\Theta''$ is either a video or a pure chord.  However, then applying Lemma \ref{light_sides} to $\Theta'$ reveals that either $\Delta_1$ or its reverse is  a near true graph, a near clip, or a near sequence, and this contradicts (ii).  This completes the proof.
\quad\quad$\Box$

\section{Graphs}
\label{graph_sec}

In this section we prove a stability lemma for trios with a graphic side, and a light side lemma for videos.   
The proof of this stability lemma has two principal parts.  The first part is a 
classification of certain small edge cuts in vertex and edge transitive graphs, and this is 
the content of the first subsection.  The second part requires an analysis of the trios which 
can arise from these, and this is done in the second subsection.

\subsection{Small Edge Cuts}

In this subsection we will prove a theorem which characterizes all edge cuts of size less than 
$2d$ (with at least one side finite) in a vertex and edge transitive graph of degree $d$.  Our techniques are the same as those 
used by Mader, Watkins, Hamidoune and others to achieve similar connectivity results in vertex transitive 
graphs.  We begin with a proposition which sharpens Observation \ref{2d_obs} and motivates this characterization.  Recall that for a subset of vertices $A$ we let $c(A)$ denote the number of edges with exactly one endpoint in $A$.

\begin{lemma}
\label{crit2cut}
Let $\Gamma = (V,E)$ be a graphic duet of degree $d$, and let 
$A \subseteq V$ be finite or cofinite.  Then
\[ \delta(A) = \frac{w^{\circ}(E)}{2} \big(2d - c(A) \big). \]
\end{lemma}

\noindent{\it Proof:}  First note that $w(\Gamma) = d w^{\circ}(E) = 2w^{\circ}(V)$.  
Let $S$ denote the set of edges with both ends in $V \setminus A$.  If $A$ is cofinite, then summing 
degrees in $V \setminus A$ gives us  $d | V \setminus A | = 2|S| + c(A)$ which implies
\[ \delta(A) = w(S) + w(\Gamma) - w(V \setminus A) = w^{\circ}(E) \big( |S|  + d  - \tfrac{d}{2} |V \setminus A| \big) = \tfrac{w^{\circ}(E)}{2} \big( 2d - c(A) \big). \]
Otherwise, $A$ is finite, and setting $R$ to be the set of edges with both ends in $A$ we have
$d|A| = 2|R| + c(A)$ (by counting degree sums in $A$).  This gives us 
\[ \delta(A)	
			=	w(A) + w(\Gamma) - w(R \cup \partial(A))	
			=	w^{\circ}(E) \big( \tfrac{d}{2} |A| + d  -  |R| - c(A) \big)	
			=	\tfrac{w^{\circ}(E)}{2} \big( 2d - c(A) \big) \]
as desired.
\quad\quad$\Box$

\bigskip

Next we introduce some graph theoretic terminology.  
The {\it girth} of a graph $\Gamma$ is the length of the shortest cycle, or $\infty$ if it has no cycle.  The \emph{line graph} of $\Gamma$, 
denoted $L(\Gamma)$, is a simple graph with vertex set $E$ with the property 
that $e,f \in E$ are adjacent (in $L(\Gamma)$) if these edges are incident with a 
common vertex (in $\Gamma$).  We define $\mathcal{LC}$ to be the class of all 
graphs which are isomorphic to $L(\Gamma)$ for some cubic graph $\Gamma$ 
with girth at least four.  If $\Gamma$ is a regular graph with girth $g$ and there exists an integer $s$ so that 
every edge of $\Gamma$ is contained in exactly $s$ cycles with length $g$, then 
we call $\Gamma$ an {\it equitable} graph.  Note that every vertex and edge transitive graph is equitable.  
Next we offer a basic proposition concerning
equitable graphs.  The proof follows from a simple case analysis, so we have only provided a sketch.

\begin{proposition}
\label{equitable_prop}
The following chart is a complete list (up to isomorphism) of all connected
equitable graphs with the indicated degree and girth parameters.
\begin{center}
\begin{tabular}{|c|c|p{2.5in}|}
\hline
Degree	&	$\mbox{ }$Girth$\mbox{ }$		&	Equitable graphs \\
\hline
3	&	3	&	Tetrahedron	\\
\hline
3	&	4	&	Cube, $K_{3,3}$	\\
\hline
4	&	3	&	Octahedron, $K_5$, members of $\mathcal{LC}$ \\
\hline
3	&	5	&	Dodecahedron, Petersen	\\
\hline
5	&	3	&	Icosahedron, $K_6$	\\
\hline
\end{tabular}
\end{center}
\end{proposition}

\noindent{\it Sketch of Proof:}
Let $\Gamma$ be an equitable graph with degree $d$, girth $g$, and
with exactly $s$ cycles of length $g$ containing each edge.  
Since every cycle which contains a vertex $u$ must contain exactly 
two edges incident with $u$, it follows that $s$ must be even whenever 
$d$ is odd.  Let us note a couple of easy facts which are helpful in the forthcoming 
analysis.  First, observe that if $g = 3$ then $s \le d-1$ and 
if $s=d-1$ we must have $\Gamma \cong K_{d+1}$.  Next suppose that $s = 2$ 
and either $d = 3$, or $g = 3$ and the set of vertices adjacent to a given vertex 
induce a cycle.  In either of these cases, we can form a polygonal surface from our 
graph by identifying each cycle of length $g$ with a regular $g$-gon.  In 
these cases, it may be helpful to argue in terms of this surface.

If $d=3$ and $g=3$, then $s=2$ and $\Gamma \cong \mathrm{Tetrahedron}$.  
If $d=3$ and $g=4$, then either $s=2$ and $\Gamma \cong \mathrm{Cube}$ or 
$s \ge 4$.  In the latter case, it follows easily from considering the 
subgraph induced by two adjacent vertices and any vertex adjancent to one
of them that $\Gamma \cong K_{3,3}$.
If $d=3$ and $g=5$, then either $s=2$ and $G \cong \mathrm{Dodecahedron}$, 
or $s \ge 4$.  In the latter case, it follows easily from considering the 
graph induced by all vertices of distance at most two from a fixed vertex
that $s=4$ and $G \cong \mathrm{Petersen}$.  If $d=4$, $g=3$, and $s=1$, then  
define a new graph $\Upsilon$ with a vertex for each triangle of $\Gamma$ and 
with two vertices adjacent if and only if the corresponding triangles share
a vertex.  Then $\Upsilon$ is cubic with girth at least four, and $\Gamma \cong L(\Upsilon)$, 
so $\Gamma$ is in $\mathcal{LC}$.  If $d=4$, $g=3$, and $s=2$ then 
$G \cong \mathrm{Octahedron}$, and if $d=4$, $g=3$ and $s=3$, then 
$G \cong K_5$.  Finally, if $d=5$ and $g=3$, then either 
$s=2$ and $G \cong \mathrm{Icosahedron}$ or $s=4$ and $G \cong K_6$.
\quad\quad$\Box$

\bigskip

Now we are ready for the main theorem of this subsection which classifies small edge cuts 
in vertex and edge transitive graphs.  The key property utilized by our argument is the 
following submodularity inequality which holds for every graph $\Gamma$ and $A,B \subseteq V(\Gamma)$
with $c(A), c(B) < \infty$.
\begin{equation}
\label{ec_submod}
c(A \cap B) + c(A \cup B) \le c(A) + c(B).
\end{equation}

\begin{theorem}
\label{cutset_classify}
Let $\Gamma = (V,E)$ be a simple connected vertex and edge-transitive graph of degree $d$ and 
let $A \subseteq V$ be finite and nonempty.  If $c(A) < 2d$ and $|A| \le \frac{1}{2}|V|$, then setting $\Gamma'$ to 
be the subgraph induced by $A$ one of the following holds.
\begin{enumerate}
\item		$\Gamma'$ is a single vertex.
\item		$\Gamma'$ is a one edge path.
\item		$\Gamma$ is a cycle and $\Gamma'$ is a path.
\item		$\Gamma$ is cubic and $\Gamma'$ is a 2-edge path.
\item		$\Gamma$ is in $\mathcal{LC}$ and $\Gamma'$ is a triangle.
\item		$\Gamma$ is isomorphic to Cube, Octahedron, 
		Dodecahedron, Icosahedron, Petersen, or $K_6$, and $\Gamma'$ is a 
		shortest cycle of $\Gamma$.
\end{enumerate}
\end{theorem}

\noindent{\it Proof:} Let $G \le Aut(\Gamma)$ act transitively on $V$ and $E$.  Now suppose (for a contradiction) that the theorem is fails for $\Gamma$, and choose a counterexample $A$ subject to the following conditions.  
\begin{enumerate}
\item[$\alpha$)] $c(A)$ is minimum.
\item[$\beta$)] $|A|$ is minimum (subject to $\alpha$).
\end{enumerate}
Note that the minimality of our counterexample implies that the graph $\Gamma'$ 
is connected.  We proceed with a sequence of claims.  

\bigskip

\noindent{(i) } $d \ge 3$, $|A| \ge 4$, and $\Gamma$ is not isomorphic
to Cube, Octahedron, Dodecahedron, Icosahedron, Petersen, or $K_6$.  In particular, 
if $d=3$ then $\Gamma$ has girth at least 6.

\smallskip

We leave the reader to verify the cases when $d \le 2$ or $|A| \le 3$ or $\Gamma$ is one of 
the graphs in the above list.  In light of Proposition \ref{equitable_prop} we may then assume
that the girth of $\Gamma$ is at least 6 in the case when $d=3$.

\bigskip

\noindent{(ii) } ${\mathit deg}_{\Gamma'}(x) > \frac{1}{2}{\mathit deg}_{\Gamma}(x)$ 
	for every $x \in A$.

\smallskip

If (ii) is false for $x$, then the theorem holds for $A \setminus \{x\}$ and we must 
have outcome 4 or 5.  Either case yields a contradiction.

\bigskip

\noindent{(iii) } If $3 \le |gA \cap A| < |A|$ for $g \in G$, then
	$|gA \cap A| = 3$ and either $d=3$ or $\Gamma \in \mathcal{LC}$.

\smallskip 

Let $k = |gA \cap A| \ge 3$ and note that 
$|V \setminus (gA \cup A)| \ge |gA \cap A| = k$.  If $c(gA \cup A) < c(A)$ then by our choice of 
counterexample, the theorem holds for either $gA \cup A$ or its complement.  Only outcomes 4 and 5 
are possible, and in either case it follows that $k=3$ so (iii) holds.  
Otherwise $c(gA \cup A) \ge c(A)$ and then by equation \ref{ec_submod} we conclude 
that $c(gA \cap A) \le c(A)$.  So, the theorem holds for $gA \cap A$ and the desired
conclusion follows immediately.  

\bigskip

\noindent{(iv) } $\Gamma$ is not cubic.

\smallskip

Suppose (for a contradiction) that $\Gamma$ is cubic.  Observe that 
by (ii) we must have $| \{ v \in V \mid {\mathit deg}_{\Gamma'}(v) = 2 \} | = c(A) < 6$.  
If the graph $\Gamma'$ has at most two vertices with degree 3, then $\Gamma$ has girth 
at most 5 and Proposition \ref{equitable_prop} gives us a contradiction to (i).  Otherwise, 
we may choose a pair of adjacent vertices $u,u' \in A$ and a pair of adjacent vertices $v,v' \in A$ so that 
${\mathit deg}_{\Gamma'}(u) = {\mathit deg}_{\Gamma'}(u') 
	= {\mathit deg}_{\Gamma'}(v) = 3$ and 
${\mathit deg}_{\Gamma'}(v') = 2$.  Since $G$ acts transitively on 
$E$ we may choose $g \in G$ so that $g\{u,u'\} = \{v,v'\}$.  Then we have 
$4 \le |gA \cap A| < |A|$ which contradicts (iii).  

\bigskip

\noindent{(v) } $\Gamma$ is not in $\mathcal{LC}$.

\smallskip

Suppose (for a contradiction) that $\Gamma \cong L(\Upsilon)$ for a cubic
graph $\Upsilon$ with girth at least four.  Observe that the vertex set of $\Upsilon$ 
corresponds to $C_3(\Gamma)$ and the edge set of $\Upsilon$ corresponds to $V$.  
It follows that $\Upsilon$ is both vertex and edge transitive.

\begin{figure}[ht]
\centerline{\includegraphics[width=2.7cm]{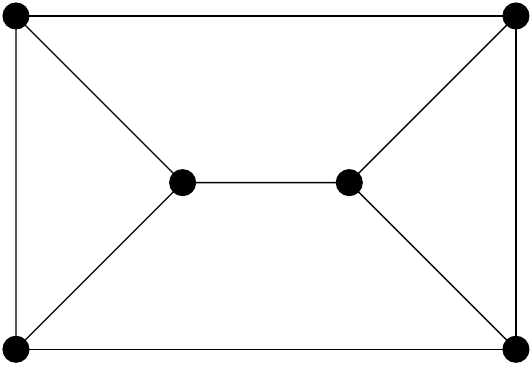}}
\caption{Prism}
\label{prism_fig}
\end {figure}

If there exists $u \in A$ with 
${\mathit deg}_{\Gamma'}(u) = 4$, then (by (ii)) we may choose $v \in A$ with 
${\mathit deg}_{\Gamma'}(v) = 3$ and $g \in G$ with $gu = v$.  However then we have 
$4 \le |gA \cap A| < |A|$ which contradicts (iii).  It follows that 
${\mathit deg}_{\Gamma'}(u) = 3$ for every $u \in A$.  Now we have $|A| = c(A) < 8$, so 
$\Gamma'$ is a cubic graph on fewer than 8 vertices, and it must then be  
isomorphic to one of $K_4$, $K_{3,3}$, or Prism (depicted in Figure \ref{prism_fig}).  
Observe that $\Gamma'$ has a 4-cycle, so $\Upsilon$ must have girth 4.  However, then  
Proposition \ref{equitable_prop} implies that $\Upsilon$ is isomorphic to either
Cube or $K_{3,3}$.  Thus $\Gamma \cong L(\mathrm{Cube})$ or 
$\Gamma \cong L(K_{3,3})$ and either possibility contradicts the structure of 
$\Gamma'$.

\bigskip

\noindent{(vi) } Let $u,v \in A$ and $g \in G$ satisfy $gu = v$ and $gA \neq A$.  Then $d$ is odd
 and further ${\mathit deg}_{\Gamma'}(u) = {\mathit deg}_{\Gamma'}(v) = \frac{d+1}{2}$.

If $u,v$ and $g$ satisfy the conditions above, then by (iii),(iv),(v) it
follows that $|A \cap gA| \le 2$ so 
both $g\{ w \in A \mid w \sim u \}$ and $\{ w \in A : w \sim v \}$ are
subsets of $\{ w \in V : w \sim v \}$ with size greater than 
$\frac{d}{2}$ and intersection of size at most $1$.  The above conclusion 
follows immediately from this.

\bigskip

If there exist $u,v \in A$ with 
${\mathit deg}_{\Gamma'}(u) \neq {\mathit deg}_{\Gamma'}(v)$, then 
choosing $g \in G$ so that $gu = v$ gives us a contradiction to (vi).  Therefore $\Gamma'$ 
is regular.  Now, $A$ cannot be a block of imprimitivity of $V$ since there are edges with both 
ends in $A$ and edges with one end in $A$ (and $G$ acts transitively on $E$).  It follows that we may choose $g \in G$ so that 
$\emptyset \neq gA \cap A \neq A$.  From here, (vi) implies that $d$ is odd and $\Gamma'$ is 
$\frac{d+1}{2}$ regular.  Now by (iv) we have $d \ge 5$ so 
$|A| = c(A) \frac{2}{d-1} < 4 \frac{d}{d-1} \le 5$.  The only possibility here is $d=5$ and $\Gamma' \cong K_4$, 
but then Proposition \ref{equitable_prop} implies that $\Gamma \cong K_6$ which contradicts (i).
\quad\quad$\Box$

\bigskip

A quick check of the statement of the previous theorem reveals that $c(A) \ge 2d - 2$ in all cases except for the
first outcome.  The next corollary follows immediately from Lemma \ref{crit2cut} and this observation.

\begin{corollary}
\label{2dm2_cor}
Let $\Gamma = (V,E)$ be a graphic duet of degree $d$.  If $A \subseteq V$ is finite or cofinite then either 
$\min \{ |A|, |V \setminus A| \} \le 1$ or $\delta(A) \le w^{\circ}(E)$.
\end{corollary}

\subsection{Stability}

The main theorem from the previous subsection will be used to characterize all maximal critical crosses $(A,B)$ of graphic duets $(V,E)$.  Here we shall determine which trios arise from such crosses by studying their behaviour under the group action.  We begin with an easy lemma concerning the transitivity of duets.  We say that a $G$-duet $\Delta = (X,Y)$ is 
\emph{incidence transitive} if $G$ acts transitively on the set $\{ (x,y) \in X \times Y \mid x \sim y \}$.

\begin{lemma}
\label{inc_trans}
If $\Delta = (X,Y)$ is a $G$-duet and $d_{\Delta}(X,Y)$, 
$d_{\Delta}(Y,X)$ are finite and relatively prime, then $\Delta$ is 
incidence transitive.
\end{lemma}

\noindent{\it Proof:}
Let $S \subseteq \{ (x,y) \in X \times Y \mid x \sim y \}$ be a 
$G$-orbit and define a new symmetric incidence relation $\sim'$
on $X \cup Y$ by the rule that $x \sim' y$ whenever $(x,y) \in S$.  
It follows immediately that  $\Delta' = (X,Y)$ equipped with the incidence 
relation $\sim'$ is a $G$-duet.  Next observe
\[ \frac{d_{\Delta}(X,Y)}{ d_{\Delta}(Y,X)} = \frac{w^{\circ}(X) }{ w^{\circ}(Y)} 
	= \frac{d_{\Delta'}(X,Y) }{ d_{\Delta'}(Y,X)}. \]
Now $d_{\Delta}(X,Y)$ and $d_{\Delta}(Y,X)$ are relatively prime, 
$d_{\Delta'}(X,Y) \le d_{\Delta}(X,Y)$, and 
$d_{\Delta'}(Y,X) \le d_{\Delta}(Y,X)$.  It follows that 
$d_{\Delta'}(X,Y) = d_{\Delta}(X,Y)$ and 
$d_{\Delta'}(Y,X) = d_{\Delta}(Y,X)$, so $\Delta = \Delta'$ and we conclude that 
$\Delta$ is incidence transitive.
\quad\quad$\Box$

\bigskip

For a graph $\Gamma = (V,E)$ an ordered pair of adjacent vertices is called an \emph{arc}, 
and we say that $\Gamma$ is \emph{arc transitive} if $Aut(\Gamma)$ acts transitively on 
the set of arcs.  

\begin{lemma}
\label{p3_trans}
Let $\Gamma = (V,E)$ be a graphic $G$-duet with degree $3$.  Then 
$G$ acts transitively on $P_3(\Gamma)$.
\end{lemma}

\noindent{\it Proof:}
It follows from the previous lemma that $\Gamma$ is $G$-incidence transitive, so 
$G$ acts transitively on the arcs of $\Gamma$.  Next we establish a correspondence 
between $P_3(\Gamma)$ and the set of arcs of $\Gamma$ 
as follows.  We associate the path with vertex sequence $x,y,z$ to the arc $(y,y')$ where 
$y'$ is unique vertex adjacent to $y$ other than $x,z$.  Since the action of $G$ on $\Gamma$ 
preserves this correspondence, it follows that $G$ acts transitively on $P_3(\Gamma)$ 
as desired.
\quad\quad$\Box$

\bigskip

If $\Gamma$ is a graph which is isomorphic to Petersen ($K_6$), 
we say that a \emph{face set for} $\Gamma$ is a subset $F$ of $C_5(\Gamma)$ 
($C_3(\Gamma)$) with the property that every edge is contained in exactly two 
members of $F$.  This set of faces then yields an embedding of $\Gamma$ in the 
projective plane as shown in Figure \ref{k6pete}.

\begin{proposition}
\label{face_orbit}
Let $\Gamma = (V,E)$ be a graphic $G$-duet with girth $g$.  Then we have
\begin{enumerate}
\item		If $\Gamma$ is isomorphic to Cube, Icosahedron, or Dodecahedron, then
		$G$ acts transitively on $C_g(\Gamma)$.
\item		If $\Gamma \cong \mathrm{Octahedron}$, then either $G$ 
		acts transitively on $C_3(\Gamma)$ or $C_3$ consists of 
		two $G$-orbits with every edge incident with exactly one 
		cycle from each orbit	.
\item		If $\Gamma$ is isomorphic to either Petersen or 
		$K_6$, then either $G$ acts transitively on $C_g(\Gamma)$ or 
		$C_g(\Gamma)$ contains two $G$-orbits, each of which is a face
		set for $\Gamma$.
\end{enumerate}
\end{proposition}

\noindent{\it Proof:} 
Let $S \subseteq C_g(\Gamma)$ be a $G$-orbit.  Since $G$ acts transitively 
on $E$, it follows that every $e \in E$ is incident with the same 
number, say $t$, of cycles in $S$.  If the degree of $\Gamma$ is odd, then 
$t$ must be even, since every such cycle incident with a vertex $v$ contains 
two edges from $\partial(v)$.  This easily implies the result.  
\quad\quad$\Box$

\begin{lemma}
\label{path_trans}
Let $\Gamma = (V,E)$ be a graph which is either a cycle or a two-way infinite path.  If $G \le {\mathit Aut}(\Gamma)$ acts transitively on $V$ and $E$, then $G$ acts transitively on $P_k(\Gamma)$ for every $k$.
\end{lemma}

\noindent{\it Proof:} 
For $i=1,2$ let $\Gamma_i$ be a  path of length $k$ in $\Gamma$ given by the vertex-edge 
sequence $v^i_1, e^i_1, v^i_2, \ldots, e^i_{k-1}, v_k^i$ (so each $e^i_j$ is an edge incident with the 
vertices $v^i_j$ and $v^i_{j+1}$).  It is sufficient to prove that there exists $g \in G$ with $g \Gamma_1 = \Gamma_2$.  
If $k$ is odd, say $k = 2j+1$, then since $G$ acts vertex transitively there 
exists $g \in G$ with $g v^1_{j+1} = g v^2_{j+1}$ (so $g$ maps the middle vertex of $\Gamma_1$ to the 
middle vertex of $\Gamma_2$) and then $g \Gamma_1 = \Gamma_2$.  If $k$ is even, say $k = 2j$ then 
we may choose $g \in G$ with $g e^1_{j-1} = e^2_{j-1}$ (so $g$ maps the middle edge of $\Gamma_1$ to the 
middle edge of $\Gamma_2$) and then $g \Gamma_1 = \Gamma_2$.  
\quad\quad$\Box$

\bigskip

We have now collected the required information concerning the actions of $G$ on various structures in 
graphic $G$-duets for our stability lemma.  However we need to take care of one more detail before we are 
ready for this proof.  In the previous subsection we saw small edge cuts arising from triangles in graphs from ${\mathcal LC}$.  However, 
this class of graphs does not appear in our definition of videos.  This disappearance is due to the following observation
showing that these trios are isomorphic to other videos.  

\begin{observation}
\label{cubiclg_equiv}
Let $\Upsilon$ be a cubic graph with girth at least $4$ and let $\Gamma = L(\Upsilon)$.  Then 
\[ \Gamma : C_3 \sim V \sim E \,\, \cong \,\, \Upsilon : V \sim E \sim P_3. \]
\end{observation}

\noindent{\it Proof:}
There are natural correspondences between $V(\Gamma)$ and $E(\Upsilon)$ and between 
$C_3(\Gamma)$ and $V(\Upsilon)$.  In addition, we may associate each path of $\Upsilon$ with 
vertex sequence $x,y,z$ to the edge of $\Gamma$ which is contained in the triangle associated
with $y$ and has one endpoint in the triangle associated with $x$ and one in the triangle associated
with $z$.  This yields the desired isomorphism.
\quad\quad$\Box$

\bigskip

We are now ready to prove the main result from this section.

\begin{lemma}[Graph Stability]
\label{graph_stability}
Every connected critical trio with a side which is a simple graph is either a video or a pure chord.
\end{lemma}

\noindent{\it Proof:}
Let $(V,E,Z)$ be a connected critical trio and assume that $\Gamma = (V,E)$ is a simple 
graph of degree $d$.  The clone relation must be trivial when restricted to $V \cup E$, and by identifying clones in $Z$ 
we may assume it is also trivial on $Z$.  Let $(A,B)$ be a cross of $(V,E)$ associated with some point in $Z$.  We 
will determine the original structure of $(V,E,Z)$ by classifying $(A,B)$ and then determining its image under 
the action of $G$.  To classify $(A,B)$, note that Lemma \ref{crit2cut} gives $c(A) < 2d$ so either $A$ or $\overline{A} = V \setminus A$ must satisfy one of the six outcomes of Theorem \ref{cutset_classify}.  
We now step through these possibilities, studying the possible orbits of $(A,B)$ for each one.  

The first outcome contradicts the assumption that $(V,E,Z)$ is connected, since in this case either 
$d(Z,V) = 1$ or $d(Z,E) = 0$.  In all the remaining cases Corollary \ref{2dm2_cor} implies that 
$\delta(A) \le w^{\circ}(E)$ so we must have $B = E \setminus N(A)$.  

For the second outcome, it is not possible for $\overline{A}$ to be the 
vertex set of a one edge path, since then $d(Z,E) = 1$ contradicting our connectivity assumption.  If 
$A$ is the vertex set of a one edge path, then since $G$ acts transitively on $E$ we have 
$(Z,V,E) \cong \Gamma : E \sim V \sim E$.  It follows from the connectivity of $(V,E,Z)$ that $|V| \ge 5$, 
so $(Z,E,V)$ is a standard video.

In case we have the third outcome, $\Gamma$ is either a cycle or a two-way 
infinite path and either $A$ is the vertex set of a path or $B$ is the edge set of a path (both if $\Gamma$ is finite).  
It then follows from Lemma \ref{path_trans} that either $(Z,V,E) \cong \Gamma : P_k \sim V \sim E$ or 
$(V,E,Z) \cong \Gamma : V \sim E \sim P_k$ for some $k$, and in either case we have that $(V,E,Z)$ is a pure chord.

For the fourth outcome, $\Gamma$ is cubic and either $A$ is the vertex set of a 2 edge path or $B$ is the edge set 
of a 2 edge path.  It then follows from Lemma \ref{p3_trans} that either $(V,E,Z) \cong \Gamma : V \sim E \sim P_3$ or 
$(Z,V,E) \cong \Gamma : P_3 \sim V \sim E$.  Again by the connectivity of $(V,E,Z)$ we must have $|V| \ge 6$, so 
$(V,E,Z)$ is a standard video.  

If we get the fifth outcome, $\Gamma$ is in ${\mathcal LC}$ and either $A$ is the vertex set of a triangle, or $B$ is 
the edge-set of a triangle.  The latter possibility cannot occur, since this would imply $d(E,Z) = 1$ contradicting our 
connectivity assumption.  To handle the former case, let us choose a cubic graph $\Upsilon$ so that 
$\Gamma \cong L(\Upsilon)$.  Now $G$ acts transitively on the vertices and edges of $\Upsilon$ so $\Upsilon$ is also a $G$-duet.  Furthermore $G$ acts transitively on the triangles of $\Gamma$ so Observation  \ref{cubiclg_equiv} implies
\[ (Z,V,E) \; \cong  \; \Gamma : C_3 \sim V \sim E  \; \cong  \; \Upsilon : V \sim E \sim P_3. \]

Finally, let us consider the sixth outcome.  In this case, Proposition \ref{face_orbit} implies that $(V,E,Z)$ is an exceptional 
video except in the special case when $\Gamma \cong \mathrm{Octahedron}$ and the action of $G$ on the triangles is
not transitive.  However, in this case we have $d(E,Z) = 1$ which contradicts our connectivity assumption.
\quad\quad$\Box$

\subsection{Light Sides of Videos}

In this subsection we will prove a simple lemma which bounds the deficiency of a video, and we 
will prove our light side lemma for videos.

\begin{lemma}[Video criticality]
\label{video_crit}
Let $\Theta = (X,Y,Z)$ be a video.  Then
\[ \delta(\Theta) \le \min \{ w^{\circ}(X), w^{\circ}(Y), w^{\circ}(Z) \}. \]
\end{lemma}

\noindent{\it Proof:}
By possibly passing to a quotient, we may assume that $\Theta$ is clone free.  
We will use Lemma \ref{crit2cut} to relate $\delta(\Theta)$ to $w^{\circ}(E)$.  
For a trio $\Theta = \Gamma : E \sim V \sim E$ based on a graph of degree $d \ge 3$  
we have $\delta(\Theta) = w^{\circ}(E) \le w^{\circ}(V)$.  
For a trio $\Theta = \Gamma : E \sim V \sim P_3$ or $\Gamma : V \sim E \sim P_3$ based on a cubic graph 
$\Gamma$ we have 
$\delta(\Theta) = \tfrac{1}{2} w^{\circ}(E) = w^{\circ}(P_3) \le w^{\circ}(E) \le  w^{\circ}(V)$.  The exceptional 
trios are straightforward to check.
\quad\quad$\Box$

\bigskip

\begin{figure}[ht]
\centerline{\includegraphics[height=10cm]{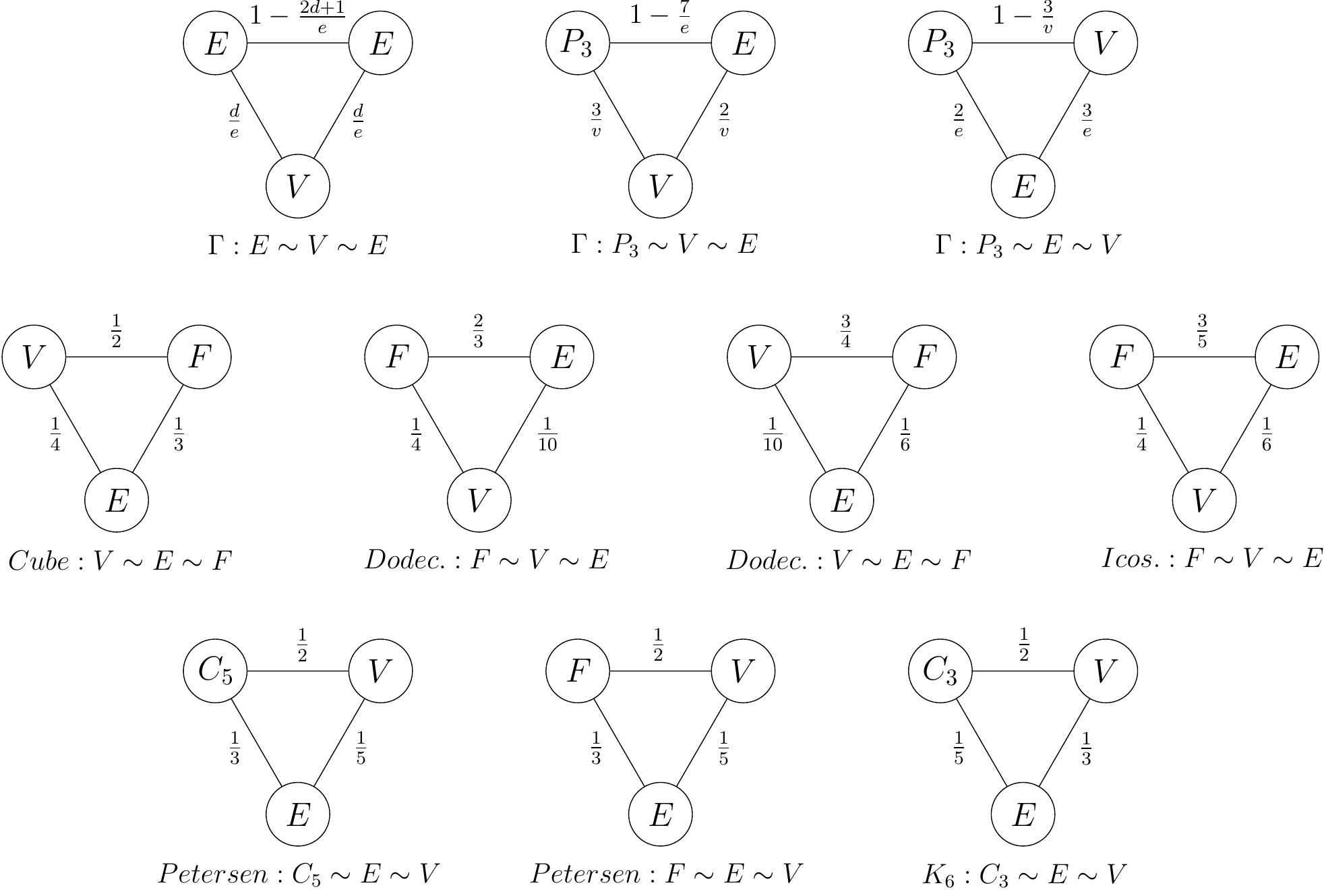}}
\caption{Densities of video sides}
\label{alldeg}
\end {figure}

For the purposes of determining the light sides of videos, it will be convenient to record the densities of the sides of 
(finite) videos.  This information appears in Figure \ref{alldeg} where the variables $d$, $e$, and $v$ denote the 
degree of $\Gamma$, $|E(\Gamma)|$, and $|V(\Gamma)|$.  

\begin{lemma}
\label{largest_video_duets}
Let $\Theta$ be a clone free video and let $(X,Y)$ be a light side of $\Theta$.  Then either 
$(X,Y)$ or $(Y,X)$ is a true graph or a clip.
\end{lemma}

\noindent{\it Proof:}
We shall repeatedly call upon the information in Figure \ref{alldeg}, and shall use $d,e,v$ as in this figure.  The proof of the lemma is 
a straightforward check for the exceptional videos.  Consider a standard video given by $(E_1, V, E_2) = \Gamma : E \sim V \sim E$.  Obviously, 
the duets $(V,E_1)$ and $(V,E_2)$ are graphs.  If the 
duet $(E_1,E_2)$ is a light side, then $1 - \frac{2d+1}{e} \le \frac{d}{e}$ and since $e = \frac{dv}{2}$ this implies $2 \le d(6-v)$ so we must 
have $v \le 5$.  However, then $\Gamma \cong K_5$ and $(E_1,E_2)$ is isomorphic to the clip $K_5(E,C_3)$.  Next consider a video given by 
$(P_3, V, E) = \Gamma : P_3 \sim V \sim E$ (where $\Gamma$ is cubic).  The side $(V,E)$ is a graph and $(P_3,V)$ is a video.  If the side 
$(P_3,E)$ is light, then $1 - \frac{7}{e} \le \frac{3}{v}$.  Now using $e = \frac{3}{2}v$ we find $v \le \frac{23}{3} < 8$.  However, then 
$\Gamma \cong K_{3,3}$, and (as seen in Figure \ref{two_videos}) we have that $(E,P_3)$ is a graph.  Finally, we consider a standard 
video of the form $(P_3,E,V) = \Gamma : P_3 \sim E \sim V$.  Both of the sides $(E,P_3)$ and $(V,E)$ are graphs.  If $(P_3,V)$ is a light side
then $1 - \frac{3}{v} \le \frac{3}{e}$.  However, then using $e = \frac{3}{2}v$ we have $v \le 5$ which contradicts the definition.
\quad\quad$\Box$

\begin{lemma}[Light Sides Video]
\label{light_sides_video}
Let $\Theta \le \Theta^*$ be trios and assume that $\Theta$ is critical and that $\Theta^*$ is a video.  
If $\Delta$ is a light side of $\Theta$, then either $\Delta$ or its reverse is either a near true graph or a near clip.  
\end{lemma}

\noindent{\it Proof:} Assume that $\Theta = (X,Y,Z; \sim)$ and $\Theta^* = (X,Y,Z; \sim^*)$ and that $\Delta = (X,Y)$.  
Since $(X,Y; \sim)$ is a light side of $\Theta$ we may assume (without loss) that $w(X,Y; \sim) \le w(Y,Z; \sim)$.  Since 
$\Theta$ is critical we must have
\begin{equation}
\label{tightside}
w(X,Y; \sim^*) - w(X,Y; \sim) < \delta(\Theta^*).
\end{equation}
Suppose (for a contradiction) that $w( X,Y ; \sim^* ) > w(Y,Z; \sim^*)$.  Let $\mcp, \mcq$ denote the clone partitions of $X,Y$ in the trio $\Theta^*$ 
and note that $w(X,Y; \sim^*)$ and $w(Y,Z; \sim^*)$ are both multiples of $w^{\circ}(\mcq)$.  Thus Lemma \ref{video_crit} implies
\[ w( X,Y; \sim) > w( X,Y; \sim^*) - \delta(\Theta^*) \ge w( X,Y; \sim^*) - w^{\circ}(\mcq) \ge w( Y,Z; \sim^*) \ge w( Y, Z ; \sim) \]
which is contradictory.  So, Lemma \ref{largest_video_duets} implies that either $(X,Y; \sim^*)$ or $(Y,X; \sim^*)$ is either a true graph or a clip.  
Equation \ref{tightside} then gives us
\[ w( X,Y; \sim) > w( X, Y; \sim^*) - \delta(\Theta^*) \ge w( X,Y ; \sim^*) - \min \{ w^{\circ}(\mcp), w^{\circ}(\mcq) \} \]
so $\Delta = (X,Y; \sim)$ is near $(X,Y; \sim^*)^{\bullet}$ relative to $(\mcp, \mcq)$ as desired.
\quad\quad$\Box$

\section{Octahedral Choruses}
\label{oct_sec}

The goal of this section will be to classify maximal critical octahedral choruses.

\subsection{Octahedral Classification}

Let $\Lambda = (X_1,\ldots,X_n)$ be an octahedral chorus with distinguished partition 
$\{ \Lambda_1, \Lambda_2, \Lambda_3 \}$.  In Section \ref{structure} we associated $\Lambda$ with 
a graph $\Gamma$ which has vertex set  $\{X_1, \ldots, X_n \}$ and an edge between $X_i$ and $X_j$ if 
$X_i$ and $X_j$ are not contained in the same member of the distinguished partition.  Furthermore, 
we split the edges of $\Gamma$ into three types: empty, partial, and full.  Here we shall refine this
division by introducing a weight function on the edges of $\Gamma$ given by the rule that the 
weight of an edge $e$ between $X_i$ and $X_j$ is given by $q(e) = q(X_i,X_j)$ (recall that $q(X_i,X_j)$ is the 
density of the duet $(X_i,X_j)$).  So, an edge $e$ is empty if $q(e) = 0$, full if $q(e) = 1$ and 
partial if $0 < q(e) < 1$.  We will also use the function $m$ given by $m(X_i) = 1$ if the member of the distinguished partition
containing $X_i$ has size two, and $m(X_i) = 2$ otherwise.  We say that the graph $\Gamma$ 
equipped with the functions $q$ and $m$ is \emph{associated} with $\Lambda$.  The following
easy observation relates the deficiency of $\Lambda$ to this graph.  Here we use $xy$ to denote
the edge between vertices $x$ and $y$.

\begin{observation}
An octahedral chorus $\Lambda$ is critical if and only if its associated graph $\Gamma$ satisfies
\[ \sum_{xy \in E(\Gamma)}  m(x) m(y) q(xy) > 6 \]
\end{observation}

\noindent{\it Proof:} This follows immediately from the definition of critical and the observation that for 
a finite group $H$ an $H$-duet $(X,Y)$ satisfies 
\[ w_H(X,Y) = d(X,Y) w^{\circ}(Y) = d(X,Y) \frac{ |H| }{ |Y| } = |H| q(X,Y). \quad\quad\Box \]

Our next proposition records some simple properties we will use frequently.  

\begin{lemma}
\label{oct_lem}
Let $\Lambda$ be an octahedral chorus and let $\Gamma$ 
be the associated graph.  Let $e,e',e'' \in E(\Gamma)$ be the edge set of a triangle.  
\begin{enumerate}
\item		If $e$ is a full edge, one of $e',e''$ is empty.
\item		If $q(e) > 0$ then $q(e') + q(e'') \le 1$.
\item		If $e,e',e''$ are nonempty then $q(e) + q(e') + q(e'') \le \frac{3}{2}$.
\item		$q(e) + q(e') + q(e'') \le 2$.
\end{enumerate}
\end{lemma}

\noindent{\it Proof:} Part 1 is immediate, parts 2 and 3 follow from Lemma \ref{nontriv_bound}, 
and part 4 follows from part 3.
\quad\quad$\Box$

\bigskip

With this, we are ready to prove our classification theorem for octahedral choruses.  

\bigskip

\noindent{\it Proof of Theorem \ref{oct_type}:}
Let $\Lambda$ be a maximal critical octahedral chorus and let $\Gamma$ be the associated weighted graph.  Note that the maximality of $\Lambda$ implies that every partial edge of $\Gamma$ must be contained in a triangle of partial edges.  Let $\zeta = \sum_{uv \in E(\Gamma)} q(uv)m(u)m(v)$, so $\zeta > 6$ by assumption.  We split into cases depending on $|V(\Gamma)|$.

\bigskip

\noindent{\it Case 1:} $|V(\Gamma)| \le 4$.

\smallskip

If $|V(\Gamma)| = 3$ then $\Gamma$ must have an empty edge as otherwise part 3 of Lemma \ref{oct_lem} 
gives $\zeta = 4 \sum_{e \in E(\Gamma)} q(e) \le 6$ which is contradictory.  However, then by maximality
$\Gamma$ must have two full edges and one empty edge (a graph which appears in Figure \ref{deg_oct}).  
Next suppose that $|V(\Gamma)| = 4$ and let $y,z \in V(\Gamma)$ be 
the two vertices with $m(y) = m(z) = 2$ and let $x,x'$ be the two vertices with $m(x) = m(x') = 1$.  If every edge is partial then 
using part 3 of  Lemma \ref{oct_lem} gives us
\[ \zeta = 2\big(q(yx) + q(zx) + q(yz)\big) + 2\big(q(yx') + q(zx') + q(yz)\big) \le 2( \tfrac{3}{2} + \tfrac{3}{2}) = 6\]
which is a contradiction.  If $yx$, $zx$ and $yz$ are partial then one of $yx'$ and $yz'$ is empty so using part 2 of Lemma 
\ref{oct_lem} we find \[ \zeta \le 4q(yz) + 2 q(xy) + 2q(xz) + 2 = 2(q(yz) + q(xy)) + 2(q(yz) + q(xz)) + 2 \le 6 \]
which is another contradiction.  Thus $\Gamma$ must have no partial edges and an easy check reveals that it must 
be one of the graphs in Figure \ref{deg_oct}.    

\bigskip

\noindent{\it Case 2:} $|V(\Gamma)| = 5$.

\smallskip

Let $x$ be the vertex with $m(x) = 2$, let $y,z,y', z'$ be the other vertices and assume that $y$ and $y'$ are 
nonadjacent.  First, observe that if all edges are partial, then summing the bound from part 3 of Lemma \ref{oct_lem} 
over all triangles gives us $\zeta \le 6$ which is a contradiction.  It follows that there must be at least one empty edge.  
Next suppose that $yz$ is empty and that all other edges except possibly $y'z'$ are nonempty.  Then using 
bounds 2 and 4 from Lemma \ref{oct_lem} we find
\begin{align*}
 \zeta = &\big( q(xy) + q(yz') \big) + \big( q(xy) + q(x z') \big) + \big( q(xz) + q(x y') \big) + \big( q(xz) + q(z y') \big) \\
 &+ \big( q(x y') + q(x z') + q( y' z' ) \big) \le  6
\end{align*}
and this is a contradiction.  It follows from this that either $\Gamma$ has an empty edge incident with $x$ or $\Gamma$ has 
two empty edges incident with some other vertex.  In either case, we may assume without loss that there is a vertex not incident with any partial edge.  We may assume (without loss) that $z'$ is such a vertex.  If all edges not incident with 
$z'$ are partial then using bound 2 from the above proposition gives us
\[ \big( q(xy) + q(yz) \big) + \big( q(xy) + q(xz) \big) + \big( q(xy') + q(y'z) \big) + \big( q(x y') + q(xz) \big) \le 4 \]
but then again $\zeta \le 6$ which gives us a contradiction.  It follows that there is at most one triangle of partial edges in $\Gamma$.  
If there is no triangle of partial edges then $\Gamma$ must be one of the graphs in Figure \ref{deg_oct} of type $-1$.  
Otherwise, applying the bound in part 3 of the above proposition to the unique triangle of partial edges reveals that the total sum of $q(uv) m(u) m(v)$ over all non partial edges must be greater than $3$, and then $\Gamma$ must be one of the graphs from Figure \ref{deg_oct} of type $2$.  

\bigskip

\noindent{\it Case 3:} $|V(\Gamma)| = 6$.

\smallskip

Call two triangles $T,T'$ of $\Gamma$ {\it opposite} if every they are vertex disjoint.  We begin by 
proving that for every pair of opposite triangles $T, T'$, there is an empty edge in $E(T) \cup E(T')$.  
To show this, let us suppose (for a contradiction) that it fails, and let $E(T) = \{ e_1, e_2, e_3\}$.  For 
$i=1,2,3$ let $T_i$ be the triangle distinct from $T$ which contains $e_i$ and let the edges of 
$T_i$ be $e_i, e_i', e_i''$.  Now using part 2 of Lemma \ref{oct_lem} gives $w(e_i') + w(e_i'') \le 1$ and part 3 of 
this proposition gives $w( E(T) ), w( E(T') ) \le \frac{3}{2}$ and this contradicts $\zeta > 6$.

Our next goal will be to prove that there is a vertex $x$ of $\Gamma$ which is incident with at least two 
empty edges but no partial edge.  To do so, let us now choose $S$ to be a minimal set of empty edges so 
that every pair of opposite triangles contains at least one.  If $|S| = 2$ then $S$ consists of two edges 
incident with a common vertex $x$ but not contained in a triangle, and in this case the vertex $x$ meets 
our goal.  Otherwise $|S| = 3$ and there exist opposite triangles $T,T'$ so that 
$S \subseteq E(T) \cup E(T')$.  If $S = E(T)$, then we may partition the edges not in $T$ into three triangles, 
and then bound 5 from Lemma \ref{oct_lem} brings us to a contradiction.  In the remaining case, there is a vertex $x$ 
incident with two edges in $S$, and every triangle containing $x$ contains an edge in $S$, so by maximality, 
$x$ is not incident with any partial edge.

Let $V(\Gamma) = \{x,x',y,y',z,z'\}$ (with $x$ as in the previous paragraph) and assume that 
$x,x'$ are nonadjacent and $y,y'$ are nonadjacent.  Suppose (for a contradiction that $x'$ is not incident with an 
empty edge.  Then using bound 2 from Lemma \ref{oct_lem} gives us
\[ \big( q(x'z) + q(zy') \big) + \big( q(x'y) + q(yz) \big) + \big( q(x' z') + q(z'y) \big) + \big( q(x'y') + q(y'z') \big) \le 4 \]
and since $x$ is incident with at least two empty edges, this is contradictory.  So, we may assume 
without loss that $x'y$ is empty.  It now follows that $y$ is not contained in any triangle of partial edges.

Based on our analysis, we now have that $\Gamma$ contains at most two triangles of partial edges, 
and if there are two, they share an edge.  If there are no partial edges, then $\zeta > 6$ implies that 
$\Gamma$ is one of the graphs in Figure \ref{all_oct} of type $-1$ or type $0$.  If there is one triangle 
of partial edges, then bound 3 from Lemma \ref{oct_lem} implies that there must be at least 5 full 
edges, and it follows that $\Gamma$ is one of the graphs in Figure \ref{all_oct} of type $1$.  Finally, 
if there are two triangles of partial edges sharing an edge, then bound 3 in Lemma \ref{oct_lem} 
implies that the sum of the weights of these partial edges is at most 3, so there must be at least 4 full 
edges, and it follows that $\Gamma$ is one of the graphs in Figure \ref{all_oct} of type $2$.
\quad\quad$\Box$

\subsection{Type 2 Classification}

In this subsection we will further analyze type 2 octahedral choruses, showing that those 
of interest to us must have type $2A$, $2B$, or $2C$.  We also record a simple lemma concerning
the deficiency of an octahedral chorus.

\begin{lemma}
\label{two_big_sides}
Let $(X,Y,Z)$ be a nontrivial maximal critical trio and assume $\delta(X,Y,Z) = w(X,Y)$.  Then 
$(X,Y,Z)$ is a pure beat relative to $(X,Y)$.
\end{lemma}

\noindent{\it Proof:} It follows from Theorem \ref{connected_bound} that $(X,Y)$ is disconnected, 
so we may let $( \mcp, \mcq )$ denote its component quotient.  Choose 
$z \in Z$ and let $(A,B)$ be the cross associated with $z$.  Then $\delta(\Theta) = \delta(A,B)$, so 
Lemma \ref{cross_in_discon} implies that $(A,B)$ is pure with respect to $(\mcp,\mcq)$.  It then follows 
from maximality that $(X,Y,Z)$ is a pure beat relative to $(X,Y)$, as desired.  
\quad\quad$\Box$

\bigskip

We are now ready to prove our classification theorem for octahedral choruses of type $2$.  

\bigskip

\noindent{\it Proof of Theorem \ref{classify_oct2}: }
Let $\Lambda$ be a maximal critical  $H$-octahedral chorus of type 2 with incidence relation $\sim$ and choose 
$X,Y,Z_1,Z_2$ so that for $i=1,2$ we have that $(X,Y)$, $(Y,Z_i)$, and $(X,Z_i)$ are sides 
of $\Lambda$ each of which is partial.  We may assume that either $\Lambda$ has rank
6 and $Z_1 \neq Z_2$ or $\Lambda$ has rank 5 and $Z_1 = Z_2$ has the property that $\{Z_i\}$ is a member of 
the distinguished partition.  For $i=1,2$ let $\Theta_i = (X,Y,Z_i; \sim )$
and note that $\Theta_i$ is a nontrivial $H$-trio.   Our definitions now imply the following 
(for both the rank 5 and rank 6 cases)
\begin{align*}
0 < \delta( \Lambda ) &= w(X,Y) + w(X,Z_1) + w(Y,Z_1) + w(X,Z_2) + w(Y,Z_2) - 2 |H| 	\\
		&= \delta( \Theta_1) + \delta( \Theta_2) - w(X,Y) 
\end{align*}
Observe first that the above equation together with the easy bound $\delta(\Theta_i) \le w(X,Y)$ (from Lemma \ref{nontriv_bound}) imply that both $\Theta_1$ and $\Theta_2$ are critical.  If $(X,Y)$ is connected, then Theorem \ref{connected_bound} implies $\delta( \Theta_i ) \le \frac{1}{2}w(X,Y)$ for $i=1,2$ which again contradicts our equation.  Thus, $(X,Y)$ is disconnected, and we will let $( \mcp, \mcq )$ be its component quotient.  Now for $i=1,2$ choose 
$z_i \in Z_i$ and let $(A_i,B_i)$ be the cross of $(X,Y)$ associated with $z_i$ in the trio $\Theta_i$, and note that the maximality of $\Lambda$ implies that 
$(A_i,B_i)$ is maximal for $i=1,2$.  It follows from Lemma \ref{cross_in_discon} that at least one of $(A_1,B_1)$ and $(A_2,B_2)$ is pure with respect to $(\mcp, \mcq)$.  If $\Lambda$ has rank 5, then $\Theta_1 = \Theta_2$  and $(A_1,B_1) = (A_2,B_2)$ is pure.  Now  
by maximality, $\Theta_1$ is a pure beat relative to $(X,Y)$, and $\Lambda$ has type $2C$.  So, we 
may now assume that $\Lambda$ has rank 6.  For $i=1,2$ if $(A_i,B_i)$ is pure, choose a maximal trio 
$\Theta_i = (X,Y,Z, \sim_i)$ with $\Theta_i \le \Theta_i^*$ and note that the maximality of $(A_i,B_i)$ implies that $\sim |_{X \cup Z_i} = \, \sim_i|_{X \cup Z_i}$ and $\sim|_{Y \cup Z_i} = \, \sim_i |_{Y \cup Z_i}$, so Lemma \ref{two_big_sides} implies that $\Theta_i^*$ is a pure chord relative to $(X,Y)$.  
If both $(A_1,B_1)$ and $(A_2,B_2)$ are pure, then the maximality of  $\Lambda$ implies $\sim|_{X \cup Y} = \; \big( \sim_1|_{X \cup Y} \big) \cap \big( \sim_2|_{X \cup Y} \big)$
and we conclude that $\Lambda$ has type $2B$.  Otherwise, we may assume (without loss) that $(A_1,B_1)$ is impure.  In this case Lemma \ref{cross_in_discon} implies that $(X,Y,Z_1; \sim)$ is an impure beat relative to $(X,Y)$.  Further, $\Theta_2 = (X,Y,Z; \sim_2)$ is a pure beat relative to 
$(X,Y)$ and by construction $\sim|_{X \cup Y} \subseteq \sim_2|_{X \cup Y}$ so $\Lambda$ is type $2A$.
\quad\quad$\Box$

\begin{lemma} 
\label{oct_ex_bound}
Let $\Lambda = (X_1,\ldots,X_n)$ be a maximal critical octahedral chorus of type $0$, $1$, or $2$ and let $(X_i,X_j)$ be a nonempty side of $\Lambda$.  Then $\delta(\Lambda) \le w(X_i,X_j)$.
\end{lemma}

\noindent{\it Proof}
If $(X_i,X_j)$ is full, then the result is immediate.  Otherwise there exists $1 \le k \le n$ so that $(X_i,X_j,X_k)$ is a nontrivial
trio for which $\delta(\Lambda) \le \delta(X_i,X_j,X_k)$.  Now the result follows from Lemma \ref{nontriv_bound}.
\quad\quad$\Box$

\section{Sequences and Near Sequences}
\label{sequence_sec}

The purpose of this section is to prove a stability lemma for near sequences and a light side 
lemma for chords.

\subsection{Sequence Stability}

In this subsection we establish some basic terminology, and then prove a stability lemma for sequences. 

Let $\Delta = (X,Y)$ be a sequence.  We define the \emph{length} of $\Delta$ to be $\ell = d(X,Y)$.  We say that two points $x,x' \in X$ (or $x,x' \in Y$) 
are \emph{consecutive} if 
$|N(x) \cap N(x')| = \ell - 1$.  It follows from our definitions (which require $2 \le \ell \le |X| - 2$) that every point in $X$ is consecutive with exactly two other points, so the graph on $X$ obtained by adding an edge between each pair of consecutive points is either a two-way infinite path or a (finite) cycle.  We say that a set $I \subseteq X$ (or $I \subseteq Y$) is an \emph{interval} if we can write $I = \{ x_1, x_2, \ldots, x_k \}$ where all these points are distinct, and $x_i$ and $x_{i+1}$ are consecutive for every $1 \le i \le k-1$.  In this case we call $x_1$ and $x_k$ the \emph{endpoints}, and we say the sequence is 
\emph{nontrivial} whenever $x_1 \neq x_k$.  Our first lemma shows that critical crosses of sequences come from intervals.

\begin{lemma}
\label{cross_seq}
Let $\Delta = (X,Y)$ be a sequence.  If $(A,B)$ is a nontrivial critical cross of $\Delta$ with $A$ finite, then $(A,B)$ is maximal and $A$, $Y \setminus B$ are intervals.
\end{lemma}

\noindent{\it Proof:}
Let $(A,B)$ be a nontrivial critical cross, and choose a maximal cross $(A^*,B^*)$ with $A^* \supseteq A$ and $B^* \supseteq B$.
Since $N(A^*)$ is finite, we may express it as a disjoint union of intervals $I_1 \cup I_2 \cup \ldots \cup I_m$.  
For each $x \in A^*$ we have that $N(x)$ is an interval, so there must exist $1 \le i \le m$ with $N(x) \subseteq I_i$.  
If $m > 1$ we may  partition $A^*$ into $\{A^*_1,A^*_2\}$ so that 
$N(A^*_1) \cap N(A^*_2) = \emptyset$.  Then Theorem \ref{connected_bound} implies
\begin{align*}
w(N(A^*)) 
	&= w(N(A^*_1)) + w(N(A^*_2)) \\
	&\ge \big( w(A^*_1) + \tfrac{1}{2}w(\Delta) \big) + \big( w(A^*_2) + \tfrac{1}{2}w(\Delta) \big)	\\
	&= w(A^*) + w(\Delta)
\end{align*}
which contradicts the criticality of $A^*$.  Thus $N(A^*) = Y \setminus B^*$ is a single interval, and then by maximality, 
$A^*$ must also be an interval.  Furthermore, if we set $k = w^{\circ}(X) = w^{\circ}(Y)$ then 
$\delta(A^*,B^*) = k  d(X,Y) + k |A^*| - k |N(A^*)| = k$.  It follows from this last equation that 
$(A,B) = (A^*,B^*)$ which completes the proof.
\quad\quad$\Box$

\bigskip

Let $G$ be a group acting on a set $X$.  We previously defined this action to be dihedral if $G$ is dihedral and the rotation subgroup of $G$ acts regularly on $X$.  We now define this action to be \emph{cyclic} if $G$ is cyclic and $G$ acts regularly on $X$, and we define it to be \emph{split} if $G$ is dihedral and $G$ acts regularly on $X$.  

Let $\Delta = (X,Y)$ be a $G$-duet which is a sequence.  The automorphism group of $\Delta$ must preserve the relation 
consecutive, and it follows easily that ${\mathit Aut}(\Delta)$ is the dihedral group ${\mathbf D}_n$ where $n = |X| \ge 4$ 
(here we permit $n = \infty$ if $X$ is infinite).  Since $\Delta$ is a $G$-duet, this gives us a homomorphism 
\[ \phi : G \rightarrow {\mathit Aut}(\Delta) \cong {\mathbf D}_n.\]
Note that the kernel of $\phi$ is given by $H = G_{ (X) } = G_{ (Y) }$.  It follows easily that the action of $G/H$ on 
$X$ and $Y$ must be either dihedral, cyclic, or split, and we say that $\Delta$ has \emph{dihedral}, \emph{cyclic}, or 
\emph{split action} accordingly.  We record this observation below.

\begin{observation}
If $(X,Y)$ is a $G$-duet which is a sequence, it has either cyclic, dihedral, or split action.
\end{observation}

We define a coset $S \in G/H$ to be a \emph{consecutive shift} if for some (and thus every) $g \in S$ we have that $gx$ and $x$ are consecutive for 
every $x \in X$.  Note that sequences with dihedral or cyclic action have consecutive shifts, but sequences with split action do not.  
Indeed, for a sequence $(X,Y)$ with split action, we may colour the points in $X$ with $\{ black, white \}$ so that consecutive points are distinct 
colour, and then $\phi(G)$ contains all rotations which map each colour class to itself, and all flips which interchange the colour classes.  
Sequences with split action must also satisfy an additional parity constraint as indicated by the following observation.

\begin{observation}
\label{alt_seq_obs}
Let $\Delta = (X,Y)$ be a  $G$-duet which is a sequence with split action.
\begin{enumerate}
\item The collection of all intervals of $Y$ of length $m$ form a single orbit under the action of $G$ if $m$ is odd, and form two orbits if $m$ is even.
\item The length of $\Delta$ is odd.
\end{enumerate}
\end{observation}

\noindent{\it Proof Sketch:} To see part 1, first suppose that $m$ is odd, and note that for any two intervals, $I,I' \subseteq Y$ with $|I| = |I'| = m$ there must exist $g \in G$ which maps the ``middle point'' of $I$ to the ``middle point'' of $I'$ and then $gI = I$.  On the other hand, if $m$ is even and $I$ is an interval with $|I| = m$ then there exists $g \in G$ which acts nontrivially on $X$ but has $gI = I$.  For part 2, note that in order to be a $G$-duet, $G$ must act transitively on $\{ N(x) \mid x \in X \}$, so $d(X,Y) = d(Y,X)$ must be odd.
\quad\quad$\Box$

\bigskip

Let $\Theta = (X,Y,Z)$ be a $G$-trio which is a pure chord.  We say that $\Theta$ has \emph{dihedral}, \emph{cyclic}, or \emph{split action} if setting $H = G_{ (X) } = G_{ (Y) } = G_{ (Z) }$ we have that the action of $G/H$ on each of $X$, $Y$, and $Z$ has the corresponding type.  It follows easily that every pure chord has either dihedral, cyclic, or split action.   

We are now ready for our stability lemma for sequences.  

\begin{lemma}[Sequence Stability]
\label{sequence_stability}
Let $\Theta$ be a maximal critical connected $G$-trio.  If $\Theta$ has a side $\Delta$ which is a sequence, then $\Theta$ is 
a pure chord.  
\end{lemma}

\noindent{\it Proof:} By identifying clones, we may assume that $\Theta = (X,Y,Z)$ is a clone free $G$-trio
and that $(X,Y)$ is a sequence.  Then choose $z \in Z$, let $(A,B)$ be the cross of $(X,Y)$ associated with $z$, and assume $A$ is finite.  If $\{ gA \mid g \in G \}$ is the set of all intervals of $X$ of size $|A|$, then $(X,Y,Z)$ is a pure chord as desired.  This completes the proof except in the case when $(X,Y)$ has split action and $|A|$ is even.  However, it then follows from Observation \ref{alt_seq_obs} that there exist distinct points $x,x' \in X$ for which $x \in gA$ if and only if $x' \in gA$ for all $g \in G$.  But then $N(x) \cap Z = N(x') \cap Z$ but $N(x) \cap Y \neq N(x') \cap Y$ and this contradicts the maximality of $\Theta$.
\quad\quad$\Box$

\subsection{Adjusting Fringed Sequences}

In this section we will prove a key lemma which allows us to adjust a fringed sequence with a special behaviour.  As a consequence of this we will derive lemmas to deal with dihedral chords of type $-1$ and type $0$.  We begin with some definitions and a significant technical point.  

Consider a $G$-duet $\Delta = (X,Y)$ which is a length $\ell$ fringed sequence relative to $( {\mathcal P}, {\mathcal Q} )$ where $J = {\mathbb Z}$ or ${\mathbb Z}/n {\mathbb Z}$ for some $n \ge \min\{ 4, \ell+1 \}$ and ${\mathcal P} = \{ P_j \mid j \in J \}$ and ${\mathcal Q} = \{ Q_j \mid j \in J \}$ satisfy
\begin{equation}
\label{examp_fringe}
 (P_i,Q_j) \mbox{ is } \left\{\begin{array}{ll}
		\mbox{partial}		&	\mbox{ if $j=i$ or $j=i+\ell-1$}	\\
		\mbox{full}	&	\mbox{ if $i+1 \le j \le i+\ell-2$}	\\
		\mbox{empty}		&	\mbox{ otherwise}
		\end{array}
		\right. 
\end{equation}
Let $H = G_{ ( \mcp ) } = G_{ (\mcq) }$ and recall that by our definition $G/H$ is dihedral, and further, the action of $G/H$ on 
$\mcp$ and $\mcq$ is dihedral.  

We now define $\Delta$ to be \emph{short} if $\ell=2$ and \emph{long} if $\ell = |J| - 1$.  Observe that the duet $(\mcp, \mcq)$ will be a sequence if and only if $\Delta$ is not long.  We will have to be mindful of this technicality since we will be interested in passing between fringed sequences and near sequences.  Also, this slight difference will force us to define a couple of terms for fringed sequence which we have already established for sequences.  

For $i,j \in J$ we say that $P_i$ and $P_j$ ($Q_i$ and $Q_j$) are \emph{consecutive} if $i = j \pm 1$.  It is straightforward to check that every automorphism of $\Delta$ must preserve this relation.  As before, we say that $S \in G/H$ is a \emph{consecutive shift} if some (and thus every) $g \in S$ we have that $gP$ and $P$ ($gQ$ and $Q$) are consecutive for every $P \in \mcp$ ($Q \in \mcq$).  It follows that there are exactly two consecutive shifts, $S$ and $S^{-1}$, one of which maps each block in $\mcp,\mcq$ to the next higher index, and one mapping each such block to the next lower index.

\bigskip

\noindent{\bf Definition:} Let $\Lambda$ be a square chorus.  We say that $\Lambda$ is \emph{small} (\emph{big}) if there exists $P \in \Lambda$ so that $(P,Q)$ is empty (full) for every $Q \in \Lambda$ which is not in the same member of the distinguished partition as $P$.

\begin{figure}[ht]
\centerline{\includegraphics[height=4.4cm]{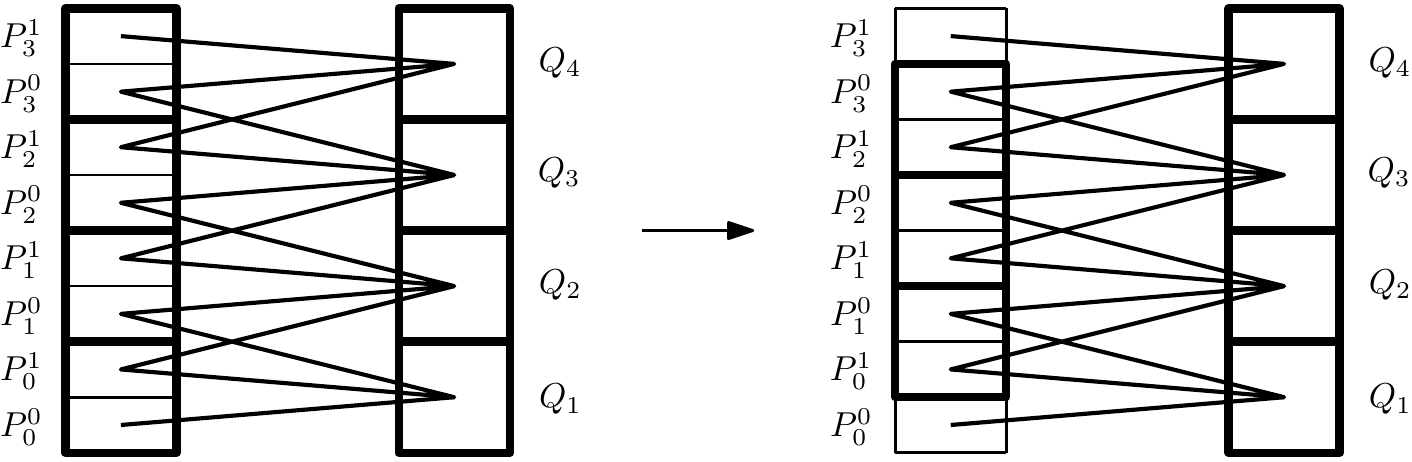}}
\caption{Adjusting a Fringed Sequence}
\label{adj_fringe}
\end {figure}

\begin{lemma}[Adjusting Fringed Sequences]
\label{adjust_fringe}
Let $(X,Y)$ be a fringed sequence of length $\ell$ relative to the systems of imprimitivity $\mcp, \mcq$ of $X,Y$ and assume it has a small associated square chorus, that $H = G_{ ( \mcp )} = G_{ (\mcq) }$, and that  $S \in G/H$ is a consecutive shift.  If $\ell=2$ then $(X,Y)$ is disconnected.  Otherwise, there exist systems of imprimitivity ${\mathcal P}'$ of $X$ and ${\mathcal Q}'$ of $Y$ so that 
\begin{itemize}
\item $G_{ ( \mcp' ) } = G_{ (\mcq') } = H$.
\item $w^{\circ}( {\mathcal P}' ) = w^{\circ}( {\mathcal Q}') = w^{\circ}( {\mathcal P} ) = w^{\circ}( {\mathcal Q} )$.  
\item $( {\mathcal P}', {\mathcal Q}' )$ is a sequence of length $\ell-1$ with dihedral action and consecutive shift $S$.
\item  $(X,Y)$ is a near sequence relative to $( {\mathcal P}', {\mathcal Q}' )$ 
\end{itemize}
\end{lemma}

\noindent{\it Proof:} Assume that $(X,Y)$ and ${\mathcal P}, {\mathcal Q}$, $H$, and $\ell$ are as given in equation \ref{examp_fringe}, let 
$\Lambda$ be an associated square chorus, and note that $\Lambda$ is nonempty by assumption.  Consider 
$(P_1,Q_1)$ and observe that by the assumption $\Lambda$ is small, there must either exist $x \in P_1$ so that $N(x) \cap Q_1 = \emptyset$ or there must exist $y \in Q_1$ so that $N(y) \cap P_1 = \emptyset$ and we shall assume the former case without loss.  It now follows from our assumptions that every $P_j$ is the disjoint union of two 
$H$-orbits which we shall denote $P_j^0$ and $P_j^1$.  We 
may assume further that $(P_j^1,Q_j)$ and $(P_j^0, Q_{j+\ell-1})$ are empty for every $j \in J$.  Now, 
${\mathcal P}' = \{ \ldots, P_0^1 \cup P_1^0, P_1^1 \cup P_2^0, P_2^1 \cup P_3^0, \ldots \}$ is a system of imprimitivity.  It follows from our construction 
that $( P_1^1 \cup P_2^0, Q_j)$ is nonempty only for $2 \le j \le \ell$.  If $\ell=2$ then $N( P_1^1 \cup P_2^0 )$ is contained in 
$Q_2$ and $(X,Y)$ is disconnected.  Otherwise, $( \mcp', \mcq)$ is a sequence of length $\ell-1$ and 
\[ w(X,Y) = 2(\ell-2)|H| + \hat{w}_H(\Lambda) > (\ell-2) w^{\circ}( \mcp ) = w( \mcp', \mcq) - w^{\circ}( \mcp' ) \]
so $(X,Y)$ is a near sequence relative to $( \mcp', \mcq)$.  The remaining properties are straightforward to verify.
\quad\quad$\Box$

\bigskip

Next we establish a result which will be used to handle dihedral chords of type $-1$.

\begin{lemma}
\label{kill-1}
There does not exist a maximal dihedral chord of type $-1$.  
\end{lemma}

\begin{figure}[ht]
\centerline{\includegraphics[height=2cm]{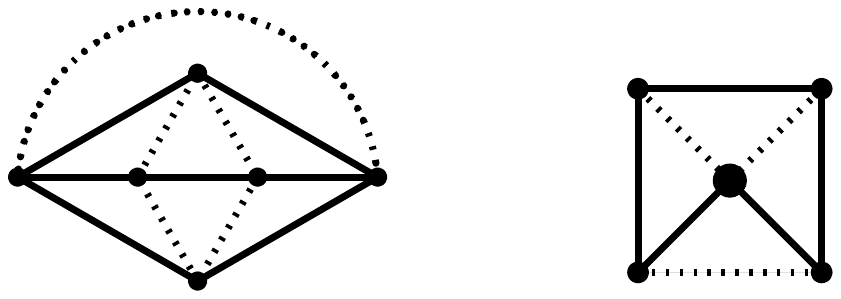}}
\caption{Two type $-1$ octahedral choruses}
\label{exc_oct}
\end {figure}

\noindent{\it Proof:} 
Suppose (for a contradiction) that $\Theta$ is a maximal dihedral chord of type $-1$, let $\Lambda$ be an associated octahedral chorus and note that our definitions imply that $\Theta$ is critical.  It follows from the definition of dihedral chord that none of the square choruses contained in $\Lambda$ can be full or empty.  This eliminates all possibilities for $\Lambda$ except those in Figure \ref{exc_oct}. 

Both of the octahedral choruses in this figure have the property that they contain two square choruses which are 
simultaneously big and small.  Therefore, we may choose a side $\Delta$ of  $\Theta$ so that $\Delta$ is a fringed sequence of length $\ell$ relative to $(\mcp,\mcq)$ and its associated $H$ square chorus $\Lambda'$ is both big and small (recall that by definition $w^{\circ}(\mcp) = w^{\circ}(\mcq) = 2|H|$).  Next apply Lemma \ref{adjust_fringe} to choose a quotient $(\mcp', \mcq')$ of $\Delta$.  Now we have 
\[ w(\Delta) = 2(\ell-2)|H| + \hat{w}_H(\Lambda') = 2(\ell-1) |H| = w( \mcp', \mcq' ) \]
so $\mcp'$ and $\mcq'$ are systems of clones in $\Delta$, and $\Delta^{\bullet}$ is a sequence.  
However, $w^{\circ}(\mcp') = w^{\circ}(\mcq') = 2|H|$, but in both of the octahedral choruses in Figure \ref{exc_oct}, one of the three associated square choruses is a rank 4 chorus which has (up to reversing) three full sides and one empty side.  It follows from this that there is a side $\Delta'$ of $\Theta$ for which either $w(\Delta')$ or $\overline{w}(\Delta')$ is an odd multiple of $|H|$, and this contradicts the maximality of $\Theta$.
\quad\quad$\Box$

\bigskip

Our next lemma uses similar logic to handle dihedral chords of type $0$.

\begin{lemma}
\label{type0alt}
A trio is a dihedral chord of type $0$ if and only if it is a a pure chord with split action.
\end{lemma}

\begin{figure}[ht]
\centerline{\includegraphics[height=4cm]{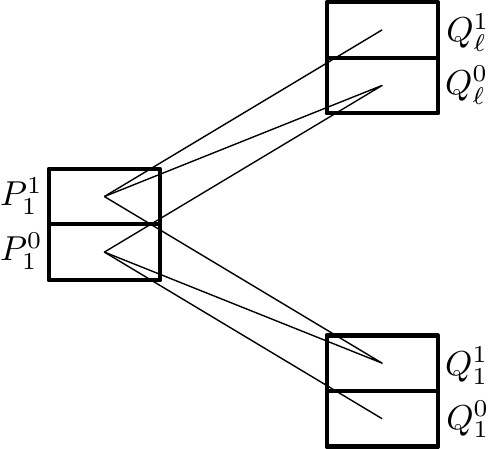}}
\caption{Adjusting a type $0$ dihedral chord}
\label{zero2pure_fig}
\end {figure}

\noindent{\it Proof:} 
First, we let $\Theta = (X,Y,Z)$ be a dihedral chord of type $0$ and show that $\Theta$ is a pure chord with split action.  
A quick check of Figure \ref{all_oct} reveals that in the graph associated with an octahedral chorus of type $0$, among 
the three squares which are associated with a square chorus, two have three full edges and one empty edge, while the other 
has one full edge and three empty edges.  It follows that we may choose a side $(X,Y)$ of $\Theta$ which is a fringed sequence 
(i.e. $w(X,Y) < \infty$) relative to $( {\mathcal P}, {\mathcal Q} )$ and has associated square chorus with 3 full and one empty
edge.  We shall assume that $\Delta = (X,Y)$ and ${\mathcal P}$, ${\mathcal Q}$, $H$, and $\ell$ are as given in equation \ref{examp_fringe}.  
By assumption, each $P_j$ ($Q_j$) is the union of two $H$-orbits which we denote $P_j^0$ and $P_j^1$ ($Q_j^0$ and $Q_j^1$).  We may further assume that for every $j \in J$ we have $(P_j^i,Q_j^k)$ is empty if and only if $i= 1$ and $k = 0$ and $(P_j^i,Q_{j+\ell-1}^k)$ is empty if and only if $i=0$ and $k=1$ as shown in Figure \ref{zero2pure_fig}.
Define ${\mathcal P}' = \{    \ldots, P_0^0, P_0^1, P_1^0, P_1^1, \ldots \}$ and ${\mathcal Q}' = \{ \ldots, Q_0^0, Q_0^1, Q_1^0, Q_1^1, \ldots \}$, and observe that ${\mathcal P}'$ and ${\mathcal Q}'$ are systems of clones, and further $( {\mathcal P}', {\mathcal Q}' )$ is a sequence with split action.  Let
$(A,B)$ be a cross of $(X,Y)$ associated with some point in $Z$, and assume (without loss) that $A$ is finite.  Now $(A,B)$ is pure with respect to 
$(\mcp', \mcq')$ so by Lemma \ref{cross_seq} there must be an interval $I \subseteq \mcp'$ of $(\mcp', \mcq')$ so that $A$ is the union of the blocks in $I$.  Furthermore, it follows from the assumption that $\Theta$ has type $0$ that $|I|$ is odd.  It then follows from Observation \ref{alt_seq_obs} that 
$\Theta$ is a pure chord with split action, as desired.

Next suppose that $\Theta = (X,Y,Z)$ is a pure chord with split action and with clone partitions $\mcp, \mcq, \mcr$ of $X,Y,Z$.  Assume that 
$w(X,Y), w(Y,Z) < \infty$ and suppose that $\mcq = \{ Q_j \mid j \in J \}$ for $J = {\mathbb Z}$ or ${\mathbb Z}/n{\mathbb Z}$ and that $Q_j$ and $Q_{j+1}$ are consecutive for every $j \in J$.  Note that $d(\mcp,\mcq)$ and $d(\mcq,\mcr)$ are odd, and $d(\mcp,\mcr)$ is either odd or $\infty$ so we must have $|J| \ge 8$.  Now define $\mcq' = \{ Q_1 \cup Q_2, Q_3 \cup Q_4, \ldots \}$, and note that $\mcq'$ is a system of imprimitivity.  Furthermore, in the duet $(\mcp, \mcq)$ every $P \in \mcp$ is joined to an interval of odd length from $\mcq$.  It follows that for every $P \in \mcp$ there is exactly one $P' \in \mcp \setminus \{P\}$ so that $P$ and $P'$ have the same neighbours in $(\mcp,\mcq')$.  Thus 
$P \cup P'$ is a block of imprimitivity and we let $\mcp'$ denote the corresponding system of imprimititivity.  Define 
$\mcr'$ on $Z$ using blocks from $\mcr$ analogously.  Note that $\Theta$ is maximal since it is a pure chord.  Therefore, $\Theta$ is a maximal critical dihedral chord relative to $\mcp', \mcq', \mcr'$ and the square choruses associated to the fringed sequences on $(X,Y)$ and $(Y,Z)$ both have the property that their associated square choruses have rank 4, no partial sides, and (up to reversing) just one empty side.  It follows from this and a check of Figure \ref{all_oct} that $\Theta$ is a dihedral chord of type $0$ relative to $\mcp', \mcq', \mcr'$, as desired.
\quad\quad$\Box$

\subsection{Light Sides}

The goal of this subsection is to prove our light side lemma for chords.  In fact, we will prove a somewhat stronger 
result which applies not only to light sides, but more generally to sides of finite weight.  
We begin with an observation which follows from a check of Figures \ref{all_oct} and \ref{deg_oct}.  

\begin{observation}
\label{big-or-small}
For every octahedral chorus of type $0$, $1$, or $2$, two of the square choruses it contains are big
and the other is small.
\end{observation}

Observation \ref{big-or-small} is especially convenient since it forces the square choruses we use to be one of 
two extreme types.  The next observation follows quite quickly from this.

\begin{observation}
\label{length-width}
Let $\Theta$ be a dihedral chord and let $\Delta$ be a side of $\Theta$ which is a fringed sequence with associated 
square chorus $\Lambda$.
\begin{enumerate}
\item If $\Delta$ is short, then $\Lambda$ is big.
\item If $\Delta$ is long, then $\Lambda$ is small.
\end{enumerate}
\end{observation}

\noindent{\it Proof:} For part 1, observe that if $\Delta$ is short and $\Lambda$ is small, then Lemma \ref{adjust_fringe} would imply that $\Delta$ is not connected, which contradicts the definition of dihedral chord.  For part 2, note that if $\Delta$ is long, then every side of $\Theta$ other than $\Delta$ or its reverse will be a short fringed sequence, which by part 1 must have an associated square chorus which is big.  Observation \ref{big-or-small} then implies that $\Lambda$ must be small.
\quad\quad$\Box$

\bigskip

The previous two observations also reveal why in the definition of dihedral chord, we could insist on the bound $|J| \ge 4$ without loss.  Were we to permit $|J| = 3$ for such a trio, it would have every side short, but then Observation \ref{big-or-small} and Lemma \ref{adjust_fringe} would force the trio to be disconnected.

\begin{lemma}
\label{light_sides_chord}
Let $\Theta = (X,Y,Z ; \sim)$ be a critical $G$-trio and assume that $\Theta \le \Theta^*$ for a song $\Theta^*$.  Let $\Delta$ be a side of $\Theta$ with $w(\Delta) < \infty$.
\begin{enumerate}
\item If $\Theta^*$ is a pure chord, then $\Delta$ is a near sequence.  
\item If $\Theta^*$ is a cyclic chord relative to $\mcp, \mcq, \mcr$ on $X, Y, Z$ then 
$\Delta$ is a near sequence relative to $(\mcp, \mcq)$, and furthermore,  $w(\mcp,\mcq) - w(\Delta) \le w^{\circ}(\mcp) - \delta(\Theta)$.
\item Assume $\Theta^*$ be a dihedral chord of type $1$ or $2$ relative to $\mcp, \mcq, \mcr$ on $X,Y,Z$ with consecutive shift $S$.   Then there exist systems of imprimitivity $\mcp', \mcq'$ so that $\Delta$ is a 
near sequence relative to $(\mcp', \mcq')$ with dihedral action and consecutive shift $S$, and furthermore, $w( \mcp', \mcq' ) - w(\Delta) \le w^{\circ}(\mcp) - \delta(\Theta)$.
\end{enumerate}
\end{lemma}

\noindent{\it Proof:} Let $\Theta^* = (X,Y,Z ; \sim^*)$, assume that $\Delta = (X,Y; \sim)$ and let $\Delta^* = (X,Y; \sim^*)$.
Then we have
\begin{equation}
\label{lsc_def_eqn}
w(\Delta^*) - w(\Delta) \le \delta(\Theta^*) - \delta(\Theta)
\end{equation}
We now split into cases depending on $\Theta^*$.

\bigskip

\noindent{\it Case 1:} $\Theta^*$ is a pure chord.

\smallskip

Let $\mcp$, $\mcq$, $\mcr$ be the clone partitions of $X$, $Y$, $Z$ in the 
trio $\Theta^*$ and note that $( \mcp, \mcq; \sim^* )$ is a sequence.  Set $k = w^{\circ}( \mcp ) = w^{\circ} ( \mcq ) = w^{\circ}( \mcr )$ and note that $k = \delta(\Theta^*)$.  Now by equation \ref{lsc_def_eqn} we have 
$w( \mcp, \mcq ) = w(\Delta^*) < w(\Delta) + k$, so $\Delta$ is a near sequence relative to $(\mcp, \mcq)$.

\bigskip

\noindent{\it Case 2:} $\Theta^*$ is a cyclic chord.

\smallskip

Choose $P \in \mcp$ let $\{ Q_1, \ldots, Q_{\ell} \} = \{ Q \in \mcq \mid Q \sim^* P \}$, and set $k = w^{\circ}(\mcp) = w^{\circ}(\mcq) = w^{\circ}(\mcr)$.  Since $\Theta^*$ is a cyclic chord, we may assume (without loss) that $(P,Q_1; \sim^*)$ is partial, and $(P, Q_i; \sim^*)$ is full for $2 \le i \le {\ell}$.  Now, $(P,Q_1; \sim^*)$ is a side of a continuation of $\Theta^*$ so we have $\delta(\Theta^*) \le w(P,Q_1; \sim^*)$ (by Lemma \ref{nontriv_bound}).  Combining this with Equation \ref{lsc_def_eqn} then yields
\begin{align*}
w( \mcp, \mcq ; \sim^*) - w(\Delta) 
	&= 	\big( w( \mcp, \mcq ; \sim^* ) - w(\Delta^*) \big) + \big( w (\Delta^*) - w(\Delta) \big) 	\\
	&\le 	k - w(P_1,Q_1; \sim^*) + \delta(\Theta^*) - \delta(\Theta)	\\
	&\le	k - \delta(\Theta).
\end{align*}
Since $\delta(\Theta) > 0$, it follows that  $\Delta$ is a near sequence relative to $(\mcp, \mcq)$, and since $k = w^{\circ}(\mcp)$ we also have $w(\mcp,\mcq) - w(\Delta) \le w^{\circ}(\mcp) - \delta(\Theta)$.

\bigskip

\noindent{\it Case 3:} $\Theta^*$ is a dihedral chord.

\smallskip

Define $H = G_{ (\mcp) } = G_{ (\mcq) } = G_{ (\mcr) }$ and note that $S \in G/H$.  It follows from our assumptions that $\Delta^*$ is a fringed sequence relative to $(\mcp, \mcq)$ with cyclic shift $S$ 
and we let $\ell$ be its length and let $\Lambda^*$ be an associated $H$-square chorus.  By assumption, $\Theta^*$ must have type $1$ or $2$ so we may choose a side $\Upsilon^*$ of $\Lambda^*$ which is partial.  Note that Lemma \ref{oct_ex_bound} implies $\delta( \Theta^* ) \le w_H( \Upsilon^* )$.  If $\Lambda^*$ is big, then $\hat{w}(\Lambda^*) \ge 2|H| + w_H(\Upsilon^*)$ so we have 
\begin{align*}
w( \mcp, \mcq ) - w(\Delta)	
	&= 	\big( w( \mcp, \mcq ) - w(\Delta^*) \big) + \big( w (\Delta^*) - w(\Delta) \big) 	\\
	&\le	w^{\circ}(\mcp) - w_H(\Upsilon^*) + \delta(\Theta^*) - \delta(\Theta)	\\
	&\le	w^{\circ}(\mcp) - \delta(\Theta).
\end{align*}
Since Observation \ref{length-width} implies that $\Delta^*$ is not long, and $\delta(\Theta) > 0$, the above equation implies that $\Delta$ is a near sequence relative to $(\mcp, \mcq)$.  Further, it inherits dihedral action and consecutive shift $S$ from $\Theta^*$.  
Otherwise $\Lambda^*$ is small, so $w_H(\Upsilon^*) \le \hat{w}(\Lambda^*) < 2|H|$ and we may apply Lemma \ref{adjust_fringe} to choose $\mcp', \mcq'$.  Now we have 
\begin{align*}
w( \mcp', \mcq' ) - w(\Delta)	
	&= 	\big( w( \mcp', \mcq' ) - w(\Delta^*) \big) + \big( w (\Delta^*) - w(\Delta) \big) 	\\
	&\le	w^{\circ}(\mcp') - w_H(\Upsilon^*) + \delta(\Theta^*) - \delta(\Theta)	\\
	&\le	w^{\circ}(\mcp') - \delta(\Theta).
\end{align*}
Since $\Delta^*$ is not short, $\Delta$ is a near sequence relative to $(\mcp', \mcq')$.  Furthermore, Lemma \ref{adjust_fringe} implies that this near sequence has dihedral action and consecutive shift $S$.  
\quad\quad$\Box$

\subsection{Crosses of Near Sequences}

In this subsection we will establish a lemma which provides strong information about critical crosses in near sequences.  
We begin with an easy lemma. 

\begin{lemma}
\label{basic_near_seq}
Let $(X,Y)$ be a duet which is a near sequence relative to $( {\mathcal P}, {\mathcal Q} )$, and let $A \subseteq X$ be finite and critical.  If $A$ is pure with respect to $\mcp$, then $A$ and $N(A)$ are unions of intervals of $\mcp$ and $\mcq$.
\end{lemma}

\noindent{\it Proof:} Observation \ref{near_obs} implies that $N(A)$ is pure with respect to $\mcq$.  Thus, 
$A' = \{ P \in \mcp \mid P \subseteq A \}$ is critical in $( \mcp, \mcq )$ and the result now follows from Lemma \ref{cross_seq}.
\quad\quad$\Box$

\bigskip

Let $(X,Y)$ be a duet, let ${\mathcal P}$, ${\mathcal Q}$ be systems of imprimitivity on $X$ and $Y$ and 
assume that $({ \mathcal P}, {\mathcal Q} )$ is a sequence.  We say that a finite set $A \subseteq X$ ($B \subseteq Y$) 
is \emph{weakly fringed} if the set of non-outer blocks of ${\mathcal P}$ (${\mathcal Q}$) forms a nontrivial interval, 
and every boundary block is an endpoint of this interval.  An infinite set is \emph{weakly fringed} if its complement is 
weakly fringed.  Finally, we say that a cross $(A,B)$ is \emph{weakly fringed} if both $A$ and $B$ are weakly fringed.
We are now ready for the main lemma from this subsection.

\begin{lemma}
\label{fringed_cross}
Let $\Delta = (X,Y)$ be a $G$-duet which is a near sequence relative to $( {\mathcal P}, {\mathcal Q} )$ and let $(A,B)$ be a maximal critical cross of $(X,Y)$.  If $A$ is not contained in any block of ${\mathcal P}$ and $B$ is not contained in any block of ${\mathcal Q}$, 
then $(A,B)$ is weakly fringed.
\end{lemma}

\noindent{\it Proof:} 
Define $\Delta' = ( {\mathcal P}, {\mathcal Q} )$, set $\ell = d( {\mathcal P}, {\mathcal Q} )$, set $k = w^{\circ}(\mcp) = w^{\circ}(\mcq)$ and let $H = G_{ (\mcp) } = G_{ (\mcq) }$.  
Now suppose (for a contradiction) that the lemma fails for $\Delta$, let $(A,B)$ be a counterexample so that 
\begin{enumerate}
\item[$\alpha$)] $\delta(A,B)$ is maximum.
\item[$\beta$)] The total number of boundary blocks in ${\mathcal P}$ and ${\mathcal Q}$ is minimum (subject to $\alpha$).
\item[$\gamma$)] $\min\{ w(A), w(B) \}$ is minimum (subject to $\alpha$, $\beta$).
\end{enumerate}
Criteria $\alpha$ and $\beta$ have a convenient behaviour in relation to purification.  Namely, if we modify $(A,B)$ by purifying at 
any boundary block $T \in \mcp \cup \mcq$, the new cross will be superior to $(A,B)$ with regard to $\alpha$ and $\beta$.  To see this, let $(A',B')$ and $(A'',B'')$ be the weak purification and the purification of $(A,B)$ at $T$.  Now Theorem \ref{purify_thm} implies $\delta(A'', B'') \ge \delta(A',B') \ge \delta(A,B)$.  Since $N(T)$ is pure 
(by Observation \ref{near_obs}) it follows that $(A',B')$ is superior to $(A,B)$ with regard to $\alpha$ and $\beta$.  Now either $(A'',B'') = (A',B')$ or 
$\delta(A'',B'') > \delta(A',B')$, and in either case $(A'',B'')$ is superior to $(A,B)$ with regard to $\alpha$ and $\beta$.

Note that by Lemma \ref{basic_near_seq} there must be at least one boundary block, and then by maximality (and Observation \ref{near_obs}) it follows that there must be a boundary block in both ${\mathcal P}$ and ${\mathcal Q}$.  
Next we establish a couple of key claims.

\bigskip

\noindent{(i)} Either there exists an interval $I$ of ${\mathcal P}$ so that $A \subseteq \cup_{P \in I} P$ and $|N_{\Delta'}(I)| = |I| + \ell - 1$ or there exists an interval $I'$ of ${\mathcal Q}$ so that $B \subseteq \cup_{Q \in I'} Q$ and $|N_{\Delta'}(I')| = |I| + \ell - 1$.  

\smallskip

We shall assume (without loss) for the purposes of this claim that $w(A) \le w(B)$.  We may also assume $w(B) < \infty$ as otherwise the claim holds 
trivially.  Choose a boundary block $P_0 \in {\mathcal P}$ and let 
$(A',B')$ be the cross obtained by purifying $(A,B)$ at $P_0$.  Note that $w(B) > \tfrac{1}{2} \overline{w}(\Delta) \ge k$ so by Theorem \ref{purify_thm} we must have $B' \neq \emptyset$.  If $B'$ is not contained in any block of ${\mathcal Q}$, then we may apply the lemma inductively to $(A',B')$, and we deduce that $(A',B')$ is weakly fringed.  But then there exists an interval $I$ of $\mcp$ with $A' \subseteq \cup_{P \in I} P$ and 
$|N_{\Delta'}(I)| = |I| + \ell - 1$, so the same holds for $A$, thus proving the claim.  Thus, we may assume that $B' \subseteq Q_0$ for some $Q_0 \in {\mathcal Q}$.  If $Q_0$ is inner with respect to $B$, then $I = {\mathcal P} \setminus N_{\Delta'}(Q_0)$ satisfies the claim for $A$, so we may assume that $Q_0$ is a boundary block of $B$.  Next consider the cross $(A'',B'')$ obtained from $(A,B)$ by weakly purifying at $P_0$ and then purifying at $Q_0$.  Since this cross is critical and $B'' = Q_0$ it follows that every $P \in {\mathcal P} \setminus N_{\Delta'}(Q_0)$ 
must be non-outer in $A''$.  However, then there must exist at least two blocks of ${\mathcal P}$ which are non-outer with respect to $A$ and are not incident with $Q_0$.  Now define $(A''',B''')$ to be the cross obtained by purifying $(A,B)$ at $Q_0$.  It follows from our assumptions that $A'''$ is not contained in any block of ${\mathcal P}$ and $B'''$ is not contained in any boundary block of ${\mathcal Q}$.  So, by induction $(A''',B''')$ is weakly fringed, and since $B \subseteq B'''$ it follows that there exists an interval $I'$ which satisfies the claim for $B$.  

\bigskip

By claim (i) we may assume that there is an interval $I$ of $\mcp$ so that $A \subseteq \cup_{P \in I} P$ and 
$|N_{\Delta'}(I)| = |I| + \ell - 1$.  Furthermore, by possibly interchanging $A,B$ we may assume that either $w(A) \le w(B)$ 
or that $\mcp$ has no inner block (if $w(B) < w(A)$ and $\mcp$ has an inner block $P$, then 
$I' = {\mathcal Q} \setminus N_{\Delta'}(P)$ is an interval of $\Delta'$ and $B \subseteq \cup_{Q \in I'} Q$).  
By choosing $I$ minimal, we may further assume that each endpoint of $I$ is a non-outer block.  

\bigskip

\noindent{(ii)} Every block in $I$ is either inner or boundary (with respect to $A$).

\smallskip

To see this, modify $(A,B)$ to form $(A',B')$ by weakly purifying at every boundary block in ${\mathcal P}$.  Then 
the set $A'$ is pure with respect to ${\mathcal P}$ and it is critical.  If $N(A') \neq Y$ then  
 Lemma \ref{basic_near_seq} implies that $A'$ is the union of the blocks in $I$, and this gives the desired result.  
 Otherwise $w(\Delta) \le \ell k$ and $w(N(A)) = ( |I| + \ell - 1) k$ so again 
 $A'$ must contain every block in $I$ and we are finished.

\bigskip

Now it will be helpful to introduce coordinates.  Assume that ${\mathcal P} = \{ P_j \mid j \in J\}$ and ${\mathcal Q} = \{ Q_j \mid j \in J\}$ 
where $J = {\mathbb Z}$ or $J = {\mathbb Z}/n {\mathbb Z}$ and $P_i \sim Q_j$  if $i \le j \le i + \ell- 1$.  Assume that $I = \{ P_1, P_2, \ldots, P_m \}$, so the blocks $Q_1, Q_2, \ldots, Q_{\ell + m - 1}$ are distinct.  Note that (ii) implies that $A \cap P_i \neq \emptyset$ for $1 \le i \le m$.  

Next we consider the case that $m \ge 3$.  Since $A \cap P_1 \neq \emptyset$ the block $Q_1$ must be either boundary or outer.  
First let us suppose that $Q_1$ is boundary and let $(A',B')$ be the cross obtained from $(A,B)$ by purifying at $Q_1$.  Note that 
our assumptions imply $A' = A \setminus P_1$.  Now, the lemma applies nontrivially to $(A',B')$ and it follows that $(A',B')$ is weakly fringed, 
so in particular $\cup_{i=3}^{\ell + m - 2} Q_i \subseteq N(A')$.  It follows that blocks $Q_3, \ldots Q_{\ell + m - 2}$ must be outer 
with respect to our original set $B$.  Next we shall show that a similar conclusion holds under the assumption that $Q_1$ is outer.  
In this case we let $(A',B')$ be a maximal cross containing $(A \setminus P_1, B \cup Q_1)$.  Note that by the maximality of $(A,B)$ 
we must again have $A' = A \setminus P_1$.  We claim that the cross $(A',B')$ is superior to $(A,B)$ with regard to our three optimization 
criteria.  It is immediate that $\delta(A',B') \ge \delta(A,B)$ and this claim holds already 
if $\delta(A',B') > \delta(A,B)$.  Otherwise we must have that $P_1$ is inner (which implies that 
$w(A) \le w(B)$) and $B' = B \cup Q_1$ (which implies that $(A',B')$ has the same total number of boundary blocks as $(A,B)$) so 
the cross $(A',B')$ improves upon $(A,B)$ with regard to our third optimization criteria.  Therefore, the lemma applies nontrivially to 
$(A',B')$ and we conclude, as before, that  $\cup_{i=3}^{\ell + m - 2} Q_i \subseteq N(A')$ so blocks $Q_3, \ldots, Q_{\ell + m - 2}$ 
must be outer with respect to $B$.  By a similar argument applied to the block $Q_{\ell + m - 1}$ we conclude that blocks 
$Q_2, \ldots, Q_{\ell + m - 3}$ must be outer with respect to $B$.  But then $B$ is weakly fringed, and by maximality $A$ is 
weakly fringed.  

Thus, we may now assume $m = 2$.  The result holds trivially whenever one of $P_1$ or $P_2$ is inner, so we may assume $P_1$ and
$P_2$ are boundary.  Now we shall consider the case $\ell \ge 3$.  Observe that if there exists $g \in G_{ P_1 } \setminus H$ then $gQ_1 = Q_{\ell}$ and 
$gQ_{\ell} = Q_1$.  Based on this we can ``split'' our duet according to the following procedure.  Define the duet $\Delta_1 = (X,Y; \sim_1)$ by the rule that $x \in X$ and $y \in Y$ satisfy $x \sim_1 y$ if $x \sim y$  and the blocks $P_i \in {\mathcal P}$ with $x \in P_i$ and $Q_j \in {\mathcal Q}$ with $y \in Q_j$ satisfy $i < j < i + \ell - 1$.  We define the duet $\Delta_2 = (X,Y; \sim_2)$ by the rule that $x \in X$ and $y \in Y$ satisfy $x \sim_2 y$ if $x \sim y$ and $x \not\sim_1 y$.  
It follows from our earlier observation that
$\Delta_1$ and $\Delta_2$ are $G$-duets and further $w(\Delta) = w(\Delta_1) + w(\Delta_2)$.  Since $B$ is not weakly fringed, we may assume (by possibly interchanging $P_1$ and $P_2$ and reindexing) that there exists $y \in B \cap (\cup_{i=2}^{\ell-1} Q_i)$.  Now set $A_1 = A \cap P_1$ and observe that 
$N(y) \cap A_1 = \emptyset$ implies that 
$w(\Delta_1) = w(N_{\Delta_1}(y)) \le (\ell-2)k - w(A_1)$.  Choose a point $x \in A_1$ and note that 
$w(\Delta_2) = w(N_{\Delta_2}(x)) \le w(N(x) \cap Q_1) + k$.  Now let $(A',B')$ be the cross obtained from $(A,B)$ by weakly purifying at $P_2$.  Then $A' = A_1 \cup P_2$ is critical, but 
\begin{align*}
w(N(A')) 
	&\ge \ell k + w(N(x) \cap Q_1) 	\\
	&\ge \ell (k-1) + w(\Delta_2)		\\
	&\ge k  + w(A_1) + w(\Delta_1) + w(\Delta_2) \\
	&= w(A') + w(\Delta) 
\end{align*}
which is a contradiction.  

Finally, we are left with the case $\ell = m = 2$.  Here we may assume that $P_1$, $P_2$, and $Q_2$ are boundary.  
Since the lemma we are establishing is unaffected by identifying clones, and 
splitting into clones, we may apply Proposition \ref{sabidussi} to assume that $G$ acts regularly on $X$ and $Y$.  

First consider the case that $H = G_{P_1}$.  In this case $\Delta_1 = (P_1,Q_1)$ and $\Delta_2 = (P_1,Q_2)$ are $H$-duets and 
$w_G(\Delta) = w_H(\Delta_1) + w_H(\Delta_2)$.  Set $A_1 = A \cap P_1$ and choose $x \in A_1$ and $y \in B \cap Q_2$.  
Since $N(y) \cap A_1 = \emptyset$ we have $w_H(\Delta_2) = w_H( N_{\Delta_2}(y) ) \le k - w_H(A_1) = k - w_G(A_1)$.  In 
addition $w_G( N_{\Delta}(A_1) \cap Q_1) \ge w_G(N_{\Delta}(x) \cap Q_1 ) = w_H(\Delta_1)$.  Now, let $(A',B')$ be the cross 
obtained from $(A,B)$ by purifying at $P_2$.  Then we have
\begin{align*}
w_G(N(A'))
	&=	2k + w_G( N_{\Delta}(A_1) \cap Q_1)	\\
	&\ge	k + w_G(A_1) + w_H(\Delta_2) + w_H(\Delta_1)	\\
	&=	w(A') + w(\Delta)
\end{align*}
and this contradicts the criticality of $A'$.

Next consider the case that $H \neq G_{P_1}$.  Since we are assuming that $G$ acts regularly, it follows that every 
$P_j$ $(Q_j)$ is the union of two $H$-orbits, which we denote by $P_j^0$ and $P_j^1$ ($Q_j^0$ and $Q_j^1$).  
Furthermore, by indexing consistently, we may assume that for $0 \le i,j \le 1$ there exist $w_{ij}$ so that 
$w_H(P_m^i, Q_m^j) = w_{ij}$.  It then follows that 
$w_H(P_1^0,Q_2^1) = w_{10}$ and $w_H(P_1^1,Q_2^1) = w_{00}$ and $w_H(P_2^0,Q_2^1) = w_{01}$ and $w_H(P_2^1,Q_2^1) = w_{11}$ as shown in the figure.

\begin{figure}[ht]
\centerline{\includegraphics[height=4cm]{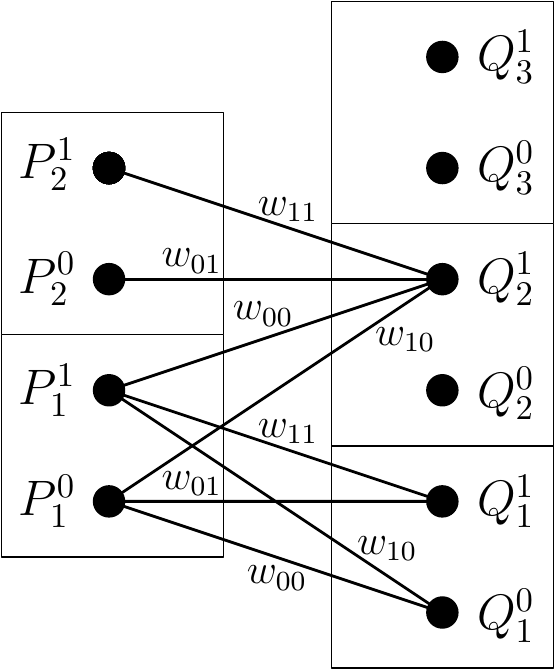}}
\label{seq22fig}
\end {figure}

Since $B$ is not weakly fringed, we may assume without loss that $B \cap Q_2^1 \neq \emptyset$.  Then by possibly reindexing, we may assume that $w_{00} + w_{10} \ge w_{01} + w_{11}$.  Choose a point $y \in B \cap Q_2^1$, and for $i=1,2$ let $A_1^i = A \cap P_1^i$.  Since $A_1^i$ is disjoint from $y$ we must have 
\begin{equation}
\label{seq22eq1}
w(A_1^0) \le |H| - w_{10} \quad\quad\mbox{and}\quad\quad w(A_1^1) \le |H| - w_{00}
\end{equation}
Furthermore, if $A_1^i$ is nonempty then we have
\begin{equation}
\label{seq22eq2}
w(N(A_1^i) \cap Q_1^j) \ge w_{ij}. 
\end{equation}
Let $(A',B')$ be the critical cross obtained from $(A,B)$ by purifying at $P_2$ (so $w(A') = 2|H| + w(A_1)$).  First suppose that $A_1^0 \neq \emptyset$ and $A_1^1 \neq \emptyset$.  In this case, averaging equation \ref{seq22eq2} gives us $w(N(A_1) \cap Q_1) \ge \frac{1}{2}( w_{00} + w_{01} + w_{10} + w_{11}) \ge w_{01} + w_{11}$.  But then using equation \ref{seq22eq1} we find
\begin{align*}
w(N(A'))
	&\ge		4|H| + w_{01} + w_{11}	\\
	&\ge		2|H| + w(A_1) + w_{00} + w_{10} + w_{01} + w_{11}	\\
	&=		w(A) + w(\Delta)
\end{align*}
which is a contradiction.  Next suppose that $A_1^1 = \emptyset$ (so $A_1^0 \neq \emptyset$).  In this case 
\begin{align*}
w(N(A'))
	&\ge		4|H| + w_{00} + w_{01}	\\
	&\ge		3|H| + w(A_1) + w_{10} +  w_{00} + w_{01}	\\
	&\ge		w(A') + w(\Delta)
\end{align*}
Finally, if $A_1^0 = \emptyset$ (so $A_1^1 \neq \emptyset$) we have
\begin{align*}
w(N(A'))
	&\ge		4|H| + w_{10} + w_{11}	\\
	&\ge		3|H| + w(A_1) + w_{00} +  w_{10} + w_{11}	\\
	&\ge		w(A) + w(\Delta)
\end{align*}
which is again a contradiction.  This completes the proof.
\quad\quad$\Box$

\subsection{Near Sequence Stability}

In this subsection we prove our stability lemma for near sequences.

\begin{lemma}[Near Sequence Stability]
\label{near_sequence_stability}
Let $\Theta$ be a maximal connected critical trio with a side which is a near sequence.  Then $\Theta$ is a pure chord, 
a cyclic chord, or a dihedral chord of type $1$, $2A$, $2B$, or $2C$.
\end{lemma}

\noindent{\it Proof:} 
Let $\Theta = (X,Y,Z)$ and assume that $\Delta = (X,Y)$ is a  near sequence
relative to $( {\mathcal P}, {\mathcal Q} )$.  Assume further that ${\mathcal P} = \{ P_j \mid j \in J \}$ and 
${\mathcal Q} = \{ Q_j \mid j \in J \}$ for $J= {\mathbb Z}$ or ${\mathbb Z} / n {\mathbb Z}$ where $P_i \sim Q_j$ if 
$i \le j \le i + \ell-1$.  Choose $z \in Z$ and let $(A,B)$ be the cross of $(X,Y)$ associated with $z$.  It now follows from 
Lemma  \ref{fringed_cross} that $(A,B)$ is weakly fringed, so we may assume that the non-outer blocks of $A$ form the 
interval $\{P_1, \ldots, P_{m} \}$.  If $(A,B)$ is pure with respect to $(\mcp,\mcq)$, then by maximality 
$\Delta^{\bullet}$ is a sequence.  Since $\Delta^{\bullet}$ is a side of $\Theta^{\bullet}$ (by maximality), 
Lemma \ref{sequence_stability} completes the proof.  It follows from this
and the maximality of $(A,B)$ that either ${\mathcal P}$ and ${\mathcal Q}$ each have exactly
one boundary block, or they each have exactly two.  Now, assume $\Theta$ is a $G$-trio and let 
$H = G_{ ( {\mathcal P} ) } = G_{ ( {\mathcal Q} ) }$.  
We now split into cases depending on the action on $( {\mathcal P}, {\mathcal Q} )$.

\bigskip

\noindent{\it Case 1:} The sequence $({\mathcal P}, {\mathcal Q} )$ has cyclic action. 

\smallskip

First suppose that $A$ and $B$ each have two boundary blocks.  Then $\Delta_1 = (P_1,Q_1)$ and 
$\Delta_{m} = (P_{m}, Q_{\ell + m-1})$ are both nontrivial $H$-duets and by maximality, we must have $w_G(\Delta) = w(\Delta_1) + w(\Delta_{m}) + (\ell-2)|H|$.  
Set $A_1 = A \cap P_1$ and $A_{m} = A \cap P_{m}$ and note that $w_G(A) = w_G(A_1) + w_G(A_{m}) + (m-2)|H|$.  
Then by Lemma \ref{nontriv_bound}, we have for $i=1$ and $i=m$ that 
$w_H(N_{\Delta_i}(A_i)) \ge w_H(\Delta_i) + w_H(A_i) - \frac{1}{2}|H|$.  But then 
\begin{align*}
w_G(N(A)) 
	&=	w_H(N_{\Delta_1}(A_1))  + w_H(N_{\Delta_2}(A_2)) + (\ell + m-3) |H|	\\
	&\ge	w_H(\Delta_1) + w_H(A_1) + w_H(\Delta_{m}) + w_H(A_{m}) + (\ell + m-4) |H| \\
	&=	w_G(A) + w_G(\Delta)
\end{align*}
which contradicts the criticality of $A$.  It follows that $A$ and $B$ must each have one boundary block.  Now it follows from 
maximality that $\Theta$ is a cyclic chord.

\bigskip

\noindent{\it Case 2:} The sequence $({\mathcal P}, {\mathcal Q} )$ has dihedral action.

\smallskip

For every $j \in J$, let $R_j$ be the set of all points in $z$ for which the cross associated with $z$ has nonempty intersection with all of the blocks $P_j, P_{j+1}, \ldots, P_{j + m - 1}$.  By maximality, every $z \in R_j$ will have $( \{z\}, P_{j+1} \cup \ldots \cup P_{j + m-2} )$ full, and by the assumption that the action is dihedral (and the assumption $(A,B)$ is impure), there must exist $z' \in R_i$ with $( \{z'\}, P_j)$ not full and $z'' \in R_j$ with 
$(\{ z''\}, P_{j+m-1} )$ not full.  It follows from this and maximality that 
$\Theta$ is a dihedral chord relative to ${\mathcal P}$, ${\mathcal Q}$, and $\{ R_j \mid j \in J \}$.  By Theorems \ref{oct_type} and 
\ref{classify_oct2} we find that the associated octahedral chorus must have type $-1$, $0$, $1$, $2A$, $2B$, or $2C$, so $\Theta$ must have one of these 
types.  Now by Lemmas \ref{kill-1} and \ref{type0alt}, we deduce that either $\Theta$ is a pure chord, or a dihedral chord of type 
$1$, $2A$, $2B$, or $2C$.  

\bigskip

\noindent{\it Case 3:} The sequence $({\mathcal P}, {\mathcal Q} )$ has split action.

\smallskip

Any two consecutive blocks of $\mcp$ or $\mcq$ form a block of imprimitivity, so it follows from the connectivity of $\Theta$ that 
neither $A$ nor $B$ can be contained in two consecutive blocks.  Thus $m \ge 3$ and since $Q_2 \cup Q_3 \ldots \cup Q_{\ell + m - 2} \subseteq N(A)$ so 
we must have $|J| \ge \ell + m \ge \ell+3$.  Suppose (for a contradiction) that $|J| = \ell +3$.  In this case $m=3$ and we must have that $P_2$ is inner, $P_1,P_3$ are boundary, and on the other side $Q_{\ell+3}$ is inner, and $Q_{\ell+2}, Q_{\ell+4} = Q_1$ boundary.  
Now define the $H$-duets $\Delta_1 = (P_1,Q_1)$ and $\Delta_2 = (P_3,Q_{\ell+2})$ and note that by this choice we must have
$w_G(\Delta) \le w_H(\Delta_1) + w_H(\Delta_2) + (\ell-2) |H|$.  Now setting $A_1 = A \cap P_1$ and $A_3 = A \cap P_3$ we have $Q_1 \not\subseteq N(A_1)$ 
and $Q_{\ell+2} \not\subseteq N(A_3)$ so Lemma \ref{nontriv_bound} implies
\begin{align*}
w_G(N(A))	&=	w_H(N_{\Delta_1}(A_1)) + w_H(N_{\Delta_2}(A_3)) + \ell |H|		\\
			&\ge	w_H(\Delta_1) + w_H(A_1) + w_H(\Delta_2) + w_H(A_3) + (\ell-1)|H|		\\
			&\ge	w_G(A) + w_G(\Delta).
\end{align*}
which contradicts the criticality of $A$.  Therefore, we must have $|J| > \ell+3$.  Since $|J|$ must be $\infty$ or even and $\ell$ is odd (by Lemma \ref{alt_seq_obs}) we conclude $|J| \ge \ell+5 \ge 8$.

Note that $N( P_1 \cup P_2) = Q_1 \cup \ldots \cup Q_{\ell+1}$ and that $(P_1 \cup P_2, Q_{2i+1} \cup Q_{2i+2})$ is full for 
$1 \le i < \frac{\ell-1}{2}$.  Now we define the systems of imprimitivity
${\mathcal P}' = \{ P_1 \cup P_2, P_3 \cup P_4, \ldots \}$ and ${\mathcal Q}' = \{ Q_1 \cup Q_2, Q_3 \cup Q_4, \ldots \}$.  It follows from our assumptions that 
$(X,Y)$ is a fringed sequence relative to $(\mcp', \mcq')$.  Now the argument from the previous case 
implies that $\Theta$ is either a pure chord or a dihedral chord of type $1$, $2A$, $2B$, or $2C$ as desired.
\quad\quad$\Box$

\section{Near Graphs and Doubling Blocks}
\label{near_graph_sec}

In this section we will prove a stability result for near true graphs and duets which have doubling blocks.  As will 
be detailed shortly, these two classes of duets are very closely related.

\subsection{Preliminaries}

Here we prove a sequence of three lemmas which will culminate with an important, albeit rather technical bound.  To motivate the first lemma in this sequence, recall that by Theorem \ref{hamidoune} every partial duet $\Delta = (X,Y)$ has a nontrivially critical block of imprimitivity $T$ which has the same deficiency as $\Delta$.  However, this lemma gives us no control over whether $T$ is included in $X$ or $Y$.  Our first lemma permits us to choose $T \subseteq Y$ under certain additional assumptions.

\begin{lemma}
Let $\Delta =(X,Y)$ be a finite partial duet of degree $d$.  Assume that there is at most one integer 
$m$ with $d \le m < |Y|$ which is not relatively prime with $|Y|$, and if such a number exists, it has 
the form $m = (s-1)t$ where $st = |Y|$ and $3t \le 2( |Y| - d +1)$.  Then there exists a nontrivially 
critical block of imprimitivity $Q \subseteq Y$ with $\delta(Q) = \delta(\Delta)$.
\end{lemma}

\noindent{\it Proof:} By Theorem \ref{hamidoune} we may choose a nontrivially critical block of imprimitivity
$P$ so that $\delta(P) = \delta(\Delta)$ and we may assume $P \subseteq X$ as otherwise we 
are done.  Let $\mcp$ be the system of imprimitivity associated with $P$ and set $m = |N(P)|$.  Now 
$w(N(P))$ is given by the equation $w(N(P)) = m w^{\circ}(Y) = m \frac{ |G| }{ |Y| }$.  However, 
by Proposition \ref{block_mult} we may choose $\ell \in {\mathbb Z}$ so that $w(N(P)) = \ell w(P) = \ell \frac{ |G| }{ |\mcp| }$.  
Since $m \ge d$ this gives us 
\begin{equation}
\label{awkward-div}
\frac{d}{|Y|} \le \frac{m}{|Y|} = \frac{w(N(P))}{|G|} = \frac{\ell}{|\mcp|} < 1.
\end{equation}
If $m$ and $|Y|$ are relatively prime, then $|\mcp| \ge |Y|$, and choosing $y \in Y$ we have
\[ \delta( \{y\} ) = w(y) + w(\Delta) - w(N(y)) = w(y) \ge w(P) \ge w(P) + w(\Delta) - w(N(P)) = \delta(P). \]
Thus, $\{y\}$ is a block of imprimitivity which satisfies the lemma.  
Otherwise our assumptions imply $m = (s-1)t$ and $|Y| = st$.  First suppose that $|\mcp| = s$ (so $\ell = s-1$).  Since every $P' \in \mcp$ satisfies
$\overline{w}(N(P')) = \frac{|G|}{s}$ we have that $\mcq = \{ Y \setminus N(P') \mid P' \in \mcp \}$ is a system of imprimitivity.  Further, 
every $Q \in \mcq$ has $\delta(Q) = \delta(P)$ so it is a nontrivially critical block of maximum deficiency. 
Thus, we may now assume $|\mcp| \neq s$.  Equation \ref{awkward-div} gives us $\frac{\ell}{|\mcp|} = \frac{m}{|Y|} = \frac{s-1}{s}$, which 
implies $|\mcp| \ge 2s$.  Now by the assumption $3t \le 2(|Y| - d + 1)$ we find
\[ \delta(P)	 =	w(P) + w(\Delta) - w(N(P))	 \le	\frac{ |G| }{2s} + \frac{ d |G| }{st } - \frac{(s-1) t |G| } {st} \le	\frac{|G|}{st} = w^{\circ}(Y). \]
However, this again implies that a single point of $Y$ is a nontrivially critical block of maximum deficiency.  
\quad\quad$\Box$

\begin{lemma}
Let $\Delta = (X,Y)$ be a duet with $4 \le |Y| \le 9$ and $\frac{|Y|}{2} < d(X,Y) < \min \{6, |Y| \}$.  
If $S \subseteq Y$ satisfies $|S| = 2$ then one of the following holds
\begin{enumerate}
\item $w(N(S)) \ge w(\Delta) + w^{\circ}(Y)$
\item $|Y| = 6$, $S$ is a block of clones, and $\Delta^{\bullet}$ is a triangle.
\end{enumerate}
\end{lemma}

\noindent{\it Proof:} Let $d = d(X,Y)$ and $k = w^{\circ}(Y)$ and observe that $w(\Delta) = dk$.  We shall assume 
that the first outcome does not hold, and prove that the second does.  By the previous lemma, 
we may choose a nontrivially critical block of imprimitivity $Q \subseteq Y$ so that 
$\delta(Q) = \delta(\Delta) \ge \delta(S) = w(S) + w(\Delta) - w(N(S)) > k$.  Since each singleton from $Y$
has deficiency $k$ it follows that $Q$ is a nontrivial block of imprimitivity.  Furthermore, since $w(N(Q))$ is a multiple
of $w(Q)$ strictly between $w(Q)$ and $w(Y) = |G|$ (by Proposition \ref{block_mult}) it must be that $(|Y|, |Q|)$ is one of $(6,2)$, $(8,2)$, or 
$(9,3)$.  First suppose $|Y| = 8$ and $|Q| = 2$.  In this case $d = 5$ so $w(\Delta) = 5k$ and 
since $w(N(Q))$ is a multiple of $2k$ we must have $w(N(Q)) = 6k$.  But then $\delta(Q) = k$ which is a contradiction.  
In the remaining cases, let  ${\mathcal Q}$ be the system 
of imprimitivity associated with $Q$, and observe that $(X,{\mathcal Q})^{\bullet}$ is a triangle.  First suppose $|Y| = 6$ 
and $|Q| = 2$.  In this case $4k = w(N(T)) \ge w(\Delta) = dk$ so we must have $d=4$.  However, then 
$w(N(Q)) = w(\Delta)$ so $Q$ is a block of clones, and the second outcome is satisfied.  In the remaining case, $|Y| = 9$ and 
$|Q| = 3$ and this implies $d=5$.  If $S$ contains points in two distinct blocks of ${\mathcal Q}$ then 
each point has a neighbourhood of weight $5k$ and these neighbourhoods overlap in a set of weight at most $3k$ 
so $w(S) \ge 7k$ which contradicts our assumption.  Otherwise, there exists $Q \in \mcq$ with $S \subseteq Q$ and 
we must have $[S] = Q$.  Now either $w(N(S)) = w(N(Q)) = 6k$ or by Theorem \ref{critical_closure} 
$w(N(S)) \ge w(S) + w(\Delta) - \frac{1}{3}w(Q) = 6k$, and since $6k = w(\Delta) + w^{\circ}(Y)$ this contradicts our assumption.
\quad\quad$\Box$

\bigskip

Now we will turn our attention to duets which have graphic quotients (the main subject of this section).  We begin
by introducing a useful bit of terminology.

\bigskip

\noindent{\bf Definition:} Let $\Delta = (X,Y)$ be a $G$-duet and let $\mcv, \mce$ be systems of imprimitivity on $X,Y$ so that 
$(\mcv, \mce)$ is a graph.  Let $E \in \mce$ and assume $E$ has ends $V_1,V_2 \in \mcv$.  We define the \emph{profile} of 
$\Delta$ \emph{with respect to} $(\mcv, \mce)$ to be $\{ w( N_{\Delta}(E) \cap V_1 ), w( N_{\Delta}(E) \cap V_2 ) \}$ (note that 
this is independent of the choice of $E$).

\begin{lemma} 
\label{technical}
Let  $\Delta = (X,Y)$ be a $G$-duet, and let $\mcv, \mce$ be systems of imprimitivity on $X,Y$.  Assume that $(\mcv,\mce)$ is a graph 
of degree $d \ge 4$ with profile $\{w_1,w_2\}$, and assume that $\mcv$ is critical.  Let $k = \frac{1}{2}w^{\circ}(\mce)$.  If $w_1 \neq w_2$ 
then $2k | w_i$ for $i=1,2$.  Otherwise $k | w_1$, and furthermore, whenever $E_1,E_2 \in \mce$ are incident with $V \in \mcv$ and 
$E_1 \neq E_2$ one of the following holds
\begin{enumerate}
\item $w( N_{\Delta}(E_1 \cup E_2) \cap V ) \ge \min \{ dk, 6k, w_1 +k \}$ 
\item The incidence geometry $( V,  \{ E \in {\mathcal E} \mid E \sim V \} )^{\bullet}$ is a triangle and $d=6$.
\end{enumerate}
\end{lemma}

\noindent{\it Proof:} Our assumptions imply $w_G^{\circ}(\mcv) = dk$ and $w_G( \mcv, \mce ) = 2dk$.  Since every block of $\mcv$ is critical, we 
have
\begin{equation}
\label{w1w2-bound}
dk < w_G( \Delta ) \le w_G ( X, \mce ) = w_1 + w_2
\end{equation}
First suppose that $G$ acts arc-transitively on $(\mcv, \mce)$ (so $w_1 = w_2$), let $V \in \mcv$, define 
$\mcq = \{ E \in \mce \mid E \sim V \}$, and let $E_1,E_2 \in \mcq$ be distinct.  Now $(V,\mcq)$ is an $H$-duet for $H = G_V$, and we have
$w_H^{\circ}(\mcq) = w_H(E_1) = \frac{1}{2}w_G(E_1) = k$.  This gives us
\[ w_1 = w_G( N(E_1) \cap V ) = w_H( N(E_1) \cap V ) = w_H(V,\mcq) = d(V,\mcq) w_H^{\circ}(\mcq) = d(V,\mcq) k \]
so $w_1$ is a multiple of $k$.  If $w_1 \ge 6k$ then the first outcome holds and we are finished.  So, by equation \ref{w1w2-bound} we may assume $d \le 9$.  
If $d=6$ and $(V, \mcq)^{\bullet}$ is a triangle, then the second outcome is satisfied.  Otherwise the previous lemma gives
\[ w_G( N_{\Delta}(E_1 \cup E_2) \cap V ) = w_H( N(E_1 \cup E_2) \cap V ) \ge w_H(V,\mcq) + w_H^{\circ}(\mcq) = w_1 + k \]
and this completes the proof for this case.

In the remaining case $G$ does not act arc-transitively on $(\mcv, \mce)$ which implies that $d$ is even.  Set $H = G_V$ as before, and observe 
that $\{ E \in \mce \mid E \sim V \}$ consists of two $H$-orbits which we denote by $\mcq_1$, $\mcq_2$.  For every $E \in \mcq_i$ we have 
$w_G( N_{\Delta}(E) \cap V ) = w_H( N_{\Delta}(E) \cap V) = w_H(V, \mcq_i)$ so we may asume (without loss) that $w_i = w_H( V, \mcq_i)$ for $i =1,2$.  
Furthermore, $w_H^{\circ}( \mcq_i ) = w_H(E) = w_G(E) = 2k$ so 
$w_i = w_H( V, \mcq_i ) = d(V, \mcq_i) w_H^{\circ}$ is a multiple of $2k$.  

If $w_1 \neq w_2$ then we have nothing left to prove, so we may assume $w_1 = w_2$.  Now let $E_1,E_2 \in \mce$ be distinct with 
$E_1, E_2 \sim V$, and suppose (without loss) that $E_1 \in \mcq_1$.  Now we have
$w_G( N_{\Delta}(E_1) \cap V) = w_H (V, \mcq_1) = w_1$ so we are finished unless $w_1 < 6k$.  
It now follows from our assumptions and equation \ref{w1w2-bound} that $w_1 = 4k$, and 
$d \in \{4,6\}$.  If $d = 4$ then $w_G(N_{\Delta}(E_1) \cap V) = dk$ and we are done.  So, we may now assume $d = 6$, which implies
$|\mcq_1| = |\mcq_2| = 3$.  It then follows from $w_1 = w_2 = 4k$ that both $(V,\mcq_1)^{\bullet}$ and $(V, \mcq_2)^{\bullet}$ are triangles.  
If $E_2 \in \mcq_1$ then $V \subseteq N_{\Delta}(E_1 \cup E_2)$ so $w_G( N_{\Delta} (E_1 \cup E_2) \cap V ) = dk$ and we are done. 
 So, we may assume $E_2 \in \mcq_2$.  Now for $i=1,2$ let $\mcp_i$ be the partition of $V$ into clones in the duet $(V, \mcq_i)$ (so $|\mcp_1| = |\mcp_2| = 3$), 
 and for $i=1,2$ choose $P_i \in \mcp_i$ to be the block of points not incident with $E_i$.  Set $T = P_1 \cap P_2$.   If $w(T) \le k$ then we have 
$w_G ( N_{\Delta}( E_1 \cup E_2 ) ) = 6k - w(T) \ge 5k$ satisfying outcome 1.  Otherwise, note that both $\mcp_1$ and $\mcp_2$ are systems 
of imprimitivity under the action of $G$, so $T$ is also a block of imprimitivity.  But then $w(T) > k$ forces $P_1 = P_2 = T$.  Thus,  
$\mcp_1 = \mcp_2$, and $( V, \{ E \in \mce \mid E \sim V \} )^{\bullet}$ is a triangle satisfying outcome 2.
\quad\quad$\Box$

\subsection{Star}

The main goal of this subsection is to prove a lemma which will be  used to handle crosses in duets with graphic quotients 
where all non-outer edges are incident with a common vertex (thus forming a star).  However, we begin by 
introducing a new definition which will capture the duets we are interested in a form which is convenient to work with.

\bigskip

\noindent{\bf Definition:} Let $(X,Y)$ be a duet and let ${\mathcal V}$, ${\mathcal E}$ be systems of imprimitivity on $X$, $Y$. 
We say that $(X,Y)$ is a \emph{weak true graph} (\emph{relative to} $({\mathcal V}, {\mathcal E})$) if the following properties hold
\begin{enumerate}
\item $( {\mathcal V}, {\mathcal E} )$ is a true graph.
\item ${\mathcal V}$ is critical.
\item If $d( {\mathcal V}, {\mathcal E} ) = 3$ then the associated profile is $\{ w^{\circ}( {\mathcal V} ) \}$.
\end{enumerate}

Our next lemma demonstrates that this definition captures the duets which will be of interest to us in this section.

\begin{lemma}
\label{weak_graph_lem}
If $\Delta$ is a near true graph then it is a weak true graph.  If $\Delta$ is a connected duet with a doubling 
block, then either $\Delta$ is a near polygon, or one of $\Delta$ or its reverse is a weak true graph.
\end{lemma}

\noindent{\it Proof:} 
If $\Delta = (X,Y)$ is a near true graph relative to $( {\mathcal V}, {\mathcal E} )$ then it follows from Observation
\ref{near_obs} that $\Delta$ is also a weak true graph relative to  $( {\mathcal V}, {\mathcal E} )$.  Next suppose that 
$\Delta$ is a connected duet with a doubling block, and let $U$ be a doubling block with $w(U)$ minimum.  
We may assume by possibly replacing $(X,Y)$ with its reverse that $U \subseteq X$ and we let ${\mathcal V}$ denote 
the corresponding system of imprimitivity.  Let $\Delta' = ( {\mathcal V}, Y )$ and note that by our assumptions 
$w(\Delta') = w(N_{\Delta'}(U)) = w(N_{\Delta}(U)) = 2 w(U)$.  Therefore
every $y \in Y$ has $w(N_{\Delta'}(y)) = 2 w(U)$ so $y$ must be incident with exactly two blocks of ${\mathcal V}$.  Let 
${\mathcal E}$ be the partition of $Y$ into clones from the duet $\Delta'$.  Then $( {\mathcal V}, {\mathcal E} )$ is 
a simple connected graph, so it must either be a polygon (in which case we are finished) or a true graph.  For the last part, 
suppose that $( {\mathcal V}, {\mathcal E} )$ is a true graph with $d( {\mathcal V}, {\mathcal E} ) = 3$ and let 
$k = \frac{1}{2} w^{\circ}( \mce ) = \frac{1}{3} w^{\circ}( \mcv )$.  Now let $E \in {\mathcal E}$ and observe that 
\[ 3k < w(\Delta) \le w(N_{\Delta}(E)) \le w(N_{\Delta'}(E)) = 6k. \]
It follows from Proposition \ref{block_mult} that $w(N_{\Delta}(E))$ must be a multiple of $2k$.  If $w(N_{\Delta}(E)) = 4k$ then 
$E$ is a doubling block which contradicts our choice of $V$.  Thus we must have $w(N_{\Delta}(E)) = 6k$ and this 
implies that property 3 holds, thus completing the proof.
\quad\quad$\Box$

\begin{lemma}
\label{pure_star}
Let $\Delta = (X,Y)$ be a weak true graph relative to $( {\mathcal V}, {\mathcal E} )$, let 
$d = d( \mcv, \mce )$, and let $(A,B)$ be a maximal critical cross of $\Delta$ with at least two 
non-outer edges.  If every non-outer edge is incident with the vertex $V$, then one of the following holds.
\begin{enumerate}
\item		$(A,B)$ is pure, ${\mathcal E}$ is critical, and $d = 3$.
\item		$[B] \neq Y$, every $y \in B$ satisfies $w(N(y) \cap V) > \frac{1}{2}w(N(y))$, 
		and there are no outer vertices in ${\mathcal V} \setminus \{V\}$.
\end{enumerate}
\end{lemma}

\noindent{\it Proof:} Let $k = \frac{1}{2}w^{\circ}( \mce ) = \frac{1}{d} w^{\circ} ( \mcv )$.  
Suppose (for a contradiction) that the lemma is false for $\Delta$ and let 
$(A,B)$ be a counterexample  with $w(B)$ maximum.  Let 
$E_1,E_2,\ldots,E_m$ be the non-outer edges and assume that $E_i$ has ends 
$V$ and $V_i$ for $1 \le i \le m$.  Set $A_i = A \cap V_i$ and $B_i = B \cap E_i$ 
for $1 \le i \le m$.  We shall assume (without loss) that $w(A_1) \le w(A_i)$ for every $1 \le i \le m$.  
So in particular, if some $V_i$ is outer, then $V_1$ is outer.  We proceed with a sequence of 
claims.  

\bigskip

\noindent{(i)} $w(V_1 \setminus A_1) > \frac{1}{2}dk$.

\smallskip

It follows from the assumption that $\mcv$ is critical that $w(\Delta) > dk$.  So, it suffices to prove that 
$w(V_1 \setminus A_1) \ge \frac{1}{2} w(\Delta)$.  Suppose (for a contradiction) that this fails.  Then for every 
$y \in B_i$ we have 
$w( N(y) \cap V_i ) \le w(V_i \setminus A_i) \le w(V_1 \setminus A_1) < \frac{1}{2}w(\Delta) = \frac{1}{2}w(N(y))$.  
But then setting $T = \{ y \in Y \mid w(N(y) \cap V) > \frac{1}{2}w(N(y)) \}$ we find that $T$ is a block of imprimitivity containing 
$B$, so $[B] \neq Y$ and outcome 2 is satisfied.  

\bigskip

\noindent{(ii)} $w(B_1 \cup B_i) > w(V_1 \setminus A_1)$ for every $2 \le i \le m$.  

\smallskip

Suppose that (ii) is violated and modify $(A,B)$ to form $(A',B')$ by purifying at every boundary 
$V_h$ where $h \in \{2, \ldots, m\} \setminus \{i\}$.  Set $B'_i = B' \cap E_i$ for $1 \le i \le m$.  
By this modification, every $h \in \{2, \ldots, m\} \setminus \{i\}$ satisfies either $B_h' = \emptyset$ or 
$V_h \subseteq N(B_h')$ so in either case $w(N(B_h') \cap V_h) \ge w(B_h')$.  By assumption 
we have $w(N(B_1') \cap V_1) = w(V_1 \setminus A_1) \ge w(B_1' \cup B_i')$, and obviously $w(N(B_i)) \ge w(\Delta)$.  
Combining these gives a contradiction to $\delta(B') > 0$.

\bigskip

\noindent{(iii)} Every $V_i$ is boundary.  Also, if $d=3$ then ${\mathcal E}$ is not critical.

\smallskip

Suppose (for a contradiction) that $V_1$ is outer (so $A_1 = \emptyset$).  Then by (ii) we have
$4k \ge w(B_1 \cup B_2) > w(V_1 \setminus A_1) = dk$ which implies $d = 3$.  
If $m = 2$ we must have $N(B_2) \neq N(E_2)$ (as otherwise $(A,B)$ is pure), but then  
$w(B_1 \cup B_2) > 3k$ implies $w(B_2) > k$, so $[B_2] = E_2$. Then Theorem \ref{critical_closure} implies
\[ w(N(B)) \ge 3k + w(N(B_2)) \ge 2\frac{1}{3}k + w(B_2) + w(\Delta) \ge w(B) + w(\Delta) \]
which is a contradiction.  If $m = 3$ and $V_2$ is also outer then $w(N(B)) \ge 6k + w(N(B_3)) 
\ge w(B) + w(\Delta)$ which is a contradiction.  Finally, if $m=3$ and $V_2,V_3$ are boundary, then we 
may purify at $V_3$ to obtain a critical cross where $V_1$ is outer and $V_2$ is boundary and this leads 
to a contradiction as before.  It follows that every $V_i$ is boundary.  For the second part, note that if 
$d=3$ and ${\mathcal E}$ is critical, then purifying at $E_1$ gives a cross $(A',B')$ which contradicts the 
choice of $(A,B)$.

\bigskip

\noindent{(iv)} Either $d = 3$ and there exists a critical system of imprimitivity ${\mathcal E}'$ 
which is a refinement of ${\mathcal E}$ with $w^{\circ}( {\mathcal E}') = k$, or $4 \le d \le 7$ 
and setting ${\mathcal E}' = {\mathcal E}$ we have that ${\mathcal E}'$ is critical.  

\smallskip

It follows from (i) and (ii)  that $4k \ge w(B_1 \cup B_2) > \frac{d}{2}k$ 
so we must have $d < 8$.  If $d \ge 4$ then we may choose $j \in \{1,2\}$ so that $w(B_j) > k$ 
and if $d = 3$ we may choose $j \in \{1,2\}$ so that $w(B_j) > \frac{3}{4}k$.  Now, (iii) implies
that every $V_i$ is boundary, so by purifying at every $V_i$ with 
$i \in \{1,2,\ldots,m \} \setminus \{j\}$ we find that $B_j$ is 
critical.  Then by applying Theorem \ref{critical_closure} we deduce that $[B_j]$ is critical.  
So, for $4 \le d \le 7$ we have ${\mathcal E}' = {\mathcal E}$ is critical and for $d=3$ 
(applying (iii)) proves the existence of a critical system of imprimitivity ${\mathcal E}'$ as claimed.

\bigskip

\noindent{(v)} $B_i \in {\mathcal E}'$ for every $1 \le i \le m$.

\smallskip

Suppose (for a contradiction) that the above property fails for $B_i$.  If $4 \le d \le 7$ then 
$E_i$ is a critical block of imprimitivity and purifying $(A,B)$ at $E_i$ results in a cross which
contradicts our choice of $(A,B)$.  If $d = 3$ then property 3 in the definition of weak true graph
together with (iii) imply that $B_i \neq E_i$.  If $B_i \not\in \mce'$ then there is a block $T \in \mce'$ 
which is boundary with respect to $(A,B)$ and purifying $(A,B)$ at $T$ gives a cross which contradicts
the choice of $(A,B)$. 

\bigskip

With this last claim in place, we are ready to complete the proof.  Now, $( \mcv, \mce' )$ is a (not necessarily simple) 
connected graph and $d' = d( \mcv, \mce' )$ satisfies $4 \le d' \le 7$.  Let $k' = \frac{1}{2} w^{\circ}( \mce' )$, and note 
that $B$ is pure with respect to $\mce'$ (by (v)) and note that $d'k' = w^{\circ}( \mcv ) = dk$.  Now let $\{ w_1, w_2 \}$ be the profile 
of $(X,Y)$ associated with $( {\mathcal V}, {\mathcal E}')$ and assume (without loss) that 
$w_1 = w(N(B_1)) \cap V_1 = w(V_1 \setminus A_1)$.  It follows from (i) that $w_1 > \frac{1}{2}d'k' \ge 2k'$ and from (ii) that 
$4k' = w(B_1 \cup B_2) > w_1$.  It now follows from Lemma \ref{technical} that $w_1 = w_2 = 3k'$ and 
from $w_1 > \frac{1}{2}d'k'$ that $d' \in \{4,5\}$.  Now applying Lemma \ref{technical} again gives us $w(N(B) \cap V) \ge 4k'$, 
but then we have
\[ w(N(B)) \ge 4k' + 3k' m \ge 6k' + 2k'm = w( X, {\mathcal E}' ) + w(B) \ge w(\Delta) + w(B) \]
which contradicts the assumption that $(A,B)$ is critical.  
\quad\quad$\Box$

\subsection{Triangle}

Our main interest in this subsection is crosses in weak graphs for which the only non-outer edges 
form a triangle.  However, before treating these, we will to establish a couple lemmas.

\begin{lemma}
\label{num2}
Let $\Delta = (X,Y)$ be a weak true graph relative to $( {\mathcal V}, {\mathcal E} )$ and assume that
$d( {\mathcal V}, {\mathcal E} ) = 3$.  If there exists a critical block $S$ contained in a vertex $V$
with $w(S) = \frac{1}{2}w(V)$ and $N(S) \neq N(V)$, then $w(N(S) \cap E) = \frac{1}{2}w(V)$ for 
every $E \in {\mathcal E}$ incident with $V$.
\end{lemma}

\noindent{\it Proof:} Let $k = \frac{1}{2}{w}^{\circ}(\mce) = \frac{1}{3} w^{\circ}(\mcv)$.  Proposition \ref{block_mult} implies that $w(N(S))$ 
is a multiple of $\frac{3}{2}k$ and $3k < w(\Delta) \le w(N(S)) < w(N(V)) = 6k$, so it must be that $w(N(S)) = \frac{9}{2}k$.  
Let $\mcv'$ be the system of imprimitivity associated with $S$ and observe that $w( \mcv', Y ) = \frac{9}{2}k$ so every $y \in Y$ is 
incident with exactly three blocks of $\mcv'$.  Now let $E \in \mce$, choose $y \in E$ and define 
$T = \{ y' \in E \mid N_{ (\mcv', Y) }(y') = N_{ (\mcv', Y) }(y) \}$.  Since $E$ is incident with exactly four blocks of $\mcv'$ and $T$ is a block 
of imprimitivity contained in $E$ we must have that $1 \le \frac{w(E)}{w(T)} \le 4$ is an integer.  
However, by construction $w( N_{ (\mcv', Y) }(T)) = \frac{9}{2}k$ so by Proposition \ref{block_mult} $w(T)$ must divide $\frac{9}{2}k$ 
and it follows that $w(T) = \frac{1}{2}k$.  Thus, the system of imprimitivity $\mce'$ associated with $T$ has $w^{\circ}( \mce' ) = \frac{1}{2}k$
and the result follows immediately from this.  
\quad\quad$\Box$

\begin{lemma}
\label{num4}
Let $\Delta = (X,Y)$ be a  weak true graph relative to $( {\mathcal V}, {\mathcal E} )$ and let 
$(A,B)$ be a maximal critical cross.  If there are exactly two non-outer vertices, then $(A,B)$ is pure.
\end{lemma}

\noindent{\it Proof:} Let $d = d(\mcv, \mce)$ and set $k = \frac{1}{2}w^{\circ}(\mce) = \frac{1}{d} w^{\circ}(\mcv)$.  
Suppose (for a contradiction) that the lemma is false for $\Delta$, and let $(A,B)$ be a counterexample for which
$w(A)$ is as large as possible.  Let 
$V_1,V_2$ be the vertices which have nonempty intersection with $A$ and set 
$A_i = A \cap V_i$ for $i=1,2$.  If $V_1,V_2$ are not incident with a common 
edge, then Theorem \ref{connected_bound} gives us 
$w(N_{\Delta}(A)) = \sum_{i=1}^2 w(N(A_i)) \ge \sum_{i=1}^2 ( w(A_i) + \tfrac{1}{2}w(\Delta) ) \ge w(A) + w(\Delta)$
which is a contradiction.  Thus $V_1,V_2$ are incident with a common edge $E$.  
If there exist non-outer edges $E_1 \neq E$ incident with $V_1$ and 
$E_2 \neq E$ incident with $V_2$, then purifying at $V_1$ results in a cross 
which contradicts our choice of $(A,B)$.  Thus, we may assume (without loss) that every edge 
incident with $V_1$ except possibly $E$ is outer.  It follows from this that 
\begin{equation}
\label{onesolidv}
w(N(A)) \ge w(N(A_2)) + 2(d-1)k.
\end{equation}
Note that $V_2$ must not be inner since $(A,B)$ is maximal and impure.  
If $w(N(A_2)) \ge w(A_2) + w(\Delta) - \frac{1}{3}w(V)$ then combining this with equation \ref{onesolidv} we find
$w(N(A)) \ge 2 (d-1) k + w(A_2) + w(\Delta) - \frac{1}{3} d k \ge w(A) + w(\Delta)$ giving us a contradiction.  Thus it follows from Corollary 
\ref{block_subset_bound2} and the assumption that $(A,B)$ is maximal that $A_2$ is a critical block with 
$w(A_2) = \frac{d}{2}k$.  If $d \ge 4$ then equation \ref{onesolidv} gives 
$w(N(A)) \ge 2(d-1) k + w(\Delta) \ge \frac{3}{2} d k + w(\Delta) \ge w(A) + w(\Delta)$ and we have another contradiction.  
So, we may assume $d = 3$ and now by Lemma \ref{num2} we have that $w(N(A_2) \cap E') = \frac{3}{2}k$ for every $E' \in \mce$ with $E' \sim V$.  
Now $V_1$ must be inner (otherwise purifying at $V_1$ gives a better cross) so
$w(N(A)) = 9k = w(A) + w( N(A_2) ) \ge w(A) + w(\Delta)$ which gives us a final contradiction.
\quad\quad$\Box$

\bigskip

\begin{figure}[ht]
\centerline{\includegraphics[height=5cm]{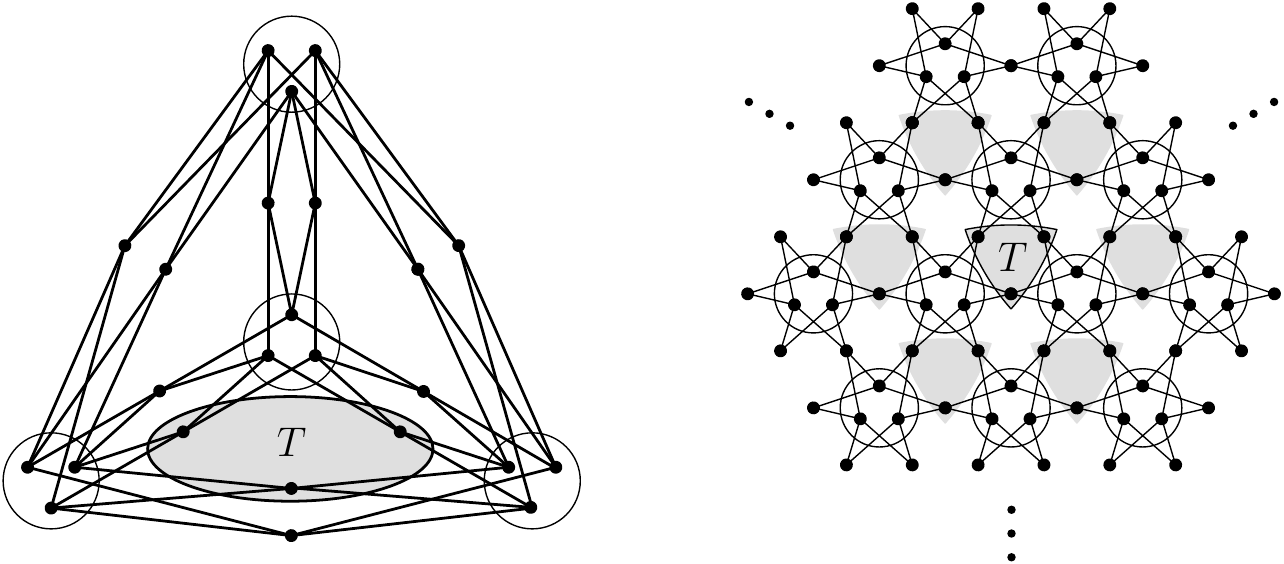}}
\caption{Critical blocks of imprimitivity in weak true graphs}
\label{triangle_block}
\end {figure}

Consider a duet $(X,Y)$ which is a weak true graph relative to $( {\mathcal V}, {\mathcal E} )$ and let $(A,B)$ be a critical cross of $(X,Y)$ for which the non-outer edges form a triangle.  It is certainly possible that $(A,B)$ is pure with respect to $( {\mathcal V}, {\mathcal E} )$ as we have already seen  examples of critical crosses in graphs whose edges form a triangle.  However, in addition to this, there are some examples of impure critical crosses for which the non-outer edges form a triangle, two of these are depicted in Figure \ref{triangle_block}.  Both of these duets are weak true graphs, where the vertex partition ${\mathcal V}$ is indicated by the circles.  In both of these cases each $V \in {\mathcal V}$ has the property that 
$(V, N(V))^{\bullet}$ is a triangle, and $T$ is a critical block of imprimitivity for which $N(T)$ is impure with respect to 
${\mathcal V}$.  Fortunately, these are in a sense the only essentially new types of critical crosses of this form, as is proved in our main lemma from this subsection.

\begin{lemma}
\label{triangle_pure}
Let $\Delta = (X,Y)$ be a weak true graph relative to $( {\mathcal V}, {\mathcal E})$, let $d = d( {\mathcal V}, {\mathcal E} )$ and $k = \frac{1}{2} w^{\circ}( {\mathcal E} )$.  If $(A,B)$ is a maximal critical impure cross so that the non-outer edges form a triangle, then $[B] \neq Y$, there are no 
	 outer vertices, and one of the following holds.
\begin{enumerate}
\item		$d=6$ and ${\mathcal E}$ is critical and 
		$(V, \{ E \in {\mathcal E} \mid E \sim V \} )^{\bullet}$ is a triangle for every $V \in {\mathcal V}$.
\item		$d=3$, there is a critical system of imprimitivity ${\mathcal E}'$ which
		is a refinement of ${\mathcal E}$ so that 
		$w^{\circ}( {\mathcal E}') = k$, and $(V, \{ E \in {\mathcal E}' \mid E \sim V \} )^{\bullet}$ is a triangle 
		for every $V \in {\mathcal V}$.
\end{enumerate}
\end{lemma}

We suppose (for a contradiction) that the lemma is false for $\Delta$ and let $(A,B)$ be a counterexample with
$w(B)$ as large as possible.  Assume that the non-outer edges form a triangle 
with vertices $V_1,V_2,V_3$ and edges $E_1,E_2,E_3$ and assume that $E_i$ is 
not incident with $V_i$ for $1 \le i \le 3$.  Also, for $1 \le i \le 3$, 
set $B_i = B \cap E_i$ and $A_i = A \cap V_i$.  We proceed with a sequence of claims.

\bigskip

\noindent{(i)} $w(V_i \setminus A_i) < w(B)$ for $1 \le i \le 3$.

\smallskip

To see this, choose $u \in B_i$ and observe that 
\[ w(B) + w(\Delta) > w(N(B)) \ge w(V_i \setminus A_i) + w(N(u)) 
	= w(V_i \setminus A_i) + w(\Delta). \]

\smallskip

\noindent{(ii)} At most one $V_i$ with $1 \le i \le 3$ is outer.

\smallskip

Suppose (for a contradiction) that $A_2 = A_3 = \emptyset$.  Then since $(A,B)$ is impure and maximal
we must have $E_1 \subseteq B$ and $V_1$ boundary.  Now by weakly purifying at $V_1$ 
we deduce that $(X \setminus (V_2 \cup V_3), E_1)$ is critical which implies that
$w(\Delta) + w(E_1) > w(V_2 \cup V_3)$ so $w(\Delta) > (2d-2)k$.  Now choose $x \in A_1$ 
and note that $B_2 \cup B_3$ is disjoint from $N(x)$ so $w(B_2 \cup B_3) < 2k$.  It follows that 
$w(B) < 4k$ and then by (i) and the assumption $A_2 = \emptyset$ we have $dk < 4k$, so 
$d = 3$.  However, then $({\mathcal V}, {\mathcal E}) \cong {\mathit Tetrahedron}$ (by Proposition \ref{equitable_prop}) 
and there are exactly two non-outer vertices in $(A,B)$ and this contradicts Lemma \ref{num4}.

\bigskip

\noindent{(iii)} If ${\mathcal E}$ is critical, then $B = E_1 \cup E_2 \cup E_3$.

\smallskip

Suppose (for a contradiction) that (iii) fails.  
If $E_i$ is boundary, then consider the cross $(A',B')$ obtained by 
purifying $(A,B)$ at $E_i$.  It follows from our choice of counterexample that 
$(A',B')$ is pure, and this implies that $V_i$ is outer with respect to $(A,B)$.  So, by 
(ii) we may assume that $V_1$ is outer, $V_2, V_3$ are boundary, $E_2, E_3$ are 
inner and $E_1$ is boundary.  If $d = 3$ then by the definition of weak true graph 
$V_3 \subseteq N(E_2)$ and this is contradictory, so $d \ge 4$.  By purifying at $E_1$ 
we have that  $(X \setminus (V_1 \cup V_2 \cup V_3), E_1 \cup E_2 \cup E_3)$
is critical.  It follows that $w(\Delta) + 6k > 3dk$.  Now choose a point $x \in A_2$ and 
note that $B_3$ must be disjoint from $N(x)$ so $2k = w(B_3) \le 2dk - w(\Delta) < (6-d)k$ which 
is a contradiction.  

\bigskip

\noindent{(iv)} If ${\mathcal E}$ is not critical then $V_i$ is boundary for $1 \le i \le 3$.

\smallskip

By (ii) we may suppose (for a contradiction) that $V_1$ is outer and $V_2,V_3$ are boundary.  
For $i=2,3$ we may purify at $V_i$ to obtain a critical cross, and it follows that $B_i$ is critical.  
Then Theorem \ref{critical_closure} implies that $[B_i]$ is critical for $i =2,3$.  Since ${\mathcal E}$ 
is not critical, it must be that $w(B_i) \le k$ for $i = 2,3$.  If $B_1$ is not critical then 
$w(N(B)) \ge w(V_1) + w(N(B_1)) \ge w(B_2 \cup B_3) + w(\Delta) + w(B_1)$ which contradicts 
the criticality of $B$.  Thus $B_1$ is critical, and by Theorem \ref{critical_closure} we have that 
$[B_1]$ is critical, and as before, this implies $w(B_1) \le k$.  However, then $w(B) \le 3k$ and
we have a  contradiction to (i).  

\bigskip

\noindent{(v)} Either ${\mathcal E}$ is not critical and there exists a critical system of imprimitivity
${\mathcal E'}$ refining ${\mathcal E}$ with either $w^{\circ}({\mathcal E}') = \frac{2}{3}k$ and 
$d = 3$ or with $w^{\circ}( {\mathcal E}') = k$ and $3 \le d \le 5$ or ${\mathcal E'} = {\mathcal E}$ is 
critical and $4 \le d \le 11$.

\smallskip

Since $w(\Delta) > dk$, by considering a point $y \in B_1$ we have that $w(N(y)) > dk$ and this implies
$\max\{ w(V_2 \setminus A_2), w(V_3 \setminus A_3) \} > \frac{1}{2}{dk}$ and it then follows from (i) that
$w(B) > \frac{1}{2}dk$.  Since $w(B) \le 3 w^{\circ}( {\mathcal E} ) = 6k$ it follows that $d \le 11$.  If $\mce$ is 
critical, then (iii) (and the definition of weak true graph) imply that $d \ge 4$ and there is nothing left to prove.
So, we may assume  ${\mathcal E}$ is not critical.  By the above we may choose $j \in \{1,2,3\}$ 
so that $w(B_j) > \frac{1}{6}dk$.  Now by (iv) $V_j$ is boundary, so by purifying at $V_j$ we deduce that $B_j$ is 
critical.  Then by Theorem \ref{critical_closure} the block $T = [B_j]$ is also critical.  Since $w(T) > \frac{1}{6}dk$ and 
$w(T)$ is a proper divisor of $2k$, it must be that either $3 \le d \le 5$ and $w(T) = k$ or $d = 3$ and $w(T) = \frac{2}{3}k$ and then
setting ${\mathcal E}'$ to be the system of imprimitivity induced by $T$ we have the desired conclusion.

\bigskip

\noindent{(vi)} $B_i \in {\mathcal E}'$ for $1 \le i \le 3$.  

\smallskip

If ${\mathcal E}' = {\mathcal E}$ then (iii) implies the result.  Otherwise ${\mathcal E}$ is not critical, so by 
(iv) every $V_i$ is boundary.  If 
there is a block of ${\mathcal E}'$ which is boundary, then purifying $(A,B)$ at this block results in a cross 
which contradicts our choice of $(A,B)$.  So each $B_i$ is a union of blocks of ${\mathcal E}'$.  If $B_i$ 
consists of more than one such block, then by purifying at $V_i$ we deduce that $B_i$ is critical, and then 
Theorem \ref{critical_closure} implies that $[B_i]$ is critical, but $[B_i] = E_i$ and this is a contradiction.

\bigskip

With this last claim in place we are ready to complete the proof.  Consider the graph $( {\mathcal V}, {\mathcal E}')$, 
let $d' = d( {\mathcal V}, {\mathcal E}')$ and set $k' = \frac{1}{2} w^{\circ}({\mathcal E}')$, and note that by our assumptions
we must have $4 \le d' \le 11$.  Let $\{ w_1, w_2 \}$ be the profile of $(X,Y)$ associated with $( {\mathcal V}, {\mathcal E}' )$, 
and note that 
\begin{equation}
\label{w1w2lb}
w_1 + w_2 = w(X, {\mathcal E}') \ge w(\Delta) > dk = d'k'.
\end{equation}
Next choose $1 \le i \le 3$ so that $E_i$ is not incident with an outer vertex and note that $N(B_i)$ must intersect 
one end of $E_i$, say $V_j$, in a set of weight $\max\{w_1, w_2\}$.  Since $V_j \not \subseteq N(B_i)$ and (by (i))
$w( N(B_i) \cap V_j ) < w(B) = 6k'$ we have the following inequality.
\begin{equation}
\label{w1w2ub}
\max\{ w_1, w_2 \} < \min\{ 6k', d'k' \}. 
\end{equation}
If $w_1 \neq w_2$, Lemma \ref{technical} implies that $d'$ is even, and $w_1,w_2$ are multiples of $2k'$.  However, then equation \ref{w1w2ub} 
implies that $\{w_1, w_2\} = \{2k', 4k'\}$ and $d' \ge 6$, but this contradicts equation \ref{w1w2lb}.  So we must have $w_1 = w_2$, and then 
Lemma \ref{technical} implies that $w_1$ is a multiple of $k'$, so equations \ref{w1w2lb} and \ref{w1w2ub} force $w_1 \in \{3k', 4k', 5k' \}$.  
If $w_1 = 5k'$ then we must have $d' \ge 6$, but then Lemma \ref{technical} implies that $w(V_1 \setminus A_1) \ge 6k' \ge w(B)$ which contradicts (i).  
Similarly, if $w_1 = 3k'$ then we must have $d' \in \{4,5\}$, and then Lemma \ref{technical} implies that $w(V_1 \setminus A_1) = w(N( B_2  \cup B_3) \cap V_1 ) \ge 4k'$
and similarly for the other two vertices.  However, then we have the contradiction $w(N(B)) \ge 12k' \ge 6k' + 6k' = w(B) + w( X, {\mathcal E}') \ge w(B) + w(\Delta)$. In the remaining case we have $w_1 = 4k'$ and note that this implies $d' \ge 5$.  If $w(V_i \setminus A_i) \ge 5k'$ for every $1 \le i \le 3$ then $w(N(B)) \ge 15k' > 6k' + 8k' = w(B) + w( X, {\mathcal E}' ) \ge w(B) + w(\Delta)$ which is contradictory.  Therefore, by Lemma \ref{technical} we deduce that $d' = 6$ and $( V, \{ E \in {\mathcal E}' \mid E \sim V \} )^{\bullet}$ is a triangle for every $V \in {\mathcal V}$.  It follows that whenever $E, E' \in {\mathcal E}'$ are incident with $V \in {\mathcal V}$ 
either $V \subseteq N(E \cup E')$ or $N(E)  \cap V = N(E') \cap V$.  Let us say that $E$ and $E'$ are \emph{clones with respect to} $V$ if the latter condition is satisfied.  Note that for every vertex $V$ there are 6 incident edges in ${\mathcal E}'$ and these are grouped into 3 doubletons of clones.  If $B_1$ and $B_2$ are not clones with respect to $V_3$, then $V_3$ is outer and we have $w(N(B)) \ge w(V_3) + 2 \cdot 4k' = 14k' \ge w(B) + w(\Delta)$ which is a contradiction.  So, we find that $B_1$ and $B_2$ are clones with respect to $V_3$, and similarly $B_2,B_3$ and $B_1,B_3$ are clones with respect to their common endpoint.  However, in this case $B = B_1 \cup B_2 \cup B_3$ is a block of imprimitivity so conclusion 1 or 2 from the lemma holds.  
\quad\quad$\Box$

\subsection{Stability}

In this section we will prove our stability lemma for near graphs and duets with doubling blocks.  This will follow easily once
we have established the following.

\begin{lemma}
\label{2and2v2}
Let $\Delta = (X,Y)$ be a  weak true graph relative to $( {\mathcal V}, {\mathcal E} )$
and let $(A,B)$ be a maximal critical cross.  If there are at least two non-outer vertices
and at least two nonadjacent non-outer edges, then $(A,B)$ is pure.
\end{lemma}

\noindent{\it Proof:} Let $d = d( \mcv, \mce )$, let $k = \frac{1}{2}w^{\circ}(\mce) = \frac{1}{d} w^{\circ}(\mcv)$, 
and note that $dk < w(\Delta)$ since $\mcv$ is critical.  We may assume $w(\Delta) < 2dk$ as otherwise 
$\mcv$ is a block of clones, and the result is immediate.  Now suppose (for a contradiction) that 
the lemma is false for $\Delta$ and choose a counterexample $(A,B)$ subject to the following conditions.
\begin{enumerate}
\item[$\alpha$)] $\delta(A,B)$ is maximum.
\item[$\beta$)] The total number of boundary blocks in $\mcv$ is minimum (subject to $\alpha$).
\end{enumerate}
As was the case in the proof of Lemma \ref{fringed_cross}, modifying $(A,B)$ by purifying at a boundary vertex $V$ 
results in a cross which is superior to $(A,B)$ with regard to $\alpha$, $\beta$.  We now proceed with a sequence of claims.

\bigskip

\noindent{(i)} There does not exist a boundary vertex $V$ so 
	that the cross obtained by purifying $(A,B)$ at $V$ is impure 
	and has three non-outer edges which form a triangle.

\smallskip

Suppose (for a contradiction) that $V$ is a vertex which violates (i), let $E_1, E_2, E_3 \in \mce$ be  the edges in the resulting 
triangle, and note that the purification of $(A,B)$ at $V$ must satisfy outcome 1 or 2 in Lemma \ref{triangle_pure}.  

First assume that outcome 1 is satisfied, and note that this implies $w(\Delta) \le w(X, \mce) = 8k$.  In 
this case $\mce$ is critical, and we now define $(A',B')$ to be the critical cross obtained from $(A,B)$ by weakly purifying 
at every boundary edge.  If $F_1, \ldots, F_{\ell}$ are the edges incident with $V$ which are not outer 
with respect to $(A,B)$ then $B' = E_1 \cup E_2 \cup E_3 \cup F_1 \ldots, F_{\ell}$ so $w(B') = 2k(3+ \ell)$.  
Furthermore, the structure of $(X,\mce)$ implies that for every $V' \in \mcv$ with $N(B') \cap V' \neq \emptyset$
we have $w( N(B') \cap V') \ge 4k$ and further, $w(N(B') \cap V') = 6k$ if $B'$ contains at least three edges 
incident with $V'$.  However, this contradicts $\delta(B') > 0$.

Next suppose that outcome 2 in Lemma \ref{triangle_pure} holds and let $\mce'$ be the corresponding 
critical system of imprimitivity on $Y$.  Now $d=3$, so Proposition \ref{equitable_prop} implies that the graph 
$( \mcv, \mce)$ is isomorphic to Tetrahedron.  It then follows from Lemma \ref{triangle_pure} and our assumptions
that there are no outer vertices with respect to $(A,B)$.  Choose a vertex $V' \in \mcv \setminus \{V\}$ so that the edge 
incident with $V$ and $V'$ is not outer, and note that $V'$ must then be incident with at least three non-outer blocks of 
$\mce'$.  Now consider modifying $(A,B)$ by repeatedly purifying at blocks of $\mce'$.  It follows from our 
assumptions that at some point in this process the vertex $V'$ will become outer.  However, if we stop this process as 
soon as the first outer vertex is encountered, the resulting cross contradicts the choice of $(A,B)$.

\bigskip

\noindent{(ii)} There exists a critical cross $(A',B')$ which is pure 
with respect to $( {\mathcal V}, {\mathcal E} )$ and has at least two inner 
vertices and at least two outer vertices.

\smallskip

If $(A,B)$ has at least two outer vertices, then let $(A',B')$ be obtained by repeatedly
weakly purifying $(A,B)$ at every boundary vertex.  It follows immediately that
$(A',B')$ has no boundary vertices, and from the maximality of $(A,B)$ that $(A',B')$ it 
has no boundary edges, so it satisfies (ii).  So, we may now assume there is at most 
one outer vertex with respect to $(A,B)$.  Then, by our assumptions, we may choose 
non-outer non-adjacent edges $E,E' \in \mce$ so that if an outer vertex exists, it is an 
endpoint of $E$.  If there exists a boundary vertex not incident with $E$ or $E'$ then 
purifying at this vertex gives a cross which contradicts $(A,B)$.  So, letting $V_1,V_2$ 
be the ends of $E$ and $V_1',V_2'$ be the ends of $E'$ we have that every non-outer edge
has both ends in $\{V_1, V_2, V_1',V_2'\}$.  

If $E,E'$ are the only non-outer edges, then by Corollary \ref{connected_set_bound} we 
have $w(N(B)) = w(N(B \cap E)) + w(N(B \cap E')) \ge w(B) + w(\Delta)$ which is 
contradictory.  Thus, we may assume that there is a non-outer edge $F$ with
ends $V_1,V_1'$.  Let $(A',B')$ be the cross obtained from $(A,B)$ 
by purifying at $V_2'$.  If $(A',B')$ is pure, then (ii) is satisfied, so we may assume it is not.  
It then follows from (i) that the non-outer edges in $(A',B')$ are precisely $E$ and $F$.  
Now by applying Lemma \ref{pure_star} we find that $V_2$ is boundary, and every 
$y \in F$ satisfies $w( N(y) \cap V_1' ) < w( N(y) \cap V_1 )$.  In this case we may purify the 
original cross $(A,B)$ at the vertex $V_2$ to obtain a new cross $(A'',B'')$.   It follows from 
a similar argument that either $(A'',B'')$ is pure (so we are done) or every $y \in F$ 
satisfies $w( N(y) \cap V_1 ) < w( N(y) \cap V_1')$.  However, this last possibility is contradictory, 
and this completes the proof of (ii).

\bigskip

\noindent{(iii)} $\Delta$ is a near true graph relative to $( {\mathcal V}, {\mathcal E} )$.

\smallskip

By (ii) we may choose a pure critical cross $(A',B')$ with at least two inner and 
two outer vertices.  If we let $A^*$ be the set of inner vertices with respect to $(A',B')$ and 
$B^*$ be the set of inner edges with respect to $(A',B')$, then 
$(A^*,B^*)$ is a critical cross in $( {\mathcal V}, {\mathcal E} )$.  Now using 
Corollary \ref{2dm2_cor} we have
\[ w( {\mathcal V}, {\mathcal E} ) - w( X, Y)	= \delta(A^*,B^*) - \delta(A',B')  < \delta(A^*) \le w^{\circ}( {\mathcal E} ) \]
and this completes the proof.

\bigskip

\noindent{(iv)} The graph $( {\mathcal V}, {\mathcal E} )$ is not isomorphic to
	${\mathit Tetrahedron}$.  

\smallskip

Suppose (iv) fails and choose two nonadjacent non-outer edges $E,E'$.  It 
follows from Lemma \ref{num4} that there is at most one outer vertex, so
$E$ and $E'$ must both be boundary, and we may assume that if an outer vertex
exists it is incident with $E$.  However, then purifying at $E$ yields a maximal 
impure cross with exactly two non-outer vertices, thus contradicting Lemma \ref{num4}.

\bigskip

\noindent{(v)} Every boundary edge is incident with every boundary vertex.

\smallskip

Suppose (for a contradiction) that $V$ is a boundary vertex and $E$ a boundary edge and 
$V \not\sim E$.  Now consider the cross $(A',B')$ obtained from $(A,B)$ by purifying at $E$.  
By construction $(A',B')$ is not pure, and by Observation \ref{near_obs} it has fewer boundary 
vertices than $(A,B)$, so it must be that $V$ is the only non-outer vertex of $(A',B')$ (otherwise 
$(A',B')$ contradicts the choice of $(A,B)$).  It follows from (iv) that $| \mcv | \ge 5$ and thus 
there must exist an edge $E'$ not incident with $V$ and not adjacent with $E$.  But then 
the cross obtained from $(A,B)$ by purifying at $V$ is a counterexample which contradicts 
our choice of $(A,B)$.  

\bigskip

\noindent{(vi)} There is exactly one boundary vertex $V$ (so every boundary 
	edge is incident with $V$).

\smallskip

If (vi) fails, then by (v) there must exist a single boundary edge $E$ so that
the boundary vertices are precisely the two ends $V_1,V_2$ of $E$.  Now every 
edge incident with $V_1$ or $V_2$ other than $E$ must be outer (were $E'$ to be an inner edge 
incident with $V_i$ for $i=1,2$ Observation \ref{near_obs} would imply that $V_i$ is outer).
But then $N( B \cap E ) \cap N( B \setminus E) = \emptyset$ and by Corollary \ref{connected_set_bound}
we find that $w(N(B)) = w(N(B \cap E)) + w(N(B \setminus E)) \ge w(B) + w(\Delta)$
which contradicts the assumption that $B$ is critical.

\bigskip

\noindent{(vii) There are at least three inner vertices.}

\smallskip

It follows from (vi) and Lemma \ref{num4} that there are at least two inner vertices.  Suppose
(for a contradiction) that there are exactly two inner vertices.  Let $(A',B')$ be the cross 
obtained by weakly purifying at the unique boundary vertex $V$ and let $(A^*,B^*)$ 
denote the cross of $( {\mathcal V}, {\mathcal E} )$ obtained by taking $A^*$ to be the set of inner vertices and 
$B^*$ to be the inner edges with respect to $(A',B')$.  Now, as before, we have that $(A^*,B^*)$ 
is critical, so by (iv) and Theorem \ref{cutset_classify} it must be that either 
$( {\mathcal V}, {\mathcal E} )$ has degree 3 and the three non-outer vertices form a path, or 
$( {\mathcal V}, {\mathcal E} )$ has degree 4 or 5 and the three non-outer vertices induce a triangle.  
A quick check reveals that the following inequality holds in all of these cases.  
\[ w( N(A) \setminus N(V) ) \ge w( A \setminus V ) + \tfrac{1}{2}dk \]
However, we must have $N( A \cap V ) \neq N(V)$ so Corollary \ref{block_subset_bound2} implies
$w(A \cap V) \ge w(\Delta) + w(A \cap V) - \frac{1}{2}dk$ and combining this with the above 
inequality yields a contradiction.

\bigskip

With this final claim in place, we can complete the proof.  There are at least three inner vertices by (vii), and 
since there are two nonadjacent non-outer edges, (vi) implies that we must also have at least three outer vertices. 
Consider the cross $(A_1,B_1)$ obtained by purifying at the 
unique boundary vertex and a cross $(A_2,B_2)$ obtained by purifying at any boundary edge, 
and note that by (iii) and Observation \ref{near_obs}, both of these crosses are pure.
Now let $A_1^*$ denote the set of inner vertices with respect to $(A_1,B_1)$ and $A_2^*$ denote the 
set of inner vertices with respect to $A_2$.  As in the previous cases, we have that 
$A_1^*$ and $A_2^*$ are critical in the duet $( {\mathcal V}, {\mathcal E} )$.  
Furthermore, by construction $A_2^* = A_1^* \setminus \{ V \}$.  By applying Theorem
\ref{cutset_classify} we deduce that 
$( {\mathcal V}, {\mathcal E} ) \cong \mathrm{Cube}$ and the set of 
non-outer vertices is the vertex set of a cycle of length four.  Since $A$ is critical the set $A' = A \cap V$ satisfies
\[ w(N(A')) + 10k \le w(N(A)) < w(A) + w(\Delta) = 9k + w(A') + w(\Delta). \]
Since $N(A') \neq N(V)$ it then follows from Corollary \ref{block_subset_bound2} and the maximality of $(A,B)$ 
that $A'$ is a critical block of imprimitivity with $w(A') = \frac{3}{2}k$.  However, then $w(A) = \frac{21}{2}k$ and 
Lemma \ref{num2} implies that $w(N(A)) = \frac{31}{2}k$ and $w(\Delta) \le w(N(A')) = \frac{9}{2}k$ and 
this gives us a final contradiction, thus completing the proof.
\quad\quad$\Box$

\begin{lemma}[Near Graph and Doubling Block Stability]
\label{near_graph_stability}
Let $\Theta$ be a maximal connected critical trio, and let $\Delta$ be a side of $\Theta$.  
\begin{enumerate}
\item If $\Delta$ is a near true graph, then $\Theta$ is a video.  
\item If $\Delta$ has a critical doubling block, then either $\Theta$ is a video or 
$\Delta$ is a near polygon.
\end{enumerate}
\end{lemma}

\noindent{\it Proof:} Let $\Theta = (X,Y,Z)$ and suppose that $\Delta = (X,Y)$.  In light of Lemma \ref{weak_graph_lem}, 
to complete the proof it suffices to show that $\Theta$ is a video under the assumption that $\Delta$ is a weak true graph.  
Choose systems of imprimitivity ${\mathcal V}$, ${\mathcal E}$ so that $\Delta$ is a weak true graph relative to
$( {\mathcal V}, {\mathcal E} )$.  Let $(A,B)$ be the cross of $\Delta$ associated with some point in $Z$ and 
note that by the connectivity of $\Theta$ this cross must have at least two non-outer vertices and at least 
two non-outer edges.  
If all non-outer edges are pairwise incident, then it follows from Lemma \ref{pure_star} or \ref{triangle_pure} (and the 
connectivity of $\Theta$) that $(A,B)$ is pure.  Otherwise, Lemma \ref{2and2v2} shows that $(A,B)$ is pure.  It now 
follows from the maximality of $\Theta$ that ${\mathcal V}$ and ${\mathcal E}$ are systems of clones and that 
$\Delta^{\bullet}$ is a true graph.  Then Lemma \ref{graph_stability} implies that $\Theta$ is a video, as desired.
\quad\quad$\Box$

\section{Clips and Near Clips}
\label{clip_sec}

In this section we establish stability lemmas for clips and near clips.

\subsection{Crosses of Clips}

In this subsection we will classify critical crosses of clips.

\begin{lemma}
\label{funnyp3}
Let $\Gamma = (V,E)$ be an arc-transitive cubic graph and let $\{A,\overline{A}\}$ be a 
partition of $V$ so that $|A| \ge 2$ and $|\overline{A}| \ge 3$ and at least one of $A,\overline{A}$ is finite.  
Then one of the following holds.
\begin{enumerate}
\item There are at least $9$ two edge paths with one end in $A$ 
	and with middle vertex in $\overline{A}$.
\item $A$ consists of two adjacent vertices.
\item $\overline{A}$ is the vertex set of a two edge path.
\item $\Gamma \cong \mathrm{Cube}$ and $A,\overline{A}$ are vertex sets of cycles of length 4.
\end{enumerate}
\end{lemma}

\noindent{\it Proof:} Suppose (for a contradiction) that the lemma is false and choose a counterexample $\{A,\overline{A}\}$ for which $c(A)$ is minimum.  The 
lemma is straightforward to verify when $|A| = 2$ or $|\overline{A}| = 3$ so we may assume $|A| \ge 3$ and $|\overline{A}| \ge 4$.  We define a two edge path to be \emph{interesting} 
if it has one end in $A$ and middle vertex in $\overline{A}$.  Now, for $0 \le i \le 3$ let $A_i = \{ x \in A : |\partial(x) \cap \partial(A)| = i \}$ and let
$\overline{A}_i = \{ y \in \overline{A} : |\partial(y) \cap \partial(A)| = i \}$.  

First suppose that there exists a vertex $v \in A_3 \cup \overline{A}_3$, and modify $\{A, \overline{A}\}$ by moving this vertex
to the other part of the partition.  This operation does not create any new interesting paths, and reduces $c(A)$, thus contradicting the choice of $\{A, \overline{A}\}$.    
Next suppose there exists a vertex $v \in A_2 \cup \overline{A}_2$ and modify $\{A, \overline{A} \}$ by moving this vertex to the 
other side of the partition.  This operation causes at least two paths to change from interesting to non-interesting, and at most two paths to change from 
non-interesting to interesting (here we use $A_3 = \overline{A}_3 = \emptyset$), so again the resulting partition contradicts our choice.  Thus, we may 
now assume $A_3 = \overline{A}_3 = A_2 = \overline{A}_2 = \emptyset$.  If $c(X) \ge 5$ then there are at least 10 interesting paths, so the first outcome holds.  
Otherwise, it follows from Lemma \ref{cutset_classify} that outcome 4 is satisfied.
\quad\quad$\Box$

\begin{lemma}
\label{p3_cross}
Let $\Gamma$ be a graphic duet with degree $3$ and let $(A,B)$ be a critical cross of $\Gamma(V,P_3)$ 
with $|A|, |B| \ge 2$.  Then $(A,B)$ is maximal and one of the following holds
\begin{enumerate}
\item $A$ consists of two adjacent vertices.
\item $\Gamma$ is isomorphic to Cube, and $A$ consists of the four vertices of a face, 
	so $[B] \neq P_3(\Gamma)$.
\end{enumerate}
\end{lemma}

\noindent{\it Proof:} 
Let $V = V(\Gamma)$ and $P = P_3(\Gamma)$, so $(V,P) = \Gamma(V,P_3)$, and 
let $k = w^{\circ}(P) = \frac{1}{3}w(V)$.  

Suppose first that $A$ is finite, let $D \subseteq P$ be the set of all two edge paths with
middle vertex in $A$, and let $D' = N(A) \setminus D$.  Now we have
\[ w(D) + k |D'| = w(D) + w(D') = w(N(A)) < w(A) + w(\Delta) = w(D) + 9k. \]
Thus, Lemma \ref{funnyp3} implies that either $A$ consists of two adjacent vertices, or $\Gamma$ is 
the cube and $A$ is the vertex set of a cycle of length four.  In either case the only way  
for $(A,B)$ to be critical is for it to be maximal, and this completes the proof.

Next suppose that $A$ is infinite, let $\overline{A} = V \setminus A$, and let $B'$ be the set of
all two edge paths with middle vertex in $\overline{A}$ but at least one endpoint in $A$.  Now 
$B \cup B'$ is contained in the set of two edge-paths with middle vertex in $\overline{A}$ so we have
\[ k |B'| + w(B) =  w(B') + w(B) \le w( \overline{A} ) < w(B) + w(\Delta) = w(B) + 9k, \]
but this contradicts the assumption $A$ is infinite and Lemma \ref{funnyp3}.  
\quad\quad$\Box$

\bigskip

For reference purposes, let record some parameters of the exceptional clips.  In the chart below, 
we have defined $k = \frac{w^{\circ}(X)}{d(X,Y)} = \frac{w^{\circ}(Y)}{d(Y,X)}$.
\begin{table}[ht]
\begin{center}
\begin{tabular}{|c|c|c|c|}
\hline
$(X,Y)$	&	$w^{\circ}(X)$	&	$w^{\circ}(Y)$	&	$\overline{w}(X,Y)$	\\
\hline
$K_5(E, C_3)$	&	$3k$		&	$3k$		&	$21k$			\\
\hline
$K_6(E, C_3)$	&	$4k$		&	$3k$		&	$48k$			\\
\hline
$\mathrm{Petersen}(E,C_5)$	&	$4k$		&	$5k$		&	$40k$		\\
\hline
$\mathrm{Dodecahedron}(V,F)$	&	$3k$		&	$5k$		&	$45k$		\\
\hline
\end{tabular}
\caption{Weights of some elements in exceptional videos}
\label{weight_video}
\end{center}
\end{table}

\begin{lemma}
\label{ec-cross}
Let $(A,B)$ be a critical cross of one of $K_5(E,C_3)$, Petersen$(E,C_5)$, or 
$K_6(E,C_3)$.  If $|A|, |B| \ge 2$ then $(A,B)$ is maximal and $A = \partial(v)$ for some vertex $v$.
\end{lemma}

\noindent{\it Proof:} 
Suppose (for a contradiction) that the lemma fails, let $V$, $E$ be the vertex and edge sets
of the graph, and let $(A,B)$ be a counterexample for which $|B|$ is minimum.  
The lemma is straightforward to verify whenever  $|B| = 2$,  
and whenever $A \subseteq \partial(v)$ for some vertex $v$ so we may further assume that $|B| \ge 3$ and 
that  $A$ contains two non-adjacent edges.  We now split into cases depending upon the input graph.

\bigskip

\begin{figure}[ht]
\centerline{\includegraphics[width=12cm]{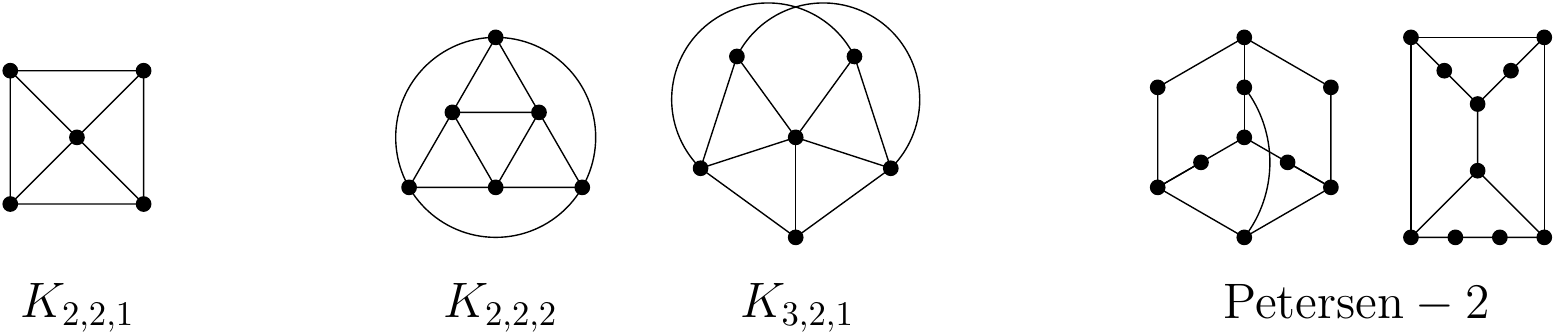}}
\caption{Subgraphs of $K_5$, $K_6$, and Petersen}
\label{clip-cross}
\end {figure}

\noindent{\sl Case 1:} $K_5(E,C_3)$

\smallskip

In this case a cross $(A',B')$ is critical if and only if $|A'| + |B'| > 7$.  If there is an edge $e \in E \setminus A$ which is 
incident with just one triangle $q \in B$, then $(A \cup \{e\}, B \setminus \{q\})$ is a superior cross.  So, we may assume no such 
edge $e$ exists.  Since $A$ contains two non-adjacent edges, the graph $(V, E \setminus A)$ must then be a subgraph of 
$K_{2,2,1}$ (as depicted on the left in Figure \ref{clip-cross}) which contains at least two triangles through each edge.  However no 
such graph exists.

\bigskip

\noindent{\sl Case 2:} $\Delta = K_6(E,C_3)$

\smallskip

In this case a cross $(A',B')$ is critical if and only if $4|A'| + 3|B'| > 48$.  It follows from the minimality of $|B|$ that there does not 
exist an edge $e \in E \setminus A$ which is incident with just one triangle in $B$.  The result 
is straightforward to verify when $(V, E \setminus A)$ contains a $K_4$ subgraph (so $A$ is contained in $\partial(u) \cup \partial(v)$ for some $u,v \in V$)
so we may assume this is not the case.  It follows from this that $A$ must either contain a three edge matching, or two vertex disjoint subgraphs 
one consisting of a triangle, the other consisting of an edge.  Therefore, $(V, E \setminus A)$ must be a subgraph of one of the two 
graphs in the middle of Figure \ref{clip-cross} which contains at least two triangles through every edge.  There is only one such graph, Octahedron, and it does 
not give a critical cross.

\bigskip

\noindent{\it Case 3:} Petersen$(E,C_5)$

\smallskip

In this case a cross $(A',B')$ is critical if and only if $4|A'| + 5|B'| > 40$.  It follows from the minimality of $|B|$ that there must 
not exist a pair of edges $e,e' \in E \setminus A$ incident with a pentagon $q \in B$ with the property that no other pentagon in $B$ 
is incident with either $e$ or $e'$.  Since $A$ contains two non-adjacent edges, the graph $(V, E \setminus A)$ must be a subgraph of 
one of the graphs on the right in Figure \ref{clip-cross} containing no pair of edges $e,e'$ as described.  There is only one such graph 
(obtained from the rightmost graph in the figure by deleting the two vertices of degree two on the bottom), and it is 
straightforward to check that does not yield a critical cross.
\quad\quad$\Box$

\begin{lemma}
\label{dodec_cross}
If $(A,B)$ is a critical cross of  Dodecahedron$(V,F)$ and $|A|, |B| \ge 2$, then $(A,B)$ is maximal and 
either $A$ consists of two adjacent vertices, or $B$ consists of two adjacent faces.
\end{lemma}

\noindent{\it Proof:} Let $V,F$ be the vertex and face sets of Dodecahedron, and note that a cross $(A',B')$ is critical if and only if 
$3|A'| + 5|B'| > 45$.  We suppose (for a contradiction) that the lemma is false, and choose a counterexample $(A,B)$ with 
$\delta(A,B)$ maximum.  The lemma is 
straightforward to verify when $|A| \le 4$ and when $|B| \le 2$ so we may further assume $|A| \ge 5$ and $|B| \ge 3$.  
Now consider a face $f \in B$ and let $f_1, \ldots, f_5$ be the faces adjacent to $f$ so that $f_i$ and $f_{i+1}$ are 
adjacent for every $1 \le i \le 5$ (working modulo 5).  If there are two vertices incident with $f$ which are not incident 
with any other member of $B$, then adding these to $A$ and removing $f$ from $B$ results in a superior cross.  
Therefore, we may assume (without loss) that $f_1, f_3 \in B$.  Now, there is a unique vertex $v$ incident with $f_2$ but 
not incident with either $f_1$ or $f_3$.  If $f_2 \not\in B$ then we may obtain a superior cross by adding $f_2$ to $B$ 
and removing $v$ from $A$.  Thus $f_1,f_2,f_3 \in B$.  If $f_4,f_5 \not\in B$ then we may obtain a superior cross by 
adding $f_4,f_5$ to $B$ and removing any vertex of $A$ incident with one of these.  Thus, $f_1, \ldots, f_5 \in B$, but 
repeating this argument results in an apparent contradiction.
\quad\quad$\Box$

\subsection{Clip Stability}

In the previous subsection we determined all critical crosses of clips.  To complete our stability lemma for clips, 
we will need to understand the images of these critical crosses under the group action.  To help appreciate a subtlety here, 
consider a $G$-duet given by $K_6(E,C_3)$.  By assumption, $G$ acts transitively on $E(K_6)$ and on $C_3(K_6)$, 
and this action preserves incidence.  However, it is not obvious, a priori, that $G$ is acting on the graph $K_6$ (i.e. acting 
on the vertex and edge sets preserving incidence).  Our first task will be to show that for each of our clips, the group $G$ is 
indeed acting on the appropriate graph.  Once we have this action in hand, it will be straightforward to determine the image 
of the crosses of interest under the action of $G$.

\begin{lemma}[Clip Stability]
\label{clip_stability}
Let $\Theta$ be a connected critical trio.  If $\Theta$ has a side which 
is a clip, then $\Theta$ is a video.
\end{lemma}

\noindent{\it Proof:} Assume $\Theta = (Z,X,Y)$ is a $G$-trio, and let $\Delta=(X,Y)$ be a side of $\Theta$ which is a clip.  
Let $z \in Z$ and let $(A,B)$ be the cross of $(X,Y)$ associated with $z$.  Now we split into cases depending on the possibilities
for $\Delta$.

\bigskip

\noindent{\it Case 1:} $\Delta \cong \Gamma(V,P_3)$ for a graphic duet $\Gamma$ of degree $3$.

\smallskip

Define two vertices of $\Gamma$ to be \emph{adjacent} if they are both incident with four members of $\Gamma(P_3)$.  These 
adjacencies correspond to $E(\Gamma)$, and it follows that $G$ is acting on the graph $\Gamma$.  It then follows 
from the assumption that $G$ acts transitively on $P_3(\Gamma)$ that $G$ is also transitive on $E(\Gamma)$.  Now 
by the connectivity of $\Theta$ and Lemma \ref{p3_cross} we find that $A$ consists of two adjacent vertices and 
$B = Y \setminus N(A)$.  Since our action is edge-transitive, we deduce that $\Theta^{\bullet} \cong \Gamma: E \sim V \sim P_3$.  

\bigskip

\noindent{\it Case 2:} $\Delta$ is isomorphic to one of $K_5(E,C_3)$, $K_6(E,C_3)$, or Petersen$(E,C_5)$.  

\smallskip

For $K_5$ and $K_6$ we define two edges to be \emph{adjacent} if there is a triangle incident with both edges, 
and for Petersen we define two edges to be \emph{adjacent} if there are four pentagons incident with 
both edges.  Maximum size sets of pairwise adjacent edges correspond to vertices, so $G$ is indeed acting 
on our graph.  Furthermore, it follows from the assumption that the action of $G$ is edge-transitive and the fact 
that our graph is not bipartite, that this action is vertex transitive.  It follows from the connectivity of $\Theta$ and Lemma 
\ref{ec-cross} that $A = \partial(v)$ for some vertex, and we must have $B = Y \setminus N(A)$ in order to be critical.  
It now follows that $\Theta^{\bullet}$ is isomorphic to one of $K_5 : V \sim E \sim C_3$ or 
$K_6: V \sim E \sim C_3$ or $\mathrm{Petersen}: V \sim E \sim C_5$ as desired.

\bigskip

\noindent{\it Case 3:} $\Delta \cong \mathrm{Dodecahedron}(V,F)$.

\smallskip

Define two vertices to be \emph{adjacent} if there are two faces incident with both vertices.  These adjacencies correspond
precisely to $E( \mathrm{Dodecahedron} )$ so $G$ is acting on $\mathrm{Dodecahedron}$, as desired.  Furthermore, $G$ acts 
transitively on the edge set (and also their duals).  Lemma \ref{dodec_cross} implies that $(A,B)$ is maximal and either $A$ 
consists of two adjacent vertices or $B$ consists of two adjacent faces.  It follows that $\Theta^{\bullet}$ must be isomorphic to a 
trio equivalent to $\mathrm{Dodecahedron}: E \sim F \sim V$ or Dodecahedron$: F \sim V \sim E$.
\quad\quad$\Box$

\subsection{Near Clip Stability}

In this subsection we prove a stability lemma for near clips.

\begin{lemma}
\label{p3_purity}
Let $(X,Y)$ be a duet which is near $\Gamma(V,P_3)$ for a graphic duet $\Gamma$ of degree $3$.
If $(A,B)$ is a maximal critical cross of $(X,Y)$ with at least two non-outer vertices and at least two non-outer paths, 
then $(A,B)$ is pure.
\end{lemma} 

\noindent{\it Proof:}
Choose block partitions $\mcv$ of $X$ and $\mcp$ of $Y$ so that $(X,Y)$ is near $\Gamma(V,P_3)$
relative to $(\mcv, \mcp)$ and let $(A,B)$ be a counterexample to the lemma chosen so that
\begin{enumerate}
\item[$\alpha$)] $\delta(A,B)$ is maximum.
\item[$\beta$)] The total number of boundary blocks in $\mcv$ and $\mcp$ is minimum (subject to $\alpha$).
\end{enumerate}
As in the previous instances, purifying $(A,B)$ at any boundary block results in a cross which will 
be superior to $(A,B)$ with regard to the above criteria.  Let $k = w^{\circ}(\mcp) = \frac{1}{3}w^{\circ}(\mcv)$ and note
that $8k < w(X,Y) \le 9k$ since $\mcp$ is critical.  We proceed with a sequence of numbered claims.

\bigskip

\noindent{(i)} There do not exist $V_1,V_2 \in \mcv$ with $A \subseteq V_1 \cup V_2$.

\smallskip

Suppose (for a contradiction) that (i) is violated and let 
$A_i = A \cap V_i$ for $i=1,2$.  If there is a boundary path not incident with $V_1$ and a boundary
path not incident with $V_2$ then by purifying at $V_1$ or $V_2$ 
yields a cross which contradicts our choice of $(A,B)$.  Thus we
may assume without loss that every boundary path is incident with
$V_1$.  Now by Corollary \ref{block_subset_bound2} we have
\[ w(N(A)) \ge 5k + w(N(A_2)) \ge \tfrac{7}{2}k + w(\Delta) + w(A_2) \ge w(A) + w(\Delta) \]
which contradicts the assumption that $A$ is critical.  

\bigskip

\noindent{(ii)} We have $w(B) < 4k$.

\smallskip

Suppose $w(B) \ge 4k$ and let $(A',B')$ be the cross obtained from $(A,B)$ by purifying at 
a boundary vertex $V$.  It follows from Theorem \ref{purify_thm} that $w(B') > k$ so $B'$ is not contained 
in a single member of $\mcp$.  It then follows from our choice of $(A,B)$ that $(A',B')$ is pure, 
and then from (i) that $\Gamma \cong \mathrm{Cube}$ and $B'$ corresponds to the set of paths contained 
in a face.  Now let $(A'',B'')$ be the weak purification of $(A,B)$ at the vertex $V$.  Since 
$k \ge \delta(A',B') \ge \delta(A'',B'') > 0$ it must be that $12k = w(A') \le w(A'') + k \le w(A) + 4k$.  
However, now we may modify $(A,B)$ by purifying at a boundary path to obtain a cross $(A''',B''')$ 
for which $w(A''') \ge w(A) -k \ge 7k$.  It follows from our choice of $(A,B)$ that $(A''',B''')$ is pure, but $A'''$ 
is not contained in two vertices, and this gives us a contradiction.  

\bigskip

\noindent{(iii)} We have $w(A) \ge 4k$

\smallskip

Otherwise, $w(A) + w(B) + w(\Delta) \le 17k \le 3k | \mcv | = w(X)$ which contradicts the assumption $(A,B)$ is critical.

\bigskip

Now consider the cross $(A',B')$ obtained from $(A,B)$ by purifying at a boundary path $P$.  It follows from 
(iii) and our choice of $(A,B)$ that $(A',B')$ is pure.  If $A'$ consists of two adjacent vertices, then 
by Theorem \ref{purify_thm} it must be that $w(A) \le 7k$ but then $3k |\mcv| = w(X) < w(A) + w(B) + w(\Delta) \le 20k$
so $|\mcv| \le 6$ which implies that $\Gamma \cong K_{3,3}$ and $B'$ corresponds to the set of four 
paths contained in a cycle of length four.  If $A'$ does not consist of two adjacent vertices, then 
by Lemma \ref{p3_cross} we find that $\Gamma \cong \mathrm{Cube}$ and again $B'$ corresponds to the 
set of four paths contained in a cycle of length four.  So, $B'$ has the same structure in both cases, and we have $w(B') = 4k$.  
Next choose $x \in A \setminus A'$ and choose 
$V \in \mcv$ with $x \in V$.  Now $w(N(V) \cap B') = 3k$ so it follows from the assumption that $w(\mcp,\mcq) - w(X,Y) < w^{\circ}(\mcp) = k$
that $w( N(x) \cap B' ) > 2k$.  But then we must have $w(B' \setminus B) \ge 2k$.  However, Theorem \ref{purify_thm} implies
$w(A \setminus A') < k$ so we have
\[ k < w(B' \setminus B) - w(A \setminus A') = \delta(A',B') - \delta(A,B) \le k - \delta(A,B) \]
which is a contradiction.  
\quad\quad$\Box$

\bigskip

Next we prove a general stability lemma which will be used to handle duets which are near one of 
 $K_5(E,C_3)$, $K_6(E,C_3)$ or $\mathrm{Petersen}(E,C_5)$.

\begin{lemma} 
\label{genlpq-stab} 
Let $(X,Y)$ be a duet which is near $\Delta$ relative to $(\mcp, \mcq)$, and assume that every critical cross $(A^*,B^*)$ of 
$(\mcp, \mcq)$ with $|A^*|, |B^*| \ge 2$ satisfies $w(A^*) = a$ and $w(B^*) = b$ and $\delta(A^*,B^*) = t$.  
If $\overline{w}(X,Y) \ge 2 w^{\circ}(\mcp) + 2 w^{\circ}(\mcq)$ and $a \ge w^{\circ}(\mcp) + 3t$ and $b \ge w^{\circ}(\mcq) + 3t$,  then every maximal critical cross $(A,B)$ of $(X,Y)$ satisfies one of:
\begin{enumerate}
\item $A$ is contained in a block of $\mcp$.
\item $B$ is contained in a block of $\mcq$.
\item $(A,B)$ is pure.
\end{enumerate}
\end{lemma}

\noindent{\it Proof:} Suppose (for a contradiction) that the lemma is false for $(X,Y)$ and choose a 
counterexample $(A,B)$ so that
\begin{enumerate}
\item[$\alpha$)] $\delta(A,B)$ is maximum.
\item[$\beta$)]  The total number of boundary blocks in $\mcp$ and $\mcq$ is minimum (subject to $\alpha$).
\end{enumerate}

First we claim that there do not exist boundary blocks $P \in \mcp$ and $Q \in \mcq$ so that purifying at $P$ yields a cross $(A',B')$ for which 
$B'$ is contained in a block of $\mcq$ and purifying at $Q$ yields a cross $(A'',B'')$ for which $A''$ is contained in a block of $\mcp$.  To verify
this, suppose otherwise, and note that Theorem \ref{purify_thm} implies $w(A) < w(A'') + w^{\circ}(\mcq) \le w^{\circ}(\mcp) + w^{\circ}(\mcq)$ and 
$w(B) < w(B') + w^{\circ}(\mcp) \le  w^{\circ}(\mcp) + w^{\circ}(\mcq)$.  However, then $0 < \delta(A,B) = w(A) + w(B) - \overline{w}(X,Y) \le 2 w^{\circ}(\mcp) + 2w^{\circ}(\mcq) - \overline{w}(X,Y)$ contradicts our assumptions.

Next suppose there exist boundary blocks $P \in \mcp$ and $Q \in \mcq$ for which $P \not\sim Q$.  Let $(A',B')$ ($(A'',B'')$) be 
obtained from $(A,B)$ by purifying at $P$ ($Q$).  Now $(A',B')$ and $(A'',B'')$ are impure, but are superior to $(A,B)$ with regard to 
criteria $\alpha$, $\beta$.  Therefore  $B'$ is contained in a block of $\mcq$ and $A''$ is contained in a block of $\mcp$, but this 
contradicts our above claim.  Therefore, every boundary block in $\mcp$ is incident with every boundary block in $\mcq$, so 
purifying at any boundary block results in a pure cross.  

By our earlier claim, we may now assume (without loss) that purifying at some block $P \in \mcp$ results in a pure cross 
$(A',B')$ for which $w(A') = a$ and $w(B') = b$.  Since $(A',B')$ is critical, our assumptions (concerning critical crosses $(A^*,B^*)$ of 
$(\mcp, \mcq)$) imply $0 < \delta(A',B') = t - ( w(\mcp,\mcq) - w(X,Y) )$, so we have 
\begin{equation}
\label{pq-close}
w(\mcp,\mcq) < w(X,Y) + t.
\end{equation}
Next define $m$ ($n$) to be the sum of the weights of the boundary blocks of $\mcp$ ($\mcq$).  Choose a point $y \in B$ so that $y$ is in 
a boundary block of $\mcq$ and note that since $y$ is incident with every boundary block of $\mcp$, equation \ref{pq-close} implies that 
$m < t$, and by a similar argument we have $n < t$.  Using $w(B) = w(B') + n$ we then find
\begin{align}
t	&>	\delta(A',B') - \delta(A,B)	\nonumber \\
	&= 	w(A') - w(A) + w(B') - w(B)	\nonumber \\
	&\ge a - w(A) - t. \label{eqn-pq-cl}
\end{align}
Next consider a cross $(A'',B'')$ obtained from $(A,B)$ by purifying at a boundary block of $\mcq$.  It follows from 
our assumptions that $(A'',B'')$ is pure, and $w(A'') < w(A') = a$, so it must be that either $A'' = \emptyset$ or $A'' \in \mcp$.  In either case 
Theorem \ref{purify_thm} implies $w(A) = w(A'') + m \le w^{\circ}(\mcp) + t$.  But then combining this with equation \ref{eqn-pq-cl} gives us  
$w^{\circ}(\mcp) + t \ge w(A) > a - 2t$ which contradicts our hypothesis.
\quad\quad$\Box$

\begin{lemma}
\label{ec-purify}
Let $(X,Y)$ be near one of $K_5(E,C_3)$, $K_6(E, C_3)$, or $\mathrm{Petersen}(E,C_5)$ relative to $(\mce, \mcc)$.  If $(A,B)$ is a maximal critical
cross of $(X,Y)$, then either $A$ is contained in a member of $\mce$ or $B$ is contained in a member of $\mcv$ or $(A,B)$ is pure.
\end{lemma}

\noindent{\it Proof:} Let $k = \frac{w^{\circ}(\mce)}{d( \mce, \mcc)} = \frac{ w^{\circ}(\mcc)}{ d(\mcc,\mce) }$ and consider a critical cross $(A^*,B^*)$ of 
$(\mce,\mcc)$ for which $|A^*|, |B^*| \ge 2$.  It follows from Lemma \ref{ec-cross} that $A^* = \partial(v)$ for a vertex $v$ so we must have one of the 
following cases.
\begin{center}
\begin{tabular}{|c|c|c|c|}
\hline
$(\mce,\mcc)$	&	$w(A^*)$	&	$w(B^*)$	&	$\delta(A^*,B^*)$	\\
\hline
$K_5(E, C_3)$	&	$12k$		&	$12k$		&	$3k$			\\
\hline
$K_6(E, C_3)$	&	$20k$		&	$30k$		&	$2k$			\\
\hline
$\mathrm{Petersen}(E,C_5)$	&	$12k$	&	$30k$		&	$2k$		\\
\hline
\end{tabular}
\end{center}
From here the result follows easily from the observation $\overline{w}(X,Y) \ge \overline{w}(\mce,\mcv)$, Table \ref{weight_video}, and 
Lemma \ref{genlpq-stab}.
\quad\quad$\Box$

\begin{lemma}
\label{dodec_purity}
Let $(X,Y)$ be near $\mathrm{Dodecahedron}(V,F)$ relative to $(\mcv, \mcf)$.  If $(A,B)$ is a maximal critical cross 
of $(X,Y)$, then either $A$ is contained in a member of ${\mathcal V}$, or  $B$ is contained in a member of ${\mathcal F}$, 
or $(A,B)$ is pure.  
\end{lemma}

\noindent{\it Proof:} Suppose (for a contradiction) that the lemma is false for $(X,Y)$, set 
 $k = \frac{w^{\circ}(\mcv)}{d(\mcv,\mcf)} = \frac{w^{\circ}(\mcf)}{d(\mcf,\mcv)}$ (so $w^{\circ}(\mcv) = 3k$ and $w^{\circ}(\mcf) = 5k$) 
 and then choose a counterexample $(A,B)$ so that
\begin{enumerate}
\item[$\alpha$)] $\delta(A,B)$ is maximum.
\item[$\beta$)] The total number of boundary blocks in $\mcv$ and $\mcf$ is minimum (subject to $\alpha$).
\end{enumerate}
Since $(A,B)$ is critical we have $w(A) + w(B) > \overline{w}(X,Y) \ge \overline{w}( \mcv, \mcf ) = 45k$ and it follows that  $\min\{w(A), w(B)\} > 20k$.  

First suppose $w(B) > 20k$ and let $(A',B')$ be the cross obtained by purifying $(A,B)$ 
at a boundary vertex.  It follows from $\alpha$ and $\beta$ that the lemma holds for $(A',B')$.  Theorem \ref{purify_thm} implies 
that $w(B') > 17k$ and since $A'$ cannot be contained in a single vertex it must be that $A'$ is the disjoint union of 
two adjacent vertices $V_1,V_2$.  If there exist boundary faces $F_1, F_2$ so that $V_i \sim F_j$ if and only if $i=j$, then by purifying at $V_1$ gives an impure cross which contradicts our choice of $(A,B)$.  So, we may 
assume (without loss) that every boundary face is incident with $V_1$.  Now by Corollary \ref{block_subset_bound2} we have
$w(N(A)) \ge 5 k + w(N(A \cap V_1)) \ge 5k + w(A \cap V_1) + w(\Delta) - \frac{3}{2}k \ge w(A) + w(\Delta)$ which
is a contradiction.  

Next suppose $w(A) > 20k$ and let $(A'',B'')$ be the cross obtained by purifying $(A,B)$ at
a boundary face.  As before, our lemma holds for $(A'',B'')$, so by Theorem \ref{purify_thm} we have
$w(A'') > 15k$ and we deduce that $B''$ is the disjoint union of two adjacent faces $F_1,F_2$.  If there exist boundary vertices
$V_1,V_2$ so that $F_i \sim V_j$ if and only if $i=j$ then purifying at $F_1$ gives us a contradiction.  So, we may assume (without loss)
that every boundary vertex is incident with $F_1$.  Now by Corollary \ref{block_subset_bound2} we have
$w(N(B)) \ge 9k + w(N(F_1 \cap B)) \ge 9k + w( F_1 \cap B ) + w(\Delta) - \frac{5}{2}k \ge w(B) + w(\Delta)$ which is a 
contradiction.
\quad\quad$\Box$

\begin{lemma}[Near clip stability]
\label{near_clip_stability}
Let $\Theta$ be a maximal connected critical trio.  If a side of $\Theta$ is a near clip, then $\Theta$ is a video.
\end{lemma}

\noindent{\it Proof:}
Assume that $\Theta = (X,Y,Z)$ and that $(X,Y)$ is a near clip relative to $(\mcp, \mcq)$.  
Let $z \in Z$ and let $(A,B)$ be the cross of $(X,Y)$ associated with $z$.  It follows from the connectivity
of $\Theta$ that $A$ is not contained in a block of $\mcp$ and $B$ is not contained in a block of $\mcq$.  
Therefore, by Lemmas \ref{p3_purity}, \ref{ec-purify}, and \ref{dodec_purity} it must be that $(A,B)$ is pure.  
However, it then follows from the maximality of $\Theta$ that 
$\mcp,\mcq$ are systems of clones, so $(X,Y)^{\bullet}$ is a clip.  Now the result follows 
by applying  Lemma \ref{clip_stability} to $\Theta^{\bullet}$.  
\quad\quad$\Box$

\section{Maximality of Songs}

In this section we will prove the easy direction of the main theorem, that every song is a maximal critical trio.
This proof will call upon a handful of lemmas which demonstrate the maximality of various type of trios.

\subsection{Impure Beats}

In this subsection we prove the following lemma concerning maximality of impure beats.

\begin{lemma}
\label{impure_beat_max}
Every impure beat with a maximal continuation in maximal.
\end{lemma}

\noindent{\it Proof:} Let $\Theta = (X,Y,Z; \sim)$ be a $G$-trio.  Assume that $\Theta$ is an impure beat relative to $(X,Y)$ with associated 
systems of imprimitivity $\mcp, \mcq, \mcr$ of $X,Y,Z$, and assume further that $\Theta$ has a maximal continuation.  
Now suppose (for a contradiction) that there exists a trio $\Theta^* = (X,Y,Z; \sim^*)$ for which $\Theta < \Theta^*$.  It follows from the maximal continuation
assumption that whenever $(P,Q,R; \sim)$ is a continuation of $\Theta$ we must have $(P,Q,R; \sim ) = (P, Q, R; \sim ^*)$.  It follows from this (and the structure of an impure beat) that $\sim^*|_{X \cup Z} = \sim|_{X \cup Z}$ and $\sim^*|_{Y \cup Z} = \sim|_{Y \cup Z}$.  Choose $z \in Z$,  let $(A,B)$ be the cross of $(X,Y)$ associated with $z$ (this is the same in either $\Theta$ or $\Theta^*$), and assume (without loss) that $A$ is finite.  Let $P_1 \in \mcp$ and $Q_1 \in \mcq$ be the unique blocks for which  $( \{z\}, P_1 )$  and $( \{z\}, Q_1 )$) are partial.  

Next consider the $G$-chorus $( \mcp, \mcq, \sim^*)$.  It follows from the assumption that $\Theta$ has a 
maximal continuation that $(\mcp, \mcq, \sim^*)$ is not a matching.  Let $\Delta = ( \mcp', \mcq', \sim^* )$ be the component of $(\mcp, \mcq, \sim^*)$ which contains $P_1$ and let $A' = \{ P \in {\mathcal P}' \mid P \subseteq A \}$ and 
$B' = \{ Q \in {\mathcal Q}' \mid Q \subseteq B \}$.  Choose a point $x \in A \cap P_1$ and observe that  
$N_{\Theta^*}(x)$ must contain a point in $Q_2$ for some $Q_2 \in \mcq \setminus \{Q_1\}$.  It follows from this that $Q_2 \in {\mathcal Q}'$ and 
$Q_2 \not\in B'$.  Now choose $P_2 \in \mcp$ to be the unique block for which $P_2 \sim Q_2$ and note that $P_2 \in A'$, so in particular 
$A' \neq \emptyset$.  Since $A$ is finite (by assumption) 
and $P_2 \subseteq A$,  the blocks of ${\mathcal P}$ and ${\mathcal Q}$ must be finite, so $\Delta$ is a  $G$-duet.
A similar argument shows that $B' \neq \emptyset$, so $(A',B')$ is a nontrivial cross in the duet $\Delta$.  Now setting $t = w^{\circ}_G( {\mathcal P} ) = w^{\circ}_G( {\mathcal Q} )$ and applying Theorem \ref{connected_bound} gives us
\[ \tfrac{1}{2} w(\Delta) \ge \delta(\Delta) \ge \delta(A',B') = w(A') + w(\Delta) - \overline{w}(B') = w(\Delta) - t. \]
This implies $2t \ge w(\Delta) = d_{\Delta}( {\mathcal P}', {\mathcal Q}' ) t$ and since $\Delta$ is connected, it follows that it is a polygon.  However, then $N_{\Delta}(P_1) = \{Q_1,Q_2\}$ is disjoint from $B'$ so 
$(A' \cup \{P_1\}, B')$ is a nontrivial cross of $\Delta$ and we have 
\[ \delta(\Delta) \ge \delta(A' \cup \{P_1\}, B') = w(\Delta) \]
which contradicts Theorem \ref{connected_bound}. 
\quad\quad$\Box$

\subsection{Octahedral Choruses}

An octahedral chorus $\Lambda = (X_1,\ldots,X_n; \sim)$ is \emph{maximal} if there does note exist another octahedral chorus 
$Lambda^* = (X_1, \ldots, X_n; \sim^*)$ with the same distinguished partition and $\Lambda < \Lambda^*$.

\begin{lemma}
\label{oct_max}
If $\Lambda = (X_1, \ldots, X_n)$ is an octahedral chorus satisfying one of the following properties, it is maximal.
\begin{enumerate}
\item If $\Lambda$ is type $2B$, or $2C$.
\item If $\Lambda$ is type $1$ or $2A$ with a maximal continuation.
\end{enumerate}
\end{lemma}

\noindent{\it Proof:} Let $x \in X_i$ and $y \in X_j$ and assume that $x \not\sim y$ and that $X_i, X_j$ are in distinct members of the distinguished partition.  It 
suffices to prove that that the following holds:

\noindent{$(\star)$} \; There exists $z$ in the ground set of $\Lambda$ so that $x \sim z \sim y$.  

If $(X_i,X_j)$ is empty, then it follows from a check of our definitions that there must exist $1 \le k \le n$ so that one of $(X_i,X_k)$, $(X_j,X_k)$ is 
full and the other is nonempty, and this implies $(\star)$.  So, we may assume $(X_i,X_j)$ is partial.  If $\Lambda$ is type $1$, then $x$ and $y$ are in a continuation of $\Lambda$, and the result follows from the maximality of this continuation.  If $\Lambda$ is type $2C$ then there exists $1 \le k \le m$ 
so that $(X_i,X_j,X_k)$ is a pure beat, and since pure beats are maximal, this again implies $(\star)$.  In the remaining cases, $\Lambda$ is either type 
$2A$ or $2B$, and we let $\Upsilon_1$ and $\Upsilon_2$ be the $H$-trios as given in the definition of these types.  If for $i=1$ or $i=2$ we have 
that $x,y$ lie in the ground set of $\Upsilon_i$ and $\Upsilon_i$ is an impure beat, then the result follows from Lemma \ref{impure_beat_max}.  Next 
suppose that there is exactly one $i \in \{1,2\}$ for which $x,y$ lie in the ground set of $\Upsilon_i$ and that $\Upsilon_i$ is a pure beat.  Setting 
$\Upsilon_i \equiv (X_i,X_j,X_k)$ we may then assume $w(X_i, X_j) = \overline{w}(X_j,X_k)$.  However, it now follows from the fact that $(X_i,X_j,X_k)$ is nontrivial and Lemma \ref{nontriv_bound} that $(\star)$ must hold (otherwise we could modify the incidence relation on $(X_i,X_j,X_k)$ to obtain a counterexample to this lemma).  In the remaining case, we may assume that $x,y$ are contained in the ground set of both $\Upsilon_1$ and $\Upsilon_2$.  Let $\sim, \sim_1, \sim_2$ denote the incidence relations
of $\Lambda$, $\Upsilon_1$, and $\Upsilon_2$.  Since $x \not\sim y$ it follows that there exists $1 \le i \le 2$ so that $x \not\sim_i y$.  However, it then follows from the fact that $\Upsilon_i$ 
is maximal that $(\star)$ holds, and this completes the proof.
\quad\quad$\Box$

\subsection{Dihedral Chords}

In this subsection we prove a maximality lemma for dihedral chords.

\begin{figure}[ht]
\centerline{\includegraphics[height=3cm]{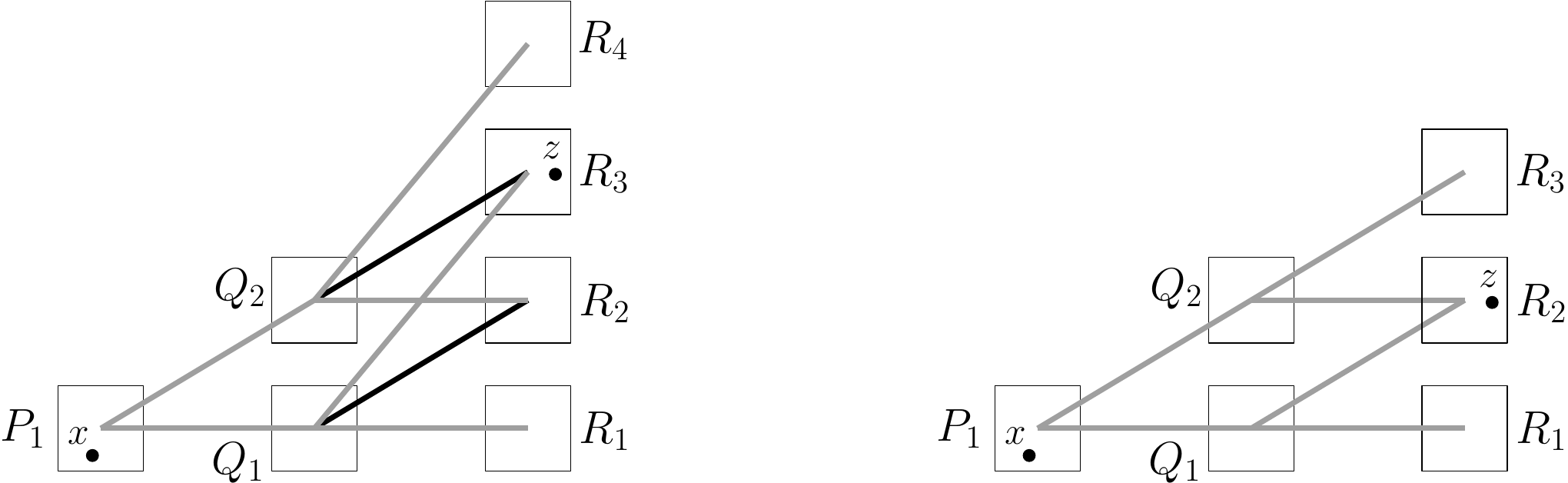}}
\caption{Parts of dihedral chords}
\label{dih-rob}
\end {figure}

\begin{lemma}
\label{chord_max}
Every trio of one of the following types is maximal.
\begin{enumerate}
\item A dihedral chord of type $0$, $2B$, or $2C$.
\item A dihedral chord of type $1$ or $2A$ with a maximal continuation.
\end{enumerate}
\end{lemma}

\noindent{\it Proof:} Let $\Theta = (X,Y,Z)$ be a dihedral chord relative to the systems of imprimitivity $\mcp, \mcq, \mcr$ on $X, Y, Z$ which satisfies one 
one of the input assumptions.  Let $H = G_{ (\mcp) } = G_{ (\mcq) } = G_{ (\mcr) }$.  Let $x \in X$ and $z \in Z$ satisfy $x \not\sim z$ and note that it suffices to prove: 

\noindent{$(\star)$} \; There exists $y$ in the ground set of $\Lambda$ so that $x \sim y \sim z$.  

If there is an associated octahedral chorus containing $x$ and $z$, then the result follows from Lemma \ref{oct_max}.  Otherwise, let $\Lambda$ be an associated octahedral chorus, and let $\Lambda_{XY}$ ($\Lambda_{YZ}$) be the square of $\Lambda$ with ground set contained in $X \cup Y$ ($Y \cup Z$).  
Choose $P \in \mcp$ with $x \in P$ and $R \in \mcr$ with $z \in R$.  First suppose that $\Lambda_{XY}$ is not small and that 
either $w(Y,Z) = \infty$ or that $(Y,Z)$ is a fringed sequence which is not short.  In this case (depicted on the left in Figure \ref{dih-rob}) we may choose $Q \in \mcq$ so that $(P,Q)$ is nonempty and $(Q,R)$ is full.  If $(P,Q)$ is 
partial, then it follows from the assumption that $\Lambda_{XY}$ is not small that we may choose $y \in Q$ so that $x \sim y$, and if $(P,Q)$ is full then we choose $y \in Q$ arbitrarily.  In either case we see that $(\star)$ holds.   
A similar argument handles the case when $\Lambda_{YZ}$ is not small and either $w(X,Y) = \infty$ or $(X,Y)$ is a fringed sequence which is not short.  
In light of Observation \ref{big-or-small} this resolves all of the cases except when both $(X,Y)$ and $(Y,Z)$ are fringed sequences which are short (and thus both $\Lambda_{XY}$ and $\Lambda_{YZ}$ are not small).  In this last case 
(depicted on the right in Figure \ref{dih-rob}) there exist $Q_1, Q_2 \in \mcq$ so that $N(P) = Q_1 \cup Q_2 = N(R)$.  Now it follows from Observation \ref{big-or-small} that $\Lambda_{XY}$ and $\Lambda_{YZ}$ are big and $w(N(x) \cap (Q_1 \cup Q_2)) > 2|H|$ and $w(N(z) \cap (Q_1 \cup Q_2)) \ge 2 |H|$.  Since 
$w(Q_1 \cup Q_2) = 4|H|$ it follows that there exists $y \in Q_1 \cup Q_2$ with $x \sim y \sim z$ as desired.
\quad\quad$\Box$

\subsection{The Proof}

In this subsection we prove a maximality lemma for videos, and then the easy direction of the main theorem.

\begin{lemma}
\label{video_max}
Every video is maximal.
\end{lemma}

\noindent{\it Proof:} For the purpose it it sufficient to consider only clone-free videos.  Let
$(X,Y,Z)$ be a clone-free video and assume (without loss) that $w(X,Y) \le w(Y,Z) \le w(Z,X)$.  Note that by Lemma 
\ref{largest_video_duets} one of $(X,Y)$ or $(Y,X)$ (and similarly one of $(Y,Z)$ or $(Z,Y)$) must either be a graph
or a clip.  To check that $(X,Y,Z)$ is maximal, it suffices to prove that for a point $z \in Z$ the cross of $(X,Y)$ associated with $z$ is 
maximal, and for a point $x \in X$ the cross of $(Y,Z)$ associated with $x$ is maximal.  This is easy to verify 
for graphs, and the other cases follow from Lemmas \ref{p3_cross}, \ref{ec-cross}, and \ref{dodec_cross}.
\quad\quad$\Box$

\bigskip

\noindent{\it Proof of the ``if'' direction of Theorem \ref{main_incidence}:} 
Let $\Theta$ be a song.  It follows from Observation \ref{2d_obs} and a routine check of our definitions that 
$\Theta$ is critical.  The fact that $\Theta$ is maximal follows from a straightforward induction and the following claim.  

\smallskip

\noindent{\it Claim: } A trio is maximal if it is one of the following. 

\smallskip

\begin{center}
\begin{tabular}{rc}
\begin{tabular}{l|l}
\quad\quad\quad\quad Structure		&		Proof	\\
\hline
Video							&		Lemma \ref{video_max}	\\
Pure beat							&		by definition	\\
Pure chord						&		straightforward	\\
Impure beat with maximal continuation	&		Lemma \ref{impure_beat_max}	\\
Cyclic chord with maximal continuation	&		straightforward	\\
Dihedral chord of type $2B$ or $2C$	&		Lemma \ref{chord_max}	\\
Dihedral chord of type $1$,$2A$ with maximal continuation	&	Lemma \ref{chord_max}
\end{tabular} 
&
\raisebox{-55pt}{\quad\quad$\Box$}
\end{tabular}
\end{center}

\section{The Structure of Critical Subset Trios}
\label{prod_again}

The goal of this section will be to prove a slight strengthening of Theorem \ref{main_subset}, which will give a necessary and sufficient 
condition for a nontrivial subset trio to be maximal and critical.

\subsection{Symmetries of Videos}

Our first step will be to record some properties of automorphism groups of videos.  Fo  standard videos which come from 
cubic graphs, the following fundamental theorem of Tutte will provide us with an additional bound.

\begin{theorem}[Tutte \cite{tutte}]
\label{tutte}
Let $\Gamma = (V,E)$ be a finite cubic arc-transitive graph and let 
$G = Aut(\Gamma)$.  Then $w_G^{\circ}(V)$ divides 48 and $w_G^{\circ}(E)$ divides 32.
\end{theorem}

\begin{figure}[ht]
\centerline{\includegraphics[height=3cm]{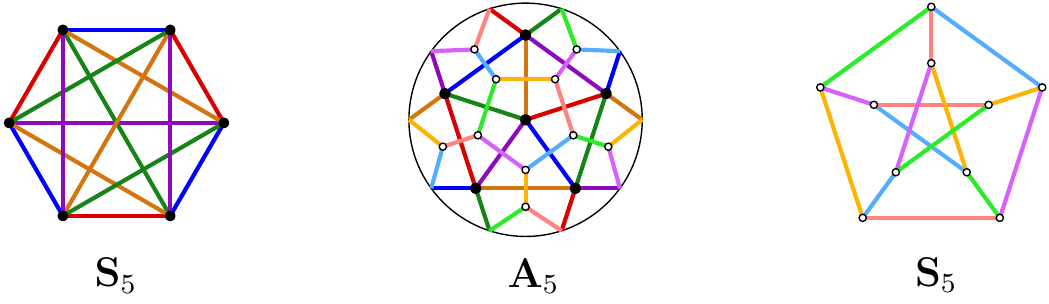}}
\caption{Automorphism groups of edge coloured graphs and maps}
\label{autgroup}
\end {figure}

In addition, we will call on some standard (and easily verified) facts concerning the automorphism groups of some 
small graphs and coloured graphs.  The most interesting of these involves the interesting action of $S_5$ on $6$ points as given 
by an outer-automorphism of ${\mathbf S}_6$.  In Figure \ref{autgroup} we have depicted $K_6$ and Petersen 
with a distinguished 5-edge colouring.  We view this edge-colouring as an equivalence relation on the edges, where two 
edges are equivalent if they have the same colour.  So, the automorphism group of this coloured graph is the subgroup of the
automorphism group of the underlying graph which preserves this equivalence relation.  In Figure \ref{autgroup}, the 
automorphism group of both the edge-coloured $K_6$ and the edge-coloured Petersen is ${\mathbf S}_5$, and each permutation
of these 5 colours corresponds uniquely to an isomorphism of the coloured graph (for Petersen this is the full automorphism group).  
In both cases, ${\mathbf S}_5$ acts transitively on the shortest cycles of the graph, but when we restrict to the subgroup 
${\mathbf A}_5$ these cycles break into two orbits, each of which yields an embedding in the projective plane.

The following lemma records the information we require concerning the automorphism groups of videos.  
For an incidence geometry $\Lambda = (X_1, \ldots, X_n)$, we say that $G \le {\mathit Aut}(\Lambda)$ is \emph{fully transitive} 
if $G$ acts transitively on $X_i$ for every $1 \le i \le n$.

\begin{proposition}
\label{video_aut}
The following chart indicates the poset of fully transitive subgroups of the 
automorphism group of each of the exceptional videos.

\begin{center}
\begin{tabular}{|c|c|}
\hline
Video		&	Fully transitive groups		 	
	\\
\hline
\raisebox{5ex}[0pt]{ $\mathrm{Cube/Octahedron} : V \sim E \sim F$ }	&
\begin{picture}(180,60)\thinlines
 \put(69,45){${\mathbf S}_4 \times {\mathbf C}_2$}
 \put(84,38){\line(-2,-3){13}}
 \put(94,38){\line(2,-3){13}}
 \put(61,3){${\mathbf S}_4$}
 \put(92,3){${\mathbf A}_4 \times {\mathbf C}_2$}
\end{picture}
	\\
\hline
\raisebox{4.0ex}[0pt]{ $\begin{array}{l}
		\mathrm{Dodecahedron} : F \sim V \sim E		\\
		\mathrm{Dodecahedron/Icosahedron} : V \sim E \sim F	\\
		\mathrm{Icosahedron} : F \sim V \sim E
			\end{array}$ }	&
\begin{picture}(180,50)\thinlines
 \put(66,35){${\mathbf A}_5 \times {\mathbf C}_2$}
 \put(88,28){\line(0,-1){12}}
 \put(83,3){${\mathbf A}_5$}
\end{picture}
	\\
\hline
$\mathrm{Petersen} : C_5 \sim E \sim V$ 	&
$ {\mathbf S}_5 $	
	\\
\hline
$\mathrm{Petersen}/K_6 : F \sim E \sim V$	&
$ {\mathbf A}_5	$	
	\\
\hline
\raisebox{4.5ex}[0pt]{ $K_6: C_3 \sim E \sim V$ }	&
\begin{picture}(180,55)\thinlines
 \put(80,40){ ${\mathbf S}_6$}
 \put(84,33){\line(-2,-3){10}}
 \put(94,33){\line(2,-3){10}}
 \put(61,3){ ${\mathbf S}_5$}
 \put(99,3){ ${\mathbf A}_6$}
\end{picture}
	\\
\hline
\end{tabular}
\end{center}

\end{proposition}

\noindent{\it Proof:} First note that the automorphism groups of the graphs Cube/Octahedron, Dodecahedron/Icosahedron, and Petersen are isomorphic 
to ${\mathbf S}_4 \times {\mathbf C}_2$, ${\mathbf A}_5 \times {\mathbf C}_2$, and ${\mathbf S}_5$.  
Every subgroup which acts transitively on the vertices and edges of these graphs (except Octahedron) must act transitively on the arcs (so must have size a multiple of $2|E|$ where $E$ is the set of edges), and it follows easily that the first two rows in the table
are valid.  For Petersen, the only proper subgroup of the automorphism group which is both vertex and edge transitive is isomorphic to ${\mathbf A}_5$.  The 5 cycles of Petersen form a single orbit under the full automorphism group, but split into two orbits of size 6 under this ${\mathbf A}_5$ subgroup, and this implies the validity of the third and fourth rows of the table.  

To verify the final row, let $E = E(K_6)$, $V= V(K_6)$, and $T = C_3(K_6)$ and let $G$ be a subgroup 
of ${\mathit Aut}(K_6 : C_3 \sim E \sim V)$ which acts transitively on 
$E$ and $V$ and $T$.  It follows from Lemma \ref{inc_trans} that $G$ acts transitively
on the arcs.  Define a relation $\sigma$ on $E$ by the rule $e \, \sigma f$ if either $e = f$ or if $e$, $f$ are nonadjacent and there 
exists an element of $G$ which transposes the ends of $e$, transposes the ends of $f$, 
and fixes the other two vertices.  Since there are 30 arcs and 20 triangles, the subgroup of $G$ which 
stabilizes an arc must have even order so there must be an element in $G$ which acts on $V$ as either a
transposition or as a permutation which consists of two disjoint transpositions.  
In the former case, it follows from edge transitivity that $G \cong {\mathbf S}_6$ and we are done.  
Thus, we may assume that $G$ contains a permutation consisting of two disjoint transpositions, so $\sigma$ is 
is nontrivial.

Since there are 15 edges, every edge $e$ must satisfy $e \,  \sigma f$ and $e \, \sigma f'$ for 
at least two distinct edges $f,f'$.  If there exist such edges $f,f'$ where $f$ and $f'$ are adjacent, 
then by composing the two corresponding permutations we find an element of $G$ which acts as a 3-cycle on $V$.  In this case it follows from the transitivity of $G$ on $T$ that $G$ contains all 3-cycles on $V$, so $G$ must be either the full automorphism group or the subgroup isomorphic to 
${\mathbf A}_6$.  In the only remaining case, $\sigma$ is an equivalence relation where each equivalence class is a perfect matching.  Following 
our earlier discussion, the subgroup of $Aut(K_6)$ which preserves the 5-edge colouring given by $\sigma$ is ${\mathbf S}_5$ and $G$ must be a subgroup 
of this.  Since the stabilizer of an arc has order at least two, $G$ must have order at least $60$.  However, we cannot have $G \cong {\mathbf A}_5$ as this would contradict the assumption that $G$ acts transitively on $T$.  Therefore, $G \cong {\mathbf S}_5$, and this completes the proof.
\quad\quad$\Box$

\subsection{Structure Theorem I}

Our goal over the next two subsections will be to derive our main structure theorem for subset trios from our structure theorem for incidence trios.  This is largely a rather straightforward exercise in specializing from group actions back to the setting of groups.  However, due to the inherent complexity of our structures, it still requires some effort.  
In this section we will introduce some helpful terminology, and then prove lemmas to handle beats, cyclic chords, and videos.  In the next subsection we will prove lemmas to take care of dihedral chords and pure chords, and then establish the structure theorem.

\bigskip

\noindent{\bf Definition:} Let $\Lambda = (X_1, \ldots, X_n)$ be a Cayley $G$-chorus and assume that the action of $G$ on $X_i$ is regular for 
$1 \le i \le n$.  Let $x_i \in X_i$ for every $1 \le i \le n$ and following Theorem \ref{sabidussi} make the following associations.
\begin{itemize}
\item Every $x \in X_i$ is associated with the unique $g \in G$ for which $g x_i = x$.
\item For every $1 \le i,j \le n$ the side $(X_i,X_j)$ is associated with the set of all group elements which are associated to a point in $N(x_i) \cap X_j$.
\end{itemize}
Now Theorem \ref{sabidussi} gives us a representation of $\Lambda$ as a Cayley chorus $\Lambda' = {\mathit CayleyC}(G; M)$ where $M$ is 
the $n \times n$ matrix with $ij$ entry equal to the subset associated with the side $(X_i,X_j)$.  We say that $\Lambda'$ is the \emph{realization of} $\Lambda$ 
\emph{using base points} $x_1, \ldots, x_n$.  

\bigskip

The following observation follows immediately from our definitions.

\begin{observation}
Let $(A,B,C)$ be a subset trio in $G$ and let $x,y,z \in G$.  The realization of ${\mathit CayleyC}(G; A,B,C)$ using
base points $(x,1),(y,2),(z,3)$ is given by
\[ {\mathit CayleyC}(G;  xAy^{-1}, yBz^{-1}, zCx^{-1}). \]
\end{observation}

We will repeatedly apply the above observation to move from a subset trio $(A,B,C)$ to a similar subset trio $(A',B',C')$ in the 
forthcoming lemmas.

\begin{lemma} 
\label{i2s_purebeat}
Let $(A,B,C)$ be a subset trio, and assume that $\Theta = {\mathit CayleyC}(G; A,B,C)$ is a pure beat relative to 
$(G \times \{ 1 \}, G \times \{2\})$.  Then there exists a proper finite subgroup $H$ so that $A$ is a left $H$-coset and 
$(A,B,C)$ is a pure beat relative to $H$.  
\end{lemma}

\noindent{\it Proof:} 
Let $\mcp$, $\mcq$ be the clone partitions of $G \times \{1\}$, $G \times \{2\}$.  Choose $x \in G$ so that $(x,1) \sim (1,2)$ and choose $z \in G$ arbitrarily.  Now let $P \in \mcp$ contain $(x,1)$, let $Q \in \mcq$ contain $(1,2)$, let $H = G_P = G_Q$, and note that $H$ is 
a proper finite subgroup of $G$.  
Next we define $\Theta' = {\mathit CayleyC}(G; A', B', C')$ be the realization of $\Theta$ using the base points $(x,1), (1,2), (z,3)$.  
Observe that by construction we have $A' = H$.  Furthermore, by the maximality of $\Theta$ we have ${\mathit stab}_L(B') = H$ 
and $C' = \overline{ (A'B')^{-1} }$.  Therefore $(A',B',C')$, and similarly $(A,B,C)$ is a pure beat relative to $H$.  Furthermore, by construction $A$ is contained in a left $H$-coset. \quad\quad$\Box$

\begin{lemma}
\label{i2s_impbeat}
Let $(A,B,C)$ be a subset trio in $G$, and assume $\Theta = {\mathit CayleyC}(G; A,B,C)$ is an impure beat relative to 
$(G \times \{1\}$, $G \times \{2\})$ with continuation $\Theta_1$, an $H$-trio.  Then there exist $A',B',C' \subseteq H$ so that
\begin{enumerate}
\item ${\mathit CayleyC}(H; A', B', C')$ is a realization of $\Theta_1$.
\item $(A,B,C)$ is an impure beat relative to $H$ with continuation $(A',B',C')$.
\item The minimal subgroup $H'$ so that $A$ is contained in a left $H'$-coset is conjugate to $H$.
\end{enumerate}
\end{lemma}

\noindent{\it Proof:} Let $\mcp$, $\mcq$, and $\mcr$ be the systems of imprimitivity on $G \times \{1\}$, $G \times \{2\}$ and $G \times \{3\}$ associated 
with this impure beat, and choose $P \in \mcp$, and $Q \in \mcq$, and $R \in \mcr$ so that $\Theta_1 = (P,Q,R)$ (so $H = G_P = G_Q = G_R$). 
Choose $x,y,z \in G$ so that $(x,1) \in P$ and $(y,2) \in Q$ and $(z,3) \in R$ and let 
$\Theta' = {\mathit CayleyC}(G; A', B', C')$ be the realization of $\Theta$ using the base points $(x,1)$, $(y,2)$, and $(z,3)$.  
As in the definition of realization, we associate each point $(u,1) \in G \times \{1\}$ with the group element $g$ for which 
$(gx,1) = (u,1)$ and we associate points in $G \times \{2\}$ and $G \times \{3\}$ with elements in $G$ using $(y,2)$ and $(z,3)$ similarly.
Note that for every block in $\mcp \cup \mcq \cup \mcr$, the set of group elements associated with this block (i.e. associated with a point in this block)  is a left $H$-coset.  

Now $H$ is precisely the set of elements associated with a point in $Q$, which gives us $\emptyset \neq A' \subset H$.  
Furthermore, since $(P,Q)$ is connected, it follows that $A'$ is not contained in a left coset of a proper subgroup of $H$.  
Since $A' = xAy^{-1}$ it follows that $H' = y^{-1}Hy$ is the unique minimal subgroup for which $A$ is contained in a left $H'$-coset.  
Now we define the following sets.  
\begin{align*}
	B_0	&= \{ g \in G \mid \mbox{$g$ is associated with a point in $N_{\Theta}(z,3) \cap Q$} \}	\\
	B_1	&= \{ g \in G \mid \mbox{$g$ is associated with a point in $N_{\Theta}(z,3) \cap (G \times \{2\} \setminus Q)$} \}	\\
	C_0	&= \{ g \in G \mid \mbox{$g$ is associated with a point in $N_{\Theta}(z,3) \cap P$} \}	\\
	C_1	&= \{ g \in G \mid \mbox{$g$ is associated with a point in $N_{\Theta}(z,3) \cap (G \times \{1\} \setminus P)$} \}.   
\end{align*}
It follows from our construction that $(B')^{-1} = B_0 \cup B_1$ and $C' = C_0 \cup C_1$.  Futhermore, we have $\emptyset \neq B_0,C_0 \subseteq H$ and $H \le {\mathit stab}_R(B_1)$, ${\mathit stab}_R(C_1)$.  This implies $H \le {\mathit stab}_L(B' \setminus H)$.
Now, for every $S \in G/H \setminus \{H\}$ there exist blocks $P' \in \mcp$ and $Q' \in \mcq$ which are associated with 
this set of elements, and these blocks satisfy $P' \sim Q'$.  Furthermore, by construction, we must have either
$S \subseteq B_1$ or $S \subseteq C_1$.  It follows from this that 
$\overline{(A'B')^{-1}} \setminus H = C \setminus H$.  
We conclude that $(A',B',C')$ and similarly $(A,B,C)$ is an impure beat relative to $H$.  Furthermore, 
$(A', B_0, C_0)$ is a continuation of $(A,B,C)$ and ${\mathit CayleyC}(H; A', B_0, C_0)$ is a realization of $\Theta_1$.
\quad\quad$\Box$

\begin{lemma}
\label{i2s_cchord}
Let $(A,B,C)$ be a subset trio in $G$ and assume $\Theta= {\mathit CayleyC}(G; A,B,C)$ is a cyclic chord with continuation $\Theta_1$, an $H$-trio.  Then there exists $A', B', C' \subseteq H$ so that 
\begin{enumerate}
\item ${\mathit CayleyC}(H; A', B', C')$ is a realization of $\Theta_1$.
\item $(A,B,C)$ is an impure cyclic chord with continuation $(A', B', C')$.
\end{enumerate}
\end{lemma}

\noindent{\it Proof:} By possibly replacing $(A,B,C)$ by a similar trio, we may assume that $w(G \times \{1\}, G \times \{2\} )$, $w(G \times \{2\}, G \times \{3\}) < \infty$.  We may choose a group $J = {\mathbb Z}$ or $J = {\mathbb  Z}/n {\mathbb Z}$ and 
systems of imprimitivity $\mcp = \{ P_j \mid j \in J \}$ of $G \times \{1\}$ and $\mcq = \{ Q_j \mid j \in J \}$ of $G \times \{2\}$ and $\mcr = \{ R_j \mid j \in J \}$ of $G \times \{3\}$ together with integers $\ell,m$ in accordance with the definition of cyclic chord and assume (without loss of generality) that 
$\Theta_1 = (P_0, Q_0, R_0)$.  Note that the subgroup $H$ as given in the statement of this lemma also aligns with that given in this 
definition.  Choose $S \in G/H$ so that every $g \in S$ satisfies $gP_j = P_{j+1}$ and $gQ_j = Q_{j+1}$ and $gR_j = R_{j+1}$ 
for every $j \in J$.  

Now choose $x,y,z \in G$ so that $(x,1) \in P_0$ and $(y,2) \in Q_0$ and $(z,3) \in R_0$ and let $\Theta' = {\mathit CayleyC}(G; A', B', C')$ be the realization of $\Theta$ using the base points $(x,1), (y,2), (z,3)$.  It follows that for every $j \in J$, the set of group elements associated with  
$P_j$ (or $Q_j$ or $R_j$) is given by $S^j$.  Therefore, $A' \cup H$ and $B' \cup H$ are $H$-stable, $\varphi_{G/H} (A')$, $\varphi_{G/H}(B')$ are basic geometric progressions with ratio $S$, and $\emptyset \neq \overline{(A'B')^{-1}} \setminus H = C' \setminus H$ and $C' \cap H \neq \emptyset$.  Thus $(A',B',C')$, and similarly $(A,B,C)$, is an impure
cyclic chord.  Furthermore, setting $A'' = A' \cap H$ and $B'' = B' \cap H$ and $C'' = C' \cap H$ we have that $(A'',B'',C'')$ is a continuation
of $(A,B,C)$ and ${\mathit CayleyC}(H; A'', B'', C'')$ is a realization of $\Theta_1$, as desired.
\quad\quad$\Box$

\begin{lemma}
\label{i2s_video}
Let $(A,B,C)$ be a subset trio in $G$ and let $\Theta = {\mathit CayleyC}(G; A, B, C)$.  If $\Theta$ is a 
video, there exists a subset trio $(A',B',C')$ similar to $(A,B,C)$ and subgroups $H \triangleleft G$ and $H \triangleleft H_1,H_2,H_3 \le G$ with 
$H_1 = {\mathit stab}_L(A') = {\mathit stab}_R(C')$, $H_2 = {\mathit stab}_R(A') = {\mathit stab}_L(B')$, and $H_3 = {\mathit stab}_R(B') = {\mathit stab}_L(C')$
satisfying one of the following
\begin{enumerate}
\item $|A'| = 2|H_2| = |B'|$ and $\delta(A',B',C') = |H_1| =  |H_3| < |H_2|$.
\item $9|H_1| = |A'| = 3|H_2|$ and $2|H_2| = |B'| = 3|H_3|$ and $\delta(A',B',C') = |H_1|$.  
Further, if $G$ is finite, then 
$[H_1 : H]$ divides $16$.
\item $4|H_1| = |A'| = 2|H_2|$ and $3|H_2| = |B'| = 2|H_3|$ and $\delta(A',B',C') = |H_1|$.  
Further, if $G$ is finite, then $[H_1 : H]$ divides $16$.
\item One of the following is satisfied.
\begin{center}
\begin{tabular}{|c|c|c|c|c|c|c|}
\hline
\begin{picture}(0,19)\end{picture} \raisebox{.95ex}{$G/H$}	
	&\raisebox{.95ex}[0pt]{$\tfrac{|G|}{|H_1|}$}	&\raisebox{.95ex}[0pt]{$\tfrac{|G|}{|H_2|}$}	
	&\raisebox{.95ex}[0pt]{$\tfrac{|G|}{|H_3|}$}	&\raisebox{.95ex}[0pt]{$\tfrac{|G|}{|A'|}$}	
	&\raisebox{.95ex}[0pt]{$\tfrac{|G|}{|B'|}$}		&\raisebox{.95ex}[0pt]{$\tfrac{|G|}{|C'|}$} 	\\
\hline
${\mathbf C_2} \times {\mathbf S_4}$, ${\mathbf C_2} \times {\mathbf A_4}$, or ${\mathbf S_4}$ 
		& $8$	& $12$	& $6$	& $4$	& $3$	& $2$	\\	
${\mathbf C_2} \times {\mathbf A_5}$ or ${\mathbf A_5}$
		& $12$	& $20$	& $30$	& $4$	& $10$	& $3/2$	\\	
${\mathbf C_2} \times {\mathbf A_5}$ or ${\mathbf A_5}$
		& $20$	& $30$	& $12$	& $10$	& $6$	& $4/3$	\\	
${\mathbf C_2} \times {\mathbf A_5}$ or ${\mathbf A_5}$
		& $20$	& $12$	& $30$	& $4$	& $6$	& $5/3$	\\	
${\mathbf S_5}$
		& $12$ & $15$ & $10$ & $3$ & $5$ & $2$ \\
${\mathbf A_5}$
		& $6$ & $15$ & $10$ & $3$ & $5$ & $2$ \\
${\mathbf S_6}$, or ${\mathbf A_6}$, or ${\mathbf S}_5$
		& $20$ & $15$ & $6$ & $5$ & $3$ & $2$ \\
\hline
\end{tabular}
\end{center}
\end{enumerate}
\end{lemma}

\noindent{\it Proof:} By possibly replacing $(A,B,C)$ by an similar trio, we may assume that $\Theta^{\bullet}$ is isomorphic to 
one of the trios listed in the definition of standard video or exceptional video.  Now for $1 \le i \le 3$ choose a subgroup $H_i \le G$ so that 
the clone partition of $G \times \{i\}$ is given by $G/H_i \times \{i\}$ for $1 \le i \le 3$.  It follows that 
$H_1 A H_2 = A$ and $H_2 B H_3 = B$ and $H_3 C H_1 = C$.  Let $H \triangleleft G$ be the kernel of the action of $G$ on $\Theta^{\bullet}$ 
and note that $H \le H_i$ for $1 \le i \le 3$.  We now split into cases depending on $\Theta^{\bullet}$.

\bigskip

\noindent{\it Case 1:} $\Theta^{\bullet}$ is isomorphic to $\Gamma : E \sim V \sim E$ for a graphic duet $\Gamma$ with degree $d \ge 3$.  

\smallskip

In this case, it follows from the degrees of the duets in $\Theta^{\bullet}$ that $|A| = 2|H_2| = |B|$ and then 
$d |H_1| = |A| = |B| = d |H_3|$.  Lemma \ref{crit2cut} gives us $\delta(A,B,C) = \delta(\Theta) = |H_1| = |H_3| < |H_2|$
so the first outcome holds.

\bigskip

\noindent{\it Case 2:} $\Theta^{\bullet}$ is isomorphic to $\Gamma : P_3 \sim V \sim E$ for graphic duet $\Gamma$ of degree $3$

\smallskip

Here it follows from the degrees of the duets in $\Theta^{\bullet}$ that $9 |H_1| = |A| = 3 |H_2|$ and 
$2|H_2| = |B| = 3|H_3|$.  Lemma \ref{crit2cut} gives us $\delta(A,B,C) = \frac{1}{2}|H_3| = |H_1|$ and Theorem \ref{tutte} implies 
that $[H_1 : H]$ divides $16$ whenever $G$ is finite.  It follows that the second outcome holds.

\bigskip

\noindent{\it Case 3:} $\Theta^{\bullet}$ is isomorphic to $\Gamma : P_3 \sim E \sim V$ for a graphic duet $\Gamma$ of degree $3$.

\smallskip

In this case it follows from the degrees of duets in $\Theta^{\bullet}$ that $4 |H_1| = |A| = 2 |H_2|$ and 
$3|H_2| = |B| = 2|H_3|$.  Lemma \ref{crit2cut} gives us $\delta(A,B,C) = \frac{1}{2}|H_2| = |H_1|$ and Theorem \ref{tutte} implies 
that $[H_1 : H]$ divides $16$ whenever $G$ is finite.  So, the third outcome is satisfied.

\bigskip

\noindent{\it Case 4:} $\Theta^{\bullet}$ is isomorphic to an exceptional video.

\smallskip

In all of these cases, we will show that one of the outcomes from the table in part 4 of the statement of the lemma 
is satisfied.
These outcomes correspond in order, to the videos 
Cube/Oct. : $V \sim E \sim F$, Dodec. : $F \sim V \sim E$, 
Dodec./Icos. : $V \sim E \sim F$, Icos. : $F \sim V \sim E$, 
Petersen: $C_5 \sim E \sim V$, Petersen/$K_6$ : $F \sim E \sim V$, 
and $K_6$ : $C_3 \sim E \sim V$.  The structure of the group $G/H$ is given by Proposition \ref{video_aut}, 
and the remaining properties are straightforward to verify (much of this information is 
encapsulated in Figure \ref{alldeg}).  
\quad\quad$\Box$

\subsection{Structure Theorem II}

In Section \ref{dihedral_sec} we introduced the notion of an octahedral configuration (recall that this is an oriented octahedron with edges labelled by subsets from a group $H$), and then in Section \ref{inc_sec} we defined a Cayley chorus.  In fact, Cayley choruses provide a natural framework to work with octahedral configurations.  To treat these in the desired manner, we introduce a natural definition.

\begin{figure}[ht]
\centerline{\includegraphics[height=5cm]{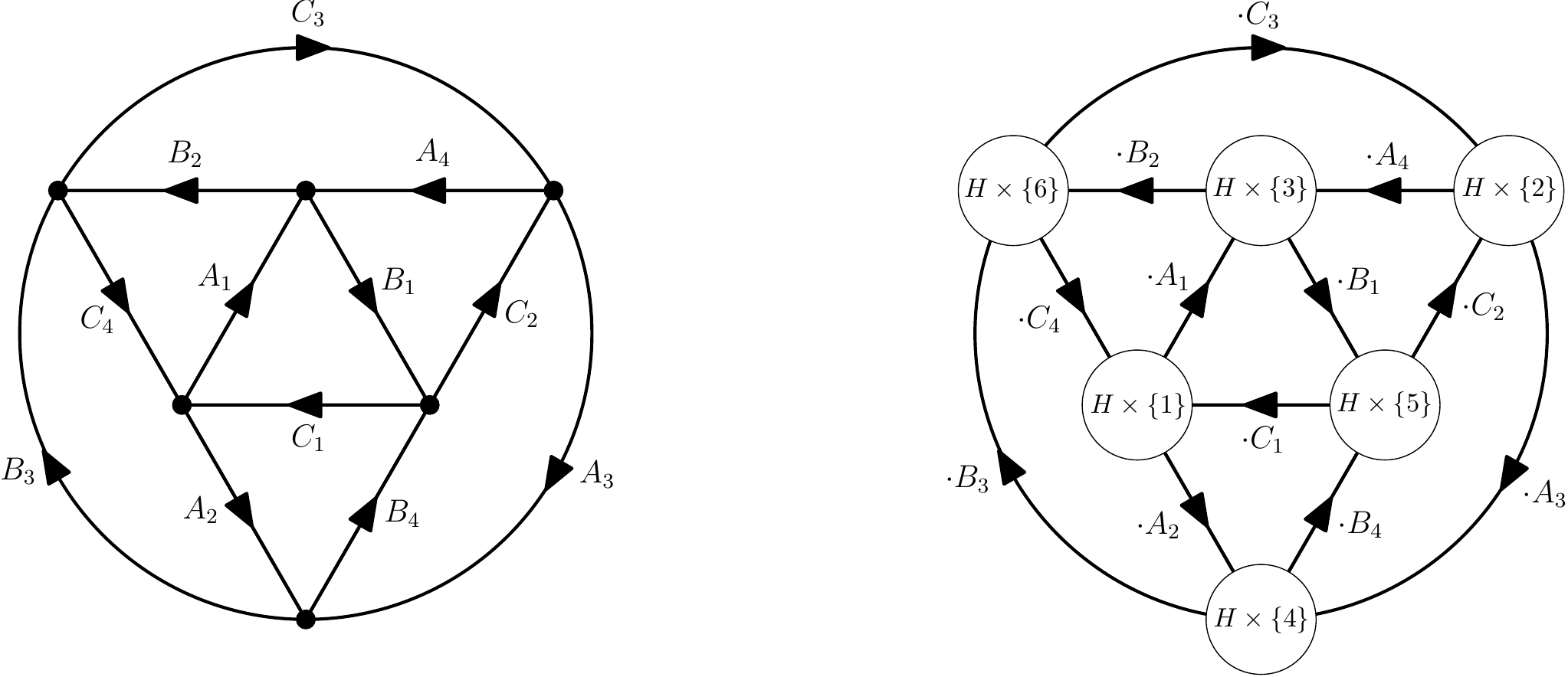}}
\caption{An octahedral configuration and a corresponding octahedral Cayley chorus}
\label{oct_cor}
\end {figure}

\noindent{\bf Definition:} Let $A_1, \ldots, A_4, B_1, \ldots, B_4, C_1, \ldots, C_4 \subseteq H$ and let $\Omega$ be an octahedral configuration 
given by these sets and the labelled graph on the left in Figure \ref{oct_cor}.  We define the \emph{Cayley octahedral chorus} 
\[ {\mathit CayleyOct}(H; A_1, \ldots, A_4, B_1, \ldots, B_4, C_1, \ldots, C_4 ) = {\mathit CayleyC}(H; M) \]
where the matrix $M$ is given below
\[ M = \left[ \begin{array}{cccccc}
	\emptyset	&	\emptyset		&	A_1		&		A_2		&		C_1^{-1}		&		C_4^{-1}	\\
	\emptyset	&	\emptyset		&	A_4		&		A_3		&		C_2^{-1}		&		C_3^{-1}	\\
	A_1^{-1}	&	A_4^{-1}		&	\emptyset	&	\emptyset		&		B_1			&		B_2		\\
	A_2^{-1}	&	A_3^{-1}		&	\emptyset	&	\emptyset		&		B_4			&		B_3		\\
	C_1		&	C_2			&	B_1^{-1}	&	B_4^{-1}		&		\emptyset		&		\emptyset	\\
	C_4		&	C_3			&	B_2^{-1}	&	B_3^{-1}		&		\emptyset		&		\emptyset		
	\end{array} \right] \]
and we equip this chorus with the following partition.
\[ \big\{ \;  \{ H \times \{1\}, H \times \{2\} \}, \;  \{ H \times \{3\}, H \times \{4\} \}, \; \{ H \times \{5\}, H \times \{6\} \} \; \big\}. \]
It follows from the definition of octahedral configuration that this chorus is collision free, and therefore it is indeed 
an octahedral chorus.  We say that this chorus \emph{corresponds} with $\Omega$, and we have depicted it on the right in 
Figure \ref{oct_cor}.  

Of course, octahedral configurations and octahedral Cayley choruses, are in essence the same.  However, 
we have definitions in place for type and continuation of an octahedral configurations which are based on subsets, and 
other definitions in place for type and continuation of an octahedral chorus which are based on properties of these 
incidence geometries.  Our next lemma shows that these definitions align properly.

\begin{lemma} 
\label{octc_lem}
Let $H$ be a finite group, let $\Omega$ be an octahedral configuration in $H$ and let $\Lambda$ be the corresponding 
Octahedral Cayley chorus.  Then we have
\begin{enumerate}
\item $\Lambda$ is not type $2C$.
\item If $\Lambda$ has type $-1$, $0$, $1$, $2A$,  or $2B$ then $\Omega$ has the same type.
\item If $\Lambda$ has type $1$ or $2A$ with continuation $\Theta$, a $K$-trio,  there are $A',B',C' \subseteq K$ so that
\begin{enumerate}
\item ${\mathit CayleyC}(K, A', B', C')$ is a realization of $\Theta$.
\item $(A',B',C')$ is a continuation of $\Omega$.
\end{enumerate}
\end{enumerate}
\end{lemma}

\noindent{\it Proof:} Let $\Omega$ be given by the labelled graph on the left in Figure \ref{oct_cor}, so $\Lambda$ is 
depicted on the right in this figure.  Note that $\Lambda$ has rank 6, so in particular $\Lambda$ is not type $2C$.  If $\Lambda$ has type $-1$ or $0$, then it follows immediately that $\Omega$ has the same type.  Next suppose that $\Lambda$ has type $1$ and assume for notational convenience that $A_1, B_1, C_1$ are proper nonempty subsets of $H$, so 
$\Theta =  (H \times \{1\}, H \times \{3\}, H \times \{5\})$ is the unique continuation of $\Lambda$.  Then $\Omega$ is type 1, and 
${\mathit CayleyC}(H; A_1,B_1,C_1)$ is a realization of $\Theta$ and $(A_1,B_1,C_1)$ is a continuation of $\Omega$, as desired.

Next assume that $\Lambda$ has type $2A$ or $2B$ and assume for notational convenience that $A_1, B_1, C_1, B_2, C_4$ are proper nonempty subsets of $H$.  
Now let $\Upsilon_1 = (H \times \{1\}, H \times \{3\}, H \times \{5\} ; \sim_1)$ and $\Upsilon_2 = (H \times \{1\}, H \times \{3\}, H \times \{6\} ; \sim_2 )$ be the $H$-trios as 
given by the definition of types $2A$ and $2B$.  For $i=1,2$ we may choose $A_1^i \subseteq H$ so that $(x,1) \in H \times \{1\}$ and $(y,3) \in H \times \{3\}$ satisfy $(x,1) \sim_i (y,3)$ if and only if $y \in x A_1^i$, and note that the definitions of types $2A$ and $2B$ (for incidence trios) imply that $A_1 = A_1^1 \cap A_1^2$.  If $\Lambda$ is type $2B$ then $\Upsilon_1$ and $\Upsilon_2$ are pure beats relative to $(H \times \{1\}, H \times \{3\})$.  Now applying Lemma \ref{i2s_purebeat} we have that there exist proper subgroups $K_1, K_2 < H$ so that $(A_1^1,B_1,C_1)$ is a pure beat relative to $K_1$ and $A_1^1$ is a left $K_1$-coset, and similarly 
$(A_1^2,B_2,C_4)$ is a pure beat relative to $K_2$ and $A_1^2$ is a left $K_2$-coset.  It follows that $\Omega$ has type $2B$ as desired.  Next suppose that $\Lambda$ has type $2A$ with continuation $\Theta$ a $K_1$-trio, and assume (without loss) that $\Upsilon_1$ is an impure beat relative to $(H \times \{1\}, H \times \{3\})$ and $\Upsilon_2$ is a pure beat relative to $(H \times \{1\}, H \times \{3\})$ (so $A_1 = A_1^1 \subseteq A_1^2$).  In this case, applying Lemma \ref{i2s_purebeat} to $\Upsilon_2$ gives us a subgroup $K_2 < H$ so that $(A_1^2, B_2, C_4)$ is a pure beat relative to $K_2$, and $A_1^2$ is a left $K_2$ coset.  Now applying Lemma \ref{i2s_impbeat} to the impure beat $(A_1, B_1, C_1)$ and the continuation $\Theta$ gives us $A',B',C' \subseteq K_1$ so that ${\mathit CayleyC}(K_1, A', B', C')$ is a realization of $\Theta$ and 
$(A',B',C')$ is a continuation of $(A_1, B_1, C_1)$ and the minimal subgroup $K_1'$ so that $A_1$ is contained in a $K_1'$ coset is conjugate to $K_1$.  Since $K_1$ and $K_1'$ are conjugate, $(A_1,B_1,C_1)$ is also an impure beat relative to $K_1'$, and since $K_1' \le K_2$ 
it follows that $\Omega$ has type $2A$ with continuation $(A',B',C')$, thus completing the proof.
\quad\quad$\Box$

\bigskip

With this lemma in place, we can now prove two of the theorems stated in Section \ref{dihedral_sec} concerning the classification of critical 
octahedral configurations.  

\bigskip

\noindent{\it Proof of Theorems \ref{oct_config_class} and \ref{oct_set_classify}:}  Let $\Omega$ be maximal critical octahedral configuration and let 
$\Lambda$ be the associated octahedral Cayley chorus.  It follows from our definitions that $\Lambda$ is also maximal and critical, and by construction $\Lambda$ 
must have rank 6.  It then follows from Theorems \ref{oct_type} and \ref{classify_oct2} that $\Lambda$ has type $-1$, $0$, $1$, $2A$, or $2B$.  Lemma 
\ref{octc_lem} now implies that $\Omega$ has type $-1$, $0$, $1$, $2A$, or $2B$ which implies both Theorems  \ref{oct_config_class} and
\ref{oct_set_classify}.
\quad\quad$\Box$

\begin{lemma} 
\label{i2s_dchord}
Let $(A,B,C)$ be a subset trio in $G$ and assume that $\Theta = {\mathit Cayley}(G; A, B, C)$ is a dihedral chord.  
Then $(A,B,C)$ is an impure dihedral chord, and further
\begin{enumerate}
\item $\Theta$ is not type $2C$.
\item If $\Theta$ has type $0$, $1$, $2A$, or $2B$ then $(A,B,C)$  has the same type.
\item If $\Theta$ has type $1$ or $2A$ with continuation $\Theta_1$, there exists a finite subgroup $K < G$ and $A',B',C' \subseteq K$ so that 
	${\mathit CayleyC}(K; A', B', C')$ is a realization of $\Theta_1$ and $(A',B',C')$ is a continuation of $(A,B,C)$.
\end{enumerate}
\end{lemma}

\noindent{\it Proof:} By replacing $(A,B,C)$ by a similar subset trio, we may choose $J = {\mathbb Z}$ or ${\mathbb Z}/n {\mathbb Z}$, systems 
of imprimitivity $\mcp = \{ P_j \mid j \in J\}$, $\mcq = \{Q_j \mid j \in J\}$, $\mcr = \{R_j \mid j \in J\}$ on $G \times \{1\}$, 
$G \times \{2\}$, $G \times \{3\}$,  a subgroup $H \triangleleft G$, and integers $\ell,m$ as in the definition of dihedral chord.  We may also assume 
(without loss) that if there is a continuation $\Theta_1$ as described in the statement of the lemma, then $\Theta_1$ has ground set a 
subset of $P_0 \cup Q_0 \cup R_0$.

By assumption the action of $G/H$ on each of $\mcp, \mcq, \mcr$ is dihedral, so we may choose a rotation $S \in G/H$
so that every $g \in S$ satisfies $gP_j = P_{j+1}$ and $gQ_j = Q_{j+1}$ and $gR_j = R_{j+1}$.  
Now choose $(x,1) \in P_0$ and $(y,2) \in Q_0$ and $(z,3) \in R_0$ and let $\Theta' = {\mathit CayleyC}(G; A', B', C')$ be the realization of $\Theta$ using these base points.  This associates each block in $\mcp \cup \mcq \cup \mcr$ with a set of group elements which consists of the union of 
two left $H$-cosets, one rotation and one flip.  Furthermore, since the action of $G$ on $\Theta$ is regular, we may refine our partitions 
$\mcp, \mcq, \mcr$ by letting $\{ P_j^0, P_j^1\}$ be a partition of $P_j$ so that the set $P_j^0$ is associated with a rotation in $G/H$ and $P_j^1$ is associated with a flip.  Define $Q_j^i$ and $R_j^i$ similarly for every $j \in J$ and $i=0,1$.  It follows immediately that 
$\{ P_j^i \mid \mbox{$j \in J$ and $i = 0,1$} \}$ and $\{ Q_j^i \mid \mbox{$j \in J$ and $i = 0,1$} \}$ and $\{ R_j^i \mid \mbox{$j \in J$ and $i = 0,1$} \}$
are systems of imprimitivity.  Furthermore, our definitions imply 
$(x,1) \in P_0^0$ and $(y,2) \in Q_0^0$ and $(z,3) \in R_0^0$.  Now we define a family of $H$-cosets as follows.
\[ \left\{ \begin{array}{c}
			A_1, A_2, A_3, A_4	
			\vspace{.1cm} \\
			B_1, B_2, B_3, B_4	
			\vspace{.1cm} \\
			C_1,C_2, C_3, C_4	
			\end{array} \right\}
		\mbox{ are associated with }
	\left\{ \begin{array}{c}
			Q_0^0, Q_0^1, Q_{\ell}^0, Q_{\ell}^1	
			\vspace{.1cm} \\
			R_0^1, R_0^1, R_m^0, R_m^1	
			\vspace{.1cm} \\
			P_0^0, P_0^1, P^0_{\ell+m}, P^1_{\ell+m} 
			\end{array} \right\}	\]
Note that this construction yields $A_1, B_1, C_1 = H$.  Next define $A^+ = A' \cup \big( \cup_{i=1}^4 A_i \big)$ and $B^+ = B' \cup \big( \cup_{i=1}^4 B_i \big)$ and $ C^+ = C' \cup \big( \cup_{i=1}^4 C_i \big)$.  It follows immediately that $A^+$, $B^+$, and $C^+$ are $H$-stable.  Furthermore, $\varphi_{G/H}(A^+)$ is a dihedral progression with ratio $S$ given by $\{H, S, S^2, \ldots, S^{\ell}\} \{H, A_2\}$.  Similarly $\varphi_{G/H}(B^+)$ is a dihedral progression with ratio $S$ and furthermore, $(\varphi_{G/H}(A^+), \varphi_{G/H}(B^+), \varphi_{G/H}(C^+))$ is a dihedral prechord with 
labelled octahedron given by the sets $A_1, \ldots, A_4, B_1, \ldots, B_4, C_1, \ldots, C_4$ and the graph on the left in Figure \ref{oct_rem}.  
It follows from this that $(A',B',C')$ is an impure dihedral chord.

Now we will follow the procedure described in Section \ref{dihedral_sec} to construct an octahedral configuration associated with $A',B',C'$.  To 
do this, we will use the variables $u,u',v,v',w,w'$ labelling the vertices of the graph on the left in Figure \ref{oct_rem}.  Let $u',v',w' = 1$ and choose 
$u \in C_2$ and $v \in A_2$ and $w \in B_2$.  Then for $1 \le i \le 4$ if the $A_i$ edge has initial vertex labelled by $x$ and terminal vertex labelled by $y$ then we set $A_i^* = x(A' \cap A_i) y^{-1}$, and define $B_i^*$ and $C_i^*$ similarly.  This gives us an octahedral configuration $\Omega$ as depicted on the right in Figure \ref{oct_rem}.  

Returning to our incidence geometry, we define $\Lambda$ to be the $H$-octahedral chorus given by
$\Lambda = ( P_0^0, P_0^1, Q_0^0, Q_0^1, R_0^0, R_0^1 )$ with distinguished 
partition $\{ \{ P_0^0, P_0^1\}$, $\{Q_0^0, Q_0^1 \}$, $\{R_0^0, R_0^1\} \}$.  So, by our assumptions, if $\Theta_1$ is a continuation of $\Theta$, 
then $\Theta_1$ is a continuation of $\Lambda$.  Now define $\Lambda'$ to be the realization of $\Lambda$ using the base points 
$(x,1)$, $(ux,1)$, $(y,2)$, $(vy,2)$, $(z,3)$, $(wz,3)$, and assume that 
\[\Lambda' = {\mathit CayleyOct}(H; A_1^{\diamond}, \ldots, A_4^{\diamond}, B_1^{\diamond}, \ldots, B_4^{\diamond}, C_1^{\diamond}, \ldots, C_4^{\diamond} ).\]
Now we have
\begin{align*}
A_1^{\diamond} &= \{ h \in H \mid (h y, 2) \sim (x,1) \} = A' \cap A_1 = A_1^*	\\
A_2^{\diamond} &= \{ h \in H \mid (h v y, 2) \sim (x,1) \} = \{ h \in H \mid hv \in A' \cap A_2 \} = (A' \cap A_2) v^{-1} = A_2^* \\
A_3^{\diamond} &= \{ h \in H \mid (h v y, 2) \sim (u x,1) \}  = \{ h \in H \mid u^{-1}hv \in A' \cap A_3 \} = u (A' \cap A_3) v^{-1} = A_3^* \\
A_4^{\diamond} &= \{ h \in H \mid (h y, 2) \sim (u x, 1) \} = \{ h \in H \mid u^{-1}h \in A' \cap A_4 \} = u (A' \cap A_4) = A_4^*
\end{align*}
Similarly we find $B_i^{\diamond} = B_i^*$ and $C_i^{\diamond} = C_i^*$ for $1 \le i \le 4$.  Therefore, $\Lambda'$ is a Cayley octahedral 
chorus corresponding to $\Omega$, and $\Lambda'$ is a realization of $\Lambda$.  The result now follows from Lemma \ref{octc_lem}.
\quad\quad$\Box$

\begin{lemma}
\label{i2s_pchord}
Let $(A,B,C)$ be a subset trio in $G$ and assume $\Theta = {\mathit CayleyC}(G; A,B,C)$ is a pure chord.
\begin{enumerate}
\item If $\Theta$ has cyclic action, then $(A,B,C)$ is a pure cyclic chord.
\item If $\Theta$ has dihedral action, then $(A,B,C)$ is a pure dihedral chord.
\item If $\Theta$ has split action, then $(A,B,C)$ is an impure dihedral chord of type $0$.
\end{enumerate}
\end{lemma}

\noindent{\it Proof:} By possibly replacing $(A,B,C)$ with a similar trio we may assume that $A$ and $B$ are finite.  Based on the definition of pure 
chord, we may choose a group $J = {\mathbb Z}$ or $J = {\mathbb Z}/n {\mathbb Z}$ and positive integers $\ell,m$ so that the clone partitions of
 $G \times \{1\}$, $G \times \{2\}$ and $G \times \{3\}$ in $\Theta$ are given by $\mcp = \{ P_j \mid j \in J \}$, $\mcq = \{ Q_j \mid j \in J \}$, and 
 $\mcr = \{ R_j \mid j \in J \}$ where $N(P_j) \cap G \times \{2\} = Q_j \cup \ldots Q_{j+\ell}$ and $N(Q_j) \cap G \times \{3\} = R_j \cup \ldots, \cup R_{j+m}$.  
 Define $H$ to be the subgroup $G_{ ( \mcp ) } = G_{ (\mcq) } = G_{ (\mcr) }$.  Next, choose $(x,1) \in P_0$ and $(y,1) \in Q_0$ and $(z,1) \in R_0$ 
 and define $\Theta' = {\mathit CayleyC}(G; A', B', C')$ to be the realization of $\Theta$ using the base points $(x,1)$, $(y,2)$, and $(z,3)$.  
We now split into cases depending on $\Theta$.  

\bigskip

\noindent{\it Case 1:} $\Theta$ is a pure chord with cyclic action.

\smallskip

In this case $G/H$ is cyclic and we may choose a consecutive shift $S \in G/H$ so that every $g \in S$ satisfies $gP_j = P_{j+1}$ and 
$gQ_j = Q_{j+1}$ and $gR_j = R_{j+1}$ for every $j \in J$.  It now follows that $S$ generates the group $G/H$ and further 
$A' = H \cup S \cup \ldots \cup S^{\ell}$ and $B' = H \cup S \cup \ldots S^m$.  It follows from our definitions that $C' = \overline{(A'B')^{-1}}$ 
is not contained in a single $H$-coset, so $(A',B',C')$, and similarly $(A,B,C)$ is a pure cyclic chord.

\bigskip

\noindent{\it Case 2:} $\Theta$ is a pure chord with split action.

\smallskip

It follows from Lemma \ref{type0alt} that $\Theta$ is a dihedral chord of type $0$.  So by Lemma \ref{i2s_dchord} we have that $(A,B,C)$ is 
an impure dihedral chord of type $0$, as desired.

\bigskip

\noindent{\it Case 3:} $\Theta$ is a pure chord with dihedral action.  

\smallskip

In this case $G/H$ is dihedral and we may choose a consecutive shift $S \in G/H$ so that every $g \in S$ satisfies $gP_j = P_{j+1}$ and 
$gQ_j = Q_{j+1}$ and $gR_j = R_{j+1}$ for every $j \in J$.  It now follows that $S$ generates the rotation subgroup of $G/H$.  Choose a flip
$F \in G/H$ so that $G_{Q_0} = H \cup F$ and then observe that $A'$ is precisely the set $(H \cup S \cup \ldots \cup S^{\ell})(H \cup F)$ so 
$\varphi_{G/H}(A')$ is a basic dihedral progression with ratio $S$.  A similar argument shows that $(B')^{-1}$ is given by 
$(H \cup F) (H \cup S^{-1} \cup \ldots \cup S^{-m}$.  So $\varphi_{G/H}(B')$ is also a basic dihedral progression with ratio $S$ and 
${\mathit stab}_L(A') = H \cup F = {\mathit stab}_R(B')$.  It follows that $(A',B',C')$, and similarly $(A,B,C)$ is a pure dihedral chord.
\quad\quad$\Box$

\bigskip

We are now ready to prove our classification theorem for maximal critical subset trios.

\begin{theorem}
\label{subset_equiv}
If $\Phi$ is a nontrival subset trio in $G$, then $\Phi$ is maximal and critical if and only if there exist a sequence of subgroups
$G = G_1 > G_2 \ldots > G_m$ and a sequence of subset trios $\Phi = \Phi_1, \Phi_2, \ldots, \Phi_m$ 
so that the following hold.
\begin{enumerate}
\item $\Phi_i$ is a subset trio in the group $G_i$ for $1 \le i \le m$.
\item $\Phi_i$ is either an impure beat, an impure cyclic chord, or an impure dihedral chord of type $1$ or $2A$  with continuation
$\Phi_{i+1}$ for $1 \le i < m$.
\item  $\Phi_m$ is either a pure beat, a pure cyclic chord, a pure dihedral chord, an impure dihedral chord of type $0$ or $2B$, or 
$\Phi_m = (A,B,C)$ where ${\mathit CayleyC}(G_m; A, B, C)$ is a video.
\end{enumerate}
\end{theorem}

\noindent{\it Proof:} For the ``only if'' direction, we let $\Phi = (A,B,C)$ be maximal and critical and define 
$\Theta = {\mathit CayleyC}(G; A,B,C)$.  It follows from Theorem \ref{main_incidence} that $\Theta$ is a song, and so the 
result follows from Lemmas  \ref{i2s_purebeat}, \ref{i2s_impbeat}, \ref{i2s_cchord}, \ref{i2s_dchord}, and \ref{i2s_pchord}, 
together with a straightforward induction.  

For the ``if'' direction, we also let $\Phi = (A,B,C)$ and define $\Theta = {\mathit CayleyC}(G; A,B,C)$.  It follows from our 
definitions and a straightforward induction that $\Theta$ is a song, so $\Theta$ is maximal and critical.  This implies 
that $(A,B,C)$ is maximal and critical and desired.
\quad\quad$\Box$

\bigskip

\noindent{\it Proof of Theorem \ref{main_subset}: } This follows from the previous theorem and Lemma \ref{i2s_video}. 
\quad\quad$\Box$

\section{Corollaries}

In this section we will prove the corollaries of the main theorem which were stated in the Introduction.

\subsection{Weak Structure}

Our first corollary will be a weak version of the structure theorem which applies to any critical pair, and has the convenience of unpacking some of our notation.  We will then use this result to derive Corollary \ref{struc_or_stable}.

\begin{corollary}[Weak Structure]
\label{weak_structure}
If $(A,B)$ is a critical pair in the group $G$, then setting $C = \overline{(AB)^{-1}}$, one of the following holds.
\begin{enumerate}
\item One of $A$, $B$, or $C$ is contained in a coset of a proper subgroup.
\item There exists $H \triangleleft G$ so that $G/H$ is a cyclic group with $|G/H| \ge 4$ for which 
 $\varphi_{G/H}(A)$ and $\varphi_{G/H}(B)$ are geometric progressions with a common ratio.  
Furthermore, one of the following holds.
	\begin{enumerate}
	\item $\delta(A,B)  \le \frac{1}{2} |H|$ and $\max\{ |AH \setminus A|, |BH \setminus B| \} + \delta(A,B) \le |H|$.
	\item $ABH = AB$ and $\delta(A,B) + |AH \setminus A| + |BH \setminus B| = |H|$.
	\end{enumerate}
\item There exists $H \triangleleft G$ so that $G/H$ is dihedral with $|G/H| \ge 8$.  In addition, there exist $A^+ \supseteq A$ and 
	$B^+ \supseteq B$ so that $A^+H = A^+$ and $B^+H = B^+$ and $\varphi_{G/H}(A^+)$, $\varphi_{G/H}(B^+)$ are dihedral 
	progressions with a common ratio. Furthermore, one of the following holds.
	\begin{enumerate}
	\item $\delta(A,B) \le \frac{1}{2} |H|$ and $\max \{ |A^+ \setminus A|, |B^+ \setminus B| \} + \delta(A,B) \le 2 |H|$.
	\item $ABH = AB$ and $|H| = |B^+ \setminus BH| = |A^+ \setminus AH| 
	= \delta(A,B) + |AH \setminus A| + |BH \setminus B|$.	
	\item $A^+B^+ = AB$ and $\delta(A,B) + |A^+ \setminus A| + |B^+ \setminus B| = 2|H|$.
	\end{enumerate}
\item There exist $A^* \supseteq A$ and $B^* \supseteq B$ so that $A^*B^* = AB$ and there exists 
$(S,T) \in \{ (A^*,B^*), (B^*,C), (C, A^*) \}$ and $H_1,H_2,H_3 \le G$ so that $H_1 S H_2 = S$ and $H_2 T H_3 = T$ 
and one of the following holds.
	\begin{enumerate}
	\item $\frac{|S|}{|H_2|} = \frac{|T|}{|H_2|} = 2$ and $\delta(A,B) + |A^* \setminus A| + |B^* \setminus B| = |H_1| = |H_3| < |H_2|$.
	\item $\Big\{ \frac{|S|}{|H_2|}, \frac{|T|}{|H_2|} \Big\} = \{2,3\}$ and  $\delta(A,B) + |A^* \setminus A| + |B^* \setminus B| = \min \big\{ | H_1|, |H_3| \big\} < |H_2|$.
	\end{enumerate}
\item There exist $H \triangleleft G$ and $H \le H_1,H_2,H_3 \le G$ so that $A^* = H_1 A H_2$ and 
$B^* = H_2 B H_3$ satisfy $A^* B^* = AB$ and $\delta(A,B) + |A^* \setminus A| + |B^* \setminus B| \le \min\{ |H_1|, |H_2|, |H_3| \}$.  Furthermore, one of the following holds.
		\begin{center}
		\begin{tabular}{|c|c|}
		\hline
		$G/H$	& $\big\{ |G|/|A^*|, |G|/|B^*|, |G|/|C| \big\}$	\\
		\hline
		${\mathbf C_2} \times {\mathbf S_4}$, ${\mathbf C_2} \times {\mathbf A_4}$, or ${\mathbf S_4}$ 		
			&	$\{2,3,4\}$	\\
		${\mathbf C_2} \times {\mathbf A_5}$ or ${\mathbf A_5}$
			& 	$\{ 3/2, 4, 10 \}$, $\{4/3, 6, 10 \}$, or $\{5/3, 4, 6 \}$ \\
		${\mathbf S_6}$, ${\mathbf A_6}$, ${\mathbf S_5}$, or ${\mathbf A_5}$ 
			&	$\{ 2, 3, 5 \}$	\\
		\hline
		\end{tabular} 
		\end{center}
\end{enumerate}
\end{corollary}

\noindent{\it Proof:} Let $(A^*,B^*,C)$ be a maximal subset trio with $A^* \supseteq A$ and $B^* \supseteq B$ and define the Cayley choruses 
$\Theta = {\mathit CayleyC}(G; A,B,C)$ and $\Theta^* = {\mathit CayleyC}(G; A^*, B^*, C)$.  We will use $\sim$ to denote the incidence relation in $\Theta$, and
$\sim^*$ for the incidence relation in $\Theta^*$.  Note that 
$\delta(A,B) = \delta(\Theta) = \delta(\Theta^*) - |A^* \setminus A| - |B^* \setminus B|$.  If $(A^*,B^*,C)$ is trivial, then the first outcome holds.  Otherwise,  
Theorem \ref{main_incidence} implies that $\Theta^*$ is a song.  We now split into cases depending on $\Theta^*$.

\bigskip

\noindent{\it Case 1:} $\Theta^*$ is a pure or impure beat.

\smallskip

In this case it follows from Lemma \ref{i2s_purebeat} or \ref{i2s_impbeat} that the first outcome holds.

\bigskip

\noindent{\it Case 2:} $\Theta^*$ is a video.

\smallskip

Lemma \ref{i2s_video} implies that either the fourth or fifth outcome holds.

\bigskip

\noindent{\it Case 3:} $\Theta^*$ is a pure chord.

\smallskip

It follows from Lemma \ref{i2s_pchord} that $(A^*,B^*,C)$ is either a pure cyclic chord, a pure dihedral chord, or an impure dihedral chord of type $0$.  If $(A^*,B^*,C)$ is a pure cyclic chord relative to $H$, then $\delta(A^*,B^*,C) = |H|$, so $A^* = AH$ and $B^* = BH$ and outcome 2b holds.
If $(A^*,B^*,C)$ is a pure dihedral chord relative to $H$, then $\delta(A^*,B^*,C) = 2|H|$ and outcome 3c holds (with $A^+ = A^*$ and $B^+ = B^*$).  In the remaining case we may assume that $(A^*,B^*,C)$ is an impure dihedral chord of type $0$ relative to $H$.  We claim that in this case, outcome 3b is satisfied.  First observe that 
$ABH = AB$ and $A^*H = A^*$ and $B^*H = B^*$ and $\delta(A^*,B^*,C) = |H|$ which implies $AH = A^*$ and $BH = B^*$.  So, to show that 3b holds, we need only find suitable sets $A^{+}$ containing $AH$ 
and $B^{+}$ containing $BH$.  Now choose sets $A^+, B^+, C^+$ in accordance with the definition of impure dihedral chord so that 
$AH \subseteq A^+$ and $BH \subseteq B^+$ and $C \subseteq C^+$ and $A^+, B^+$ are dihedral progressions with a common ratio.  It follows from the definition of a type $0$ that the ends $(A_1,A_2,A_3,A_4)$ of $A^+$ satisfy one of the following
\begin{itemize}
\item $| \{ 1 \le i \le 4 \mid A_i \subseteq AH \} | = 3$ and $| \{1 \le i \le 4 \mid A_i \cap AH = \emptyset \} | = 1$
\item $| \{ 1 \le i \le 4 \mid A_i \subseteq AH \} | = 1$ and $| \{1 \le i \le 4 \mid A_i \cap AH = \emptyset \} | = 3$
\end{itemize}
In the former case $A^+$ contains $AH$ as desired.  In the latter we may shorten $A^+$ by replacing it with either
$A^+ \setminus (A_1 \cup A_2)$ or $A^+ \setminus (A_3 \cup A_4)$ to obtain a dihedral progression containing $AH$ as desired.  
A similar argument for $B^+$ shows that outcome 3b holds.

\bigskip

\noindent{\it Case 4:} $\Theta^*$ is a cyclic chord.

\smallskip

In this case we will show that outcome 2a is satisfied.  Suppose that $\Theta^*$ is a cyclic chord relative to the systems of imprimitivity $\mcp, \mcq, \mcr$ on $G \times \{1\}$, $G \times \{2\}$, $G \times \{3\}$, 
and let $H = G_{ (\mcp) } = G_{ (\mcq) } = G_{ (\mcr) }$.  It follows from the structure of a cyclic chord that $G/H$ is cyclic, so we may choose a consecutive shift $S \in G/H$.  Since $\Theta^*$ has a continuation which is a nontrivial $H$-trio with the same deficiency as $\Theta^*$, Lemma \ref{nontriv_bound} gives us $\delta(A,B) \le \delta(\Theta^*) \le \frac{1}{2}|H|$.  
Now by Lemma \ref{light_sides_chord} we find that $(G \times \{1\}, G \times \{2\} ; \sim )$ is a near sequence relative to $(\mcp, \mcq)$ and 
$w( \mcp, \mcq ; \sim ) - w(G \times \{1\}, G \times \{2\} ; \sim ) \le |H| - \delta(\Theta)$.  It follows from this that we may choose integers $\ell < m$ so that the set $A^{+} = S^{\ell} \cup S^{\ell + 1} \ldots \cup S^{m}$ satisfies $A \subseteq A^{+}$ and $|A^{+} \setminus A| \le |H| - \delta(A,B)$.  However, we must then have $A^{+} = AH$, which implies $|AH \setminus A| \le |H| - \delta(A,B)$.  A similar argument for the set $B$ shows that outcome 2a holds, as desired.

\bigskip

\noindent{\it Case 5:} $\Theta^*$ is a dihedral chord of type $1$ or $2$.

\smallskip

In this case we will show that outcome 3a is satisfied.  Suppose that $\Theta^*$ is a dihedral chord relative to the systems of imprimitivity $\mcp, \mcq, \mcr$ on $G \times \{1\}$, $G \times \{2\}$, $G \times \{3\}$, 
and set $H = G_{ (\mcp) } = G_{ (\mcq) } = G_{ (\mcr) }$.  By assumption we have that $G/H$ is dihedral, so we may choose a consecutive shift $S \in G/H$.  Observe that $\Theta^*$ cannot be a dihedral chord of type 2c, and if 
it has type 2b then $\delta(\Theta^*) = |K|$ for a proper subgroup $K < H$ so $\delta(\Theta^*) \le \frac{1}{2}|H|$.  Otherwise, $\Theta^*$ is either type $1$ or type $2a$, so it has a continuation which is a nontrivial $K$-trio for some $K \le H$ with the same deficiency as $\Theta^*$.  It follows from this and Lemma \ref{nontriv_bound} that $\delta(A,B) \le \delta(\Theta^*) \le \frac{1}{2}|H|$.
  
Now apply Lemma \ref{light_sides_chord} to the side $(G \times \{1\}, G \times \{2\} ; \sim )$ of $\Theta$.  This gives us systems of imprimitivity $\mcp', \mcq'$ of $G \times \{1\}$, $G \times \{2\}$ so that 
 $(G \times \{1\}, G \times \{2\} ; \sim )$ is a near sequence relative to $( \mcp', \mcq' )$ with dihedral action and consecutive shift $S$, and furthermore $w( \mcp', \mcq' ; \sim ) - w(G \times \{1\}, G \times \{2\} ; \sim ) \le 2|H| - \delta(\Theta)$.  It follows from this that we may choose integers $\ell < m$ and a flip $F \in G/H$ so that the set $A^{+} = (S^{\ell} \cup S^{\ell+1} \ldots \cup S^m)(H \cup F)$ satisfies $A \subseteq A^{+}$ and 
 $|A^{+} \setminus A| \le 2|H| - \delta(A,B)$.  A similar argument
for the set $B$ shows that outcome 3a holds, as desired.
\quad\quad$\Box$

\bigskip

\noindent{\it Proof of Corollary \ref{struc_or_stable}:} 
If $(A,B)$ is not critical then 3 holds for $H_1 = H_2 = H_3 = \{1\}$.  Otherwise, we apply Corollary
\ref{weak_structure}.  If we get outcome 1 from this corollary, then 1 holds for $(A,B)$.  If we get outcomes
2a or 3a then 2 holds  for $(A,B)$.  Outcomes 2b and 3b give us a subgroup $H \triangleleft G$ so that 
setting $H_1 = H_2 = H_3 = H$ we have that 3 holds.  If we get outcome 3c, then there exist $A^+, B^+ \subseteq G$ 
with $A \subseteq A^+$ and $B \subseteq B^+$ so that $(A^+,B^+,C)$ is a pure dihedral chord, so the subgroups 
$H_1 = {\mathit stab}_L(A^+)$, and
$H_2 = {\mathit stab}_R(A^+) = {\mathit stab}_L(B^+)$ and $H_3 = {\mathit stab}_R(B^+)$ all satisfy 
$[H_i : H] = 2$ and it follows that  3 holds.  Finally, outcomes 4 and 5 give us subgroups $H_1, H_2, H_3$ 
which satisfy 3. 
 \quad\quad$\Box$

\subsection{Deficiency vs. Disparity}

For a pair $(A,B)$ of subsets of $G$, we define the \emph{disparity} of $(A,B)$ to be $\nu(A,B) = \big| |A| - |B| \big|$.  In this section we will 
consider pairs $(A,B)$ for which $\delta(A,B) > \nu(A,B)$, or equivalently, $|AB| < 2 \min \{ |A|, |B| \}$ (note that any critical pair of the form $(A,A)$ satisfies this inequality).  Pairs which have deficiency greater than discrepancy can only form videos in a particular manner, as we will 
see in the next lemma, and this fact will allow us to simplify our weak structure theorem in this case.

\begin{lemma}
\label{video_disc}
Let $(A,B,C)$ be a subset trio in $G$ and assume $\Theta = {\mathit CayleyC}(G; A,B,C)$ is a video.  If $\delta(\Theta) > \nu(A,B)$  
then $\Theta \cong \Gamma : E \sim V \sim E$ for a graphic duet $\Gamma$.  Therefore,
there exists $H \le G$ so that $H = {\mathit stab}_R(A) = {\mathit stab}_L(B)$ satisfies $|A| = 2 |H| = |B|$ and 
$\delta(A,B,C) = | {\mathit stab}_L(A)| = | {\mathit stab}_R(B)| < |H|$.
\end{lemma}

\noindent{\it Proof:} Let $(X,Y,Z)$ be an arbitrary video and assume $|w(X,Y) - w(Y,Z)| < \delta(X,Y,Z)$.  If $w(X,Y) \neq w(Y,Z)$, then follows from Lemma \ref{video_crit} that $| w(X,Y) - w(Y,Z) | \ge w^{\circ}(Y) \ge \delta(\Theta)$ which is a contradiction.  So, it must be that $w(X,Y) = w(Y,Z)$.  A straightforward analysis
using Figure \ref{alldeg} implies that this is only possible when $(X,Y,Z) \cong \Gamma : E \sim V \sim E$ for a graphic duet $\Gamma$.

Based on this analysis, we find that $\Theta \cong \Gamma : E \sim V \sim E$ for a graphic duet $\Gamma$ with degree at least three.  If the clone 
partition of  $G \times \{i\}$ (in the trio $\Theta$) is given by $G/H_i \times \{i\}$ for $1 \le i \le 3$ then $H_1 = {\mathit stab}_L(A)$ and ${\mathit stab}_R(A) = H_2 = {\mathit stab}_L(B)$ and $H_3 = {\mathit stab}_R(B)$.  Further, $|A| = 2|H_2| = |B|$ and $\delta(\Theta) = \delta(A,B,C) = |H_1| = |H_3| < |H_2|$ as desired. 
\quad\quad$\Box$

\begin{lemma}
\label{conj-no-vid}
Let $A,B \subseteq G$ be finite and assume $B = x^{-1}Ax$ for some $x \in G$ and that $(A,B)$ is critical.  If $C = \overline{(AB)^{-1}}$ and $(A^*,B^*,C)$ is  a maximal critical trio with $A \subseteq A^*$ and $B \subseteq B^*$, then ${\mathit CayleyC}(G; A^*,B^*,C)$ is not a video.
\end{lemma}

\noindent{\it Proof:} Suppose (for a contradiction) that $\Theta = {\mathit CayleyC}(G; A^*, B^*, C)$ is a video.  Since 
$\delta(A,B) > 0 = \nu(A,B)$ it follows that $\delta(\Theta) = \delta(A^*,B^*) > \nu(A^*,B^*)$.  So, by the previous lemma the subgroups 
$H_1 = {\mathit stab}_L(A^*)$, $H_2 = {\mathit stab}_R(A^*) = {\mathit stab}_L(B^*)$, $H_3 = {\mathit stab}_R(B^*)$ satisfy 
$|A^*| = 2|H_2| = |B^*|$ and $\delta(\Theta) = |H_1| = |H_3| < |H_2|$.  

Now, on one hand, we have 
\[ 2 |A^* \setminus A| = |A^* \setminus A| + |B^* \setminus B| = \delta(A^*,B^*,C) - \delta(A,B,C) < |H_1|\]
so $A$ is very efficiently contained in $A^*$ which has left stabilizer $H_1$.  On the other hand, the set
$B^* = H_2 B$ contains $B = x^{-1}Ax$ so $2|H_2| = | H_2 x^{-1} A x | = | x H_2 x^{-1} A |$, which means that 
the subgroup $K = x H_2 x^{-1}$ satisfies $|KA| = 2|K| = 2 |H_2|$.  Thus, $A$ is also efficiently contained in 
$KA$ which is left stabilized by $K$.  

Let $d = |A^*| / |H_1|$ (so $d$ is the number of right $H_1$-cosets contained in $A^*$) and note that $d \ge 3$ by assumption.  
Since $|A^* \setminus A| < \frac{1}{2} |H_1|$ it follows that we may choose $y \in A$ so that $|H_1 y \cap A| > (1-\frac{1}{2d})|H_1| \ge \frac{2}{3} |H_1|$.
Since $|KA| = 2|K|$ it follows that we may choose $z \in A$ so that $|Kz \cap H_1 y \cap A| > \frac{1}{3}|H_1|$.  However, then 
$Kz \cap H_1 y$ must be a right $L$-coset for a subgroup $L \le H_1 \cap K$ with $[H_1 : L] \le 2$.  

Now consider the set $A' = LA$.  Since $L \le H_1$ and $H_1 = {\mathit stab}_L(A^*)$ it must be that $A' \subseteq A^*$ and since 
$|A^* \setminus A| < |L|$ we must have $A' = A^*$.  Note that this implies $|A'| = 2|H_2| = 2|K|$.  However, $L \le K$ so $A' \subseteq KA$ and since $|A'| = 2|K|$ we must have $A' = KA$.  However, this contradicts ${\mathit stab}_L(A^*) = H_1$, thus completing the proof.
\quad\quad$\Box$

\begin{corollary}[Deficiency vs.\ Disparity]
\label{def_vs_disp}
Let $A,B \subseteq G$ be finite and nonempty, and assume $|AB| < 2 \min \{ |A|, |B| \}$.  Let $x \in A$ and $y \in B$ and assume that 
$H \le G$ is the unique minimal subgroup for which $A \subseteq xH$.  Then $B \subseteq Hy$, and one of the following holds.
\begin{enumerate}
\item Either $AB = xHy$ or there exists $K < H$ and $z \in H$ so that $x( H \setminus zK )y \subseteq AB$.
\item There exists $K \triangleleft H$ so that $H/K$ is a cyclic group with $|H/K| \ge 4$ so that  $\varphi_{H/K}(x^{-1}A)$ and $\varphi_{H/K}(By^{-1})$ are geometric progressions with a common ratio.   Furthermore, one of the following holds.
	\begin{enumerate}
	\item $\delta(A,B)  \le \frac{1}{2} |K|$ and $\delta(A,B) + \max\{ |AK \setminus A|, |KB \setminus B| \} \le |K|$.
	\item $AKB = AB$ and $\delta(A,B) + |AK \setminus A| + |KB \setminus B| = |K|$.
	\end{enumerate}
\item There exists $K \triangleleft H$ so that $H/K$ is dihedral with $|H/K| \ge 8$ and sets $A^+ \supseteq A$ and $B^+ \supseteq B$ so that
$A^+K = A^+$ and $KB^+ = B^+$ and $\varphi_{H/K}(x^{-1}A^+)$, $\varphi_{H/K}(B^+y^{-1})$ are dihedral progressions with a common ratio. Furthermore, one of the following holds.
	\begin{enumerate}
	\item $\delta(A,B) \le \frac{1}{2} |K|$ and $\delta(A,B) + \max \{ |A^+ \setminus A|, |B^+ \setminus B| \} \le 2 |K|$.
	\item $AKB = AB$ and $|K| = |B^+ \setminus KB| = |A^+ \setminus AK| 
	= \delta(A,B) + |AK \setminus A| + |KB \setminus B|$.	
	\item $A^+B^+ = AB$ and $\delta(A,B) + |A^+ \setminus A| + |B^+ \setminus B| = 2|K|$.
	\end{enumerate}
\item There exist $K_1, K_2, K_3 \le H$ so that $K_1 A K_2 B K_3 = AB$ and furthermore
	\begin{enumerate}
	\item $|AK_2| = 2|K_2| = |K_2 B|$
	\item $\delta(A,B) + |AK_2 \setminus A| + |K_2 B \setminus B| = |K_1| = |K_3| < |K_2|$.
	\end{enumerate}
\end{enumerate}
Furthermore, if $A = B$, then $xH = Hx$ and outcome 4 may be excluded.
\end{corollary}

\noindent{\it Proof:} First we prove that $B$ is contained in a single right $H$-coset.  Suppose (for a contradiction) that $y,z \in B$ satisfy $Hy \cap Hz = \emptyset$.  Then we have $|AB| \ge |A \{y,z \}| = 2 |A|$ which is contradictory.  Thus we have $A \subseteq xH$ and $B \subseteq Hy$.  Furthermore, by a similar argument, $H$ must be the minimal subgroup so that $B \subseteq Hy$.  Now define 
the sets $A' = x^{-1}A$ and $B' = By^{-1}$.  It follows that $(A',B')$ is a critical pair in $H$.  Furthermore, neither $A'$ nor $B'$ is contained in a coset of 
a proper subgroup of $H$.  

Next define $C = H \setminus (A'B')^{-1}$ and note that we may assume $C \neq \emptyset$ as otherwise the first outcome is satisfied.  
Let $(A^*,B^*,C)$ be a maximal trio in $H$ with $A' \subseteq A^*$ and $B' \subseteq B^*$ and define $\Theta = {\mathit CayleyC}(H; A^*, B^*, C)$.    
Now the remaining part of the argument follows the proof of Corollary \ref{weak_structure} with a few exceptions which we now detail.

First, if $\Theta$ is a pure or impure beat, then $C$ must be contained in a coset of a proper subgroup, giving us the first outcome.  
Second, if $\Theta$ is a video, then the bound $\nu(A',B') < \delta(A',B')$  implies $\nu(A^*,B^*) < \delta(A^*,B^*) = \delta(\Theta)$, so Lemma \ref{video_disc} implies the existence of subgroups $K_1,K_2,K_3$ as given in the final outcome.  

Finally, if $A = B$ then we may take $x = y$, so $A \subseteq xH$ and $A \subseteq Hx$.  However, this implies 
$A \subseteq xH \cap Hx = x( H \cap x^{-1} H x )$ and then by the choice of $H$ we must have $x^{-1}Hx = H$.  Furthermore, 
the sets $A' = x^{-1}A$ and $B' = A x^{-1}$ now satisfy $B' = x A' x^{-1}$ so Lemma \ref{conj-no-vid} implies that $\Theta$ cannot be a video, 
which excludes the fourth outcome.
\quad\quad$\Box$

\bigskip

\noindent{\it Proof of Corollary \ref{a_squared}:} Applying Corollary \ref{def_vs_disp} to the pair $(A,A)$ implies that $xH=Hx$ and one of the first three 
outcomes of this corollary holds.  The second and third outcomes of this corollary immediately yield the second and third outcomes for Corollary \ref{a_squared}.  If we get the first outcome of Corollary \ref{def_vs_disp} and $AA = xHx$ then since $|A| \ge 2$ the subgroup $K = \{1\} < H$ has the desired properties.  Otherwise this first outcome gives us a subgroup $K$ and some conjugate of $K$ has the desired properties.
\quad\quad$\Box$

\subsection{Generalizing Kneser's Theorem}

In this subsection we prove Corollaries \ref{big_cor} and \ref{cor_basic} which generalize of Kneser's addition theorem.  We 
will deduce these as corollaries of Theorem \ref{subset_equiv} which allows us to work primarily in the context of subsets of 
groups (instead of incidence geometry), however this theorem still has a conclusion concerning videos which will require us
to consider these particular incidence geometries.  

The corollaries we are proving here concern the notion of conjugate stability (recall that for a subgroup $H \le G,$ a set $A \subseteq G$ is $H$-conj-stable if for every $x \in A$ there exists $y \in G$ so that $x(y^{-1}Hy) \subseteq A$).  For a subgroup $H \le G$ we define a trio $(A,B,C)$ 
to be $H$-\emph{controlled} if $A,B,C$ are all $H$-conj-stable and $\delta(A,B,C) \le |H|$.  We say that $(A,B,C)$ is \emph{controlled} 
if it is $H$-controlled for some $H \le G$.  The main result of this subsection (from which we will derive Corollaries \ref{big_cor} and \ref{cor_basic}) 
is as follows.

\begin{theorem}
\label{controlled}
Every maximal subset trio is controlled.
\end{theorem}

Although the notion of conjugate stability is much weaker than that of stability, it does have some convenient properties, as indicated by the following straightforward observation.

\begin{observation}
\label{conj-stab-obs}
Let $A,B \subseteq G$ and let $g \in G$.
\begin{enumerate}
\item If $H \le {\mathit stab}_L(A)$ or $H \le {\mathit stab}_R(A)$, then $A$ is $H$-conj-stable.
\item If $A$ is $H$-conj-stable, then $AB$, $BA$, $gA$, $Ag$, and $A^{-1}$ are all $H$-conj-stable.
\item If $(A,B,C)$ and $(A',B',C')$ are similar subset trios, and $A,B,C$ are $H$-conj-stable, then $A',B',C'$ are $H$-conj-stable.
\end{enumerate}
\end{observation}

Our proof will call upon a lemma concerning videos which we give next.

\begin{lemma}
\label{control_video}
If $(A,B,C)$ is a subset trio and $\Theta = {\mathit CayleyC}(G; A, B, C)$ is a video, then $(A,B,C)$ is controlled.
\end{lemma}

\noindent{\it Proof:} By possibly replacing $(A,B,C)$ by a similar trio, we may assume that 
$\Theta^{\bullet}$ is isomorphic to one of the trios given in the definition of standard and exceptional video.  Choose 
subgroups $H_i \le G$ for $1 \le i \le 3$ so that $G/H_i \times \{i\}$ is the clone partition of $G \times \{i\}$ given by $\Theta$.  Now, Lemma 
\ref{crit2cut} relates the deficiency of $\Theta$ with the number of edges in certain edge-cuts.  
Two of the edge-cuts appearing in our structures have size $2d-2$ where $d$ is the degree of the graph $\Gamma$, while 
all of the others have size $2d-1$.  We first handle these two special cases, and 
then give a general argument for the others. 

\bigskip

\noindent{\it Case 1:} $\Theta^{\bullet} \cong \Gamma : E \sim V \sim E$ for a graphic duet $\Gamma$.

\smallskip

In this case $H_1$ and $H_3$ are subgroups stabilizing edges of $\Gamma$ so $H_1$ and $H_3$ are conjugate.  Since 
$H_1 \le {\mathit stab}_L(A)$ and $H_3 \le {\mathit stab}_R(B)$ and $H_1 \le {\mathit stab}_R(C)$ we have that $A,B,C$ are 
all $H_1$-conj-stable.  It follows from Lemma \ref{crit2cut} that $\delta(A,B,C) = \delta(\Theta) = |H_1|$, so $(A,B,C)$ is 
$H_1$-controlled.

\bigskip

\noindent{\it Case 2:} $\Theta^{\bullet} \cong \mathrm{Cube/Octahedron}: V \sim E \sim F$.

\smallskip

We may assume (without loss) that $\Theta^{\bullet} \cong \mathrm{Cube}: V \sim E \sim F$.  In this case, it follows immediately 
that $A$ and $B$ are $H_2$-conj-stable since $H_2 \le {\mathit stab}_R(A)$ and $H_2 \le {\mathit stab}_L(B)$.  To see 
that $C$ is also $H_2$-conj-stable, we shall construct another incidence geometry.  For this purpose, we define 
the relation $x \antipode y$ if $x$ and $y$ are antipodal faces in the cube.  
 
\begin{figure}[ht]
\centerline{\includegraphics[width=6cm]{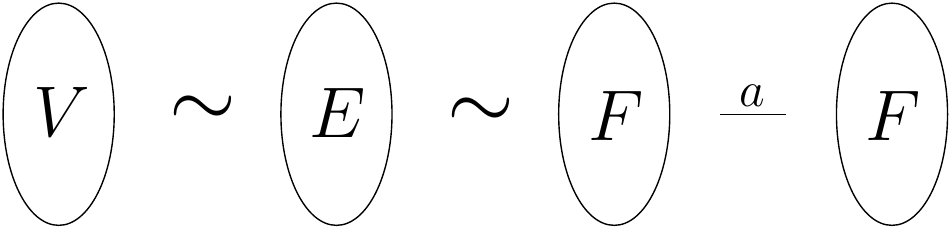}}
\caption{An incidence geometry based on Cube}
\label{cube_inc_geom}
\end {figure}

Let $V = V( \mathrm{Cube} )$ and $E = E( \mathrm{Cube} )$ and $F = F( \mathrm{Cube} )$.  Next define an incidence geometry $\Lambda = (V, E, F \times \{1\}, F \times \{2\})$ with incidence given by the following 
rules for every $v \in V$ and $e\in E$ and $(x,1) \in F \times \{1\}$ and $(y,2) \in F \times \{2\}$.
\begin{itemize}
\item $v$ and $e$ are incident in $\Lambda$ if they are incident in Cube.
\item $e$ and $(x,1)$ are incident in $\Lambda$ if $e,x$ are incident in Cube.
\item $(x,1)$ and $(y,2)$ are incident in $\Lambda$ if $x \antipode y$.   
\item there are no incidences between $V$ and $F \times \{1\} \cup F \times \{2\}$ or between $E$ and $F \times \{2\}$.
\end{itemize}
We already have a realization of a part of $\Lambda$ using the subgroups $H_1, H_2, H_3$ and the sets $A$ and $B$.  By applying Theorem \ref{sabidussi} we may extend this by choosing 
a set $D \subseteq G$ with $H_3 D H_3 = D$ so that the chorus ${\mathit CayleyC}(G; M)$ given by the matrix
\[ M = \left[ \begin{array}{cccc}
  			\emptyset		&A			& \emptyset	& \emptyset	\\
			A^{-1}		&\emptyset	& B			& \emptyset	\\
			\emptyset		& B^{-1}		& \emptyset	& D			\\
			\emptyset		& \emptyset	& D^{-1}		& \emptyset 
			\end{array} \right] \]
has quotient $( G/H_1 \times \{1\} , G/H_2 \times \{2\} , G/H_3 \times \{3\}, G/H_3 \times \{4\} )$ strongly isomorphic to $\Lambda$.  It follows from 
this construction that $ABD = C$.  Since $A$ is $H_2$-conj-stable, $C$ must also be $H_2$-conj-stable, and since 
$\delta(A,B,C) = \delta(\Theta) = |H_2|$ we have that $(A,B,C)$ is $H_2$-controlled.

\bigskip

\noindent{\it Case 3:} $\Theta^{\bullet}$ is not isomorphic to one of the videos from Case 1 or Case 2.

\smallskip

In all of these cases, one of the subgroups, say $H_i$, corresponds to the stabilizer of an edge, and another, say $H_j$ corresponds to the stabilizer of a vertex.  Choose $H \le G$ to be the stabilizer of an arc.  Then, $H$ is conjugate to a subgroup of $H_i$ and to a subgroup of $H_j$.  Since each of the three sets $A,B,C$ is either right or left stabilized by one of the subgroups $H_i$ or $H_j$, it follows that $A,B,C$ are all $H$-conj-stable.  
Now, the stabilizer of an arc will either be equal to the stabilizer of the corresponding edge, or a subgroup of this group of index two, so we have $|H| \ge \frac{1}{2}|H_i|$.  Then by Lemma \ref{crit2cut} we have $\delta(A,B,C) = \delta(\Theta) = \frac{1}{2}|H_i| \ge |H|$ which shows that 
$(A,B,C)$ is $H$-controlled.
\quad\quad$\Box$

\bigskip

We are now ready for the main result of this subsection.

\bigskip

\noindent{\it Proof of Theorem \ref{controlled}:}  Let $\Phi$ be a maximal subset trio in $G$.  If $\delta(\Phi) \le 0$ then 
$\Phi$ is $\{1\}$-controlled, so we may assume $\Phi$ is critical.  If $\Phi$ is trivial, then it is similar to 
$(G,G,\emptyset)$ so it is $G$-controlled.  Otherwise, $\Phi$ is a nontrivial maximal critical trio, so we may apply Theorem
\ref{subset_equiv} to obtain sequence of subset trios $\Phi_1, \ldots, \Phi_m$.  The desired result now follows from 
Lemma \ref{control_video}, the following claim, and a straightforward induction.

\bigskip

\noindent{\it Claim:} Let $(A,B,C)$ be a subset trio in $G$.
\begin{enumerate}
\item If $(A,B,C)$ is a pure beat, pure cyclic chord, or pure dihedral chord, then $(A,B,C)$ is controlled.  
\item If $(A,B,C)$ is an impure beat or an impure cyclic chord with a continuation which is $K$-controlled, 
		then $(A,B,C)$ is $K$-controlled.  
\item If $(A,B,C)$ is an impure dihedral chord of type $0$ or $2B$, then it is controlled.  If it has type $1$ or $2A$ with 
	 a continuation which is $K$-controlled, then $(A,B,C)$ is $K$-controlled.
\end{enumerate}

So to complete the proof, we need only to establish the claim, which we do in parts.

\bigskip

\noindent{Proof of 1. } If $(A,B,C)$ is either a pure beat relative to $H$, or a pure cyclic chord relative to $H$, 
then it follows from Observation \ref{conj-stab-obs} and our definitions that $(A,B,C)$ is $H$-controlled.  If $(A,B,C)$ is 
a pure dihedral chord relative to $H$, then there exist subgroups $H_1, H_2, H_3$ 
where each $H_i$ has the form $H_i = H \cup F$ for a flip $F \in G/H$ 
so that $H_1 A H_2 = A$ and $H_2 B H_3 = B$ and $H_3 C H_1 = C$.  Now, in the dihedral group ${\mathbf D}_n$, the 
subgroups of the form $\{1,f\}$ where $f$ is a flip will form a single conjugacy class if $n$ is odd, and form two conjugacy 
classes if $n$ is even.  It follows that we may choose $1 \le i < j \le 3$ so that $H_i$ and $H_j$ are conjugate.  Now since
each of our sets $A,B,C$ is either left or right stabilized by $H_i$ or $H_j$ it follows that $A,B,C$ are all $H_i$-conj-stable, 
and since $\delta(A,B,C) = 2|H|$, we have that $(A,B,C)$ is $H_i$-controlled, as desired.

\bigskip

\noindent{Proof of 2. } Let $(A,B,C)$ be an impure beat or impure cyclic chord relative to $H$ with a $K$-controlled 
continuation (so $K < H$).  By possibly replacing $(A,B,C)$ with a similar trio, we may assume that
$(A \cap H, B \cap H, C \cap H)$ is a continuation which is a $K$-controlled $H$-trio.  Now $A \setminus H$ and $B \setminus H$ and $C \setminus H$
are all unions of either left or right $H$-cosets, so it follows immediately from this and the assumption that $(A',B',C')$ is 
$K$-controlled that $(A,B,C)$ is $K$-controlled.

\bigskip

\noindent{Proof of 3. } Let $(A,B,C)$ be an impure dihedral chord relative to $H$ and assume it has type $0$, $1$, $2A$, or $2B$.
If $(A,B,C)$ is type $0$, then $A,B,C$ are all $H$-stable and $\delta(A,B,C) = |H|$ so $(A,B,C)$ is $H$-controlled.  Next suppose 
that $(A,B,C)$ is type $1$ and that $(A',B',C')$ is a continuation of $(A,B,C)$ (so $(A',B',C')$ is an $H$-trio) which is $K$-controlled.
In this case, there exist $R, S, T \in G/H$ so that $A \setminus R, B \setminus S, C \setminus T$ are $H$-stable and 
$(A \cap R, B \cap S, C \cap T)$ is similar to $(A',B',C')$ when viewed as subset trios of $G$.  It follows immediately that 
$A,B,C$ are $K$-conj-stable, and thus $(A,B,C)$ is $K$-controlled.  

In the remaining cases, $(A,B,C)$ is an impure dihedral 
chord of type $2A$ or $2B$ and we may assume that $\Omega$ is an associated octahedral configuration and that 
the subsets $A',B',C', B'', C'' \subseteq H$ and subgroups $K',K'' < H$ are as given in the definitions of these types.  It follows that
for every $S \in G/H$ the set $A \cap S$ ($B \cap S$, $C \cap S$) is either empty, equal to $S$, or conjugate to $A'$ 
(to one of $B', B''$, or to one of $C', C''$). If $(A,B,C)$ has type $2A$, then $(A',B',C')$ is an impure beat relative to $K'$ with a continuation which is $K$-controlled, so we have that $A',B',C'$ are all $K$-conj-stable and $\delta(A',B',C') = \delta(A,B,C) \le |K|$.  Since $K' \le K''$ and for $x \in A$ we have that $(xK'', B'', C'')$ is a pure beat relative to $K''$ it must be that $B''$ and $C''$ are also $K$-conj-stable, and it follows that 
$(A,B,C)$ is $K$-controlled.  
If $(A,B,C)$ has type $2B$, then $(xK',B',C')$ is a pure beat relative to $K'$ so these sets 
are $K'$-conj-stable, and $(xK'',B'',B'')$ is a pure beat relative to $K''$ so these sets are $K''$-conj-stable, and then for the subgroup
$K = K' \cap K''$ we find that $A,B,C$ are all $K$-conj-stable, and $\delta(A,B,C) = |A'| = |K|$ so $(A,B,C)$ is $K$-controlled, as desired.
\quad\quad$\Box$

\bigskip

\noindent{\it Proof of Corollay \ref{big_cor}:} 
Let $A,B \subseteq G$ be finite and nonempty, let $C = \overline{(AB)^{-1}}$ and choose $A \subseteq A^*$ and $B \subseteq B^*$ so that 
$(A^*,B^*,C)$ is a maximal subset trio.  By the previous theorem there exists a subgroup $H$ so that $A^*,B^*,C$ are $H$-conj-stable and 
$|A^*| + |B^*| - |AB| = \delta(A^*,B^*,C) \le |H|$.  Since $A^*$ is $H$-conj-stable and $AB = A^*B^*$, the set $AB$ is also $H$-conj-stable, 
and this completes the proof.
\quad\quad$\Box$

\bigskip

\noindent{\it Proof of Corollary \ref{cor_basic}:} 
This follows immediately from Corollary \ref{big_cor}.
\quad\quad$\Box$

\section*{Appendix}

Here we prove a result mentioned in Section \ref{intro_sec} which is of interest in conjunction with a theorem of Freiman.  

\begin{proposition*}
Let $B \subseteq G$ and $K \le G$ be finite and assume $KBK = B$ and $|B| = 2|K|$ and $|B^2| < 2|B|$.  Then there exist 
$L \le K \le H \le G$ and $x \in G$ so that the following hold.
\begin{enumerate}
\item  $B \subseteq xH = Hx$
\item $[K:L] \le 2$
\item $L \triangleleft H$ and $H/L$ is either cyclic or dihedral.
\end{enumerate}
\end{proposition*}

\noindent{\it Proof:} Let $x \in B$ and let $H$  be the minimal subgroup for which $B \subseteq xH$.  Suppose (for a contradiction) that 
$B \not\subseteq Hx$ and choose $y \in B \setminus Hx$.  Now $xHx \cap xHy = \emptyset$ so $Bx \cap By = \emptyset$, but then $|B^2| \ge |B \{x, y \}| \ge 2|B|$ gives us a contradiction.  Thus $B \subseteq Hx \cap xH = x( x^{-1}Hx \cap H )$ from which it follows that $x^{-1}Hx = H$ and thus 1 is satisfied.  
If $K \neq {\mathit stab}_L(B)$ then $B$ is a single left $H$-coset and setting $L=K$ satisfies the proposition.  So, we may assume 
$K = {\mathit stab}_L(B)$, and by a similar argument $K = {\mathit stab}_R(B)$.  

Now define $\Delta = {\mathit CayleyC}(G; B)$ and note that $K = {\mathit stab}_L(B) = {\mathit stab}_R(B)$, so we have that 
$G/K \times \{1\}$ and $G/K \times \{2\}$ are the clone partitions of $G \times \{1\}$ and $G \times \{2\}$.  Consider the sets $P = H \times \{1\}$ and 
$Q = Hx \times \{2\}$ and note that our construction implies that $(P,Q)$ is a connected $H$-duet.  Furthermore, it follows from the stability properties of $B$ 
that $(P,Q)^{\bullet}$ is a polygon.  Define $L = \{ g \in H \mid \mbox{$g$ fixes each point in $(P,Q)^{\bullet}$} \}$ and note that by construction 
$L \triangleleft H$.  Now $ K \times \{1\}$ is a single point in $(P,Q)^{\bullet}$ with stabilizer $K$, and it follows from the structure of a polygon that 
$[K:L] \le 2$.  The action of $H$ on $(P,Q)^{\bullet}$ gives us an isomorphism from $H/L$ to a subgroup of the automorphism group of a polygon 
(which is dihedral), and this yields the final part of the proposition.
\quad\quad$\Box$

\end{document}